\documentclass[12pt]{article}
\usepackage[authoryear,round]{natbib}
\usepackage{graphicx}
\usepackage{latexsym}
\usepackage[left=1.5in,top=1.5in,right=1.5in,bottom=1.5in]{geometry}\usepackage{amsmath}
\usepackage{booktabs}
\usepackage{color}
\usepackage{tikz}
\usetikzlibrary{shapes}
\usetikzlibrary{arrows}
\definecolor{darkblue}{rgb}{0,0.4,0.9}
\definecolor{gray10}{rgb}{0.1,0.1,0.1}
\definecolor{gray20}{rgb}{0.2,0.2,0.2}
\definecolor{gray30}{rgb}{0.3,0.3,0.3}
\definecolor{gray40}{rgb}{0.4,0.4,0.4}
\definecolor{gray60}{rgb}{0.6,0.6,0.6}
\definecolor{gray80}{rgb}{0.8,0.8,0.8}
\definecolor{gray90}{rgb}{0.9,0.9,.9}
\definecolor{gray95}{rgb}{0.95,0.95,.95}
\definecolor{gray96}{rgb}{0.96,0.96,.96}
\definecolor{lgreen} {RGB}{180,210,100}
\definecolor{dblue}  {RGB}{20,66,129}
\definecolor{ddblue} {RGB}{11,36,69}
\definecolor{lred}   {RGB}{220,0,0}
\definecolor{nred}   {RGB}{224,0,0}
\definecolor{norange}{RGB}{230,120,20}
\definecolor{nyellow}{RGB}{255,221,0}
\definecolor{ngreen} {RGB}{98,158,31}
\definecolor{dgreen} {RGB}{78,138,21}
\definecolor{nblue}  {RGB}{28,130,185}
\definecolor{jblue}  {RGB}{20,50,100}
\definecolor{nnyellow}{RGB}{235,200,0}
\definecolor{purple}{RGB}{150, 0, 120}
\definecolor{sgGreen} {RGB}{20, 180, 50}
\definecolor{revised}{rgb}{0,0,0.9}
\newtheorem{definition}{Definition}

\newtheorem{theorem}{Theorem}

\newtheorem{lemma}{Lemma}

\newcommand{\nl}{\newline}
\newcommand{\pl}{\parallel}
\newcommand{\openr}{\hbox{${\rm I\kern-.2em R}$}}
\newcommand{\openn}{\hbox{${\rm I\kern-.2em N}$}}
\newcommand{\argmax}[1]{\underset{#1}{\operatorname{argmax}}\;}

\usepackage{authblk}

\title{Higher Order Targeted Maximum Likelihood Estimation}

\author{Mark van der Laan}
\author{Zeyi Wang}
\author{Lars van der Laan}
\affil{University of California, Berkeley}
 
\begin{document}
\maketitle
 
\begin{abstract}
Asymptotic efficiency of targeted maximum likelihood estimators (TMLE)  of target features of the data distribution relies on a second order remainder being asymptotically negligible \citep{vanderLaan&Rubin06}. 
In previous work we proposed a nonparametric MLE termed Highly Adaptive Lasso (HAL) which parametrizes the relevant functional of the data distribution in terms of a multivariate real valued cadlag function that is assumed to have finite variation norm 
\citep{Benkeser&vanderLaan16,vanderLaan15,vanderLaan17}. We showed that the HAL-MLE converges in Kullback-Leibler dissimilarity at a rate $n$-1/3 up till log-n factors \citep{Bibaut&vanderLaan19}. Therefore, by using HAL as initial density estimator in the TMLE, the resulting HAL-TMLE is an asymptotically efficient estimator for realistic statistical models only assuming that the relevant nuisance functions of the data density are cadlag and have finite variation norm \citep{vanderLaan17}. However, in finite samples, the second order remainder can dominate the sampling distribution so that inference based on asymptotic normality would be anti-conservative.

 In this article we propose a new higher order (say $k$-th order) TMLE, generalizing the regular (first order) TMLE, which is like a regular TMLE targeting sequentially defined data-adaptive higher order HAL-regularized TMLE-fluctuations of the target parameter.  We prove that it satisfies an exact linear expansion, in terms of the efficient influence functions  of the sequentially defined higher order  fluctuations of the target parameter, with a remainder that is a $k+1$-th order remainder, and a HAL-regularization-bias term that is controlled by setting the $L_1$-norm in the HAL-MLE. We show that this HAL-regularization-bias term is, even without undersmoothing,  guaranteed to be of small enough order. As a consequence, this $k$-th order TMLE allows statistical inference only relying on the $k+1$-th order remainder being negligible. 
 
We also provide a finite sample rational for the higher order TMLE that demonstrates that it will be superior to the first order TMLE by (iteratively) locally (around the initial estimator) minimizing (in the choice of initial estimator) the exact finite sample remainder of the first order TMLE.
The second order TMLE is demonstrated for nonparametric estimation of the integrated squared density and for the treatment specific mean outcome. We also provide an initial  simulation study for the second order TMLE in these two examples confirming the theoretical analysis. 
\end{abstract}


\section{Introduction}
Consider the problem of statistical estimation of a real valued target feature $\Psi(P_0)$ of the probability distribution $P_0$ based on observing $n$ independent and identically distributed copies  $O_1,\ldots,O_n$ from $P_0$, and knowing that $P_0$ is an element of a specified  set ${\cal M}$ of possible probability distributions. 
 This set ${\cal M}$ of possible distributions of $O$  is called the statistical model for the data distribution $P_0$. This article generalizes immediately to  Euclidean valued target parameters. 
 
We consider the case that the target feature is a pathwise differentiable mapping $\Psi:{\cal M}\rightarrow\openr$ from the set ${\cal M}$ of possible data distributions to the real line.  That is, for a collection of paths $\{P_{\delta,h}:\delta\}\subset {\cal M}$ through $P$ with score $h\in {\cal H}$ at $\delta=\delta_0\equiv 0$, we have $
\frac{d}{d\delta_0}\Psi(P_{\delta_0,h})=P D^{(1)}_P h$, where $D^{(1)}_P$ is an element of the tangent space $T_P({\cal M})$ at $P$. Here we use the notation $Pf\equiv \int f(o)dP(o)$ for the expectation of $f(O)$ under $P$. The tangent space $T_P({\cal M})$ is defined as the closure of the linear span of all the scores $h\in {\cal H}$ in the Hilbert space $L^2_0(P)$, consisting of all functions of $O$ with mean zero and finite variance, endowed with inner product $\langle h_1,h_2\rangle_P=P h_1h_2$. The unique element $D^{(1)}_P$ is called the canonical gradient of the pathwise derivative of $\Psi$ at $P$.
Let $P_n$ be the empirical probability measure of $O_1,\ldots,O_n$, and we use the notation $Pf\equiv \int f(o)dP(o)$.

Efficiency theory for this statistical estimation problem teaches us that an estimator $\hat{\Psi}(P_n)$ is asymptotically efficient at $P_0$  if and only if it is asymptotically linear with influence curve  equal to the canonical gradient $D^{(1)}_{P_0}$ of the pathwise derivative of the target parameter \citep{bkrw1997}: $\hat{\Psi}(P_n)-\Psi(P_0)=P_n D^{(1)}_{P_0}+o_P(n^{-1/2})$.
A TMLE is a two stage substitution estimator $\Psi(P_n^*)$ that first obtains an initial estimator $P_n^0\in {\cal M}$ of the data distribution $P_0$ and subsequently computes the MLE along a parametric path $\{\tilde{P}^{(1)}(P_n^0,\epsilon):\epsilon\}\subset {\cal M}$ through the initial estimator $P_n^0$ at $\epsilon=0$, chosen so that the linear span of the scores of $\epsilon$ at $\epsilon=0$  spans the canonical gradient $D^{(1)}_{P_0}$ \citep{vanderLaan&Rubin06,vanderLaan08, vanderLaan&Rose11,vanderLaan&Rose18}.
Such a path is often called a least favorable parametric submodel.
One then defines $P_n^*=\tilde{P}_n^{(1)}(P_n^0)=\tilde{P}^{(1)}(P_n^0,\epsilon_n^{(1)})$ as this TMLE-update of $P_n^0$, where $\epsilon_n^{(1)}=\arg\min_{\epsilon} P_n L(\tilde{P}^{(1)}(P_n^0,\epsilon))$ and $L(P)=-\log p$ is the log-likelihood loss. 

 We defined a universal least favorable path $\tilde{P}_n^{(1)}(P,\epsilon)$ for a one dimensional target estimand as a one-dimensional path through $P$ so that the score at any $\epsilon$, not just at $\epsilon =0$, equals the canonical gradient $D^{(1)}_{\tilde{P}^{(1)}(P_n^0,\epsilon)}$ of $\Psi$ at $\tilde{P}^{(1)}(P_n^0,\epsilon)$ \citep{vanderLaan&Gruber15,vanderLaan&Rose18,Cai&vanderLaan19}. 
We also defined a universal least favorable path for a multidimensional target estimand as a one dimensional (data dependent) path for which  the derivative at  $\epsilon$ of the log-likelihood equals the Euclidean norm of the empirical mean of vector efficient influence curve at $\tilde{P}^{(1)}(P_n^0,\epsilon)$ \citep{vanderLaan&Gruber15}. Both paths are constructed by locally tracking the one dimensional and multidimensional local least favorable path, respectively. A one-step TMLE using a universal least favorable path solves the canonical gradient equation $P_n D^{(1)}_{\tilde{P}^{(1)}(P_n^0,\epsilon_n^{(1)})}=0$ exactly. A TMLE update defined by the MLE for a local least favorable path will only solve this equation up till a second order term, so that, in that case, iteration of the TMLE update procedure might need to be employed for a few times till the equation is solved at a desired level. Thus, using a universal least favorable path in the TMLE makes the TMLE a one-step TMLE, thereby making the TMLE more robust in finite samples than an iterative TMLE.

Let $R^{(1)}(P,P_0)\equiv \Psi(P)-\Psi(P_0)-(P-P_0) D^{(1)}_P$, where $(P-P_0)D^{(1)}_P=-P_ 0 D^{(1)}_P$, be the so called exact second order remainder for target parameter $\Psi$. 
By definition of $R^{(1)}(P,P_0)$, a TMLE $P_n^*$ solving $P_n D^{(1)}_{P_n^*}=0$  allows an exact expansion of the form:
\[
\Psi(P_n^*)-\Psi(P_0)=(P_n-P_0)D^{(1)}_{P_n^*}+R^{(1)}(P_n^*,P_0).\]
 This remainder $R^{(1)}(P,P_0)$ is generally a second order difference in $p$ and $p_0$, and by computing it in a particular estimation problem, one can confirm the precise structure of this exact remainder. 

In previous work we proposed a general MLE of $P_0$ or relevant functionals thereof, which we named the Highly Adaptive Lasso (HAL), or HAL-MLE \citep{vanderLaan15,Benkeser&vanderLaan16,vanderLaan17}. 
The HAL-MLE of the data density $p_0$  involves first parametrizing $p=p_{\theta}$ by a multivariate real valued cadlag function $\theta(p): [0,\tau]\subset\openr^d \rightarrow\openr$, and computing $\theta_n^C=\arg\min_{\theta,\pl \theta\pl_v\leq C} P_n L(P_{\theta})$ the MLE over all cadlag functions in the parameter space with (sectional) variation norm $\pl \theta\pl_v\leq C$ bounded by a constant $C$. Due to our general representation of a multivariate real valued cadlag function $\theta=\int_{[0,\tau]}\phi_x(u) d\theta(u)$ as an infinite linear combination of tensor products of zero order spline basis functions $\phi_x(u)=I(x\geq u)$ at knot-points $u\in [0,\tau]$, and defining the variation norm $\pl \theta\pl_v\equiv
\int_{[0,\tau]} \mid d\theta(u)\mid$ \citep{Gill&vanderLaan&Wellner95,vanderLaan15,vanderLaan17},
it follows that one can compute this HAL-MLE by maximizing over a linear combination of a large set of spline basis functions $\phi_x(u_j)$  indexed by knot points $u_j$ under the constraint that the $L_1$-norm of the coefficient vector is bounded by $C$. As a consequence, HAL-MLEs can generally be implemented with available Lasso implementations such as glmnet$()$ in $R$.
HAL-MLE can also be separately computed for different functionals of $P_0$ such as univariate conditional densities, as long as these functional parameters can be identified as the minimizer of a risk function. 

We have shown that the HAL-MLE converges in loss based dissimilarity $d_0(\tilde{P}_n,P_0)\equiv P_0 L(\tilde{P}_n)-P_0L(P_0)$ at minimal at a rate $n^{-2/3}(\log n)^d$, even when the parameter space only assumes cadlag and finite variation norm \citep{Bibaut&vanderLaan19}. By using Cauchy-Schwarz inequality, the exact remainder $R^{(1)}(P,P_0)$ can generally be bounded in terms of $d_0(P,P_0)$, assuming that $p_0$ is uniformly bounded away from zero  on the relevant support for $\Psi$ (i.e., the so called positivity assumption). This implies that a TMLE using HAL-MLE $\tilde{P}_n$ as initial estimator $P_n^0$ will be asymptotically efficient: $R^{(1)}(\tilde{P}^{(1)}(\tilde{P}_n,\epsilon_n),P_0)=O_P(n^{-2/3}(\log n)^d)$. The only other condition that is needed  for asymptotic efficiency is that $\{D^{(1)}_P:P\in {\cal M}\}$ is a $P_0$-Donsker class, and that is generally implied by these HAL-models making the cadlag and bounded variation assumption: the class of multivariate real valued cadlag functions with  a universal bound on their variation norm is a uniform Donsker class with a nice entropy  \citep{vanderVaart&Wellner96,
Bibaut&vanderLaan19}.
One can optimally select the variation norm for the HAL-MLE with cross-validation, and, by including various HAL-MLEs (indexed by tuning parameters such as the maximal level of interaction of basis functions, or prescreening) in the library of a super-learner, the super-learner will achieve at minimal the same rate of convergence as the best choice among these HAL-MLEs
\citep{vanderLaan&Dudoit03,vanderVaart&Dudoit&vanderLaan06,vanderLaan&Dudoit&vanderVaart06},

Therefore, in great generality, such an HAL-MLE-based TMLE is asymptotically efficient estimator of pathwise differentiable target parameters for realistic statistical models \citep{vanderLaan17}. Wald type confidence intervals are computed as $\Psi(P_n^*)\pm 1.96\sigma_n/n^{1/2}$, where $\sigma_n^2$ is an estimator of the asymptotic variance $\sigma^2_0=P_0 \{D^{(1)}_{P_0}\}^2$. providing asymptotically valid confidence intervals.

However, the performance of the estimator and corresponding confidence intervals heavily rely on $R^{(1)}(P_n^*,P_0)$ being small, at least w.r.t. the leading term $(P_n-P_0)D^{(1)}_{P_n^*}$. Therefore, the curse of dimensionality still plays a fundamental role in finite sample inference, even though the HAL-TMLE was able to deal with the {\em asymptotic} curse of dimensionality by having a rate of convergence that is hardly affected by the dimension $d$.

The literature on higher order efficient estimation aims to address this by aiming to construct estimators  that allow for inference that only  assumes that a higher order remainder (higher than second order) is negligible  \citep{levit1975,ibragimov1981,pfanzagl1982,pfanzagl1985,bickel1982annals,robins2008imsnotes,
robins2009metrika,li2011statsprobletters,vandervaart2014statscience}.  The above references aim to achieve this with an higher order extension of the one-step estimator, and have focussed on the second order extension, which  involves adding a second order U-statistic of an approximate second order efficient influence function to an initial estimator of the target estimand. An excellent review of this  general higher-order extension of the one-step  estimator is provided in \cite{vandervaart2014statscience}.  
We developed an analogue second order TMLE \citep{diaz2016ijb,Carone2014techreport,Caroneetal17,vanderLaan&Rose18}, which uses the same approximate higher order influence functions to propose extra fluctuation parameters for the regular TMLE-update based on the local least favorable path. \cite{pfanzagl1982} highlighted the appeal of devising an higher order efficient estimator through updating in the model space, even decades before the actual development of TMLE, realizing that such an approach would be significantly more robust.

The current second order efficiency theory and its construction of corresponding second order estimators relies on the target parameter being second order pathwise differentiable with a second order canonical gradient, and, in general, higher order efficiency relies on higher order pathwise differentiability with corresponding higher order canonical gradients. Unfortunately, almost all first order pathwise differentiable target parameters of interest in realistic statistical models are not second order pathwise differentiable, let alone, higher order pathwise differentiable
(e.g, \citep{robins2008imsnotes}). Efforts to still utilize this theory to construct second order one-step estimator  or second order TMLE are therefore only resulting in limited practical and theoretical improvements. 

These methods essentially assume a sieve (e.g., family of parametric models) and for each sieve employ the second order efficient estimators based on the sieve-specific second order efficient influence function, even though this second order efficient influence function does not exist in the limit (i.e., are infinite). By necessity, these methods have to trade off  the increasing (as complexity of sieve grows) variance of the second order sieve based estimators as estimators of sieve projected target parameters with the bias due to the gap between the sieve and true statistical model. Even though these approaches based on the higher order efficiency theory can result in finite sample gains, they generally do not result in meaningful gains in rate of convergence for the exact remainder. In addition, they can easily become unstable. For example, for the nonparametric estimation of the treatment specific mean, the second order canonical gradient relies on inverse weighting by a marginal density of the baseline covariate vector.


In this article, we follow a different approach within the TMLE framework and establish a new higher order TMLE for (first order) pathwise differentiable target parameters that essentially achieves the goals of the traditional higher order efficiency theory, but without any need to assume the too stringent forms of higher order pathwise differentiability. Contrary to the traditional higher order efficiency theory, the expansion is not higher order in terms of U-statistics (which thereby jump from a remainder that is $O_P(n^{-1/2})$ in one step to a remainder that is $O_P(n^{-1})$ etc), but is higher order in terms of differences of the (initial) density estimator and true density. Our new $k$-th order TMLE can be straightforwardly defined, relying on pathwise differentiability of a sequentially and recursively defined set of higher order fluctuations of the target parameter, always fluctuating in direction of the canonical gradient.of the previously defined fluctuation of the target parameter.

 The key idea is to replace the initial estimator $P_n^0$ in the first order TMLE $\Psi(\tilde{P}^{(1)}(P_n^0,\epsilon_n))$ by a TMLE $P_n^*$ of $\Psi_n^{(1)}(P_0)\equiv \Psi(\tilde{P}^{(1)}(P_0,\epsilon_n(P_0))$, i.e. of the target parameter of $P_0$ defined as the one obtained by replacing the initial by its oracle choice $P_0$, where, for now we act as if $\Psi_n^{(1)}$ is pathwise differentiable, and we discuss this subtlety in the next paragraph.
In this manner, a second order TMLE of $\Psi(P_0)$ is just a first order TMLE that uses as initial estimator a TMLE $P_n^*=\tilde{P}_n^{(2)}(P_n^0)$ of $\Psi_n^{(1)}(P_0)$ that is fully tailored to its purpose/role in estimation of $\Psi(P_0)$. Since $P_n^*=\tilde{P}_n^{(2)}(P_n^0)$ is its own TMLE update of an initial estimator $P_n^0$, we could now  also define that initial estimator $P_n^0$ as a TMLE $P_n^*$ of  $\Psi_n^{(2)}(P_0)\equiv \Psi_n^{(1)}(\tilde{P}_n^{(2)}(P_0))=\Psi(\tilde{P}_n^{(1)}\tilde{P}_n^{(2)}(P_0))$, thereby tailor it for its purpose in the final TMLE $\Psi(\tilde{P}_n^{(1)}\tilde{P}_n^{(2)}(P_n^*))$ of $\Psi(P_0)$. Iterating this strategy results in the  definition of a $k$-th order TMLE for general integer $k=1,\ldots$, up till the modification discussed in next paragraph.

The second order TMLE relies on pathwise differentiability of $\Psi_n^{(1)}:{\cal M}\rightarrow\openr$ at $P$ that can occur as an update of the initial estimator $P_n^0$, and thereby the existence of its (first order) canonical gradient, which we term second order canonical gradient (even though it is not the canonical gradient of a second order pathwise derivative as in current higher order efficiency theory). It happens to be the case this parameter $\Psi_n^{(1)}(P)$  is smooth in $P$ up till the dependence of the MLE $\epsilon_n^{(1)}(P)$ on $P$. That is,  $P\rightarrow\epsilon_n^{(1)}(P)=\arg\max_{\epsilon}P_n \log \tilde{p}_n^{(1)}(p,\epsilon)$ is not pathwise differentiable at a $P=P_1$ due to $P_n$ not being absolutely continuous w.r.t. $P_1$: $dP_n/dP_1$ does not exist for most $P_1$.  Therefore, we replace the empirical mean in the log-likelihood by the expectation w.r.t. an HAL-MLE $\tilde{P}_n$. The fact that $d\tilde{P}_n/dP$ exists for all $P$ that can occur as an initial or a higher order TMLE-update of the initial estimator guarantees that $\Psi_n^{(1)}$ is pathwise differentiable at $P$ and thus has a canonical gradient $D^{(2)}_{n,P}$.
In our treatment specific mean example, we can set the marginal distribution of the covariates under $\tilde{P}_n$ equal to the empirical distribution, since the initial estimator $P_n^0$ uses as marginal covariate distribution the empirical  and its TMLE-updates are not affecting this marginal empirical distribution (since it is an NPMLE itself). Therefore, contrary to the second order estimators referenced above, our second order TMLE completely avoids estimation of a marginal density of a high dimensional covariate vector in this example. 

By undersmoothing this HAL-MLE $\tilde{P}_n$, this HAL-regularized MLE $\epsilon_{\tilde{P}_n}^{(1)}(P)$ (i.e., a plug-in estimator plugging in $\tilde{P}_n$) will still behave as the efficient regular MLE  $\epsilon_{P_n}^{(1)}(P)$ \citep{vanderLaan19efficient,vanderLaan&Benkeser&Cai19}. (Our results show that this is only needed at $P=P_0$. )
 From an intuitive point of view, this is due to $\tilde{P}_n$ being a nonparametric MLE itself, thereby allowing the MLE to increase the likelihood in the direction of the empirical measure as far as needed to make sure that the score $\tilde{P}_n D^{(1)}_{\tilde{P}^{(1)}(P,\epsilon_n^{(1)}(P))}$ of $\epsilon_n^{(1)}(P)$ approximates the score $P_nD^{(1)}_{\tilde{P}^{(1)}(P,\epsilon_{P_n}^{(1)}(P))}$ of the regular $\epsilon_{P_n}^{(1)}(P)$. Formally, as shown in \citep{vanderLaan19efficient,vanderLaan&Benkeser&Cai19}, this is true since an HAL-MLE solves uniformly a large class of empirical score equations, so that it solves also the linear span of these score equations that best approximates a desired score, and this approximation is tuned by the increasing the $L_1$-norm of the HAL-MLE.

The basic formal idea behind our analysis of this higher order HAL-regularized TMLE is the following. A TMLE $P_n^*$ of a target parameter such as $\Psi_n^{(1)}(P_0)$ solving $\tilde{P}_n D^{(2)}_{n,P_n^*}=0$, using the above mentioned HAL-MLE-regularization in the MLE-update, sets the directional derivative of $P\rightarrow \Psi_n^{(1)}(P)$ at $P_n^*$ in the direction $P_n^*-P_0$, equal to $(P_n^*-P_0)D^{(2)}_{P_n^*}=(\tilde{P}_n-P_0)D^{(2)}_{P_n^*}$, which thus equals  $o_P(n^{-1/2})$ noise.  
 In addition, since
 \[
 \Psi_n^{(1)}(P)-\Psi_n^{(1)}(P_0)=\tilde{P}_nD^{(1)}_{P_0}+ (\tilde{P}_n-P_0)\{D^{(1)}_{\tilde{P}_n^{(1)}(P)}-D^{(1)}_{P_0}\}+R^{(1)}(\tilde{P}_n^{(1)}(P),P_0),\]
 the directional derivative of $\Psi_n^{(1)}(P)$ at $P_n^*$ equals  the directional derivative at $P_n^*$ of the  exact total remainder  $P\rightarrow \bar{R}^{(1)}(\tilde{P}_n^{(1)}(P),P_0)$ for that target parameter defined by
\[
\bar{R}^{(1)}(\tilde{P}_n^{(1)}(P),P_0)=(\tilde{P}_n-P_0)\{D^{(1)}_{\tilde{P}_n^{(1)}(P)}-D^{(1)}_{P_0}\}+R^{(1)}(\tilde{P}_n^{(1)}(P),P_0).\]
Therefore, by using a TMLE $P_n^*$ of $\Psi_n^{(1)}(P_0)$ we are  locally optimizing the exact total remainder $\bar{R}^{(1)}(\tilde{P}_n^{(1)}(P),P_0)$ in $P$ up till noise, which we formally show to actually be true, and enforce that the relevant derivative at $P_n^*$ in the direction of $(P_n^*-P_0)$ is $o_P(n^{-1/2})$. (By replacing $P_0$ by $\tilde{P}_n$ in this argument, the directional derivative is exactly equal to zero, allowing us to obtain particularly nice exact expansions for the higher order TMLE, as we will show. ) The latter will form the fundamental ingredient to show that the targeted exact remainder $R^{(1)}(\tilde{P}_n^{(1)}(P_n^*),P_0)$ behaves as an empirical mean of the second order efficient influence function $D^{(2)}_{P_0})$ and a third order difference. By iteration, this will allow us to establish that  $k$-th order HAL-regularized TMLE behaves as a sum of empirical means of mean zero higher order canonical gradients and a $k+1$-th order difference. 

Our exact expansion for this $k$-th order HAL-regularized TMLE also includes a HAL-regularization bias term including as leading term $(\tilde{P}_n-P_n)D^{(1)}_{\tilde{P}_n^{(1)}(P_0)}$, which is fully controlled by undersmoothing HAL. However, we will show that by using empirical TMLE-updates in the $k$-th order TMLE, still using the same least favorable paths with canonical gradients of the HAL-regularized fluctuations of the target parameter, thereby depending on the HAL-MLE $\tilde{P}_n$, this bias term reduces to a negligible term theoretically. Indeed our practical studies confirm that undersmoothing is not even needed anymore for this empirical $k$-th order TMLE. So our final recommendation is to use the $k$-th order empirical TMLE. Nonetheless, the theoretical understanding (including our exact expansion) of the latter $k$-th order  empirical TMLE requires first understanding the HAL-regularized $k$-th order TMLE. 

\subsection{Organization of article}
In Section \ref{section2} we define the HAL-regularized $k$-th order TMLE, and some important versions of its implementations that guarantee that it solves the HAL-regularized higher order efficient influence equations $\tilde{P}_n D^{(j)}_{n,\tilde{P}_n^{(j)}(P)}=0$ at user supplied precision.
In Section \ref{section3} we show that the second  order TMLE step $P_n^*=\tilde{P}_n^{(2)}(P)$ minimizes the total remainder $P\rightarrow \bar{R}^{(1)}(\tilde{P}_n^{(1)}(P),P_0)$  for the first order TMLE under the constraint that one improves log-likelihood, thereby providing the finite sample rational of the second order TMLE. Since the same applies to the higher order TMLE-updates, but now w.r.t. the total remainder of the parameter it is targeting, this provides the general finite sample rational  for the HAL-regularized $k$-th order TMLE (and thereby also for the empirical $k$-th order TMLE. In Section \ref{section4} we provide  exact $k+1$-th order expansions for the HAL-regularized $k$-th order TMLE in terms of empirical means of higher order canonical gradients, an HAL-regularization bias term, and a final exact remainder $R_n^{(k)}(\tilde{P}_n^{(k)}(P_n^0),P_0)$, which is further expressed in terms of the exact remainder $R_n^{(k)}(\tilde{P}_n^{(k)}(P),\tilde{P}_n)$ at the HAL-MLE $\tilde{P}_n$. It is shown that the latter equals an iteratively targeted $R_n^{(1)}(\tilde{P}_n^{(1)}\ldots \tilde{P}_n^{(k)}(P),\tilde{P}_n)$, where each $\tilde{P}_n^{(j)}(P)$ actually targets this remainder, setting the stage for establishing it is a $k+1$-th order difference.
In Section \ref{sectionemphotmle} we provide the resulting exact expansion for the empirical $k$-th order TMLE and demonstrate that it essentially removes the HAL-regularization bias term for the HAL-regularized $k$-th order TMLE, while still obtaining the same $k$-th order exact expansion. 
The next two Sections \ref{section5} and \ref{section6} focus on establishing that the $k$-th order remainder at $\tilde{P}_n$ is indeed a $k+1$-th order difference.
Firstly, in  Section \ref{section5} we formally define generalized higher order differences $R_n(P,P_0)$ and provide the fundamental understanding of how using a TMLE $P_n^*$ of $P\rightarrow R_n(P,\tilde{P}_n)$  increases the order of difference of $R_n(P_n^*,\tilde{P}_n)$ by setting its directional derivative at $P_n^*$ in the direction $P_n^*-\tilde{P}_n$ equal to zero.
Sequential application of this result allows us in  Section \ref{section6} to establish that this exact remainder $R_n^{(k)}(\tilde{P}_n^{(k)}(P),\tilde{P}_n)$  is a generalized  $k+1$-th order difference of $P$ and $\tilde{P}_n$, up till the precision at which the sequentially defined TMLE solves the  higher order HAL-regularized canonical gradient score equations. 

In Section \ref{section8} 
we demonstrate how our exact expansion can be used to obtain $k$-th order confidence intervals based on the  $k$-th order TMLE,  which take into account the exact expansion for the $k$-th order TMLE up till the $k+1$-th order remainder. 
In Section \ref{section9} we provide formal 
algebra for sequentially  analytically computing the higher order canonical gradients.
In Section \ref{section10} we   study the second order TMLE for nonparametric estimation of the integral of the square of the density. 
In Section \ref{section11} we study the second order TMLE for nonparametric estimation of the treatment specific mean (i.e., average causal effect of treatment on binary outcome).
In Section \ref{sec:ave_dens} we provide simulation results for the second order TMLE for nonparametric estimation of the integrated square density. 
 In Section \ref{section11b} we show simulation results for the second order TMLE for nonparametric estimation of the treatment specific mean. 
We conclude with a discussion in Section \ref{section12}. 

\subsection{Appendix}
Various proofs of lemmas and theorems in the main article are deferred to the Appendix. 
In addition, we present a number of basic investigations and extensions  in the Appendix. 
Specifically, in Appendix \ref{section7} we show that by targeted undersmoothing the HAL-MLE $\tilde{P}_n$ in the HAL-regularized $k$-th order TMLE  we will not only have that the HAL-MLE regularized score equations $\tilde{P}_n D^{(j)}_{n,\tilde{P}_n^{(j)}(P)}=0$ are solved, but  that also the empirical mean $P_n D^{(j)}_{n,\tilde{P}_n^{(j)}(P_0)}\approx 0$ up till the desired precision.
In Appendix \ref{sectionfirstgradient} we demonstrate how one can compute a canonical gradient of a pathwise differentiable target parameter in terms of a linear combination of HAL-basis functions of the tangent space, defined by applying the inverse of the symmetric covariance matrix of the vector of basis functions to a  numerical pathwise derivative along the path with score the basis function or covariance of the basis function with an initial  gradient. The inverse can be determined with the Choleski decomposition, or one can approximate it with  $L_1$-penalized linear least squares regression.
 This insight allows one to compute canonical gradients without relying on analytics and Hilbert space theory, and without relying on having a closed form representation of the canonical gradient. 
In Appendix \ref{section9b} we apply this result to compute the higher order canonical gradients as well (which are just canonical gradients of $\Psi_n^{(k)}(P)$). Computing the numerical pathwise derivative now involves computing the score of the least favorable path under a fluctuation of the initial along a  path with score the basis function, across all basis function. We then show how to recursively program the $K$-th  order TMLE, by describing how to compute the $k+1$-th order TMLE at initial $P$ in terms of applying the $k$-th order TMLE and some extra calculations.

In Appendix \ref{sectionthotmle} we make the observation that our analysis of the HAL-regularized  higher order TMLE applies to any $\tilde{P}_n$ and $P_n^0$, thereby allowing us to target these estimators towards the goals of making the HAL-regularized higher order TMLE behave exactly as using a regular TMLE update for each $\tilde{P}_n^{(j)}(P)$ (i.e. maximizing empirical likelihood instead of HAL-regularized likelihood). We present an iterative HAL-regularized higher order TMLE that also targets the HAL-MLE and explain how it is motivated by making  the HAL-MLE regularized $\epsilon_{\tilde{P}_n}^{(j)}(P_0)$ behave as the regular empirical MLE $\epsilon_{P_n}^{(j)}(P_0)$. In Appendix  \ref{sectioncvhotmle} we propose a variation of the iterative higher order TMLE by targeting the HAL-MLE trained on training sample based on validation sample, thereby obtaining a cross-validated/cross-fitted (iterative) higher order TMLE. We show that such a cross-fitted higher order TMLE  allows for asymptotic normality under significantly weaker conditions than a regular higher order TMLE, analogue to the CV-TMLE being more robust than the regular TMLE.

\section{Definition of the  HAL-regularized $k$-th order TMLE}\label{section2}





\subsection{The HAL-MLE}
Let $\tilde{P}_n\in {\cal M}$ be an HAL-MLE of $P_0$, where we are reminded that an HAL-MLE is an MLE over the parameter space implied by the model under the restriction that the variation norm is bounded by a constant, which itself will be data adaptively selected.
To emphasize its dependence on a variation norm bound $C$, we will also use the notation $\tilde{P}_n^{C}$ to denote the $C$-specific HAL-MLE.  Let $C_{n,cv}$ denote the cross-validation selector, while $C_n$ denotes the selector used for $\tilde{P}_n=\tilde{P}_n^{C_n}$.
$\tilde{P}_n$ will replace $P_n$ in the MLE-update steps along least favorable paths. $P_n$ is known to be an NPMLE of $P_0$, so an HAL-MLE $\tilde{P}_n$ is a natural smooth model-based  analogue of $P_n$, especially, considering that  an HAL-MLE, even without undersmoothing, uniformly solves the scores of all its non-zero coefficients at rate $O_P(n^{-2/3})$ \citep{vanderLaan19efficient,vanderLaan&Benkeser&Cai19}. Implementations are available for HAL-MLE of conditional probabilities (logistic regression), conditional densities or hazards, and conditional means, among others (e.g., HAL9001 in R, based on glmnet). Our results for the HAL-regularized  higher order TMLE shows that  $\tilde{P}_n$ only needs to  imitate the NPMLE $P_n$  in the sense that $\epsilon_{\tilde{P}_n}^{(j)}(P_0)\approx \epsilon_{P_n}^{(j)}(P_0)$, specifically for $j=1$.

\subsection{Sequential definition of $k$ target parameters: Defining a next target parameter of data distribution as the previous target parameter applied to  the TMLE-update of the data distribution}

{\bf Defining first order TMLE targeting $\Psi(P)$:}
Let $\{\tilde{P}^{(1)}(P,\epsilon):\epsilon\}\subset {\cal M}$ be a local least favorable parametric submodel through $P$ with the property that the one-step TMLE update $
\tilde{P}_n^{(1)}(P)=\tilde{P}^{(1)}(P,\epsilon_{n}^{(1)}(P))$ using as initial ''estimator'' or off-set $P$ satisfies $\tilde{P}_n D^{(1)}_{\tilde{P}_n^{(1)}(P)}=o_P(n^{-1/2})$ at user supplied precision so that its value is negligible for inference.  This should hold at $P=P_0$, and, at TMLE updates $P=\tilde{P}_n^{(2)}\ldots\tilde{P}_n^{(k)}(P_n^0)$ of $P_n^0$ as defined below. 
At minimal, this requirement relies on  the submodel $\tilde{P}^{(1)}(P,\epsilon)$ having a score at $\epsilon=0$ that spans the canonical gradient $D^{(1)}_{P}$.
If the submodel is a universal least favorable submodel, then one can define $\epsilon_{n}^{(1)}(P)=\arg\min_{\epsilon}\tilde{P}_n L(\tilde{P}^{(1)}(P,\epsilon))$ as the MLE of the smoothed log-likelihood, which would imply 
$\tilde{P}_n D^{(1)}_{\tilde{P}_n^{(1)}(P)}=0$ exactly. If the submodel is only a local least favorable path, then defining $\epsilon_n^{(1)}(P)$ as the MLE as above still suffices at $P=P_0$, but to have the desired behavior at $P$ that can occur as $P_n^0$ or its TMLE-update, one might have to make $P_n^0$ itself a targeted estimator (as in our iterative HAL-regularized higher order TMLE defined below).

An important option is to simply define $\epsilon_n^{(1)}(P)$ as the solution of $\tilde{P}_n D_{\tilde{P}^{(1)}(P,\epsilon)}=0$, in which case we still guarantee zero contributions from these TMLE-score values, while it still closely resembles the MLE update.

 We recommend to use a local least favorable submodel with a univariate $\epsilon$, either using the MLE $\epsilon_n^{(1)}(P)$ or defining it as solution of $\tilde{P}_n D^{(1)}_{\tilde{P}_n^{(1)}(P,\epsilon)}=0$. This simplifies the calculations of  the subsequent higher order canonical gradients $D^{(j)}_{n,P}$ below.  We will use  the notation $\tilde{P}_n^{(1)}(P)$ for this general possible choice satisfying $\tilde{P}_n D^{(1)}_{\tilde{P}_n^{(1)}(P_0)}=o_P(n^{-1/2})$. 

{\bf Defining the first order fluctuation of target parameter, $\Psi_n^{(1)}(P)$:}
We now define a new target parameter $\Psi_n^{(1)}:{\cal M}\rightarrow \openr$ by $\Psi_n^{(1)}(P)=\Psi(\tilde{P}_n^{(1)}(P))$. Due to the regularization of $\epsilon_n^{(1)}(P)$ with $\tilde{P}_n$ this target parameter will be pathwise differentiable as well. We note that $P\rightarrow\Psi(\tilde{P}^{(1)}(P,\epsilon))$ is generally already pathwise differentiable at $P$ for any $\epsilon$, so that the only reason for $P\rightarrow \Psi(\tilde{P}^{(1)}_n(P,\epsilon_{P_n}^{(1)}(P))$ not being pathwise differentiable is due to $\epsilon_{P_n}(P)$ not being pathwise differentiable in $P$ due to using $P_n$. However, $\epsilon_n^{(1)}(P)=\epsilon_{\tilde{P}_n}^{(1)}(P)$ will be pathwise differentiable at $P$ as long as $d\tilde{P}_n/dP$ exists: see Appendix \ref{AppendixA}.
 
The parameter $\Psi_n^{(1)}(P)$ is a data dependent parameter due to its dependence on $\epsilon_n^{(1)}(P)$, but, beyond this random coefficient, it is a fixed target parameter of $P$. In fact, suppose we replace $\epsilon_n^{(1)}(P)$ by $\epsilon_0^{(1)}(P)=\arg\min_{\epsilon}P_0 L(\tilde{P}^{(1)}(P,\epsilon))$. Then the analogue target parameter $\Psi_0^{(1)}(P)=\Psi(\tilde{P}^{(1)}_0(P))$ with $\tilde{P}^{(1)}_0(P)=\tilde{P}^{(1)}(P,\epsilon_0^{(1)}(P))$ is now a fixed parameter and is pathwise differentiable. Note that we could also use the notation $\Psi_{\tilde{P}_n}^{(1)}(P)$ to emphasize that its dependence on data is through $\tilde{P}_n$.

{\bf Defining the second order fluctuation of target parameter,  $\Psi_n^{(2)}(P)$:}
Let $D^{(2)}_{n,P}$ be the canonical gradient of the pathwise derivative of $\Psi_n^{(1)}$ at $P$ and let $R_n^{(2)}(P,P_0)\equiv \Psi_n^{(1)}(P)-\Psi_n^{(1)}(P_0)+P_0D^{(2)}_{n,P}$ be the exact remainder. We note that $D^{(2)}_{n,P}$ is a random canonical gradient through $\tilde{P}_n$, and it represents a plug-in estimator of the ''oracle'' canonical gradient $D^{(2)}_{0,P}$ of $\Psi_0^{(1)}$ obtained by replacing $\epsilon_0^{(1)}(P)$ by $\epsilon_n^{(1)}(P)$.
Let $\{\tilde{P}_n^{(2)}(P,\epsilon):\epsilon\}$ be a least favorable submodel through off-set $P$ targeting $\Psi_n^{(1)}(P_0)$, and let $\tilde{P}_n^{(2)}(P)=\tilde{P}_n^{(2)}(P,\epsilon_n^{(2)}(P))$ be the TMLE update using again the regularized $\epsilon_n^{(2)}(P)=\arg\min_{\epsilon}\tilde{P}_n L(\tilde{P}_n^{(2)}(P,\epsilon))$ or solution of $\tilde{P}_n D^{(2)}_{\tilde{P}_n^{(2)}(P,\epsilon)}=0$. The local least favorable submodel is now data dependent due to $D^{(2)}_{n,P}$ being data dependent.
We can now define ${\Psi}_n^{(2)}:{\cal M}\rightarrow \openr$ as $\Psi_n^{(2)}(P)=\Psi_n^{(1)}(\tilde{P}_n^{(2)}(P))=\Psi(\tilde{P}_n^{(1)}\tilde{P}_n^{(2)}(P))$.

{\bf Iterating this process up till the $k$-th order  fluctuation of target parameter, $\Psi_n^{(k-1)}(P)$:}
In general,  sequentially, for $j=1,\ldots,k-1$, given $\Psi_n^{(j-1)}$ and its TMLE update $\tilde{P}_n^{(j)}(P)$, we 1) define $\Psi_n^{(j)}:{\cal M}\rightarrow\openr$ by $\Psi_n^{(j)}(P)=\Psi_n^{(j-1)}(\tilde{P}_n^{(j)}(P))$; 2) determine its canonical gradient $D^{(j+1)}_{n,P}$ at $P$; 3) the TMLE update $\tilde{P}_n^{(j+1)}(P)$ using a regularized $\epsilon_n^{(j)}(P)$ tailored to solve $\tilde{P}_n D^{(j+1)}_{n,\tilde{P}_n^{(j+1)}(P)}\approx 0$, and its exact remainder $R_n^{(j+1)}(P,P_0)=\Psi_n^{(j)}(P)-\Psi_n^{(j)}(P_0)+P_0 D^{(j+1)}_{n,P}$.


To summarize, for the $k=1$-th order TMLE, we just define $\Psi$; $D^{(1)}_P$; $\tilde{P}_n^{(1)}(P)$; $R^{(1)}(P,P_0)$. For the $k=2$-th order TMLE, we also define  $\Psi_n^{(1)}$; $D^{(2)}_{n,P}$; $\tilde{P}_n^{(2)}(P)$; $R_n^{(2)}(P,P_0)$. And so on.
We also note that $\Psi_n^{(j)}$, $D^{(j)}_{n,P}$ and $R_n^{(j)}(P,P_0)$ are random through $\tilde{P}_n$, and specifically, through its dependence on $\epsilon_n^{(1)}$,\ldots,$\epsilon_n^{(j)}$ in the definition of the TMLE updates.

\subsection{Sequentially defining the HAL-regularized TMLE-updates, using the previous TMLE as initial for the next TMLE, starting at the $k$-th TMLE.}
We are now ready to define the HAL-regularized $k$-th order TMLE.
Let $P_n^0$ be an initial estimator of $P_0$. One option is to select 
$P_n^0=\tilde{P}_n^{C_{n,cv}}$ as the HAL-MLE with the variation norm selected with the cross-validation selector (while $\tilde{P}_n$ might use an undersmoothed $C_n$). Alternatively, $P_n^0$ is a super-learner based on a library including $\tilde{P}_n^{C_{n,cv}}$ as a candidate algorithm.
Let $P_n^{k,(k),*}=\tilde{P}^{(k)}_n(P_n^0)$ be the  HAL-regularized TMLE targeting $\Psi_n^{(k-1)}(P_0)$ using as initial estimator $P_n^0$.
We then define ${P}_n^{k,(k-1),*}=\tilde{P}_n^{(k-1)}(P_n^{k,(k),*})$ as the regularized-TMLE targeting
$\Psi_n^{(k-2)}(P_0)$ using as initial estimator the previous TMLE $P_n^{k,(k),*}$.
Sequentially, we define $P_n^{k,(j),*}=\tilde{P}_n^{(j)}(P_n^{k,(j+1),*})$, $j=k-1,\ldots,1$, till the final $P_n^{k,(1),*}=\tilde{P}_n^{(1)}(P_n^{k,(2),*})$.
This final TMLE $P_n^{k,(1),*}$ is the $k$-th order TMLE of $P_0$ targeting $\Psi(P_0)$,  and $\Psi(P_n^{k,(1),*})$ is the resulting HAL-regularized $k$-th order plug-in TMLE of $\Psi(P_0)$.

\subsection{$j+1$-th order TMLE is obtained by replacing the initial estimator in the $j$-th order TMLE by the TMLE of $\Psi_n^{(j)}(P_0)$.}
Note that with $P_n^*=\tilde{P}_n^{(2)}(P_n^0)$, we have that the HAL-regularized second order TMLE is given by $\Psi_n^{(1)}(P_n^*)=\Psi(\tilde{P}_n^{(1)}(P_n^*))=\Psi(\tilde{P}_n^{(1)}\tilde{P}_n^{(2)}(P_n^0))$. In general, the HAL-regularized $j+1$-th order TMLE is given by the TMLE $\Psi_n^{(j)}(P_n^*)=\Psi(\tilde{P}_n^{(1)}\ldots\tilde{P}_n^{(j+1)}(P_n^0))$ of $\Psi_n^{(j)}(P_0)$ using $P_n^*=\tilde{P}_n^{(j+1)}(P_n^0)$. So, one can obtain the HAL-regularized $k$-th order TMLE by sequentially determining the HAL-regularized $j$-th order TMLE, $j=1,\ldots,k$,  each time replacing the initial estimator $P_n^0$ in the $j$-th order TMLE by the TMLE of $\Psi_n^{(j)}(P_0)$: 1) First order TMLE is the TMLE $\Psi(\tilde{P}_n^{(1)}(P_n^0))$ of $\Psi(P_0)$; 2) Second order TMLE replaces $P_n^0$ by $\tilde{P}_n^{(2)}(P_n^0)$ and therefore is the TMLE $\Psi_n^{(1)}(\tilde{P}_n^{(2)}(P_n^0))$ of $\Psi_n^{(1)}(P_0)$; and so on. One can iterate this till one ends up at the desired order $k$ of the HAL-regularized  TMLE.


\subsection{Using multivariate least favorable path separately targeting different components of $P$}
Suppose that $P=(Q_m(P):m=1,\ldots,M)$ for variation independent parameters $Q_m(P)$. For example, $Q_m(P)$ could be a conditional density of $X_m$, given $X_1,\ldots,X_{m-1}$, while $O=(X_1,\ldots,X_M)$. Then, $D^{(j)}_{n,P}=\sum_m D^{(j)}_{n,P,m}$ has a corresponding sum decomposition so that $D^{(j)}_{n,P,m}$ equals the canonical gradient of $\Psi_n^{(j-1)}(P)$ when only viewed as function of $Q_m(P)$ while fixing $Q_{l}(P)$ for $l\not =m$.
In the TMLE-updates $\tilde{P}_n^{(j)}(P)$ one can target each equation $\tilde{P}_nD^{(j)}_{n,P,m}\approx 0$ so that all of them are approximately equal to $0$, which requires using a multidimensional $\epsilon$ in the local least favorable submodel or the corresponding universal path that tracks this multivariate local least favorable submodel locally and iteratively. 
In our treatment specific mean example the definition of $\tilde{P}_n^{(2)}(P)$ involves separately targeting the treatment mechanism and outcome regression based on local or universal least favorable paths through these nuisance parameters.

\subsection{Iterative  $k$-th order TMLE  based on local least favorable paths to guarantee that all regularized higher order efficient score equations are solved exactly.}
Consider the case that one uses MLE-updates with local least favorable paths so that $\tilde{P}_n D^{(j)}_{n,\tilde{P}_n^{(j)}(P)}$ will not be exactly solved at the relevant $P$, $j=k,\ldots,1$. 
Let $P_n^{k,(1),*}(P)$ be the $k$-th order TMLE that uses as initial estimator $P$.
Let $P_n^0$ be the initial estimator we start out with. Let
$P_n^{k,(1),1}=P_n^{k,(1),*}(P_n^0)$ be the first $k$-th order TMLE. We then compute the second $k$-th order TMLE
$P_n^{k,(1),2}=P_n^{k,(1),*}(P_n^{k,(1),1})$ that uses the previous $k$-th order TMLE $P_n^{k,(1),1}$ as initial estimator. We can iterate this as $P_n^{k,(1),l}=P_n^{k,(1),*}(P_n^{k,(1),l-1})$, $l=1,2,\ldots$. Note that the HAL-regularized log-likelihood $\tilde{P}_n L(P_n^{k,(1),l})$ is increasing in $l$, due to each TMLE update $\tilde{P}_n^{(j)}(P)$ increasing the regularized log-likelihood relative to its initial $P$, $j=k,\ldots,1$. Therefore, this $k$-th order iterative TMLE $P_n^{k,(1),l}$ will converge as $l$ increases at which point the log-likelihood plateaus. The user can iterate this till step $l$ at which  $\tilde{P}_n D^{(j)}_{P_n^{k,(j),*}}\approx 0$ is solved at the desired precision for all $j$, where $P_n^{k,(j),*}=\tilde{P}_n^{(j)}\ldots\tilde{P}_n^{(k)}(P_n^{k,(1),l})$ is the final $j$-th order TMLE-update at that $l$-th step.
This algorithm is just the analogue of the  iterative first order TMLE, but now applied to   the $k$-th order TMLE.

In Appendix \ref{sectionthotmle} we augment this iterative  HAL-regularized higher order TMLE with simultaneous iterative targeting of the HAL-MLE $\tilde{P}_n$, so that in the limit it is also guaranteed that $P_n D^{(j)}_{P_n^{k,(j),*}}=0$ for all $j=1,\ldots,k$.

\subsection{A  $k$-th order TMLE that sequentially optimizes along the universal least favorable paths till log-likelihood is maximized}
We can also carry out a universal least favorable path analogue of the first order TMLE.
Namely, replace in the definition of the $k$-th order TMLE $\tilde{P}_n^{(j)}(P,\epsilon_n^{(j)}(P))$ by $\tilde{P}_n^{(j)}(P,\delta)$ for a small $\delta$ chosen in the direction in which the log-likelihood $\tilde{P}_n L(\tilde{P}_n^{(j)}(P,\epsilon))$ increases, $j=k,\ldots,1$. This results in a $\delta$-restricted $k$-th order TMLE $P_n^{k,1,\delta}(P_n^0)$, in which each TMLE-update only moved by an amount $\delta$. We can then make $P_n^{k,1,\delta}(P_n^0)$ be the new initial estimator, and thereby define 
$P_n^{k,1,2\delta}(P_n^0)$ as the $\delta$-restricted $k$-th order TMLE that uses $P_n^{k,1,\delta}(P_n^0)$ as initial estimator and only uses small $\delta$-updates in all TMLE update steps. As before, the $\tilde{P}_n$-log-likelihood increases at each step so that we can iterate this till $\tilde{P}_n L(P_n^{k,1,l\delta})$ is maximized. If at a particular step $l$, one of the $\tilde{P}_n^{(j)}$ updates are not able to increase its log-likelihood anymore along a small $\delta$-step, then one just sets $\delta=0$ for that particular update (but it might kick back in at the next step $l+1$).
Again, as mentioned above, this can be further extended with a universal path for $\tilde{P}_n$ so that in the end $\tilde{P}_n$ is targeted towards solving the empirical higher order efficient score equations as well (Section \ref{sectionthotmle}).

We conclude that, we can use a  local least favorable path $\tilde{P}_n^{(j)}(P,\epsilon)$ for all $j$, thereby simplifying the calculation of the higher order canonical gradients $D^{(j)}_{n,P}$, and still achieve $\tilde{P}_n D^{(j)}_{P_n^{k,(j),*}}=0$ exactly or up till user supplied precision.

Of course, neither one of these two iterative  HAL-regularized $k$-th order TMLEs described above are needed if we define $\tilde{P}_n^{(j)}(P)=\tilde{P}_n^{(j)}(P,\epsilon_n^{(j)}(P))$ with $\epsilon_n^{(j)}(P)$ be the solution of the score equation $\tilde{P}_n D^{(j)}_{n,\tilde{P}_n^{(j)}(P,\epsilon)}=0$. In fact, the latter definition of $\tilde{P}_n^{(j)}(P)$ has advantages by guaranteeing an exact $k$-th order remainder.

 \section{The HAL-regularized second  order TMLE  minimizes the total remainder w.r.t. choice of initial for first order TMLE under constraint that one improves log-likelihood}\label{section3}
 This section focusses on the second order TMLE, but it generalizes to the statement that the $j$-th order TMLE $\tilde{P}_n^{(j)}(P)$ minimizes the exact total remainder $P\rightarrow \bar{R}_{n}^{(j-1)}(\tilde{P}_n^{(j-1)}(P),P_0)$
 for  the TMLE expansion \[
 \Psi_n^{(j-2)}(\tilde{P}_n^{(j-1)}(P))-\Psi_n^{(j-2)}(\tilde{P}_n^{(j-1)}(P_0))=\tilde{P}_n D^{(j-1)}_{n,P_0}+\bar{R}_n^{(j-1)}(\tilde{P}_n^{(j-1)}(P),P_0),\] under the constraint that one improves the log-likelihood, $j=2,3,\ldots$. Since the exact decomposition of the $k$-th order TMLE presented in the next section is decomposed into a sum of  $ \Psi_n^{(j-2)}(\tilde{P}_n^{(j-1)}(P))-\Psi_n^{(j-2)}(\tilde{P}_n^{(j-1)}(P_0))$ at $P=\tilde{P}_n^{(j)}\ldots\tilde{P}_n^{(k)}(P_n^0)$, this would then show that our $k$-th order TMLE uses a  $P$ in this contribution
$ \Psi_n^{(j-2)}(\tilde{P}_n^{(j-1)}(P))-\Psi_n^{(j-2)}(\tilde{P}_n^{(j-1)}(P_0))$ that optimizes this difference, for each $j$.

In particular, this provides the rational of using a TMLE $\tilde{P}_n^{(2)}(P_n^0)$ of $\Psi_n^{(1)}(P_0)$ as initial estimator in the first order TMLE, resulting in the definition of the second order TMLE. 

Consider  the first order TMLE $\Psi(P_n^{1,*})$ with $P_n^{1,*}=\tilde{P}_n^{(1)}(P^0)$ for an initial estimator $P^0$, satisfying $\tilde{P}_n D^{(1)}_{P_n^{1,*}}=0$. This TMLE satisfies
 \[
\Psi(\tilde{P}_n^{(1)}(P^0))-\Psi(P_0)=\tilde{P}_n D^{(1)}_{P_0}+\bar{R}^{(1)}(\tilde{P}_n^{(1)}(P^0),P_0),\]
where
\begin{equation}\label{3a}
\bar{R}^{(1)}(\tilde{P}_n^{(1)}(P^0),P_0)=(\tilde{P}_n-P_0)\{D^{(1)}_{\tilde{P}_n^{(1)}(P^0)}-D^{(1)}_{P_0}\}+R^{(1)}(\tilde{P}_n^{(1)}(P^0),P_0).\end{equation}

Our goal is to select $P^0$ so that it makes this exact total remainder $\bar{R}^{(1)}(\tilde{P}_n^{(1)}(P^0),P_0)$ as small as possible.  For that purpose we view $P\rightarrow \bar{R}^{(1)}(\tilde{P}_n^{(1)}(P),P_0)$ as a function in the choice of initial $P$.
We can view computing the TMLE $P_n^{2,*}$ of $\Psi_n^{(1)}(P_0)$ as running a steepest descent algorithm for minimizing the exact total remainder $\bar{R}^{(1)}(\tilde{P}_n^{(1)}(P),P_0)$ in $P$, with the only twist that it stops the moment that it cannot move in the direction of $P_0$ anymore, as measured by not being able to increase the likelihood anymore.
We will show that now.

\subsection{The TMLE $P_n^{(2),*}$ follows the same path as the oracle steepest descent algorithm  minimizing the exact total remainder.}
Suppose our goal is to minimize $P\rightarrow \{\bar{R}^{(1)}(\tilde{P}_n^{(1)}(P),P_0)\}^2$.
That is,  by (\ref{3a}), we want to  minimize
\[
P\rightarrow f_{n,0}(P)\equiv \left\{\Psi(\tilde{P}_n^{(1)}(P))-\Psi(P_0)-\tilde{P}_nD^{(1)}_{P_0}\right\}^2.\]
Let's imagine that $f_{n,0}(P)$ is a known mapping in $P$.
To minimize this function, one starts with an initial $P_n^0$, construct a path $\tilde{P}^{(2)}(P_n^0,\delta)$ with score at $\delta=0$ the canonical gradient of $f_{n,0}(P)$ at $P=P_n^0$, and use this local path iteratively by always moving with small amounts $\delta$ and in the direction of the canonical gradient at the current $P$. Such an algorithm can be called a steepest gradient descent algorithm.
One would run this steepest gradient descent algorithm till the canonical gradient of this criterion $f_{n,0}(P)$ is small enough. 

So we would  compute the canonical gradient of the pathwise derivative of $f_{n,0}$ at $P$.
The pathwise derivative $\frac{d}{d\delta_0}f_{n,0}(P_{\delta_0,h})$ is given by
\[
2\left(\Psi(\tilde{P}_n^{(1)}(P))-\Psi(P_0)-\tilde{P}_n D^{(1)}_{P_0}\right)\frac{d}{d\delta_0}\Psi_n^{(1)}(\tilde{P}_n^{(1)}(P_{\delta_0,h})).\]
Recall that $D^{(2)}_{n,P}$ is the canonical gradient of  the pathwise derivative of $\Psi_n^{(2)}(P)=\Psi_n^{(1)}(\tilde{P}_n^{(1)}(P))$ at $P$. This proves the following lemma. 
\begin{lemma}
The canonical gradient of criterion $P\rightarrow f_{n,0}(P)\equiv \{\bar{R}^{(1)}(\tilde{P}_n^{(1)}(P),P_0)\}^2$ at $P$ is given by
\[
D_{f_{n,0},P}=\left\{ \Psi(\tilde{P}_n^{(1)}(P))-\Psi(P_0)-\tilde{P}_n D^{(1)}_{P_0}\right \}D^{(2)}_{n,P},\]
where $D^{(2)}_{n,P}$ is the canonical gradient of $\Psi_n^{(1)}$ at $P$.
\end{lemma}
Thus, this canonical gradient is just a constant times $D^{(2)}_{n,P}$. Therefore a locally steepest descent path  such as $p_{\epsilon}=(1+\epsilon D_{f_{n,0},P})p$ through a density $p$ for the purpose of minimizing $f_{n,0}(P)$ would be identical to $p_{\epsilon}=(1+\epsilon D^{(2)}_{n,P})p$.
Therefore, a  local steepest descent path for minimizing $f_{n,0}(P)$ in $P$ is identical to the least favorable path used by the TMLE of $\Psi_n^{(1)}(P_0)$. In particular, the universal steepest path defined by locally tracking the local steepest path for minimizing $f_{n,0}(P)$ is identical to the universal least favorable path for the TMLE of $\Psi_n^{(1)}(P_0)$.
A steepest descent  algorithm for minimizing $f_{n,0}(P)$ would stop at $P=P_n^*$ when $D^{(2)}_{n,P_n^*} =0$ at which point it has reached a local minimum in neighborhood of $P_n^0$. However, the TMLE $P_n^{(2),*}$ of $\Psi_n^{(1)}(P_0)$ would stop before that, namely it would stop when the log-likelihood along this path is not increasing anymore. At that point, $\tilde{P}_nD^{(2)}_{n,P_n^{2,*}}=0$.

Thus the {\em only} difference between the TMLE of $\Psi^{(2)}_n(P_0)$ and the oracle steepest descent algorithm for minimizing the exact total remainder $f_{n,0}(P)=(\bar{R}^{(1)} (\tilde{P}_n^{(1)}(P),P_0)^2$ is that the TMLE stops earlier, while the oracle algorithm would proceed along the same path in a direction that is not informed by data anymore.
So the TMLE stops earlier due to not being able to increase the likelihood anymore, at which point there is no further guidance on how to minimize the exact remainder.

\subsection{The TMLE also moves in the same direction as the steepest descent algorithm for minimizing the exact total remainder.}
The canonical gradients of $f_{n,0}(P)$ and $\Psi_n^{(1)}(P)$ differ by a scalar $C_n$. Therefore, if $C_n>0$, then the steepest descent algorithm for $f_{n,0}(P)$ would choose small steps $\epsilon<0$, while if $C_n<0$, it would choose small steps $\epsilon>0$. 
The TMLE chooses small steps $\epsilon$ whose sign is driven by the requirement to increase the log-likelihood. Therefore, we want to know if the path in the direction of increasing the log-likelihood agrees with the direction needed for minimizing the exact remainder $f_{n,0}(P)$.
This scalar $C_n$ is given by:
\begin{eqnarray*}
C_n&=&\Psi(\tilde{P}_n^{(1)}(P))-\Psi(P_0)-\tilde{P}_n D^{(1)}_{P_0}\\
&=&
\{\Psi(\tilde{P}_n^{(1)}(P))-\Psi_n^{(2)}(P_0)\}+\Psi_n^{(2)}(P_0)-\Psi(P_0)-\tilde{P}_n D^{(1)}_{P_0}\\
&\approx &\Psi^{(2)}_n(P)-\Psi_n^{(2)}(P_0)-\tilde{P}_n D^{(1)}_{P_0}+O_P(n^{-1}),
\end{eqnarray*}
 where we use that
\[
\Psi(\tilde{P}_n^{(1)}(P_0))-\Psi(P_0)=(\tilde{P}_n-P_0)D^{(1)}_{\tilde{P}_n^{(1)}(P_0)}+R^{(1)}(\tilde{P}_n^{(1)}(P_0),P_0);\]
the exact remainder  $R^{(1)}$ on right-hand side is $O_P(n^{-1})$ due to $\tilde{P}_n^{(1)}(P_0)$ being a TMLE using as initial the true $P_0$ and $\tilde{P}_n$ being appropriately undersmoothed (Section \ref{section7} and Lemma \ref{lemmaexactepsn1});  $(\tilde{P}_n-P_0)D^{(1)}_{\tilde{P}_n^{(1)}(P_0)}=\tilde{P}_n D^{(1)}_{P_0}+O_P(n^{-1})$, because of the same reason.

So we conclude that the scalar $C_n$ equals $\Psi(\tilde{P}_n^{(1)}(P))-\Psi(\tilde{P}_n^{(1)}(P_0)$ plus a random error term $O_P(n^{-1/2})$ that is approximately $n^{-1/2}N(0,\sigma^2_0=P_0\{D^{(1)}_{P_0}\}^2)$. 
The steepest descent algorithm for  minimizing $f_{n,0}(P)$ would move in the direction of  the negative of its canonical gradient, $-D_{f_{n0},P}$, since the goal is to minimize $f_{n,0}(P)$, and the gradient provides the direction of steepest increase. Suppose $\Psi_n^{(2)}(P)-\Psi_n^{(2)}(P_0)>0$.
Then, $D_{f_{n,0},P}=C_n D^{(2)}_{n,P}$ with $C_n>0$. So the steepest descent applied to $f_{n0}$ would move in direction $-D^{(2)}_{n,P}$, and thus chooses $\epsilon<0$. Due to the least favorable path be the path of maximal square change in target estimand, our TMLE algorithm wants to move $\Psi_n^{(2)}(P)$ downwards towards the truth $\Psi_n^{(2)}(P_0)$, so it would use $\epsilon<0$ in the local least favorable path. So in this case both the TMLE and the steepest descent algorithm for $f_{n,0}$  move in the same direction $\epsilon<0$.
Suppose now that $\Psi_n^{(2)}(P)-\Psi_n^{(2)}(P_0)<0$. Increasing the likelihood along the least favorable path corresponds to moving $\Psi_n^{(2)}(P)$ upwards, which corresponds with setting $\epsilon>0$.
The steepest descent applied to $f_{n0}$ would now have a canonical gradient that is $-C_nD^{(2)}_{n,P}$ with $C_n>0$. However, it wants to minimize $\bar{R}_{2n,0}$, so it would select $\epsilon>0$.  So again, the TMLE follows the  steepest descent algorithm applied to $f_{n,0}(P)$, not only by using the same path but also by  moving in the same direction along the path (the direction that decreases the exact remainder $f_{n,0}(P)$.

\subsection{Demonstration of targeting the total remainder}

\begin{figure}[ht!]
\includegraphics[scale = .6]{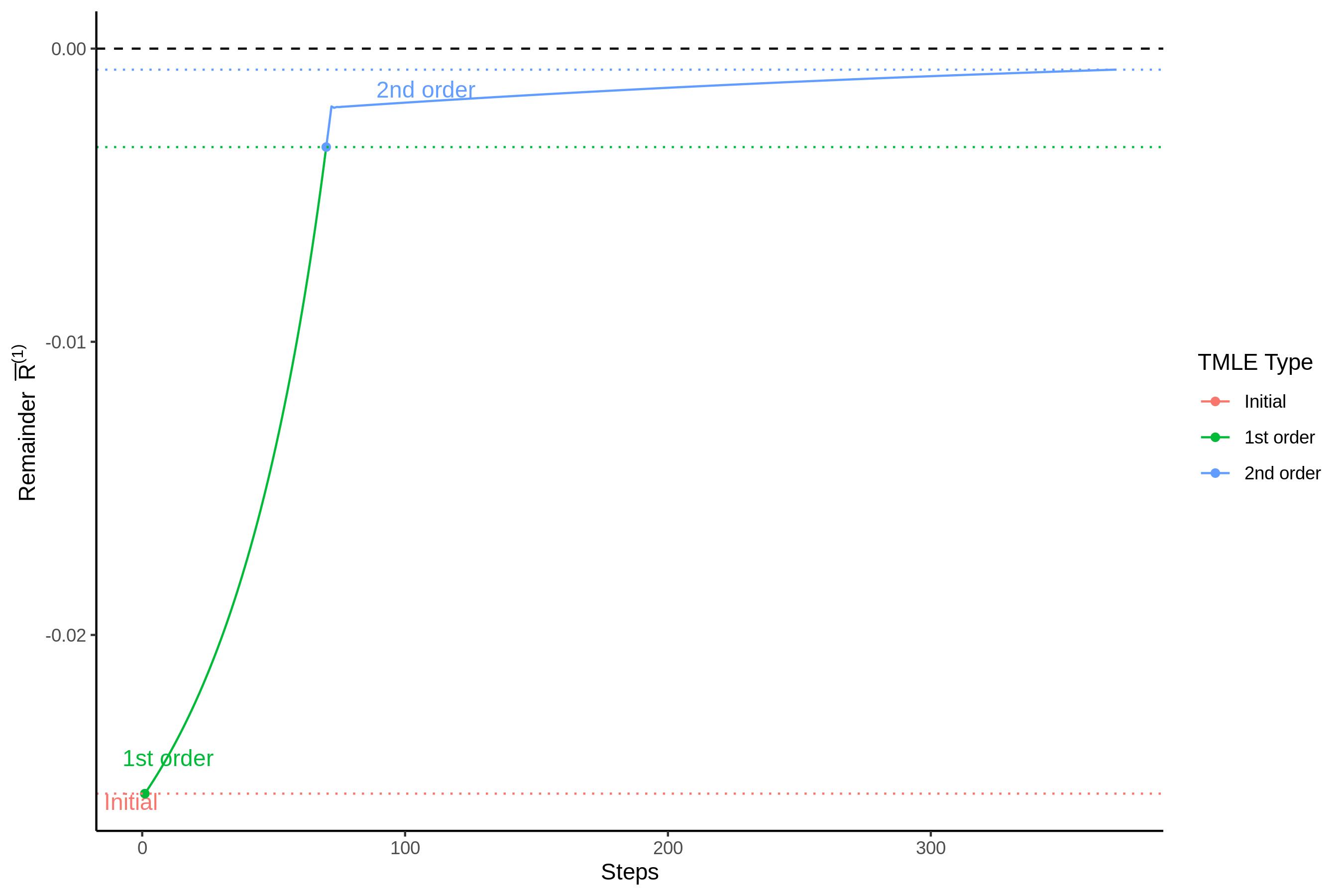}
  \caption{The impact of second order TMLE steps. The exact total remainder $\bar{R}^{(1)}(\tilde P_n^{(1)}(P), P_0)$ of first order TMLE is controlled due to the second order updates $P_n^0 \mapsto P_n^{2, *}$. }
\label{fig:track}
\end{figure}

To demonstrate the direction of steepest descent searches and the effect of second order updates $P = P_n^{(2), *}$ in minimizing the total remainder $\bar{R}^{(1)}(\tilde P_n^{(1)}(P), P_0)$, we simulate discrete univariate data and estimate the integrated square of density according to Section \ref{sec:ave_dens} with $K=4$, $n=500$, and the initial $p_n^0$ is biased by adding a same mass of $0.06$ to all the supports of the empirical pmf $p_n$ and scaling to sum 1. The first order update is iterated following the locally least favorable path with $|\epsilon^{(1)}_n| < 0.1$ at each iteration till $| \tilde{P}_n D^{(1)}_{P^{1, *}_n}| < 1/n$. The second order updates then follow with a fixed step size of $d\epsilon^{(2)} = 0.01$ along the universal least favorable path till $| \tilde{P}_n D^{(2)}_{n, P^{(2), *}_n}| < 1/n$. Figure \ref{fig:track} illustrates the process where the exact total remainder, $\bar{R}^{(1)}(\tilde P_n^{(1)}(P), P_0)$, left by the first order TMLE at $P=P_n^0$, is minimized due to the second order updates $P_n^0 \mapsto P_n^{2, *}$.

\subsection{Summary}
We conclude that the second order TMLE $\tilde{P}_n^{(2)}(P_n^0)$ represents a steepest descent algorithm that uses the  same path as the oracle steepest descent algorithm for minimizing the square of the exact (unknown) total remainder $P\rightarrow \bar{R}^{(1)}(\tilde{P}_n^{(1)}(P),P_0)$, starting at $P_n^0$. However, the TMLE stops the moment the log-likelihood is not improving anymore. Interestingly, that moment corresponds exactly with the moment the TMLE would not be able to know the signal of the scalar $C_n$ anymore in the canonical gradient of the square  $f_{n,0}(P)$ of the exact total remainder, and thus not know anymore in  what direction a steepest descent algorithm for $f_{n,0}(P)$ would move.

If the steepest descent algorithm for minimizing $f_{n,0}$ goes beyond this point,  it will be moving in random directions that switch from negative to positive. 
However, if one truly has available $f_{n,0}$ one would still know how to proceed, potentially finding a minimum of $f_{n,0}(P)$ at a $P$ that is far from $P_0$. Therefore, the TMLE stops at the right moment for the purpose of minimizing $f_{n,0}(P)$ when one also cares about fitting $P_0$, and thereby finding a local minimum in close neighborhood of the desired minimum  $P_0$.

Since the total remainder achieves its desired global minimum at $P_0$, this restriction that it should stop once the likelihood cannot be improved anymore is completely sensible: without the log-likelihood, there is no criterion anymore to decide in what direction to move along the least favorable path. That is, the TMLE $\tilde{P}_n^{(2)}(P_n^0)$ is doing the optimal possible job in locally minimizing the square of the exact total remainder at initial $P_n^0$.

\section{Exact $k+1$-th order expansions for the HAL-regularized $k$-th order TMLE}\label{section4}

\subsection{Exact expansion for HAL-regularized $k$-th order TMLE in terms of $R_n^{(k)}$ at $P_0$.}
The following theorem provides an exact  $k+1$-th order expansion for the HAL-regularized $k$-th order TMLE. 
\begin{theorem}\label{theorem1}
\begin{eqnarray*}
\Psi(P_n^{k,(1),*})-\Psi(P_0)
&=&\Psi(\tilde{P}_n^{(1)}\circ\ldots\circ \tilde{P}_n^{(k)}(P_n^0))-\Psi(P_0)\\
&=&\sum_{j=0}^{k-2} (P_n-P_0)D^{(j+1)}_{n,\tilde{P}_n^{(j+1)}(P_0)}+R_n^{(j+1)}(\tilde{P}_n^{(j+1)}(P_0),P_0)\\
&&+(P_n-P_0)D^{(k)}_{n,\tilde{P}_n^{(k)}(P_n^0)}+R_{n}^{(k)}(\tilde{P}_n^{(k)}(P_n^0),P_0)\\
&&-\sum_{j=0}^{k-2}P_n D^{(j+1)}_{n,\tilde{P}_n^{(j+1}(P_0)}-P_n D^{(k)}_{n,\tilde{P}_n^{(k)}(P_n^0)}
.
\end{eqnarray*}
Note that one can replace any $P_n D^{(j+1)}_{n,\tilde{P}_n^{(j+1)}(P_0)}$ by $\tilde{P}_nD^{(j+1)}_{n,\tilde{P}_n^{(j+1)}(P_0)}$, $j=1,\ldots,k-1$.

Suppose we use a $k$-th order TMLE that satisfies $E_n\equiv -\sum_{j=1}^{k		}\tilde{P}_n D^{(j)}_{n,\tilde{P}_n^{(j)}(P_0)}=0$ exactly (which can be arranged by defining $\epsilon_n^{(j)}(P)$ as a solution of the corresponding efficient score equation $\tilde{P}_n D^{(j)}_{n,\tilde{P}_n^{(j)}(P,\epsilon)}=0$).

Then, we can also represent this as:
\begin{eqnarray*}
\Psi(P_n^{k,(1),*})-\Psi(P_0)&=& (P_n-P_0)D^{(1)}_{\tilde{P}_n^{(1)}(P_0)}+R^{(1)}(\tilde{P}_n^{(1)}(P_0),P_0)\\
&&+\sum_{j=1}^{k-2} (\tilde{P}_n-P_0)D^{(j+1)}_{n,\tilde{P}_n^{(j+1)}(P_0)}+R_n^{(j+1)}(\tilde{P}_n^{(j+1)}(P_0),P_0)\\
&&+(\tilde{P}_n-P_0)D^{(k)}_{n,P_n^{(k)}(P_n^0)}+R_n^{(k)}(P_n^{(k)}(P_n^0),P_0)\\
&&+(\tilde{P}_n-P_n) D^{(1)}_{\tilde{P}_n^{(1)}(P_0)}.
\end{eqnarray*}
\end{theorem}
This result follows from the following expansion of the $k$-th order TMLE:
\begin{eqnarray*}
\Psi(P_n^{k,(1),*})-\Psi(P_0)&=&
\Psi(\tilde{P}_n^{(1)}\ldots\tilde{P}_n^{(k-1)} P_n^{(k)}P_n^0)-\Psi(P_0)\\
&=&
\sum_{j=0}^{k-2}\Psi_n^{(j)}(\tilde{P}_n^{(j+1)}(P_0))-\Psi_n^{(j)}(P_0)\\
&&+
\Psi_n^{(k-1)}(P_n^{(k)}(P_n^0))-\Psi_n^{(k-1)}(P_0),
\end{eqnarray*}
and using the exact expansion for each of these TMLE differences:
\[
\Psi_n^{(j)}(\tilde{P}_n^{(j+1)}(P_0))-\Psi_n^{(j)}(P_0)=-P_0D^{(j+1)}_{n,\tilde{P}_n^{(j+1)}(P_0)}+R_n^{(j+1)}(\tilde{P}_n^{(j+1)}(P_0),P_0), 
\]
for $j=0,\ldots,k-2$, and
\[
\Psi_n^{(k-1)}(P_n^{(k)}(P_n^0))-\Psi_n^{(k-1)}(P_0)=-P_0D^{(k)}_{n,P_n^{(k)}(P_n^0)}+R_{20}^{(k)}(P_n^{(k)}(P_n^0),P_0).\]
Further details are presented in Appendix \ref{AppendixB}.
This exact expansion of the $k$-th order TMLE can also be understood as implied by
\[
\Psi(P_n^{k,(1),*})-\Psi(P_0)=\Psi_n^{(k-1)}(P_0)-\Psi(P_0)+
\Psi_n^{(k-1)}(\tilde{P}_n^{(k)}(P_n^0))-\Psi_n^{(k-1)}(P_0)
,\]
where the second difference can be replaced by the exact expansion for the TMLE of $\Psi_n^{(k-1)}(P_0)$, and the first term equals a plug-in parametric (sequential and HAL-MLE regularized) MLE minus truth according a correctly specified parametric  model (with parameters $\epsilon_n^{(j)}$, $j=k,\ldots,1$).

   \subsection{Relating exact remainder $R_n^{(k)}(P,P_0)$ at true data distribution to exact remainder  $R_n^{(k)}(P,\tilde{P}_n)$ at HAL-MLE.}
  It appears that the most natural results are obtained by studying the remainder $R_n^{(k)}(P,\tilde{P}_n)$ at $\tilde{P}_n$ (instead of $P_0$), due to the HAL-regularized TMLE procedure being aimed at maximizing the likelihood at $\tilde{P}_n$. The key property of $\tilde{P}_n$ that makes the results nice is that $\tilde{P}_n^{(j)}(\tilde{P}_n)=\tilde{P}_n$ for all $j=1,\ldots,k-1$, and possibly for $j=k$ (if it uses $\tilde{P}_n$).
  The following lemma relates $R_n^{(j)}(P,\tilde{P}_n)$ to $R_n^{(j)}(P,P_0)$.
  
 \begin{lemma}\label{lemma3}
We have for $j=2,\ldots,k$
\[
\begin{array}{l}
R_n^{(j)}(\tilde{P}_n^{(j)}(P),P_0)=R_n^{(j)}(\tilde{P}_n^{(j)}(P),\tilde{P}_n)-R_n^{(j)}(\tilde{P}_n^{(j)}(P_0),\tilde{P}_n)\\
+R_n^{(j)}(\tilde{P}_n^{(j)}(P_0),P_0)
+(\tilde{P}_n-P_0) \{D^{(j)}_{n,\tilde{P}_n^{(j)}(P_0)}-D^{(j)}_{n,\tilde{P}_n^{(j)}(P)} \}.
\end{array}
\]
\end{lemma}
{\bf Proof:}
Let $P_n^{j,*}=\tilde{P}_n^{(j)}(P)$. We have 
\[
\begin{array}{l}
\Psi_n^{(j-1)}(P_n^{j,*})-\Psi_n^{(j-1)}(P_0)=\Psi_n^{(j-1)}(P_n^{j,*})-\Psi_n^{(j-1)}(\tilde{P}_n)\\
\hfill -
\left( \Psi_n^{(j-1)}(P_0)-\Psi_n^{(j-1)}(\tilde{P}_n)\right) \\
=R_n^{(j)}(P_n^{j,*},\tilde{P}_n)-\tilde{P}_n  D^{(j)}_{P_n^{j,*}}-(\Psi_n^{(j-1)}(P_0)-\Psi_n^{(j-1)}(\tilde{P}_n)),
\end{array}
\]
where we used that $\Psi_n^{(j-1)}(P)-\Psi_n^{(j-1)}(\tilde{P}_n)=-\tilde{P}_n D^{(j)}_{P}+R_n^{(j)}(P,\tilde{P}_n)$.
Now, we note that 
\[
\begin{array}{l}
\Psi_n^{(j-1)}(P_0)-\Psi_n^{(j-1)}(\tilde{P}_n)=(\Psi_n^{(j-1)}(P_0)-\Psi_n^{(j-1)}(\tilde{P}_n^{(j)}(P_0))\\
\hfill +\Psi_n^{(j-1)}(\tilde{P}_n^{(j)}(P_0))-\Psi_n^{(j-1)}(\tilde{P}_n)\\
=(\Psi_n^{(j-1)}(P_0)-\Psi_n^{(j-1)}(\tilde{P}_n^{(j)}(P_0))+R_n^{(j)}(\tilde{P}_n^{(j)}(P_0),\tilde{P}_n)-\tilde{P}_n D^{(j)}_{n,\tilde{P}_n^{(j)}(P_0)}.
\end{array}
\]
We have \[
\Psi_n^{(j-1)}(\tilde{P}_n^{(j)}(P_0))-\Psi_n^{(j-1)}(P_0) =-P_0D^{(j)}_{n,\tilde{P}_n^{(j)}(P_0)}+R_n^{(j)}(\tilde{P}_n^{(j)}(P_0),P_0).\]
Thus, combining the above expressions yields
\[
\begin{array}{l}
\Psi_n^{(j-1)}(P_n^{j,*})-\Psi_n^{(j-1)}(P_0)= R_n^{(j)}(P_n^{j,*},\tilde{P}_n)-R_n^{(j)}(\tilde{P}_n^{(j)}(P_0),\tilde{P}_n)\\
-\tilde{P}_n  D^{(j)}_{P_n^{j,*}}
+\tilde{P}_n D^{(j)}_{n,\tilde{P}_n^{(j)}(P_0)}+R_n^{(j)}(\tilde{P}_n^{(j)}(P_0),P_0)
-P_0D^{(j)}_{n,\tilde{P}_n^{(j)}(P_0)}
.
\end{array}
\]
However, we also know that \[
\Psi_n^{(j-1)}(P_n^{j,*})-\Psi_n^{(j-1)}(P_0)=-P_0D^{(j)}_{P_n^{j,*}}+R_n^{(j)}(P_n^{j,*},P_0).\]
Setting the two expressions equal  to each other gives
\[
\begin{array}{l}
R_n^{(j)}(P_n^{j,*},P_0)=
R_n^{(j)}(P_n^{j,*},\tilde{P}_n)-R_n^{(j)}(\tilde{P}_n^{(j)}(P_0),\tilde{P}_n)
-\tilde{P}_n  D^{(j)}_{P_n^{j,*}}
+\tilde{P}_n D^{(j)}_{n,\tilde{P}_n^{(j)}(P_0)}\\
\hfill +R_n^{(j)}(\tilde{P}_n^{(j)}(P_0),P_0)
-P_0D^{(j)}_{n,\tilde{P}_n^{(j)}(P_0)}
+P_0D^{(j)}_{P_n^{j,*}}\\
=R_n^{(j)}(P_n^{j,*},\tilde{P}_n)-R_n^{(j)}(\tilde{P}_n^{(j)}(P_0),\tilde{P}_n)+
R_n^{(j)}(\tilde{P}_n^{(j)}(P_0),P_0)\\
\hfill +P_0\{  D^{(j)}_{P_n^{j,*}}-D^{(j)}_{n,\tilde{P}_n^{(j)}(P_0)}\}
-\tilde{P}_n \{ D^{(j)}_{P_n^{j,*}}- D^{(j)}_{n,\tilde{P}_n^{(j)}(P_0)} \}\\
=R_n^{(j)}(P_n^{j,*},\tilde{P}_n)-R_n^{(j)}(\tilde{P}_n^{(j)}(P_0),\tilde{P}_n)+
R_n^{(j)}(\tilde{P}_n^{(j)}(P_0),P_0)\\
\hfill +(P_0-\tilde{P}_n)\{  D^{(j)}_{P_n^{j,*}}-D^{(j)}_{n,\tilde{P}_n^{(j)}(P_0)}\}.
\end{array}
\]
Finally, we write the last term as $(\tilde{P}_n-P_0)\{  D^{(j)}_{n,\tilde{P}_n^{(j)}(P_0)}-
 D^{(j)}_{P_n^{j,*}}\}$.
$\Box$

\subsection{Representing $k$-th exact remainder $R_n^{(k)}(P,\tilde{P}_n)$ at HAL-MLE as a sequentially targeted exact remainder $R_n^{(1)}$ of target estimand.}
The following lemma establishes that $R_n^{(k)}(\tilde{P}_n^{(k)}(P),\tilde{P}_n)=R^{(1)}(\tilde{P}_n^{(1)}\ldots\tilde{P}_n^{(k)}(P),\tilde{P}_n)$, under the assumption that  $\tilde{P}_n D^{(1)}_{n,\tilde{P}_n^{(1)}(P)}=0$ and $\tilde{P}_n D^{(k)}_{n,\tilde{P}_n^{(k)}(P)}=0$.

\begin{lemma}\label{lemma2}
For $j=1,\ldots,k$, we have
\[
R_n^{(j)}(\tilde{P}_n^{(j)}(P),\tilde{P}_n)=R^{(1)}(\tilde{P}_n^{(1)}\ldots \tilde{P}_n^{(j)}(P),\tilde{P}_n)
-\tilde{P}_nD^{(1)}_{\tilde{P}_n^{(1)}\ldots \tilde{P}_n^{(j)}(P)}+
\tilde{P}_nD^{(j)}_{n,\tilde{P}_n^{(j)}(P)}.\]
As a consequence, for $j=0,\ldots,k-1$, we have
\[
\begin{array}{l}
R_n^{(j+1)}(\tilde{P}_n^{(j+1)}(P),\tilde{P}_n)=R_n^{(j)}(\tilde{P}_n^{(j)}(\tilde{P}_n^{(j+1)}(P) ) ,\tilde{P}_n)
-\tilde{P}_nD^{(j)}_{n,\tilde{P}_n^{(j)}\tilde{P}_n^{(j+1)}(P)}+
\tilde{P}_nD^{(j+1)}_{\tilde{P}_n^{(j+1)}(P)}.
\end{array}
\]
We also note
 \begin{eqnarray*}
 \Psi_n^{(j)}(P)-\Psi_n^{(j)}(\tilde{P}_n)&=&
\Psi_n^{(j-1)}(\tilde{P}_n^{(j)}(P))-\Psi_n^{(j-1)}(\tilde{P}_n)\\
&=&-\tilde{P}_n D^{(j)}_{n,\tilde{P}_n^{(j)}(P)}+R_n^{(j)}(\tilde{P}_n^{(j)}(P),\tilde{P}_n)\\
&\approx& R_n^{(j)}(\tilde{P}_n^{(j)}(P),\tilde{P}_n),\end{eqnarray*}
where latter is an exact equality if  $\tilde{P}_n^{(j)}(P)$ solves $\tilde{P}_n D^{(j)}_{n,\tilde{P}_n^{(j)}(P)}=0$ exactly.
Thus, a TMLE $\tilde{P}_n^{(j+1)}(P_n^0)$ targeting $\Psi_n^{(j)}(P_0)$ is equivalent with a TMLE targeting $R_n^{(j)}(\tilde{P}_n^{(j)}(P_0),\tilde{P}_n)$.

\end{lemma}
{\bf Proof of Lemma \ref{lemma2}:}
By definition of the exact remainder $R^{(1)}(P,P_1)$ for general $(P,P_1)\in {\cal M}^2$, we have
\[
R^{(1)}(\tilde{P}_n^{(1)}(P),\tilde{P}_n)=\Psi(\tilde{P}_n^{(1)}(P))-\Psi(\tilde{P}_n)+\tilde{P}_n D^{(1)}_{\tilde{P}_n^{(1)}(P)}.\]
In other words,  for any $P$
\begin{equation}\label{a1}
 \Psi(\tilde{P}_n^{(1)}(P))-\Psi(\tilde{P}_n)=R^{(1)}(\tilde{P}_n^{(1)}(P),\tilde{P}_n)-\tilde{P}_n D^{(1)}_{\tilde{P}_n^{(1)}(P)}.\end{equation}
By definition of $R^{(2)}(\tilde{P}_n^{(2)}(P),\tilde{P}_n)$, we also have \[
\begin{array}{l}
R_n^{(2)}(\tilde{P}_n^{(2)}(P),\tilde{P}_n)=\Psi(\tilde{P}_n^{(1)}\tilde{P}_n^{(2)}(P))-\Psi(\tilde{P}_n)+\tilde{P}_nD^{(2)}_{n,\tilde{P}_n^{(2)}(P)}.
\end{array}
\]
Applying the expression for (\ref{a1}) with $P=\tilde{P}_n^{(2)}(P)$ gives
\[
\begin{array}{l}
R_n^{(2)}(\tilde{P}_n^{(2)}(P),\tilde{P}_n)=
R^{(1)}(\tilde{P}_n^{(1)}\tilde{P}_n^{(2)}(P),\tilde{P}_n)-\tilde{P}_nD^{(1)}_{\tilde{P}_n^{(1)}\tilde{P}_n^{(2)}(P)}+\tilde{P}_n D^{(2)}_{n,\tilde{P}_n^{(2)}P}.
\end{array}
\]

Similarly, we have \[
\begin{array}{l}
R_n^{(3)}(\tilde{P}_n^{(3)}(P),\tilde{P}_n)=\Psi(\tilde{P}_n^{(1)}\tilde{P}_n^{(2)}\tilde{P}_n^{(3)}(P))-\Psi(\tilde{P}_n)+\tilde{P}_nD^{(3)}_{n,\tilde{P}_n^{(3)}(P)}\\
=\Psi_n^{(1)}(\tilde{P}_n^{(2)}\tilde{P}_n^{(3)}(P))-\Psi_n^{(1)}(\tilde{P}_n)+ \tilde{P}_nD^{(3)}_{n,\tilde{P}_n^{(3)}(P)}\\
=-\tilde{P}_nD^{(2)}_{\tilde{P}_n^{(2)}\tilde{P}_n^{(3)}(P)}+R_n^{(2)}(\tilde{P}_n^{(2)}\tilde{P}_n^{(3)}(P),\tilde{P}_n)+
\tilde{P}_nD^{(3)}_{n,\tilde{P}_n^{(3)}(P)}\\
= -\tilde{P}_nD^{(2)}_{\tilde{P}_n^{(2)}\tilde{P}_n^{(3)}(P)}
+R^{(1)}(\tilde{P}_n^{(1)}\tilde{P}_n^{(2)}\tilde{P}_n^{(3)}(P),\tilde{P}_n)-\tilde{P}_nD^{(1)}_{\tilde{P}_n^{(1)}\tilde{P}_n^{(2)}\tilde{P}_n^{(3)}(P)}\\
\hfill +\tilde{P}_n D^{(2)}_{n,\tilde{P}_n^{(2)}\tilde{P}_n^{(3)}(P)}
+\tilde{P}_nD^{(3)}_{n,\tilde{P}_n^{(3)}(P)}\\
=R^{(1)}(\tilde{P}_n^{(1)}\tilde{P}_n^{(2)}\tilde{P}_n^{(3)}(P),\tilde{P}_n)
-\tilde{P}_nD^{(1)}_{\tilde{P}_n^{(1)}\tilde{P}_n^{(2)}\tilde{P}_n^{(3)}(P)}+
\tilde{P}_nD^{(3)}_{\tilde{P}_n^{(3)}(P)}.
\end{array}
\]

So in general
\[
R_n^{(j)}(\tilde{P}_n^{(j)}(P),\tilde{P}_n)=R^{(1)}(\tilde{P}_n^{(1)}\ldots \tilde{P}_n^{(j)}(P),\tilde{P}_n)
-\tilde{P}_nD^{(1)}_{\tilde{P}_n^{(1)}\ldots \tilde{P}_n^{(j)}(P)}+
\tilde{P}_nD^{(j)}_{n,\tilde{P}_n^{(j)}(P)}.\]
This proves the first statement.
As a consequence, this implies
\[
\begin{array}{l}
R_n^{(j+1)}(\tilde{P}_n^{(j+1)}(P),\tilde{P}_n)
=R^{(1)}(\tilde{P}_n^{(1)}\ldots \tilde{P}_n^{(j+1)}(P),\tilde{P}_n) -\tilde{P}_nD^{(1)}_{\tilde{P}_n^{(1)}\ldots \tilde{P}_n^{(j+1)}(P)}\\
\hfill +
\tilde{P}_nD^{(j+1)}_{\tilde{P}_n^{(j+1)}(P)}\\
=R^{(1)}(\tilde{P}_n^{(1)}\ldots \tilde{P}_n^{(j)}(\tilde{P}_n^{(j+1)}(P)),\tilde{P}_n)
-\tilde{P}_nD^{(1)}_{\tilde{P}_n^{(1)}\ldots \tilde{P}_n^{(j+1)}(P)}+
\tilde{P}_nD^{(j+1)}_{\tilde{P}_n^{(j+1)}(P)}\\
=R_n^{(j)}(\tilde{P}_n^{(j)}(\tilde{P}_n^{(j+1)}(P) ) ,\tilde{P}_n)+\tilde{P}_nD^{(1)}_{\tilde{P}_n^{(1)}\ldots \tilde{P}_n^{(j)}\tilde{P}_n^{(j+1)}(P)} -\tilde{P}_nD^{(j)}_{n,\tilde{P}_n^{(j)}\tilde{P}_n^{(j+1)}(P)}\\
\hfill -\tilde{P}_nD^{(1)}_{\tilde{P}_n^{(1)}\ldots \tilde{P}_n^{(j+1)}(P)}+
\tilde{P}_nD^{(j+1)}_{\tilde{P}_n^{(j+1)}(P)}\\
=R_n^{(j)}(\tilde{P}_n^{(j)}(\tilde{P}_n^{(j+1)}(P) ) ,\tilde{P}_n)
-\tilde{P}_nD^{(j)}_{n,\tilde{P}_n^{(j)}\tilde{P}_n^{(j+1)}(P)}+
\tilde{P}_nD^{(j+1)}_{\tilde{P}_n^{(j+1)}(P)}.
\end{array}
\]

This proves the second statement and thereby the lemma.  $\Box$

 \subsection{Exact expansion for the $k$-th order TMLE in terms of exact $k+1$-th order remainder at the HAL-MLE}
 The last Lemma \ref{lemma3} gives us an expression for $R_n^{(k)}({P}_n^{(k)}(P_n^0),P_0)$. Substituting this into  our exact expansion  for the $k$-th order TMLE of Theorem \ref{theorem1} provides the following exact expansion for the $k$-th order TMLE. 
 \begin{theorem}\label{theorem2}
  We have
  \begin{eqnarray*}
\Psi(P_n^{k,(1),*})-\Psi(P_0)
&=&\Psi(\tilde{P}_n^{(1)}\circ\ldots\circ \tilde{P}_n^{(k)}(P_n^0))-\Psi(P_0)\\
&=&\sum_{j=1}^{k}\left\{ (\tilde{P}_n-P_0)D^{(j)}_{n,\tilde{P}_n^{(j)}(P_0)}+R_n^{(j)}(\tilde{P}_n^{(j)}(P_0),P_0)\right\} \\
&& +R_n^{(k)}(\tilde{P}_n^{(k)}(P_n^0),\tilde{P}_n)-R_n^{(k)}(\tilde{P}_n^{(k)}(P_0),\tilde{P}_n)\\
&&-\sum_{j=1}^{k}\tilde{P}_n D^{(j)}_{n,\tilde{P}_n^{(j}(P_0)}.
\end{eqnarray*}
For example, if we use a $k$-th order TMLE that satisfies $\tilde{P}_n D^{(j)}_{n,\tilde{P}_n^{(j)}(P_0)}=0$ exactly for $j=1,\ldots,k$, then
\begin{eqnarray*}
\Psi(P_n^{k,(1),*})-\Psi(P_0)&=&(P_n-P_0)D^{(1)}_{\tilde{P}_n^{(1)}(P_0)}+R^{(1)}(\tilde{P}_n^{(1)}(P_0),P_0)\\
&&+\sum_{j=2}^{k}\left\{ (\tilde{P}_n-P_0)D^{(j)}_{n,\tilde{P}_n^{(j)}(P_0)}+R_n^{(j)}(\tilde{P}_n^{(j)}(P_0),P_0)\right\} \\
&& +R_n^{(k)}(P_n^{(k)}(P_n^0),\tilde{P}_n)-R_n^{(k)}(P_n^{(k)}(P_0),\tilde{P}_n)\\
&&+(\tilde{P}_n-P_n)D^{(1)}_{\tilde{P}_n^{(1)}(P_0)},
\end{eqnarray*}
and, by above lemma,  \[
R_n^{(k)}(P_n^{(k)}(P),\tilde{P}_n)=R^{(1)}(\tilde{P}_n^{(1)}\ldots P_n^{(k)}(P),\tilde{P}_n).\]
\end{theorem}
The proof of this theorem is an immediate corollary of Theorem \ref{theorem1}, Lemma \ref{lemma3} and Lemma \ref{lemma2}.
In the next section we prove that $R^{(1)}(\tilde{P}_n^{(1)}\ldots P_n^{(k)}(P),\tilde{P}_n)$
is a $k+1$-th order difference.
The next subsection bounds the undersmoothing term $(\tilde{P}_n-P_n)D^{(1)}_{\tilde{P}_n^{(1)}(P_0)}$ showing that it is bounded in terms of $(\tilde{P}_n-P_n)D^{(1)}_{P_0}$.

\subsection{Exact bound for the HAL-regularization bias term.}
The exact expansion has a bias term due to using $\epsilon_n^{(1)}(P_0)=\epsilon_{\tilde{P}_n}^{(1)}(P_0)$ instead of  the standard MLE $\epsilon_{P_n}^{(1)}(P_0)$ (where the latter  is trivially $O_P(n^{-1/2})$ and asymptotically linear) in the definition of the $k$-th order TMLE, given by $(\tilde{P}_n-P_n)D^{(1)}_{\tilde{P}_n^{(1)}(P_0)}$, and the smaller order analogues $(\tilde{P}_n-P_n)D^{(j)}_{n,\tilde{P}_n^{(j)}(P_0)}$, $j=2,\ldots,k$.

Consider the case that $\tilde{P}_n D^{(1)}_{\tilde{P}_n^{(1)}(P_0)}=0$.
The following lemma provides a finite sample bound on $\epsilon_n^{(1)}(P_0)$ and thereby for the term $P_n D^{(1)}_{\tilde{P}_n^{(1)}(P_0)}=-(\tilde{P}_n-P_n)D^{(1)}_{\tilde{P}_n^{(1)}(P_0)}$ in the exact expansion. It shows that it is bounded by $(\tilde{P}_n-P_n)D^{(1)}_{P_0}$, so that it is all about a plug-in HAL-MLE $\tilde{P}_nf$ of a mean $P_0f$ being asymptotically equivalent with the sample mean $P_n f$ for a fixed function $f=D_{P_0}^{(1)}$. This is verifiable based on data $P_n$ at any given $f$, so that one can tune $L_1$-norm $C$ in $\tilde{P}_n$ accordingly, and, we also have the option to target $\tilde{P}_n$ towards this mean $P_0D^{(1)}_{P_0}$ as in Section \ref{sectionthotmle}. 

\begin{lemma}\label{lemmaexactepsn1}
Let $I_{\xi_n,P_0}\equiv 
-\frac{d}{d\xi_n}\tilde{P}_n D^{(1)}_{\tilde{P}^{(1)}(P_0,\xi_n)}$ for $\xi_n$ sandwiched by $0$ and $\epsilon_{\tilde{P}_n}^{(1)}(P_0)$, and notice that $I_{\xi_n,P_0}\rightarrow_p I_{P_0}=P_0\{D^{(1)}_{P_0}\}^2$.
We have
\[
\epsilon_{\tilde{P}_n}^{(1)}(P_0)=I_{\xi_n,P_0}^{-1}(\tilde{P}_n-P_0)D^{(1)}_{P_0}\]
for some $\xi_{1n}$ sandwiched by $0$ and $\epsilon_{\tilde{P}_n}^{(1)}(P_0)$. 
In particular,
\[
\epsilon_{\tilde{P}_n}^{(1)}(P_0)=I_{\xi_n,P_0}^{-1}(\tilde{P}_n-P_n)D^{(1)}_{P_0}+I_{\xi_n,P_0}^{-1}(P_n-P_0)D^{(1)}_{P_0}.\]
Under the assumption that $P_n D^{(1)}_{\tilde{P}^{(1)}(P_0,\epsilon)}$ is differentiable in $\epsilon$, and using that $P_n D^{(1)}_{\tilde{P}^{(1)}(P_0,\epsilon_{P_n}(P_0)}=0$ and $\tilde{P}_n D^{(1)}_{\tilde{P}^{(1)}(P_0,\epsilon_n(P_0)}=0$, 
this bound also implies that \[
P_n D^{(1)}_{\tilde{P}^{(1)}(P_0,\epsilon_n^{(1))}(P_0))}=P_n \frac{d}{d\xi_n}D^{(1)}_{\tilde{P}^{(1)}(P_0,\xi_{1n})} (\epsilon_{\tilde{P}_n}^{(1)}(P_0)-\epsilon_{P_n}^{(1)}(P_0)) \]
for some $\xi_{1n}$ sandwiched by $0$ and $\epsilon_{\tilde{P}_n}^{(1)}(P_0)$. 
Here the left-hand side equals $(\tilde{P}_n-P_n) D^{(1)}_{\tilde{P}^{(1)}(P_0,\epsilon_n^{(1)}(P_0)}$.
\end{lemma}
Under weak regularity conditions, we also have
\[
\epsilon_{P_n}^{(1)}(P_0)=I_{n,P_0}^{-1}(P_n-P_0)D^{(1)})_{P_0},\]
where $I_{n,P_0}\equiv 
-\frac{d}{d\xi_{1n}}{P}_n D^{(1)}_{\tilde{P}^{(1)}(P_0,\xi_{1n})}$ for a $\xi_{1n}$ sandwiched by $0$ and $\epsilon_{{P}_n}^{(1)}(P_0)$.
Thus, then
\[
\epsilon_{\tilde{P}_n}^{(1)}(P_0)-\epsilon_{P_n}^{(1)}(P_0)=I_{\xi_n,P_0}^{-1}(\tilde{P}_n-P_n)D^{(1)}_{P_0}+(I_{\xi_n,P_0}^{-1}-I_{n,P_0}^{-1})(P_n-P_0)D^{(1)}_{P_0}.\]
Thus, the lemma above proves then that 
\[
P_n D^{(1)}_{\tilde{P}^{(1)}(P_0,\epsilon_n^{(1)}(P_0))}=O_P((\tilde{P}_n-P_n)D^{(1)}_{P_0})+
O_P((I_{\xi_n,P_0}^{-1}-I_{n,P_0}^{-1})(P_n-P_0)D^{(1)}_{P_0}),\]
with well understood constants.
The last term is typically $O_P(n^{-1})$, so that  the key in bounding $
P_n D^{(1)}_{\tilde{P}^{(1)}(P_0,\epsilon_n^{(1))}(P_0)}$ is to  make 
  $(\tilde{P}_n-P_n)D^{(1)}_{P_0}$ small, which is controlled by undersmoothing the HAL-MLE.
  
{\bf Proof:}
 For notational convenience, let $\epsilon_n=\epsilon_{\tilde{P}_n}^{(1)}(P_0)$.
Let $U(\epsilon_n,\tilde{P}_n)=\tilde{P}_n D^{(1)}_{\tilde{P}^{(1)}(P_0,\epsilon_n)}$ and note that
$U(\tilde{P}_n,\epsilon_n)=0$ and $U(P_0,\epsilon_0=0)=0$.
We have 
\[
U(\epsilon_n,\tilde{P}_n)-U(\epsilon_0,\tilde{P}_n)=-U(\epsilon_0,\tilde{P}_n-P_0)=
-(\tilde{P}_n-P_0)D^{(1)}_{P_0}.\]
Using an exact tailor expansion gives that the left-hand side equals
$\frac{d}{d\xi_n}U(\xi_n,P_0) \epsilon_n$ for some $\xi_n\in [0,\epsilon_n]$ or $\xi_n\in [\epsilon_n,0]$. Let $I_{\xi_n,P_0}=-\frac{d}{d\xi_n}U(\xi_n,P_0)$. 
Thus, we have
\[
\epsilon_n=I_{\xi_n,P_0}^{-1}(\tilde{P}_n-P_0)D^{(1)}_{P_0}\mbox{ $\Box$}.\]

\subsection{Discussion of expansion of HAL-regularized  $k$-th order TMLE of Theorem \ref{theorem2}}

{\bf Exact remainders $R^{(j)}(\tilde{P}_n^{(j)}(P_0),P_0)=O_P(n^{-1})$:}
The HAL-MLE $\tilde{P}_n$ generally satisfies  that $(\tilde{P}_n-P_n)D^{(1)}_{P_0}=O_P(n^{-2/3})$ (Appendix \ref{section7}).
Then, it follows that $R_n^{(1)}(\tilde{P}_n^{(1)}(P_0),P_0)$ is $O_P(n^{-1})$, due to $\tilde{P}_n^{(1)}(P_0)$ being an MLE for a one-dimensional correctly specified one dimensional parametric model and thereby that  $d_0(\tilde{P}_n^{(1)}(P_0),P_0)=O_P(\{ \epsilon_n^{(1)}(P_0)\}^2)=O_P(n^{-1})$ (by Lemma \ref{lemmaexactepsn1}).
Similarly, this applies to $R^{(j)}(\tilde{P}_n^{(j)}(P_0),P_0)$ for $j=2,\ldots,k$. It is not just that  these terms are small in rate, these terms also have a constant in front of this rate that is unaffected by the curse of dimensionality.

{\bf Canonical gradients $D^{(j)}_{n,\tilde{P}_n^{(j})(P_0)}$ are almost fixed mean zero independent random variables:}
Note also that each of these empirical process terms have little bias since the random function $D^{(j)}_{n,\tilde{P}_n^{(j)}(P_0)}$ is only random through TMLE-updates $\tilde{P}_n^{(1)}\ldots\tilde{P}_n^{(j)}(P_0)$ that start with an initial being the true $P_0$. Thus, the data dependence is only through a finite dimensional vector of well behaved MLE coefficients $\epsilon_n^{(j)}$ at previous TMLE updates of $P_0$. In addition, $P_0 D^{(j)}_{n,P_0}=0$ while $\tilde{p}_n^{(j)}(P_0)-p_0=O_P(n^{-1/2})$ under $(\tilde{P}_n-P_n)D^{(1)}_{P_0}=O_P(n^{-1/2})$ as studied in previous subsection.  

{\bf Empirical means of higher order efficient scores are solved by some undersmoothing of HAL-MLE:}
Suppose that $\tilde{P}_n D^{(j)}_{n,\tilde{P}_n^{(j)}(P)}=0$ for all $P$, including $\tilde{P}_nD^{(j)}_{n,\tilde{P}_n^{(j)}(P_0)}=0$, $j=1,2,\ldots,k$.
The  HAL-regularization bias terms $(\tilde{P}_n-P_n) D^{(j)}_{n,\tilde{P}_n^{(j)}(P_0)}$ are dominated by $(\tilde{P}_n-P_n)D^{(1)}_{P_0}$, so that some undersmoothing $\tilde{P}_n$  suffices.  

{\bf Empirical process terms $(\tilde{P}_n-P_0)D^{(j+1)}_{n,\tilde{P}_n^{(j)}(P_0)}$
are $j+1$-th order differences:}
We  will show in Section \ref{section6} that $(\tilde{P}_n-P_0)D^{(j+1)}_{n,\tilde{P}_n^{(j)}(P_0)}$ represents a  $j+1$-th order difference of $\tilde{P}_n$ and $\tilde{P}_n^{(j)}(P_0)$. 
So it appears that these higher order empirical mean terms resemble the performance of higher order $U$-statistics, by 1) being linear in $(\tilde{P}_n-P_0)\approx (P_n-P_0)$ and 2) representing a $j+1$-th order difference.  Moreover, by undersmoothing $\tilde{P}_n$ they behave as empirical means of independent mean zero random variables, not affected by the curse of dimensionality.

{\bf $k$-th exact remainder $R_n^{(k)}(P_n^{(k)}(P),\tilde{P}_n)$ is a $k+1$-th order difference:}
We prove in Sections \ref{section5} and \ref{section6} that $R_n^{(k)}(P_n^{(k)}(P),\tilde{P}_n)=R^{(1)}(\tilde{P}_n^{(1)}\ldots\tilde{P}_n^{(k)}(P),\tilde{P}_n)$ is a $k+1$-th order difference of $P$ and $\tilde{P}_n$, assuming $\tilde{P}_n D^{(j)}_{\tilde{P}_n^{(j)}(P)}=0$ for all $j$.
This remainder is the only term in the exact expansion for the $k$-th order TMLE that is affected by the curse of dimensionality: for example, suppose we use $P_n^0=\tilde{P}_n^{C_{n,cv}}$, if the dimension of $O$ is large, then the loss based dissimilarity $d_0(\tilde{P}_n^{C_{n,cv}},P_0)=O_P(n^{-2/3}(\log n)^d)$ will suffer from  a large constant in front of the rate. In these cases, having a remainder that is a $k+1$-th order difference will make a big difference in finite sample performance.
This gain is also reflected by the fact, generally, $R_n^{(k)}(\tilde{P}_n^{(k)}(P_n^0),P_0)=o_P(n^{-1/2})$ if $\pl p_n^0-p_0\pl_{\infty}=o_P(n^{-1/(2(k+1))})$.

\section{The higher order TMLE that uses empirical TMLE-updates}\label{sectionemphotmle}
Consider the case that $\tilde{P}_n^{(j)}(P)=\tilde{P}(P,\tilde{\epsilon}_n^{(j)}(P))$, where $\tilde{\epsilon}_n^{(j)}(P)$ is the solution of $\tilde{P}_n D^{(j)}_{n,\tilde{P}_n^{(j)}(P,\epsilon)}=0$, $j=1,\ldots,k$. 
Let $P_n^{(j)}(P)=\tilde{P}(P,\epsilon_n^{(j)}(P))$, where $\epsilon_n^{(j)}(P)$ is the solution of ${P}_n D^{(j)}_{n,\tilde{P}_n^{(j)}(P,\epsilon)}=0$, $j=1,\ldots,k$. In other words, $P_n^{(j)}(P)$ represents the regular empirical TMLE-updates, still using the canonical gradients $D^{(j)}_{n,P}$  that depend on the HAL-MLE $\tilde{P}_n$. 

In this section we consider the empirical $k$-th order TMLE defined by $\Psi(P_n^{(1)}\ldots P_n^{(k)}(P_n^0))$, and contrast it to the HAL-regularized $k$-th order TMLE $\Psi(\tilde{P}_n^{(1)}\ldots\tilde{P}_n^{(k)}(P_n^0))$. Specifically, we will show that it reduces the HAL-regularization bias term  by an order of difference, thereby further minimizing the need for undersmoothing the HAL-MLE.

 One way to understand this bias reduction is the following. 
 If we start the regularized higher order TMLE with initial $P_n^0=\tilde{P}_n$, then it returns $\tilde{P}_n$ itself.
This shows that the the HAL-MLE $\tilde{P}_n$ is itself a  regularized $k$-th order TMLE
satisfying  this exact expansion with the only  bias term $\sum_{j=1}^k (\tilde{P}_n-P_n)D^{(j)}_{n,\tilde{P}_n^{(j)}(P_0)}$, whose dominating term is given by $(\tilde{P}_n-P_n)D^{(1)}_{\tilde{P}_n^{(1)}(P_0)}$. On the other hand, if we plug $\tilde{P}_n$ as initial in the empirical higher order TMLE then it will sequentially target $\tilde{P}_n$ (and thus removes bias at each step) towards the target parameters $\Psi_n^{(k-1)}(P_0)$, $\ldots$, $\Psi(P_0)$, respectively. This must thus mean that it is removing bias from  $\sum_j (\tilde{P}_n-P_n)D^{(j)}_{n,\tilde{P}_n^{(j)}(P_0)}$.

In the appendix, using the above logic to understand the bias reduction due to using empirical TMLE updates at each step $j$,  we establish an exact expansion for the empirical $k$-th order TMLE presented in the next theorem.

\begin{theorem}\label{theoremempiricalhotmle}
Let $P_n^{(j)}(P)=\tilde{P}(P,\epsilon_n^{(j)}(P))$, where $\epsilon_n^{(j)}(P)$ is the solution of ${P}_n D^{(j)}_{n,\tilde{P}_n^{(j)}(P,\epsilon)}=0$, $j=1,\ldots,k$. In other words, $P_n^{(j)}(P)$ represents the regular empirical TMLE-updates, still using the canonical gradients $D^{(j)}_{n,P}$ that depend on the HAL-MLE $\tilde{P}_n$.

The empirical $k$-th order TMLE satisfies the following exact expansion:
\[
\begin{array}{l}
\Psi(P_n^{(1)}\ldots P_n^{(k)}(P_n^0))-\Psi(P_0)\\
=
\sum_{j=1}^{k}\left\{ (P_n-P_0)D^{(j)}_{n,{P}_n^{(j)}(P_0)}+R^{(j)}_n({P}_n^{(j)}(P_0),P_0)\right\}
\\
+ R^{(k)}_n({P}_n^{(k)}(P_n^0),\tilde{P}_n)-R^{(k)}_n({P}_n^{(k)}(P_0),\tilde{P}_n)\\
+\sum_{j=1}^k(\tilde{P}_n-P_n)\{D^{(j)}_{n,{P}_n^{(j)}(P_0)} -D^{(j)}_{n,{P}_n^{(j)}P_n^{(j+1)}\ldots P_n^{(k)}(P_n^0)}\}\\
+\sum_{j=1}^{k-1} \left\{R_n^{(j)}(P_n^{(j)}\ldots P_n^{(k)}P_n^0,\tilde{P}_n)-R_n^{(j)}({P}_n^{(j)}(P_0),\tilde{P}_n)\right\}\\
-\sum_{j=1}^{k-1}\left\{R_n^{(j)}(\tilde{P}_n^{(j)}P_n^{(j+1)}\ldots P_n^{(k)}P_n^0,\tilde{P}_n)-R_n^{(j)}(\tilde{P}_n^{(j)}(P_0),\tilde{P}_n)\right\}\\
\end{array}
\]
\end{theorem}
{\bf Comparison with exact expansion of HAL-regularized $k$-th order TMLE:}
Consider the case that $\tilde{P}_n^{(j)}(P)=\tilde{P}(P,\tilde{\epsilon}_n^{(j)}(P))$, where $\tilde{\epsilon}_n^{(j)}(P)$ is the solution of $\tilde{P}_n D^{(j)}_{n,\tilde{P}_n^{(j)}(P,\epsilon)}=0$, $j=1,\ldots,k$. 
Comparing this exact expansion with the exact expansion for the regularized $k$-th order TMLE, it follows that the first two lines represent the same leading expansion but with $\tilde{P}_n^{(j)}$ replaced by the empirical TMLE updates $P_n^{(j)}$. This represents an improvement since the exact remainders $R^{(j)}_n({P}_n^{(j)}(P_0),P_0)$ are now $O_P(n^{-1})$ without need for undersmoothing: i.e., $P_n^{(j)}(P_0)$ is a simple parametric standard/empirical MLE according to a correctly specified one-dimensional model.

The last three rows in the exact expansion for the empirical $k$-th order TMLE need to  be compared with the HAL-regularization bias term of the regularized $k$-th order TMLE given by $\sum_{j=1}^k (\tilde{P}_n-P_n)D^{(j)}_{n,\tilde{P}_n^{(j)}(P_0)}$. The first row of the last three rows is directly an improvement of the latter term, by subtracting from $D^{(j)}_{n,\tilde{P}_n^{(j)}(P_0)}$ its estimate. 

In addition, 
\[
\sum_{j=1}^{k-1} \left\{R_n^{(j)}(P_n^{(j)}\ldots P_n^{(k)}P_n^0,\tilde{P}_n)-R_n^{(j)}(\tilde{P}_n^{(j)}P_n^{(j+1)}\ldots P_n^{(k)}P_n^0,\tilde{P}_n)-\right\}\]
 is at minimal a third order difference. The worst case contribution from the latter sum  would come from the first term $R^{(1)}(P_n^{(1)}P,P_0)-R^{(1)}(\tilde{P}_n^{(1)}P,P_0)$, for $P=P_n^{(1)}\ldots P_n^{(k)}(P_n^0)$, which equals 
\[
\begin{array}{l}
R^{(1)}(\tilde{P}(P,\epsilon_n^{(1)}(P_n^0)),P_0)-R^{(1)}(\tilde{P}(P,\tilde{\epsilon}_n^{(1)}(P)),P_0)\\
=\frac{d}{d\xi}R^{(1)}(\tilde{P}(P,\xi),P_0)(\epsilon_n^{(1)}(P)-\tilde{\epsilon}_n^{(1)}(P)\end{array}
\]
at some intermediate point $\xi$ between $\epsilon_n^{(1)}(P)$ and $\tilde{\epsilon}_n^{(1)}(P)$.  This behaves as $d_0(\tilde{P}_n^1(P),P_0)^{1/2}(\tilde{P}_n-P_n)D^{(1)}_{\tilde{P}_n^{(1)}(P)}$, which shows  that also  this HAL-regularlization bias term represents  a third order difference instead of second order difference as with the regularized $k$-th order TMLE. This does not take into account that $R_n^{(1)}(\tilde{P}_n^{(1)}(P),P_0)$ at $P=P_n^{(2)}\ldots P_n^{(k)}(P_n^0)$ should be a better behaved remainder than at (e.g.)  $P=P_n^0$, due to the sequential targeting targeting this remainder. Therefore, in practice we expect an additional reduction due to using a $k$-th order TMLE with $k\geq 2$ with increasing benefit as $k$ increases.  
 
The last two rows in the exact expansion actually represent this third order term minus 
\[
\sum_j\{
R_n^{(j)}(\tilde{P}_n^{(j)}(P_0),\tilde{P}_n)-R_n^{(j)}(P_n^{(j)}(P_0),\tilde{P}_n)\},\] suggesting that it represents an even smaller term. 
Therefore, we conclude that the exact expansion for the empirical $k$-th order TMLE is superior to the exact expansion for the regularized $k$-th order TMLE. Based on these theoretical considerations, and our practical experience in simulations, contrary to the HAL-regularized $k$-th order TLME, we conclude that  the empirical $k$-th order TMLE does not require any undersmoothing for $\tilde{P}_n$ anymore.

{\bf ZEYI: Can you add here your graph showing that empirical second order TMLE is as good as undersmoothed HAL regularized second order TMLE, making the point I make above.}

\subsection{Understanding that HAL-regularization bias term is  a third order difference.}
The leading term in the undersmoothing term of the exact expansion for the $k$-th order TMLE is given by 
$(\tilde{P}_n-P_n)\{D^{(1)}_{n,{P}_n^{(1)}(P_0)} -D^{(1)}_{n,{P}_n^{(1)}P_n^{(2)}\ldots P_n^{(k)}(P_n^0)}\}$, which is of the  form $(\tilde{P}_n-P_n)f_n$, where $\pl f_n\pl_{P_0}=O_P(d_0^{1/2}(P_n^0,P_0))$. In this subsection we demonstrate that this is a third order difference. 
Let $\tilde{f}_n=f_n/\pl f_n\pl_{P_0}$ which is a function with a bounded $\pl \cdot\pl_{P_0}$-norm. Then,  
\[
(\tilde{P}_n-P_n)f_n=\pl f_n\pl_{P_0} (\tilde{P}_n-P_n)S_{\tilde{P}_n,\tilde{f}_n},\]
 where $S_{\tilde{P}_n,\tilde{f}_n}=\tilde{f}_n-\tilde{P}_n \tilde{f}_n$ is a score at $\tilde{P}_n$. Since $\tilde{P}_n$ is an MLE it solves a set of score equations $P_n S_{\tilde{P}_n,j}=0$ for $j=1,\ldots,J$, where this set increases as the $L_1$-norm of the HAL-MLE increases. Let $\sum_j \beta_n(j)S_{\tilde{P}_n,j}$ be the projection of $S_{\tilde{P}_n,\tilde{f}_n}$ onto the linear span of these scores in $L^2(P_0)$. 
 We have
 \[ (\tilde{P}_n-P_n)f_n=\pl f_n\pl_{P_0} (\tilde{P}_n-P_n)\left\{ S_{\tilde{P}_n,\tilde{f}_n} 
 -\sum_j\beta_n(j)S_{\tilde{P}_n,j}\right\}.\]
 One can further write $(\tilde{P}_n-P_n)=(\tilde{P}_n-P_0)-(P_n-P_0)$.
 The contribution from $(\tilde{P}_n-P_0)$ is a third order difference involving the $L^2(P_0)$ rate of $f_n$; the rate $d_0(\tilde{P}_n,P_0)^{1/2}$ and the rate of the oracle approximation of $S_{\tilde{P}_n,\tilde{f}_n}$ by the linear combination of scores solved by $\tilde{P}_n$. 
 If we set $P_n^0$ equal to an HAL-MLE, then this will be $O_P(n^{-2/3}(\log n)^{k_1})$ times the $L^2(P_0)$-norm of the oracle approximation of $S_{\tilde{P}_n,\tilde{f}_n}$. 
Thus this term can be expected to converge to zero almost as fast as $n^{-1}$. 
In fact, it is better than this third order difference since $(\tilde{P}_n-P_n)$ will behave as a nicer less biased term than $d_0^{1/2}(\tilde{P}_n,P_0)$.

\section{Definition of higher order differences, and fundamental understanding of how a TMLE targeting a remainder increases its order of difference}\label{section5}
In this section we will show that a $j$-th order remainder $R_n(P,P_0)$ can be represented as a $j+1$-th order remainder $R_n^+(P,P_0)$ plus a scaled derivative $1/j \frac{d}{dP}R_n(P,P_0)(P-P_0)$. In the next section we will apply this result at  $P=P_n^*$ and $P_0=\tilde{P}_n$, where $P_n^*$ is the TMLE of the target the parameter $P\rightarrow R_n(P,\tilde{P}_n)$ at $P=P_0$, so that this scaled directional derivative term will be equal to 0. 

\subsection{Defining a $k$-th order polynomial difference, and generalized $k$-th order difference}
\begin{definition}\label{def5a}
Consider a remainder term $R_n:{\cal M}^2\rightarrow\openr$  that, when evaluated at $(P,P_0)$, is a $k$-th order polynomial in differences of the form $H_{1,n,P_0,j}(P)-H_{1,n,P_0,j}(P_0)$ and $H_{2,n,P_0,j}(P)-H_{2,n,P_0,j}(P_0)$ of $P$ and $P_0$ defined as follows: for a $(P,P_0)\in {\cal M}^2$ 
\[
\begin{array}{l}
R_n(P,P_0)=\int \prod_{j=1}^{k}(H_{1,n,P_0,j}(P)-H_{1,n,P_0,j}(P_0)) d\mu_{1,n,P_0}\\
\hfill +
\int \prod_{j=1}^k (H_{2,n,P_0,j}(P)-H_{2,n,P_0,j}(P_0))d\mu_{2,n,P_0},\end{array}
\]
where $\mu_{1,n,P_0}$ and $\mu_{2,n,P_0}$ are measures, possibly dependent on $P_0$ and the data $P_n$.
We will refer to such a term $R_n(P,P_0)$ as a $k$-th order polynomial difference in $P$ and $P_0$.

If either $\mu_{1,n,P,P_0}$ or $\mu_{2,n,P,P_0}$ (or both)  also depend on $P$, then we refer to $R_n(P,P_0)$ as a $k$-th order difference in $P$ and $P_0$.
This definition of $k$-th order polynomial difference and $k$-th order differences generalizes immediately to sums of more than 2 integrals, each integral representing a $k$-th order difference or polynomial difference.

Finally, we define a generalized $k$-th order difference as a sum of $l$-th order differences for $l=k,\ldots,m$ for some $m>k$.
\end{definition}

Thus, for a $k$-th order polynomial difference we assume that the two measures do not depend on $P$ itself, but it can depend on the data $P_n$ and $P_0$.
We also allow that the functional mapping $P_1\rightarrow H_{1,n,P_0,j}(P_1)$ can depend on the data $P_n$
and on $P_0$, but, again, it can not be indexed by $P$ itself.
For notational convenience, we might now and then suppress the dependence of $H_{1,n,P_0,j}$ and $H_{2,n,P_0,j}$ on the data and $P_0$ and simply denote them with $H_{1,j}$ and $H_{2,j}$.
By having  $H_{1,j}=H_{1j^{\prime}}$ for $j\not =j^{\prime}$ the definition of a $k$-th order difference includes differences to a power  higher than $1$. For example, the definition of a third order polynomial difference  may include terms such as $\int (H_1(P)-H_1(P_0))^2(H_2(P)-H_2(P_0)) d\mu_{1,n,P_0}+\int (H_2(P)-H_2(P_0))^3 d\mu_{2,n,P_0}$.
However, note that the definitions of a $k$-th order polynomial difference and $k$-th order difference assume that the first term  is of same order as the second term. A generalized $k$-th order difference is a sum of a $k$-th order difference, plus additional $l$-th order differences with $l>k$.

\subsection{A $k$-th order difference is sum of $k$-th and $k+1$-th order polynomial difference}
The following lemma shows that a $k$-th order difference  is  a sum of a $k$-th and $k+1$-th order polynomial difference. 
\begin{lemma}\label{lemma6a}
Consider a $k$-th order difference $R_n:{\cal M}^2\rightarrow\openr$ and assume that its integration measures are dominated by a measure not depending on $P$. That is, in Definition \ref{def5a} of a $k$-th order polynomial difference involving a sum of  2 terms we assume that $\mu_{1,n,P,P_0}$ and $\mu_{2,n,P,P_0}$
are dominated w.r.t. a measure  $\tilde{\mu}_{1,n,P_0}$ and $\tilde{\mu}_{2,n,P_0}$ not depending on $P$.
Then,   $R_n:{\cal M}^2\rightarrow\openr$ can be decomposed into a sum of a $k$-th order polynomial difference and a $k+1$-th order polynomial difference. 
\end{lemma}
{\bf Proof:}
For simplicity, suppose that
\[
R_n(P,P_0)=\int \prod_{j=1}^{k}(H_{1,j}(P)-H_{1,j}(P_0)) d\mu_{1,n,P,P_0}.\]
Assume $d\mu_{1,n,P,P_0}=\tilde{H}_{1,n,P,P_0} d\tilde{\mu}_{1,n,P_0}$.
Then,
\[
\begin{array}{l}
R_n(P,P_0)=\int \prod_{j=1}^{k}(H_{1,j}(P)-H_{1,j}(P_0)) \tilde{H}_{1,n,P,P_0}d\tilde{\mu}_{1,n,P_0}\\
=\int \prod_{j=1}^{k}(H_{1,j}(P)-H_{1,j}(P_0)) \left (\tilde{H}_{1,n,P,P_0}-\tilde{H}_{1,n,P_0,P_0}\right)
d\tilde{\mu}_{1,n,P_0}\\
+\int \prod_{j=1}^{k}(H_{1,j}(P)-H_{1,j}(P_0))\tilde{H}_{1,n,P_0,P_0}
d\tilde{\mu}_{1,n,P_0}.
\end{array}
\]
The first term is a $k+1$-th order polynomial difference, which follows by defining $H_{1,k+1}(P)=\tilde{H}_{1,n,P,P_0}$, and defining its measure as $d\tilde{\mu}_{1,n,P_0}$.
The second term is a $k$-th order polynomial difference w.r.t. measure $d\tilde{\mu}_{1,n,P_0}$.
$\Box$

\subsection{A TMLE of target estimand sets its directional derivative at the TMLE in direction of the true data distribution equal to mean zero noise}
The $k$-th order TMLE involves targeting the remainder terms by computing a TMLE of the remainder term itself. In the next subsection we determine what happens to a $k$-th order difference $R_n(P,P_0)$ when its directional derivative at $P$ in direction $P-P_0$ equals zero. Here we first make the observation that the TMLE indeed controls this directional derivative.

The next lemma provides this particular perspective on the TMLE,  i.e., that it optimizes its target under the constraint that it needs to move in the direction of the true $P_0$.

\begin{lemma}
Let $P_n^*\in {\cal M}$ be a TMLE of a pathwise differentiable target parameter $R:{\cal M}\rightarrow\openr$ with canonical gradient $D^*_{P}$ and exact remainder $R_2(P,P_0)\equiv R(P)-R(P_0)+P_0 D^*_P$,  so that $P_n D^*_{P_n^*}=0$. 

Suppose that $R()$ and $R_2()$ can be extended so that $R_2(P,P+\delta(P-P_0))$, $R(P+\delta(P-P_0))$ are well defined for small $\delta$ around $0$, and so that the extension still satisfies $R_2(P,P+\delta(P-P_0))=R(P)-R(P+\delta(P-P_0))+(P+\delta(P-P_0))D^*_P$, and  $\frac{d}{d\delta_0}R_2(P,P+\delta_0 (P-P_0))=0$.

Then, 
\[
\frac{d}{dP_n^*}R(P_n^*)(P_n^*-P_0)=(P_n-P_0)D^*_{P_n^*}.\]
That is, the directional derivative of $R$ at $P_n^*$ in direction of $(P_n^*-P_0)$ is given by $(P_n-P_0)D^*_{P_n^*}$, and thus, behaves as an empirical mean $P_nD^*_{P_0}+o_P(n^{-1/2})$ of mean zero independent identically distributed random variables, under weak regularity conditions. 

Similarly, if $P_n^*$ is a regularized TMLE in which the empirical log-likelihood $P_n L(P)$ is replaced by an HAL-regularized $\tilde{P}_n L(P)$ for an HAL-MLE $\tilde{P}_n$  so that $\tilde{P}_n D^*_{P_n^*}=0$, then, 
\[
\frac{d}{dP_n^*}R(P_n^*)(P_n^*-P_0)=(\tilde{P}_n-P_0)D^*_{P_n^*}.\]
\end{lemma}
{\bf Proof:}
Firstly, we consider the case that the model ${\cal M}$ is convex, in which case the continuous extension condition is not needed.
For a path $P_{\delta,h}=P+\delta \int_{\cdot} hdP+o(\delta)$ with score $h$, we have
$\frac{d}{d\delta_0}R(P_{\delta_0,h})=P D^*_P h$ at $\delta_0=0$. This proves that the directional derivative of $R$ at $P$ in direction $H$ is given by $\int D^*_P dH$.
Consider now a path $p_{\delta,h}=(1+\delta h)p$ with $h=(p-p_0)/p$, which corresponds with 
$p+\delta(p-p_0)$, and due to convexity of the model, this is a valid path. 
Thus, the directional derivative of $R$ at $P_n^*$ in direction $H=\int_{\cdot}(p_n^*-p_0)/p dP=(P_n^*-P_0)(\cdot)$ is given by
\[
P_n^* D^*_{P_n^*}(p_n^*-p_0)/p_n^*=(P_n^*-P_0)D^*_{P_n^*}=-P_0D^*_{P_n^*}.\]
 Due to $P_n D^*_{P_n^*}=0$, we have $-P_0 D^*_{P_n^*}=(P_n-P_0)D^*_{P_n^*}$, thereby establishing the desired result. 
Similarly, this applies to the regularized TMLE using $\tilde{P}_n$ instead of $P_n$.

Let's now consider the more general case that the model is not necessarily convex. 
We have for any pair $(P_n^*,P_0)$ that $R(P_n^*)-R(P_0)=-P_0 D^*_{P_n^*}+R_2(P_n^*,P_0)$, for the exact remainder $R_2()$. Set $P_0=P_n^*+\delta (P_n^*-P_0)$, assuming that $R_2(P,P_0)$ and $R(P_0)$ can be extended to such elements. Then, applying the identity gives
\[
\begin{array}{l}
\Psi(P_n^*)-\Psi(P_n^*+\delta(P_n^*-P_0))=-(P_n^*+\delta(P_n^*-P_0))D^*_{P_n^*}\\+R(P_n^*,P_n^*+\delta (P_n^*-P_0))\\
=-\delta (P_n^*-P_0)D^*_{P_n^*}+ R(P_n^*,P_n^*+\delta (P_n^*-P_0))\\
=\delta P_0 D^*_{P_n^*}+ R(P_n^*,P_n^*+\delta (P_n^*-P_0)).\end{array}
\]
Thus,
\[
\begin{array}{l}
-\frac{d}{d\delta_0}\Psi(P_n^*+\delta_0(P_n^*-P_0))=P_0D^*_{P_n^*}+\lim_{\delta\rightarrow 0}\delta^{-1}R(P_n^*,P_n^*+\delta (P_n^*-P_0))\\
=P_0 D^*_{P_n^*}
=-(P_n-P_0)D^*_{P_n^*}.
\end{array}
\]
This proves the result. 
$\Box$

\subsection{Representation of a $k$-th order polynomial difference as a sum of a $k+1$-th order difference and its scaled directional derivative.}

We have the following lemma providing a representation of a $k$-th order polynomial difference $R_n(P,P_0)$ as a sum of the directional derivative at $P$ in direction $P-P_0$ and a $k+1$-th order difference. At a $P=P_n^*$ being a TMLE at which the directional derivative is approximately equal to an empirical mean of mean zero independent random variables, this then shows that a $k$-th order difference $R_n(P_n^*,P_0)$ at this TMLE $P_n^*$ will behave as a $k+1$-th order difference plus a well understood empirical mean. 

\begin{lemma}\label{lemma5}
Suppose that $R_n:{\cal M}^2\rightarrow\openr$ is a $k$-th order polynomial difference, and for simplicity, consider the case that it has a single term 
$R_n(P,P_0)=\int \prod_{j=1}^k(H_{1,j}(P)-H_{1,j}(P_0)) d\mu_{1,n,P_0}$.
Then, at any $P_n^*\in {\cal M}$, we have
\begin{eqnarray*}
R_n(P_n^*,P_0)&=&\frac{1}{k}\sum_{j=1}^k \int \prod_{i\not =j}(H_{1,i}(P_n^*)-H_{1,i}(P_0))R_{2,P_n^*,H_{1,j}()}(P_n^*,P_0) d\mu_{1,n,P_0} \\
&&+k^{-1}\frac{d}{dP_n^*}R_n(P_n^*,P_0)(P_n^*-P_0)\\
&\equiv&
R_n^{+}(P_n^*,P_0)+k^{-1}\frac{d}{dP_n^*}R_n(P_n^*,P_0)(P_n^*-P_0),
\end{eqnarray*}
where
\begin{eqnarray*}
R_{2,P,H_{1,j}()}(P,P_0)&=&H_{1,j}(P)-H_{1,j}(P_0)-\frac{d}{dP}H_{1,j}(P)(P-P_0)  .
\end{eqnarray*}

More generally, if $R_n(P_n^*,P_0)$ is a $k$-th order polynomial difference, then 
$R_n(P_n^*,P_0)=R_n^{+}(P_n^*,P_0)+k^{-1}\frac{d}{dP_n^*}R_n(P_n^*,P_0)(P_n^*-P_0)$,
where $R_n^+(P_n^*,P_0)$ is defined as above. Under weak regularity conditions specified in Lemma \ref{lemma5b} we have that $R_n^+(P_n^*,P_0)$  is a generalized $k+1$-th order difference, and, specifically, a sum of a $k+1$-th and $k+2$-th order polynomial difference.
\end{lemma}
Note that $R_{2,P,H_{1,j}()}(P,P_0)$ 
 is the second order remainder for a first order Tailor expansion at $P$ of $H_{1,j}(P_0)-H_{1,j}(P)$:\[-R_{2,P,H_{1,j}()}(P,P_0)=H_{1,j}(P_0)-H_{1,j}(P)-\frac{d}{dP}H_{1,j}(P)(P_0-P).\]

{\bf Proof of Lemma \ref{lemma5}:}
Firstly, we note that by the product rule of differentiation:					
\[\begin{array}{l}\frac{d}{dP^*_n}R_n(P_n^*,P_0)(P_n^*-P_0)\\
=\sum_{j=1}^{k}\int \prod_{i\not =j}(H_{1,i}(P_n^*)-H_{1,i}(P_0))\frac{d}{dP^*}H_{1,j}(P^*_n)(P^*_n-P_0) d\mu_{1,n,P_0}.
\end{array}
\]
Now note that \[\begin{array}{l}R_n(P,P_0)=\frac{1}{k}\sum_{j=1}^k \int \prod_{i\not =j}(H_{1,i}(P)-H_{1,i}(P_0))(H_{1,j}(P)-H_{1,j}(P_0)) d\mu_{1,n,P_0}\\
=
\frac{1}{k}\sum_{j=1}^k \int \prod_{i\not =j}(H_{1,i}(P)-H_{1,i}(P_0))R_{2,P,H_{1,j}()}(P,P_0)d\mu_{1,n,P_0}\\
+\frac{1}{k}\sum_{j=1}^k \int \prod_{i\not =j}(H_{1,i}(P)-H_{1,i}(P_0))\frac{d}{dP}H_{1,j}(P)(P-P_0)   d\mu_{1,n,P_0}\\
=R_n^+(P,P_0) 
+\frac{1}{k} \frac{d}{dP}R_n(P,P_0)(P-P_0).\end{array}\]
The following Lemma \ref{lemma5b} shows that $R_n^+(P,P_0)$ is indeed a generalized $k$-th order difference, and, specifically, a sum of a $k+1$ and $k+2$-th order polynomial difference. This completes the proof. 
$\Box$

\begin{lemma}\label{lemma5b}
Under weak regularity conditions, if $R_n(P_n^*,P_0)$ is a $k$-th order polynomial difference, then   we have that ${R}_n^+(P_n^*,P_0)$ is a  sum of a $k+1$-th and $k+2$-th order polynomial difference. 
In particular, this holds if $R_{2,P,H_{1,j}}(P,P_0)$ can be represented as  \[
R_{2,P,H_{1,j}()}(P,P_0)(x)=\int \prod_{l=1}^2 (H_{1,l,j}(P)-H_{1,l,j}(P_0))(x,y)  \bar{H}_{1,j,P,P_0}(x,y)d\mu_{j,P_0,x}(y)\]
for a measure $d\mu_{j,P_0,x}(y)=I(x=y)C_{1,j,P_0}(x)+d\mu_{j,P_0}^c(y)$, some function $C_{1,j,P_0}(x)$, and a continuous measure $d\mu_{j,P_0}^c$.
In that case, \[
\begin{array}{l}
R_{2,P,H_{1,j}()}(P,P_0)=C_{1,j,P_0}(x)\prod_{l=1}^2 (H_{1,l,j}(P)-H_{1,l,j}(P_0))(x,x) \bar{H}_{1,j,P,P_0}(x,x)\\
+\int   \prod_{l=1}^2 (H_{1,l,j}(P)-H_{1,l,j}(P_0))(x,y)  \bar{H}_{1,j,P,P_0}(x,y) d\mu_{j,P_0}^c(y).
\end{array}
\]
\end{lemma}
The proof of Lemma \ref{lemma5b} is presented in the Appendix.

\subsection{Representing a generalized $k$-th order difference $R_n()$ as a sum of a unique generalized $k+1$-th order difference $R_n^+$ and the scaled directional derivative.}
A typical remainder $R_n^{(1)}(\tilde{P}_n^{(1)}(P),P_0)$ is a second order generalized difference in $P$ and $P_0$, and can be represented as a second order polynomial difference plus a generalized third order difference: say $R_n^{(1)}(\tilde{P}_n^{(1)}(P),P_0)=R_{n,2}^p(\tilde{P}_n^{(1)}(P),P_0)+R_{n,3}(\tilde{P}_n^{(1)}(P),P_0)$. Let $P_n^*$ be a TMLE so that $d/dP_n^* R_n^{(1)}(\tilde{P}_n^{(1)}(P_n^*),P_0)(P_n^*-P_0)$ reduces to a nice  empirical mean term, say $e_n$. This implies that the directional derivative of $R_{n,2}^p$ at $P_n^*$ equals $e_n$ minus the directional derivative of $R_{n,3}$.
Applying  Lemma \ref{lemma5} to $R_{n,2}^p$  yields \[
R_{n,2}^p=R_{n,2}^{p,+}+1/2 e_n-1/2\frac{d}{dP_n^*}R_{n,3}(\tilde{P}_n^{(1)}(P_n^*),P_0)(P_n^*-P_0),\]
and thereby
\[
R_n=R_{n,2}^{p,+}+R_{n,3}-1/2\frac{d}{dP_n^*}R_{n,3}(\tilde{P}_n^{(1)}(P_n^*),P_0)(P_n^*-P_0)+1/2 e_n.\]
So $R_n$ is represented by a third order difference $R_{n,2}^{p,+}$, a generalized third order difference $R_{n,3}$, and a scaled directional derivative of the latter, plus the $e_n$-noise term controlled by the TMLE $P_n^*$. The lemma below states that the directional derivative of a generalized third order difference is a generalized third order difference itself. Thus, we have succeeded in representing $R_n$ as a generalized third order difference $R_n^+$ plus an empirical mean $e_n$-noise term controlled by the TMLE.

For the $k$-th order TMLE we want to be able to iterate this process by expressing $P_n^*=\tilde{P}_n^{(2)}(P_n^0)$, replacing $P_n^0$ by $P$ so that we obtain a new representation of the left-over remainder  $R_n^+$ in  $P$ and $P_0$, and using a next TMLE $P_n^*$ that targets this representation $R_n^+$  as a parameter in $P$. So as in previous step, we want to show that this generalized third order difference in $P$ and $P_0$ equals a generalized fourth order difference in $P$ and $P_0$ plus the scaled directional  derivative of this generalized third order difference at $P=P_n^*$ in direction $P-P_0$.

  Therefore, it is crucial that we can generalize the above Lemma \ref{lemma5} (providing a representation for a $k$-th order {\em  polynomial} difference as a sum of a generalized $k+1$-th order difference and the scaled directional derivative) to a representation of a {\em generalized} $k$-th order difference as a sum of a generalized $k+1$-th order difference plus its scaled directional derivative, where the latter would be controlled by the TMLE. 
  The following lemma establishes the desired result, and provides the building blocks and reasoning towards that result.

\begin{lemma}\label{lemma6}
We have the following facts:
\begin{itemize}
\item
By Lemma \ref{lemma6a} a $k$-th order difference $R_n:{\cal M}^2\rightarrow\openr$ can be decomposed into a sum of a $k$-th order polynomial difference and a $k+1$-th order polynomial difference. 
\item As a consequence, a generalized $k$-th order difference can be decomposed as a sum of a $k$-th order polynomial difference and $l$-th order polynomial differences, $l=k+1,\ldots,m$ for some $m>k$. The $k$-th order polynomial difference component is the natural target for the $R_n^+$-operation in the previous lemma.
\item Under a weak regularity condition, we have the following: If $R_n:{\cal M}^2\rightarrow\openr$ is a $k$-th order polynomial difference, then $\dot{R}_n:{\cal M}^2\rightarrow\openr$ defined by $\dot{R}_n(P,P_0)\equiv
\frac{d}{dP}R_n(P,P_0)(P-P_0)$ is a sum of a $k$-th and $k+1$-th order polynomial difference. 

In particular,  given $R_n(P,P_0)$ is defined in terms of $H_{1,j}$ and $H_{2,j}$ as in our definition of a $k$-th order polynomial difference, it suffices to assume that $\frac{d}{dP}H_{1,j}(P)(P-P_0)(x)=K_{1,j,P,P_0}(x)(p-p_0)(x)+\int K_{1,j,P,P_0,x}(o)(p-p_0)(o)d\mu(o)$ for certain functions $K_{1,j,P,P_0}(x)$ and kernel $(x,o)\rightarrow K_{1,j,P,P_0}(x,o)$, and similarly for $H_{2,j}$.
\item As a consequence, if $R_n:{\cal M}^2\rightarrow\openr$ is a generalized $k$-th order difference, then $\dot{R}_n:{\cal M}^2\rightarrow\openr$ is a generalized $k$-th order difference: that is, in particular, it can be written as sum of $l$-th order polynomial differences, $l=k,\ldots,m$ for some $m>k$.
\item
Lemma \ref{lemma5} established: 
If $R_n(P_n^*,P_0)$ is a $k$-th order polynomial difference, then 
$R_n(P_n^*,P_0)=R_n^{+}(P_n^*,P_0)+k^{-1}\frac{d}{dP_n^*}R_n(P_n^*,P_0)(P_n^*-P_0)$, where
$R_n^{+}$ is the uniquely defined derivative reduction of $R_n$, and $R_n^+$ is a sum of a $k+1$-th and $k+2$-th  order polynomial difference. 
\item  Suppose now that $R_n:{\cal M}^2\rightarrow\openr$ (is a generalized $k$-th order difference and thus) equals a sum of a $k$-th order polynomial difference $R_{n,k}^p:{\cal M}^2\rightarrow\openr$ and a generalized $k+1$-th order difference $R_{n,k+1}:{\cal M}^2\rightarrow\openr$. Let $R_{n,k}^{p,+}:{\cal M}^2\rightarrow\openr$ be the derivative reduction of $R_{n,k}^p$ defined by Lemma \ref{lemma5}. Then, it follows that
\[
\begin{array}{l}
R_n(P_n^*,P_0)=R_{n,k}^{p,+}(P_n^*,P_0)+R_{n,k+1}(P_n^*,P_0)+
k^{-1}\frac{d}{dP_n^*}R_n(P_n^*,P_0)(P_n^*-P_0)\\
\hfill -k^{-1}\frac{d}{dP_n^*}R_{n,k+1}(P_n^*,P_0)(P_n^*-P_0).\end{array}
\]
By the above stated facts, $R_{n,k}^{p,+}(P,P_0)$ is  a sum of a $k+1$ and $k+2$-th order polynomial difference; $k^{-1}d/dP R_{n,k+1}(P,P_0)(P-P_0)$ is a sum of a $k+1$ and $k+2$-th order polynomial difference.  Therefore, it follows that 
\[
R_n=R_n^{+}+\frac{1}{k}\dot{R}_n,\]
where $R_n^{+}:{\cal M}^2\rightarrow\openr$  is a generalized $k+1$-th order difference.
\end{itemize}
\end{lemma}
{\bf Proof of Lemma  \ref{lemma6}:}
The third statement is shown by Lemma \ref{lemma6b} in the Appendix. The last statement is shown by Lemma \ref{lemma6c} in the Appendix, but is also straightforwardly implied by previous facts.
$\Box$

\subsection{Defining the generalized $k+1$-th order difference $R_n^+$ of a generalized $k$-th order difference $R_n$.}
The previous lemma allows us now to generalize the definition of $R_n^+$. 
\begin{definition}
The previous lemma, defines, for a  given generalized $k$-th order difference $R_n:{\cal M}^2\rightarrow\openr$, a unique
generalized $k+1$-th order difference $R_n^+:{\cal M}^2\rightarrow\openr$ such that $R_n=R_n^++\frac{1}{k}\dot{R}_n$, and $\dot{R}_n:{\cal M}^2\rightarrow\openr$ was previously defined by $\dot{R}_n(P,P_0)=\frac{d}{dP}R_n(P,P_0)(P-P_0)$.
\end{definition}

\section{Establishing that the $k$-th exact remainder $R_n^{(k)}(P_n^{(k)}(P),\tilde{P}_n)$ at HAL-MLE  is a generalized $k+1$-th order difference}\label{section6}
The advantage of applying the $R_n=R_n^+ +\dot{R}_n$ decomposition,  presented in the previous Lemma \ref{lemma6}, to $P\rightarrow R_n^{(j)}(\tilde{P}_n^{(j+1)}(P),\tilde{P}_n)$, instead of $P\rightarrow R_n^{(j)}(\tilde{P}_n^{(j+1)}(P),P_0)$,  is that the directional derivative $\frac{d}{dP_n^*}R_n^{(j)}(\tilde{P}_n^{(j+1)}(P),\tilde{P}_n)(P_n^*-\tilde{P}_n)$ at $P_n^*=\tilde{P}_n^{(j)}(P)$ is exactly equal to $\tilde{P}_n D^{(j+1)}_{n,P_n^*}=0$. This was our motivation to presenting the exact remainder $R_n^{(k)}(\tilde{P}_n^{(k)}(P),P_0)$ in terms of $R_n^{(k)}(\tilde{P}_n^{(k)}(\cdot),\tilde{P}_n)$ (Lemma \ref{lemma3}).

Since the  $R_n^+$-transformation of a generalized $j$-th order difference $R_n$ is itself a member of the family of generalized $j+1$-th order differences, at a next step, one can apply this same operation to $R_n^+$ (to obtain $R_n^{+,+}$), and, thereby generate a generalized higher and higher order difference, up till the scaled directional derivative terms controlled by the sequentially defined TMLEs. Specifically,  given $R_n^{(j)}(P,\tilde{P}_n)$ is a generalized $j+1$-th order difference, we define $R_n^{(j+1)}(P,\tilde{P}_n)=
R_n^{(j)}(\tilde{P}_n^{(j+1)}(P),\tilde{P}_n)$ and recognize that this is now a $j+1$-th order difference in $P$ and $\tilde{P}_n$, using that $\tilde{P}_n^{(j+1)}(\tilde{P}_n)=\tilde{P}_n$.
Suppose now that we define a TMLE $P_n^*$ that targets $P\rightarrow R_n^{(j+1)}(P,\tilde{P}_n)$ so that the directional derivative at $P=P_n^*$ in direction $P_n^*-\tilde{P}_n$ equals zero. 
Then, $R_n^{(j+1)}(P_n^*,\tilde{P}_n)= R_n^{(j+1),+}(P_n^*,\tilde{P}_n)$ is a generalized $j+2$-th order difference in $P_n^*$ and $\tilde{P}_n$. These steps  can then be iterated till the desired $R_n^{(k)}(\tilde{P}_n^{(k)}(P),\tilde{P}_n)$.

\subsection{$R_n^{(j)}(\tilde{P}_n^{(j)}(P),\tilde{P}_n)$ is a generalized $j+1$-th order difference} 
The above arguments are formalized by the following lemma.
\begin{lemma}\label{lemmaa2}
Suppose  for a given $j$, we have $R_n^{(j)}(\tilde{P}_n^{(j)}(P),\tilde{P}_n)$ is a generalized $j+1$-th order difference.
 Let $P_n^{*}=\tilde{P}_n^{(j+1)}(P)$ be a TMLE of $\Psi_n^{(j)}(P_0)=\Psi_n^{(j-1)}(\tilde{P}_n^{(j)}(P_0))$ so that $\tilde{P}_n D^{(j+1)}_{P_n^*}\approx 0$. Then, we  have
 \[
\frac{d}{dP_n^*}\Psi_n^{(j-1)}(\tilde{P}_n^{(j)}(P_n^*))(P_n^*-\tilde{P}_n)=\tilde{P}_n D^{(j+1)}_{n,P_n^*}\approx 0,\]
 and
\[ \frac{d}{dP_n^*} R_n^{(j)}(\tilde{P}_n^{(j)}(P_n^*),\tilde{P}_n)(P_n^*-\tilde{P}_n)=\tilde{P}_n D^{(j+1)}_{n,P_n^*}\approx 0.\]
 By Lemma \ref{lemma6}, this gives
 \[
 R_n^{(j)}(\tilde{P}_n^{(j)}(P_n^*),\tilde{P}_n)={R}_n^{(j),+}(\tilde{P}_n^{(j)}(P_n^*),\tilde{P}_n)+ \frac{1}{j+1}\tilde{P}_n D^{(j+1)}_{n,P_n^*},\]
 where ${R}_n^{(j),+}(\tilde{P}_n^{(j)}(P_n^*),\tilde{P}_n)$ is  generalized $j+2$-th order difference.
By Lemma \ref{lemma2} 
\[
\begin{array}{l}
R_n^{(j+1)}(\tilde{P}_n^{(j+1)}(P),\tilde{P}_n)=R_n^{(j)}(\tilde{P}_n^{(j)}(\tilde{P}_n^{(j+1)}(P) ) ,\tilde{P}_n)
-\tilde{P}_nD^{(j)}_{n,\tilde{P}_n^{(j)}\tilde{P}_n^{(j+1)}(P)}\\
\hfill +
\tilde{P}_nD^{(j+1)}_{\tilde{P}_n^{(j+1)}(P)}.
\end{array}
\]
so that 
\[
\begin{array}{l}
R_n^{(j+1)}(\tilde{P}_n^{(j+1)}(P),\tilde{P}_n) 
+\tilde{P}_nD^{(j)}_{n,\tilde{P}_n^{(j)}\tilde{P}_n^{(j+1)}(P)}-
\tilde{P}_nD^{(j+1)}_{\tilde{P}_n^{(j+1)}(P)}\\
={R}_n^{(j),+}(\tilde{P}_n^{(j)}\tilde{P}_n^{(j+1)}(P),\tilde{P}_n)
+ \frac{1}{j+1}\tilde{P}_n D^{(j+1)}_{n,\tilde{P}_n^{(j+1)}(P)},
\end{array}
\]
and thus
\[
\begin{array}{l}
R_n^{(j+1)}(P_n^*,\tilde{P}_n)={R}_n^{(j),+}(\tilde{P}_n^{(j)}(P_n^*),\tilde{P}_n)
-\tilde{P}_nD^{(j)}_{n,\tilde{P}_n^{(j)}\tilde{P}_n^{(j+1)}(P)}+
\tilde{P}_nD^{(j+1)}_{\tilde{P}_n^{(j+1)}(P)}\\
\hfill 
+ \frac{1}{j+1}\tilde{P}_n D^{(j+1)}_{n,\tilde{P}_n^{(j+1)}(P)}.
\end{array}
\]
Thus, if the HAL-regularized TMLE scores are equal to zero, and $P_n^*=\tilde{P}_n^{(j+1)}(P)$, then
\[
\begin{array}{l}
R_n^{(j)}(\tilde{P}_n^{(j)}(P_n^*),\tilde{P}_n)=
R_n^{(j+1)}(P_n^*,\tilde{P}_n)=R_n^{(j),+}(\tilde{P}_n^{(j)}(P_n^*),\tilde{P}_n)\end{array}
\]
is a generalized $j+2$-th order difference.
By Lemma \ref{lemma2} this  then also shows that $R^{(1)}(\tilde{P}_n^{(1)}\ldots \tilde{P}_n^{(j+1)}(P),\tilde{P}_n)$ is a $j+2$-th order difference.
It follows that $R_n^{(j+1)}(\tilde{P}_n^{(j+1)}(P),\tilde{P}_n)$ is also a $j+2$-th order difference in $P$ and $\tilde{P}_n$. 

Therefore, by induction, if the HAL-regularized TMLE score equations $\tilde{P}_n D^{(j)}_{n,\tilde{P}_n^{(j)}(P)}=0$ are solved exactly,  this proves that if $R^{(1)}(\tilde{P}_n^{(1)}(P),\tilde{P}_n)$ is a generalized second order difference of $P$ and $\tilde{P}_n$, then \[
R_n^{(j)}(\tilde{P}_n^{(j)}(\tilde{P}_n^{(j+1)}(P),\tilde{P}_n)=R^{(1)}(\tilde{P}_n^{(1)}\ldots\tilde{P}_n^{(j)}(\tilde{P}_n^{(j+1)}(P)),\tilde{P}_n)\]
 is a generalized $j+2$-th order difference of $P$ and $\tilde{P}_n$, for $j=2,\ldots,k$. 
 In particular, $R_n^{(k)}(P_n^{(k)}(P),\tilde{P}_n)$ is a generalized $k+1$-th order difference of $P$ and $\tilde{P}_n$.
   \end{lemma}
   The statements in the lemma also represent the proof of these statements. 
  In the Appendix we show  that the propagation of the HAL-regularized TMLE score values, in case they are not exactly equal to zero, during the subsequent TMLE derivative reductions (i.e., $R_n=R_n^++\dot{R}_n$-operation) is controlled and generally bounded by the original TMLE score values.

 \subsection{HAL-MLE empirical process terms of $j$-th order efficient influence function represents a generalized $j$-th order difference.}
 A remaining question is the order of difference of the HAL-empirical process terms
$(\tilde{P}_n- P_0)D^{(j) }_{n,\tilde{P}_n^{(j)}(P_0)}$ of $\tilde{P}_n$ and $P_0$, even though we will show that it approaches $O_P(n^{-1})$ under an undersmoothing condition on $\tilde{P}_n$. One would conjecture that it is a $j$-th order difference. This would then teach us that even when using a non-undersmoothed HAL-MLE this term will become negligible for small $j$. 
The next lemma establishes this. 

\begin{lemma}\label{lemmahalempprocess}
For $j=2,\ldots,k$, we have
\[
\begin{array}{l}
(\tilde{P}_n-P_0)D^{(j)}_{n,\tilde{P}_n^{(j)}(P_0)}=
R_n^{(j)}(\tilde{P}_n^{(j)}(P_0),\tilde{P}_n)
-R_n^{(j-1)}(\tilde{P}_n^{(j-1)}(P_0),\tilde{P}_n)\\
\hfill -R_n^{(j)}(\tilde{P}_n^{(j)}(P_0),P_0)
+
\tilde{P}_n D^{(j-1)}_{n,\tilde{P}_n^{(j-1)}(P_0)}.
\end{array}
\]
In particular, for $j=2,\ldots,k$,
\[
\begin{array}{l}
(\tilde{P}_n-P_0)D^{(j)}_{n,\tilde{P}_n^{(j)}(P_0)} +R_n^{(j)}(\tilde{P}_n^{(j)}(P_0),P_0)    =
\tilde{P}_n D^{(j-1)}_{n,\tilde{P}_n^{(j-1)}(P_0)}\\
\hfill +R_n^{(j)}(\tilde{P}_n^{(j)}(P_0),\tilde{P}_n)
-R_n^{(j-1)}(\tilde{P}_n^{(j-1)}(P_0),\tilde{P}_n)
.
\end{array}
\]
Thus, if the TMLE $\tilde{P}_n^{(j-1)}(P_0)$ solves it target equation $\tilde{P}_n D^{(j-1)}_{n,\tilde{P}_n^{(j-1)}(P_0)}=0$, then the HAL-empirical process term behaves as difference of $j+1$-th and $j$-th order remainder of $P_0$ and $\tilde{P}_n$, up till a perfect remainder $R_n^{(j)}(\tilde{P}_n^{(j)}(P_0),P_0)$:
\[
\begin{array}{l}
(\tilde{P}_n-P_0)D^{(j)}_{n,\tilde{P}_n^{(j)}(P_0)}=
R_n^{(j)}(\tilde{P}_n^{(j)}(P_0),\tilde{P}_n)
-R_n^{(j-1)}(\tilde{P}_n^{(j-1)}(P_0),\tilde{P}_n)\\
\hfill -R_n^{(j)}(\tilde{P}_n^{(j)}(P_0),P_0).
\end{array}
\]
\end{lemma}
{\bf Proof:}
We have
 \[
\begin{array}{l}
\Psi_n^{(j-1)}(\tilde{P}_n^{(j)}(P_0))-\Psi_n^{(j-1)}(P_0)=
\Psi_n^{(j-1)}(\tilde{P}_n^{(j)}(P_0))-\Psi_n^{(j-1)}(\tilde{P}_n)\\
\hfill -(\Psi_n^{(j-1)}(P_0)-\Psi_n^{(j-1)}(\tilde{P}_n))\\
=-\tilde{P}_n D^{(j)}_{n,\tilde{P}_n^{(j)}(P_0)}+R_n^{(j)}(\tilde{P}_n^{(j)}(P_0),\tilde{P}_n)-
(\Psi_n^{(j-2)}(\tilde{P}_n^{(j-1)}(P_0))-\Psi_n^{(j-2)}(\tilde{P}_n) )\\
=-\tilde{P}_n D^{(j)}_{n,\tilde{P}_n^{(j)}(P_0)}+R_n^{(j)}(\tilde{P}_n^{(j)}(P_0),\tilde{P}_n)
+\tilde{P}_n D^{(j-1)}_{n,\tilde{P}_n^{(j-1)}(P_0)}\\
\hfill -R_n^{(j-1)}(\tilde{P}_n^{(j-1)}(P_0),\tilde{P}_n).
\end{array}
\]
Thus, for $j=2,\ldots,k$
\[
\begin{array}{l}
\Psi_n^{(j-1)}(\tilde{P}_n^{(j)}(P_0))-\Psi_n^{(j-1)}(P_0)=R_n^{(j)}(\tilde{P}_n^{(j)}(P_0),\tilde{P}_n)
-R_n^{(j-1)}(\tilde{P}_n^{(j-1)}(P_0),\tilde{P}_n)\\
\hfill +
\tilde{P}_n D^{(j-1)}_{n,\tilde{P}_n^{(j-1)}(P_0)} -\tilde{P}_n D^{(j)}_{n,\tilde{P}_n^{(j)}(P_0)}.
\end{array}
\]
We also have 
\[
\Psi_n^{(j-1)}(\tilde{P}_n^{(j)}(P_0))-\Psi_n^{(j-1)}(P_0)=-P_0D^{(j)}_{n,\tilde{P}_n^{(j)}(P_0)}+R_n^{(j)}(\tilde{P}_n^{(j)}(P_0),P_0).\]
Setting the expressions equal and solving for $-P_0 D^{(j)}_{n,\tilde{P}_n^{(j)}(P_0)}$ yields
\[
\begin{array}{l}
-P_0D^{(j)}_{n,\tilde{P}_n^{(j)}(P_0)}=
-R_n^{(j)}(\tilde{P}_n^{(j)}(P_0),P_0) +
R_n^{(j)}(\tilde{P}_n^{(j)}(P_0),\tilde{P}_n)\\
-R_n^{(j-1)}(\tilde{P}_n^{(j-1)}(P_0),\tilde{P}_n)
+
\tilde{P}_n D^{(j-1}_{n,\tilde{P}_n^{(j-1)}(P_0)} -\tilde{P}_n D^{(j)}_{n,\tilde{P}_n^{(j)}(P_0)}.
\end{array}
\]
Thus,
\[
\begin{array}{l}
(\tilde{P}_n-P_0)D^{(j)}_{n,\tilde{P}_n^{(j)}(P_0)}=
R_n^{(j)}(\tilde{P}_n^{(j)}(P_0),\tilde{P}_n)
-R_n^{(j-1)}(\tilde{P}_n^{(j-1)}(P_0),\tilde{P}_n)\\
\hfill -R_n^{(j)}(\tilde{P}_n^{(j)}(P_0),P_0)
+
\tilde{P}_n D^{(j-1)}_{n,\tilde{P}_n^{(j-1)}(P_0)}.
\end{array}
\]
This proves the lemma. $\Box$


\section{Higher order inference for $k$-th order TMLE}\label{section8}
The following also applies to the empirical $k$-th order TMLE, with the advantage that the HAL-regularization bias term $E_n$ is of smaller order than for the HAL-regularized $k$-th order TMLE.
Let $E_n\equiv -\sum_{j=1}^{k}(\tilde{P}_n-P_n) D^{(j)}_{n,\tilde{P}_n^{(j)}(P_0)}$ be the HAL-regularization bias term due to  not exactly solving the empirical efficient influence curve equations for the $k$ TMLEs $\tilde{P}_n^{(j)}(P_0)$ using as initial $P_0$ of $\Psi_n^{(j)}(P_0)$, $j=1,\ldots,k$. As discussed earlier, under some undersmoothing this term is negligible.
We also showed that  under some undersmoothing $\sum_{j=1}^{k-1}R_n^{(j)}(P_n^{(j)}(P_0),P_0)=O_P(n^{-1})$ is negligible. Even when not controlled at this rate $n^{-1})$, but at a lower rate,  it concerns a one-dimensional approximation, so that this term is not affected by curse of dimensionality.

In addition, Lemma \ref{lemmaa2}  showed that $R_n^{(k)}(\tilde{P}_n^{(k)}(P),\tilde{P}_n)$ is a $k+1$-th order difference under the assumption that $\tilde{P}_n D^{(j)}_{\tilde{P}_n^{(j)}(P_0)}=0$ for $j=1,\ldots,k$ (and small otherwise). 
 Therefore, the exact expansion of Theorem \ref{theorem2} for the HAL-regularized $k$-th order TMLE $P_n^{k,(1),*}$ provides  the following higher order approximation:
\[
\Psi(P_n^{k,(1),*})-\Psi(P_0)=
\sum_{j=1}^{k} (P_n-P_0)D^{(j)}_{n,P_n^{(j)}(P_0)}+r_k(n),
\]
where $r_k(n)$ represents a term only affected by the curse of dimensionality through  a $k+1$-th order difference $R_n^{(k)}(\tilde{P}_n^{(k)}(P_n^0),\tilde{P}_n)-R_n^{(k)}(\tilde{P}_n^{(k)}(P_0),\tilde{P}_n)$ (in addition to  some small terms close to $O_P(n^{-1})$-terms).
Let \[
\bar{D}^{k}_n\equiv \sum_{j=1}^{k} D^{(j)}_{n,{P}_n^{k,(j),*}}.\] This represents the estimate of the influence curve $\bar{D}^k_{n,P_0}\equiv \sum_{j=1}^{k} D^{(j)}_{n,P_n^{(j)}(P_0)}$.
The leading term in this estimated influence curve is $D^{(1)}_{P_n^{k,(1),*}}$, since the the other ones are higher order differences themselves, so that $\bar{D}^{k}_n$ will converge to the efficient influence curve $D^{(1)}_{P_0}$.

Let
\[
\sigma^2_n\equiv \frac{1}{n}\sum_{i=1}^n \{\bar{D}^k_n(O_i)\}^2\]
be the sample variance of this estimated influence curve $\bar{D}^k_n(O_i)$, $i=1,\ldots,n$.
This estimate of the variance takes into account the higher order expansion. 
A corresponding $0.95$-confidence interval is then given by
\[
\Psi(P_n^{k,(1),*})\pm 1.96 \sigma_n/n^{1/2}.\]
It is clear that this is a confidence interval that takes into account the exact expansion  up till its $k+1$-th order remainder. 

\paragraph{Remark: Robust estimation of variance of higher order TMLE}
For the sake of the variance estimation $\sigma_n$ we recommend using a $P_n^0$ different from the HAL-MLE  $\tilde{P}_n^{C_{n,v}}$ so that the empirical variances of $D^{(j)}_{n}$ are not underestimated due to the TMLE updates being close to $\tilde{P}_n$. In this manner, we do a better job estimating the desired efficient influence curve $D^{(j)}_{n,\tilde{P}_n^{(j)}(P_0)}$ (which has a difference $\tilde{P}_n-\tilde{P}_n^{(j)}(P_0)$, while $\tilde{P}_n-\tilde{P}_n^{(j)}(P_n^{(j-1),*})$ will be non reflective of this difference when $P_n^0=\tilde{P}_n^{C_{n,cv}}$).
Alternatively, and generally  preferable, one computes a $V$-fold cross-validated variance $\sigma^2_{n,cv}$ obtained by computing the sample variance over the validation sample  of an estimator of $\bar{D}^k_{n,P_0}$ based on plugging in estimators of $P_0$ based on the training sample, averaged across the $V$ sample splits.

\section{Algebra for sequential computation of higher order canonical gradients}\label{section9}


\subsection{Sequential computation of the  canonical gradients in terms of adjoint of score operator}
Suppose we already computed the canonical gradient  $D^{(1)}_P$ of $\Psi$.
That is, we have for all paths $P_{\delta_0,h}$ with score $h$
\[
\frac{d}{d\delta_0}\Psi(P_{\delta_0,h})=P D^{(1)}_P h.\]
Consider now the next parameter $\Psi_n^{(1)}(P)=\Psi(\tilde{P}_n^{(1)}(P))$.
We have
\begin{eqnarray*}
\frac{d}{d\delta_0}\Psi_n^{(1)}(P_{\delta_0,h})&=&\frac{d}{d\delta_0}\Psi(\tilde{P}_n^{(1)}(P_{\delta_0,h})).
\end{eqnarray*}
Consider the score operator $A_{n,P}^{(1)}: T_P\subset L^2_0(P)\rightarrow T_{\tilde{P}_n^{(1)}(P)}\subset L^2_0(\tilde{P}_n^{(1)}(P))$ defined by $A_{n,P}^{(1)}(h)=\frac{d}{d\delta_0}\log \tilde{p}_n^{(1)}(p_{\delta_0,h})$
mapping the score of the original path $p_{\delta,h}$ into the score $A_{n,P}^{(1)}(h)$ of the path $\tilde{p}_n^{(1)}(p_{\delta,h})$ through $\tilde{p}_n^{(1)}(p)$.  Here, for any $P\in {\cal M}$,  $T_P\subset L^2_0(P)$ is the  tangent space of our model ${\cal M}$ at $P$.
Since $\tilde{p}_n^{(1)}(p_{\delta,h})$ is a path through $\tilde{p}_n^{(1)}(p)$ with score $A_{n,P}^{(1)}(h)$, by definition of the pathwise derivative of $\Psi$ at $\tilde{P}_n^{(1)}(P)$, we have
\[
\frac{d}{d\delta_0}\Psi(\tilde{P}_n^{(1)}(P_{\delta_0,h}))=
\langle D^{(1)}_{\tilde{P}_n^{(1)}(P))},A_{n,P}^{(1)}(h)\rangle_{\tilde{P}_n^{(1)}(P)}.\]
Let $A_{n,P}^{(1),\top}:T_{\tilde{P}_n^{(1)}(P)}\rightarrow T_P\subset L^2_0(P)$ be the adjoint of $A_{n,P}$ defined by, for any $h\in T_P$ and $V\in T_{\tilde{P}_n^{(1)}(P))}$, 
\[
\langle A_{n,P}^{(1)}(h),V\rangle_{\tilde{P}_n^{(1)}(P)}=\langle h,A_{n,P}^{(1),\top}(V)\rangle_P.\]
Then, 
\[
\frac{d}{d\delta_0}\Psi(\tilde{P}_n^{(1)}(P_{\delta_0,h})=
\langle A_{n,P}^{(1),\top}D^{(1)}_{\tilde{P}_n^{(1)}(P))},h\rangle_P.\]
This proves that the canonical gradient of $\Psi_n^{(1)}$  at $P$ is given by
 \[
 D^{(2)}_{n,P}= A_{n,P}^{(1),\top}D^{(1)}_{\tilde{P}_n^{(1)}(P)}\in T_P\subset L^2_0(P).\]
 That is, we expressed the next canonical gradient in terms of previous canonical gradient and the adjoint of the score operator $A_{n,P}^{(1)}$.

Consider now the target parameter $\Psi_n^{(2)}:{\cal M}\rightarrow\openr$ defined by $\Psi_n^{(2)}(P)=\Psi_n^{(1)}(\tilde{P}_n^{(2)}(P))$.
Let $A_{n,P}^{(2)}:T_P\subset L^2_0(P)\rightarrow T_{\tilde{P}_n^{(2)}(P)}\subset L^2_0(\tilde{P}_n^{(2)}(P))$ be defined by $A_{n,P}^{(2)}(h)=\frac{d}{d\delta_0}\log \tilde{p}_n^{(2)}(p_{\delta_0,h})$.
Since $\tilde{P}_n^{(2)}(P_{\delta,h})$ is a path through $\tilde{P}_n^{(2)}(P)$ with score $A_{n,P}^{(2)}(h)$, we have, by definition of pathwise differentiability of $\Psi_n^{(1)}$ at $\tilde{P}_n^{(2)}(P)$, 
\[
\frac{d}{d\delta_0}\Psi_n^{(1)}(\tilde{P}_n^{(2)}(P_{\delta_0,h}))=
\langle D^{(2)}_{\tilde{P}_n^{(2)}(P)},A_{n,P}^{(2)}(h)\rangle_{\tilde{P}_n^{(2)}(P)}.\]
Let $A_{n,P}^{(2),\top}:T_{\tilde{P}_n^{(2)}(P)}\rightarrow T_P$ be the adjoint of $A_{n,P}^{(2)}$.
Then, it follows that the canonical gradient of $\Psi_n^{(2)}$ at $P$ is given by
\[
D^{(3)}_{n,P}=A_{n,P}^{(2),\top}D^{(2)}_{\tilde{P}_n^{(2)}(P)}\in T_P\subset L^2_0(P).\]

In general, supposes that we computed the canonical gradient $D^{(j)}_{n,P}\in T_P$ for target parameter $\Psi_n^{(j-1)}:{\cal M}\rightarrow\openr$ and determined a corresponding $\tilde{P}_n^{(j)}(P)$ TMLE mapping targeting $\Psi_n^{(j-1)}(P_0)$. Our goal is to determine the canonical gradient of 
$\Psi_n^{(j)}(P)=\Psi_n^{(j-1)}(\tilde{P}_n^{(j)}(P))$ at $P$.
Let $A_{n,P}^{(j)}:T_P\rightarrow T_{\tilde{P}_n^{(j)}(P)}$ be defined by $A_{n,P}^{(j)}(h)=\frac{d}{d\delta_0}\log \tilde{p}_n^{(j)}(p_{\delta_0,h})$.
Let $A_{n,P}^{(j),\top}:T_{\tilde{P}_n^{(j)}(P)}\rightarrow T_P$ be its adjoint defined by $\langle A_{n,P}^{(j)}(h),v\rangle_{\tilde{P}_n^{(j)}(P)}=\langle h,A_{n,P}^{(j),\top}(v)\rangle_P$. 
By pathwise differentiability of $\Psi_n^{(j-1)}$ at $\tilde{P}_n^{(j)}(P)$, we have 
\[
\frac{d}{d\delta_0}\Psi_n^{(j-1)}(\tilde{P}_n^{(j)}(P_{\delta_0,h}))=
\langle D^{(j)}_{n,\tilde{P}_n^{(j)}(P)},A_{n,P}^{(j)}(h)\rangle_{\tilde{P}_n^{(j)}(P)}=
\langle A_{n,P}^{(j),\top}D^{(j)}_{n,\tilde{P}_n^{(j)}(P)},h\rangle_P.\]
Thus, the canonical gradient of $\Psi_n^{(j)}$ at $P$ is given by 
\[
D^{(j+1)}_{n,P}= A_{n,P}^{(j),\top}D^{(j)}_{n,\tilde{P}_n^{(j)}(P)}\in T_P.\]

So we conclude that one can compute the canonical gradients sequentially, starting with $D^{(1)}_P$ and $\tilde{P}_n^{(1)}(P)$.  
Given, one has computed $D^{(j)}_{n,P}$ and $\tilde{P}_n^{(j)}(P)$, one computes 1) the score operator
$A_{n,P}^{(j)}(h)=\frac{d}{d\delta_0}\log \tilde{p}_n^{(j)}(p_{\delta_0,h})$; 2) its adjoint $A_{n,P}^{(j),\top}:T_{\tilde{P}_n^{(j)}(P)}\rightarrow T_P$; and one computes 3) $D^{(j+1)}_{n,P}= A_{n,P}^{(j),\top}D^{(j)}_{n,\tilde{P}_n^{(j)}(P)}\in T_P$, and 4) one determines $\tilde{P}_n^{(j+1)}(P)$ based on a local least favorable path through $P$ with score $D^{(j+1)}_{n,P}$. 
The first step is the hard work. The second step is a matter of Fubini's theorem, changing order of integration, and is generally not hard. The third step is just a matter of plugging in, and reorganizing the terms so it has a convenient form. The fourth step is generally trivial since one has already determined the type of paths one uses. 

\subsection{Computation of the adjoint of score operator}
So the most  important task is to determine  the score operator $A_{n,P}^{(j)}:T_P\rightarrow T_{\tilde{P}_n^{(j)}(P)}$ and its adjoint $A_{n,P}^{(j),\top}:T_{\tilde{P}_n^{(j)}(P)}\rightarrow T_P$.
Therefore, in this subsection we study this task. We focus on nonparametric models ${\cal M}$ and thereby nonparametric tangent space $T_P({\cal M})=L^2_0(P)$ for all $P$. For models that have a non-saturated tangent space, we can first consider the case that the model is nonparametric, and compute the resulting canonical gradients $D^{(j)}_{n,P,nonp}$. Subsequently, one projects the canonical gradient $D^{(j)}_{n,P,nonp}$  for the nonparametric model onto the tangent space $T_P({\cal M})$ to obtain the desired canonical gradient $D^{(j)}_{n,P}$. Therefore, the task addressed in this subsection for the nonparametric model is equally relevant for any non-saturated model.


\begin{lemma}\label{lemmaformulacangradient} \ \nl
{\bf Setting:}
Suppose that the model ${\cal M}$ is nonparametric; $\Psi:{\cal M}\rightarrow\openr$ has canonical gradient $D^{(1)}_P$ at $P$; $\tilde{p}_n^{(1)}(p,\epsilon)=(1+\epsilon D^{(1)}_P)p$; $\tilde{p}_n^{(1)}(p)=(1+\epsilon_n^{(1)}(P) D^{(1)}_P)p$, where
$\epsilon_n^{(1)}(P)=\arg\min_{\epsilon}\tilde{P}_n L(\tilde{p}_n^{(1)}(p,\epsilon))$, and $L(P)=-\log p$.
Let $j\in \{1,\ldots,k-1\}$ be given and suppose that we already determined
$\tilde{p}_n^{(j)}(p)=(1+\epsilon_n^{(j)}(P)D_{n,P}^{(j)})p$, where $D_{n,P}^{(j)}$ is the canonical gradient of $\Psi_n^{(j-1)}:{\cal M}\rightarrow\openr$ at $P$ defined by $\Psi_n^{(j-1)}(P)=\Psi(\tilde{P}_n^{(1)}\ldots\tilde{P}_n^{(j-1)}(P))$;  $\epsilon_n^{(j)}(P)=\arg\min_{\epsilon}\tilde{P}_n L(\tilde{p}_n^{(j)}(p,\epsilon))$; 
$\tilde{p}_n^{(j)}(p,\epsilon)=(1+\epsilon D_{n,P}^{(j)})p$. We want to compute the canonical gradient $D^{(j+1)}_{n,P}$ of $\Psi_n^{(j)}$ at $P$.

{\bf Definitions:}
Accordingly,  the class of paths through $p$ is given by $p_{\delta,h}=(1+\delta h)p$, with $h\in L^2_0(P)$.
Let  $A_{\epsilon_n^{(j)},P}:L^2_0(P)\rightarrow \openr$ be defined by 
$A_{\epsilon_n^{(j)},P}(h)=\frac{d}{d\delta_0}\epsilon_n^{(j)}(p_{\delta_0,h})$ as  the pathwise derivative of $\epsilon_n^{(j)}:{\cal M}\rightarrow\openr$ at $P$. Let $D_{\epsilon_n^{(j)},P}\in L^2_0(P)$ be the canonical gradient of $\epsilon_n^{(j)}$ at $P$ so that $A_{\epsilon_n^{(j)},P}(h)=\langle D_{\epsilon_n^{(j)},P},h\rangle_P$.
In addition, let $A_{D_n^{(j)},P}:L^2_0(P)\rightarrow L^2(P)$ be defined by
 $A_{D_n^{(j)},P }(h)=\frac{d}{d\delta_0}D_{n,P_{\delta_0,h}}^{(j)}$, i.e.,  the pathwise derivative at $P$ of  $P\rightarrow D_{n,P}^{(j)}$  considered as mapping from ${\cal M}$ to $L^2(P)$. 
 Let $D_{n,p}^{(j)}=D_{n,P}^{(j)}$, but viewed as a functional of the density $p$ instead of measure $P$.
Let $c_{n,P}^{j,-1}\equiv  \left\{\tilde{P}_n \frac{D_{n,P}^{(j),2} }{(1+\epsilon_n^{(j)}(P)D_{n,P}^{(j)})^2}\right\}^{-1} $.

{\bf Score operator $A_{n,P}^{(j)}$:}
We have
\begin{eqnarray*}
A_{n,P}^{(j)}(h)&\equiv& \frac{d}{d\delta_0}\tilde{p}_n^{(j)}(p_{\delta_0,h})/\tilde{p}_n^{(j)}(p)\\
&=& h+\frac{A_{\epsilon_n^{(j)},P}(h)D_{n,P}^{(j)}p}{\tilde{p}_n^{(j)}(p)} +\epsilon_n^{(j)}(P)  A_{D_n^{(j)},P}(h) \frac{p}{\tilde{p}_n^{(j)}(p)},
\end{eqnarray*}
where
\begin{eqnarray*}
A_{\epsilon_n^{(j)},P}(h)
&=&c_{n,P}^{j,-1}P\left(  \frac{p_n}{p} A_{D_n^{(j)},P}(h)\frac{1}{(1+\epsilon_n^{(j)}(P)D_{n,P}^{(j)})^2}\right)\\
\end{eqnarray*}
Let  $A_{D_n^{(j)},P}^{\top}:L^2(P)\rightarrow L^2_0(P)$ be the adjoint of $A_{D_n^{(j)},P}:L^2_0(P)\rightarrow L^2(P)$ so that
\[
\langle A_{D_n^{(j)},P}(h),v\rangle_P=\langle h,A_{D_n^{(j)},P}^{\top}\rangle_P.\]

If $P\not =\tilde{P}_n$, the canonical gradient $D_{\epsilon_n^{(j)},P}$ of $A_{\epsilon_n^{(j)},P}:L^2_0(P)\rightarrow \openr$ is given by:
\[
D_{\epsilon_n^{(j)},P}=c_{n,P}^{j,-1} A_{D_n^{(j)},P}^{\top}\left(\frac{d\tilde{P}_n}{dP} \frac{1}{(1+\epsilon_n^{(j)}(P)D_{n,P}^{(j)})^2}\right).\]
If $P=\tilde{P}_n$, then $A_{\epsilon_n^{(j)},P}(h)=c_{n,P}^{j,-1}P A_{D_n^{(j)},P}(h)=
c_{n,P}^{j,-1}P\frac{d}{d\delta_0}D^{(j)}_{n,P_{\delta_0,h}}=-c_{n,P}^{j,-1}P D_{n,P} h$ so that  $D_{\epsilon_n^{(j)},P}=-\{\tilde{P}_nD^{(j),2}_{n,\tilde{P}_n}\}^{-1} D^{(j)}_{n,\tilde{P}_n}$.

{\bf Adjoint $A_{n,p}^{(j),\top}$ of score operator:}
Let $A_{n,P}^{(j),\top}:L^2_0(\tilde{P}_n^{(j)}(P))\rightarrow L^2_0(P)$ be the adjoint of $A_{n,P}^{(j)}:L^2_0(P)\rightarrow L^2_0(\tilde{P}_n^{(j)}(P))$.
Then,
\[
\begin{array}{l}
A_{n,P}^{(j),\top}(V)= V(1+\epsilon_n^{(j)}(P)D_{n,P}^{(j)})
 +D_{\epsilon_n^{(j)},P}E_P \{V D_{n,P}^{(j)} \}
 + \epsilon_n^{(j)}(P) A_{D_n^{(j)},P}^{\top}\left( V\right).
 \end{array}
 \]

 
 {\bf Next canonical gradient $D^{(j+1)}_{n,P}$:} 
The next canonical gradient is thus given by:
\[
\begin{array}{l}
D^{(j+1)}_{n,P}=A_{n,P}^{(j),\top}\left(D^{(j)}_{n,\tilde{P}_n^{(j)}(P)} \right)\\
 =
 D^{(j)}_{n,\tilde{P}_n^{(j)}(P)}(1+\epsilon_n^{(j)}(P)D_{n,P}^{(j)} )
 +D_{\epsilon_n^{(j)},P}E_P \{D^{(j)}_{n,\tilde{P}_n^{(j)}(P)}D^{(j)}_{n,P} \}\\
\hfill  + \epsilon_n^{(j)}(P) A_{D_n^{(j)},P}^{\top}\left( D^{(j)}_{n,\tilde{P}_n^{(j)}(P)}\right).
 \end{array}
 \]
 Finally, for a model in which the tangent space $T_P$ is not saturated, one can still first compute these canonical gradients $D^{(j)}_{n,P,nonp}$ for the nonparametric model with the above formulas, and subsequently project them onto $T_P$. That is, 
 \[
 D^{(j)}_{n,P}=\Pi(D^{(j)}_{n,P,nonp}\mid T_P).\]
\end{lemma}

Notice that indeed $A_{n,P}^{(j),\top}(V)\in L^2_0(P)$ for a $V\in L^2_0(\tilde{P}_n^{(j)}(P))$.
In addition, we can now explicitly verify that $D^{(j)}_{n,\tilde{P}_n}=0$ as predicted by our theory.
That is,  if $P=\tilde{P}_n$, we have $A_{n,P}^{(j),\top}(V)=V-\{\tilde{P}_n D_{n,\tilde{P}_n}^{(j)2}\}^{-1} D_{n,\tilde{P}_n}^{(j)} E_{\tilde{P}_n} V D_{n,\tilde{P}_n}^{(j)}$.
 If we apply this to $V=D_{n,\tilde{P}_n}^{(j)}$ as in the formula for $D^{(j+1)}_{\tilde{P}_n}$ we obtain
 \[
D^{(j+1)}_{n,\tilde{P}_n}= D^{(j)}_{n,\tilde{P}_n}-  \{\tilde{P}_nD_{n,\tilde{P}_n}^{(j)2} \}^{-1} D_{n,\tilde{P}_n}^{(j)} E_{\tilde{P}_n} \{ D_{n,\tilde{P}_n}^{(j)} \}^2=0.\]

\paragraph{Roadmap for computing next canonical gradient:}
Thus, computing the next canonical gradient is reduced to computing the adjoint  $A_{D_n^{(j)},P}^{\top}:L^2(P)\rightarrow L^2_0(P)$  of $A_{D_n^{(j)},P}:L^2_0(P)\rightarrow L^2(P)$ defined by the directional derivative $A_{D_n^{(j)}}(h)=\frac{d}{dp}D_{n,p}^{(j)}(h p)$ of $D_n^{(j)}$ at $P$. Given this adjoint $A_{D_n^{(j)}}^{\top}$ the canonical gradient $D_{\epsilon_n^{(j)},P}$ follows and thereby also $D^{(j+1)}_{n,P}$, by the above displayed formulas for $D_{\epsilon_n^{(j)},P}$ and $D^{(j+1)}_{n,P}$. From this it also follows that the sequential computation of $D^{(j)}_{n,P}$ mostly requires computing $A_{D_n^{(j)},P}$ for the sequentially defined $D^{(j)}_{n,P}$.

 \paragraph{Analogue formula when $\epsilon_n^{(j)}(P)$ is defined as solution of TMLE score equation.}
  If we define $\epsilon_n^{(j)}(P)$ as the solution of $\tilde{P}_n D^{(j)}_{n,\tilde{P}_n^{(j)}(P,\epsilon)}=0$ instead of defining it as an MLE as in the lemma above, then we have that the pathwise derivative $\frac{d}{d\delta_0}\epsilon_n^{(j)}(P_{\delta_0,h})$ is given by
 \begin{eqnarray*}
 A_{\epsilon_n^{(j)},P}(h)&=&
 d_P^{(j),-1}
\tilde{P}_n \left(  \frac{d}{d\tilde{p}_n^{(j)}(p)} D^{(j)}_{n,\tilde{p}_n^{(j)}(p)} \left(\tilde{p}_n^{(j)}(p) h\right)\right) \\
&& +d_P^{(j),-1}\epsilon_n^{(j)}(P) \tilde{P}_n \left( \frac{d}{d\tilde{p}_n^{(j)}(p)}D^{(j)}_{n,\tilde{p}_n^{(j)}(p)}\left( p \frac{d}{dp}D^{(j)}_{n,p}(hp) \right)\right ),
\end{eqnarray*}
  where
  \[
  d_P^{(j)}\equiv 
   -\left\{\frac{d}{d\epsilon_n^{(j)}(P)} \tilde{P}_n D^{(j)}_{ \tilde{P}_n^{(j)}(P,\epsilon_n^{(j)}(P))} \right\}.\]
  Let $D_{\epsilon_n^{(j)},P}$ be the canonical gradient again of $A_{\epsilon_n^{(j)},P}$ as can be determined from this expression, so that $A_{\epsilon_n^{(j)},P}(h)=\langle D_{\epsilon_n^{(j)},P},h\rangle_P$. With this modification of the definition of $D_{\epsilon_n^{(j)},P}$, we now still have that the next canonical gradient $D^{(j+1)}_{n,P}$ is given by the above general formula:
 \[
\begin{array}{l}
D^{(j+1)}_{n,P}=A_{n,P}^{(j),\top}(D^{(j)}_{n,\tilde{P}_n^{(j)}(P)} )\\
 =
 D^{(j)}_{n,\tilde{P}_n^{(j)}(P)}(1+\epsilon_n^{(j)}(P)D_{n,P}^{(j)} )
 +D_{\epsilon_n^{(j)},P}E_P \{D^{(j)}_{n,\tilde{P}_n^{(j)}(P)}D^{(j)}_{n,P} \}\\
\hfill  + \epsilon_n^{(j)}(P) A_{D_n^{(j)},P}^{\top}\left( D^{(j)}_{n,\tilde{P}_n^{(j)}(P)}\right).
 \end{array}
 \] 
 The proof of the lemma  is presented in the Appendix \ref{AppendixF}.

 \section{Example I: Nonparametric estimation of Integrated  square of density}\label{section10}

{\bf Defining the statistical estimation problem:}
  Suppose we observe $n$ i.i.d. copies of $O\sim P_0$ and let the statistical model ${\cal M}$ be locally saturated.  Suppose that ${\cal M}$ is dominated by a measure $\mu$, and let $p=dP/d\mu$ for $P\in {\cal M}$. Typically, $\mu$ is the Lebesgue measure.  Let $\Psi:{\cal M}\rightarrow\openr$ be defined by $\Psi(P)=\int p^2(o)d\mu(o)$.  
Consider the collection of paths  $\{p_{\delta,h}=(1+\delta h)p:\delta\}$ through $p$ with score  $h\in L^2_0(P)$ at $\delta=0$.
 $\Psi$ is pathwise differentiable at $P$ along this class of paths with canonical gradient $D^{(1)}_{P}=2p(O)-2\Psi(P)$. We have
   $R^{(1)}(P,P_0)\equiv \Psi(P)-\Psi(P_0)+P_0D^{(1)}_P=-\int (p-p_0)^2 d\mu$. 
   This statistical estimation problem has been viewed as s challenging theoretical estimation problem (for example, \citep{Bickel&Ritov88,Gine08}). In particular, the second order remainder is non-forgiving by causing a negative bias that is immensely affected by the curse of dimensionality, making it particular suitable for our second order TMLE. 
   
   {\bf HAL-MLE:}
   Let $\tilde{p}_n$ be an HAL-MLE of the true density $p_0$. For example, if $O$ is one-dimensional on  interval $[0,\tau]$ one could model the hazard $\lambda(o)=p(o)/\int_o^{\tau}p(x)dx=\exp(\int_{[0,\tau]}\phi_{u}(o) dF(u))$ with the exponential link applied to a linear combination of the HAL-spline basis functions $\phi_{u}(o)=I(o\geq u)$ with knot-point $u$, where $F$ represents any cadlag function with bounded variation norm $\int_{[0,\tau]}\mid dF(u)\mid$ (with its measure extended to be defined on the edges, so that it corresponds with the sectional variation norm in \citep{vanderLaan15}). Given a finite set of knot points this model is approximated by 
   $\exp(\sum_j \beta_j \phi_{u_j})$. The HAL-MLE is then obtained by maximizing the likelihood under the constraint that $\sum_j\mid \beta_j\mid\leq C$, where $C$ is the variation norm bound. If $O$ is multidimensional then one can factorize the density in conditional univariate densities and apply this same modeling strategy  to each conditional density separately, where one uses either a separate $C$ for each conditional density or one can define a single $C$ as the $L^1$-norm of the stacked vector of coefficients across the different conditional densities. The variation norm $C\geq C_{n,cv}$ represents a  tuning parameter  for the $k$-th order TMLE that measures the degree of undersmoothing used in the HAL-MLE $\tilde{P}_n$. For an implementation of such an HAL-MLE in this particular estimation problem, we refer to \citep{Cai&vanderLaan20}.
   
{\bf Local least favorable path for first order TMLE-update:}   Let $\tilde{p}_n^{(1)}(p)=(1+\epsilon_n^{(1)}(P)D^{(1)}_P)p$ be the HAL-regularized TMLE update of $p$ using the local least favorable path $\tilde{p}_n^{(1)}(p,\epsilon)=(1+\epsilon D^{(1)}_P)p$ and MLE $\epsilon_n^{(1)}(P)=\arg\max_{\epsilon}\tilde{P}_n \log \tilde{p}_n^{(1)}(p,\epsilon)$ of its unknown coefficient $\epsilon$, where we assume that the $P$ is chosen close enough to $P_0$ or $\tilde{P}_n$ so that $\epsilon_n^{(1)}(p)$ is small enough so that $\tilde{p}_n^{(1)}(p,\epsilon_n^{(1)}(p))$ is a density. This defines the second order target estimand $\Psi_n^{(1)}:{\cal M}\rightarrow\openr$ defined by $\Psi_n^{(1)}(P)=\Psi(\tilde{P}_n^{(1)}(P))$.

  To define the second order TMLE, we need to compute the  canonical gradient $D^{(2)}_{n,P}$ of $\Psi_n^{(1)}(P)=\Psi(\tilde{P}_n^{(1)}(P))$, which we obtain by applying Lemma \ref{lemmaformulacangradient}. 
  \begin{lemma}\label{lemmaex1a}
  Define the scalars
\begin{eqnarray*}
c_{n,P}^{(1)}&=&\tilde{P}_n \frac{D^{(1),2}_P}{(1+\epsilon_n^{(1)}(P)D^{(1)}_P)^2}\\
d_P&\equiv&  E_P \{D^{(1)}_{\tilde{P}_n^{(1)}(P)}D^{(1)}_{P} \} c_{n,P}^{(1),-1}.\end{eqnarray*}
  The second order canonical gradient $D^{(2)}_{n,P}$ is given by
  \[
\begin{array}{l}
D^{(2)}_{n,P}=A_{n,P}^{(1),\top}(D^{(1)}_{\tilde{P}_n^{(1)}(P)} )\\
 =
 D^{(1)}_{\tilde{P}_n^{(1)}(P)}(1+\epsilon_n^{(1)}(P)D_{P}^{(1)} )\\
 +2d_P
 \left\{ \tilde{p}_n \frac{1}{(1+\epsilon_n^{(1)}(P)D^{(1)}_P)^2} -
2p \tilde{P}_n \frac{1}{(1+\epsilon_n^{(1)}(P)D^{(1)}_P)^2} \right\}\\
-2d_P  P \left\{ \tilde{p}_n\frac{1}{(1+\epsilon_n^{(1)}(P)D^{(1)}_P)^2}\right\}\\
+4d_P \Psi(P)
\tilde{P}_n  \frac{1}{(1+\epsilon_n^{(1)}(P)D^{(1)}_P)^2} \\
  +2 \epsilon_n^{(1)}(P) p\left(D^{(1)}_{\tilde{P}_n^{(1)}(P)}-2E_P D^{(1)}_{\tilde{P}_n^{(1)}(P)}\right)\\
-2\epsilon_n^{(1)}(P) P\left\{ p\left( D^{(1)}_{\tilde{P}_n^{(1)}(P)}-2E_P D^{(1)}_{\tilde{P}_n^{(1)}(P)}\right)\right\}.  
  \end{array}
 \]
 \end{lemma}
 Note that at $P=\tilde{P}_n$ we have $d_P=1$, and \[
 \begin{array}{l}
 D^{(2)}_{n,\tilde{P}_n}=D^{(1)}_{\tilde{P}_n^{(1)}(P)} -2\tilde{p}_n
 -2\Psi(\tilde{P}_n)+4\Psi(\tilde{P}_n)\\
 =2\tilde{p}_n-2\Psi(\tilde{P}_n)-2\tilde{p}_n+2\Psi(\tilde{P}_n)=0,
 \end{array}
 \]
 as predicted by theory.

\subsection{The second order TMLE.}
Let $\tilde{p}_n^{(2)}(p,\epsilon)=(1+\epsilon D^{(2)}_{n,p})p$, and
let $\tilde{P}_n^{(2)}(P)$ be a $\tilde{P}_n$-based TMLE update based on the universal least favorable path implied by this local least favorable path so that 
$\tilde{P}_n D^{(2)}_{n,\tilde{P}_n^{(2)}(P)}=0$ exactly.

Let $P_n^0$ be an initial estimator such as HAL-MLE  $\tilde{P}_n^{C_{n,cv}}$ in which the $L_1$-norm is selected with cross-validation.
We then compute $\tilde{p}_n^{(2)}(p_n^0)$, which is the TMLE of $\Psi_n^{(1)}(P_0)$. Subsequently,  we compute the TMLE $\tilde{p}_n^{(1)}\tilde{p}_n^{(2)}(p_n^0)$ targeting $\Psi(P_0)$ that uses $\tilde{p}_n^{(2)}(p_n^0)$ as initial estimator.  If the first order TMLE-update $\epsilon_n^{(1)}(\tilde{p}_n^{(2)}(p_n^0))$ is too large to preserve the density condition, then we make a smaller step $\epsilon$, and update $p_n^0$ with $p_n^{1}=\tilde{p}_n^{(1)}(\tilde{p}_n^{(2)}(p_n^0))$, redo $\tilde{p}_n^{(2)}(p_n^1)$ with this initial estimator $p_n^1$, compute $\tilde{p}_n^{(1)}(p_n^1)$ (possibly with small fixed $\epsilon$-step again), and  and iterate this a  few times till the full MLE $\epsilon_n^{(1)}$ works. Asymptotically, with probability tending to 1, $\tilde{p}_n^{(1)}\tilde{p}_n^{(2)}(p_n^0)$ is already well defined at the first time.

The second order TMLE of $\Psi(P_0)$ is given by
$\Psi(\tilde{P}_n^{(1)}\tilde{P}_n^{(2)}(P_n^0))$. We can iterate this by replacing the initial estimator $P_n^0$ by $\tilde{P}_n^{(1)}\tilde{P}_n^{(2)}(P_n^0)$, and thereby use the iterative second order TMLE. In the sequel let $P_n^0$ either be the original initial estimator or a few times updated second order TMLE. We always have $\tilde{P}_n D^{(2)}_{\tilde{P}_n^{(2)}(P_n^0)}=0$, and, by using a few times iterated second order TMLE we also have that
  $\tilde{P}_n D^{(1)}_{\tilde{P}_n^{(1)}\tilde{P}_n^{(2)}(P_n^0)}\approx 0$ at user supplied precision. 

\subsection{Exact expansion of the second order TMLE}
By our general exact expansion for the second order TMLE we have
\[
\begin{array}{l}
\Psi(\tilde{P}_n^{(1)}\tilde{P}_n^{(2)}(P_n^0))-\Psi(P_0)=
(P_n-P_0)D^{(1)}_{\tilde{P}_n^{(1)}(P_0)}+R^{(1)}(\tilde{P}_n^{(1)}(P_0),P_0)\\
\hfill +(P_n-P_0)D^{(2)}_{\tilde{P}_n^{(2)}(P_0)}+R_n^{(2)}(\tilde{P}_n^{(2)}(P_0),P_0)\\
\hfill +R_n^{(2)}(P_n^{(2)}(P_n^0),\tilde{P}_n)-R_n^{(2)}(P_n^{(2)}(P_0),\tilde{P}_n)\\
\hfill -P_n D^{(1)}_{\tilde{P}_n^{(1)}(P_0)}-P_n D^{(2)}_{\tilde{P}_n^{(2)}(P_0)}.
\end{array}
\]
We have $R^{(1)}(\tilde{P}_n^{(1)}(P_0),P_0)=-\int (\tilde{p}_n^{(1)}(p_0)-p_0)^2 dx$.
By Lemma \ref{lemma2}, we also have that \[
\begin{array}{l}
R_n^{(2)}(\tilde{P}_n^{(2)}(P_n^0),\tilde{P}_n)=R^{(1)}(\tilde{P}_n^{(1)}\tilde{P}_n^{(2)}(P_n^0),\tilde{P}_n)-\tilde{P}_n D^{(1)}_{\tilde{P}_n^{(1)}\tilde{P}_n^{(2)}(P_n^0)}+\tilde{P}_n D^{(2)}_{\tilde{P}_n^{(2)}(P_n^0)}.
\end{array}
\]
Due to using a TMLE $\tilde{P}_n^2(P_n^0)$ that solves its target score equation $\tilde{P}_n D^{(2)}_{\tilde{P}_n^{(2)}(P_n^0)}=0$, 
and by using an iterative second order TMLE, we also have that $\tilde{P}_n D^{(1)}_{\tilde{P}_n^{(1)}\tilde{P}_n^{(2)}(P_n^0)}\approx 0$. 
So we have that  \[
\begin{array}{l}
R_n^{(2)}(\tilde{P}_n^{(2)}(P_n^0),\tilde{P}_n)=R^{(1)}(\tilde{P}_n^{(1)}\tilde{P}_n^{(2)}(P_n^0),\tilde{P}_n).\end{array}
\]
Again, by Lemma \ref{lemma2}, 
\[
\begin{array}{l}
R_n^{(2)}(\tilde{P}_n^{(2)}(P_0),\tilde{P}_n)=R^{(1)}(\tilde{P}_n^{(1)}\tilde{P}_n^{(2)}(P_0),\tilde{P}_n)-\tilde{P}_n D^{(1)}_{\tilde{P}_n^{(1)}\tilde{P}_n^{(2)}(P_0)}+\tilde{P}_n D^{(2)}_{\tilde{P}_n^{(2)}(P_0)}.
\end{array}
\]
By definition of $\tilde{P}_n^{(2)}(P)$ as a TMLE exactly solving its target score equation, we also have 
$\tilde{P}_n D^{(2)}_{\tilde{P}_n^{(2)}(P_0)}=0$.

So we obtain
\[
\begin{array}{l}
\Psi(\tilde{P}_n^{(1)}\tilde{P}_n^{(2)}(P_n^0))-\Psi(P_0)=
(P_n-P_0)D^{(1)}_{\tilde{P}_n^{(1)}(P_0)}+R^{(1)}(\tilde{P}_n^{(1)}(P_0),P_0)\\
\hfill +(P_n-P_0)D^{(2)}_{\tilde{P}_n^{(2)}(P_0)}+R_n^{(2)}(\tilde{P}_n^{(2)}(P_0),P_0)\\
\hfill+R^{(1)}(\tilde{P}_n^{(1)}\tilde{P}_n^{(2)}(P_n^0),\tilde{P}_n)-R^{(1)}(\tilde{P}_n^{(1)}\tilde{P}_n^{(2)}(P_0),\tilde{P}_n)\\
\hfill +\tilde{P}_n D^{(1)}_{\tilde{P}_n^{(1)}\tilde{P}_n^{(2)}(P_0)} 
 -P_n D^{(1)}_{\tilde{P}_n^{(1)}(P_0)}-P_n D^{(2)}_{\tilde{P}_n^{(2)}(P_0)}.
\end{array}
\]
Lemma \ref{lemmaexactepsn1}  shows that the empirical mean of the efficient scores  at $P_0$ (i.e. of one-dimensional correct parametric models) are bounded by $(\tilde{P}_n-P_n)D^{(1)}_{P_0}$ and $O_P(n^{-1})$, so that  by some undersmoothing of $\tilde{P}_n$ these are nicely controlled.


Since $R^{(1)}(P,P_0)$ and $R_n^{(2)}(P,P_0)$ are quadratic differences in $p-p_0$, and $\tilde{P}_n^{(1)}(P_0)$ and $\tilde{P}_n^{(2)}(P_0)$ are standard MLEs  in simple correct parametric models, we have \[
R^{(1)}(\tilde{P}_n^{(1)}(P_0),P_0)=O_P(\{\epsilon_n^{(1)}(P_0)\}^2)= O_P(\tilde{P}_n-P_n)^2D^{(1)}_{P_0}))+O_P(n^{-1}),\]
 and  $R_n^{(2)}(P_n^{(2)}(P_0),P_0)=O_P((\tilde{P}_n-P_n)^2D^{(2)}_{n,P_0})+O_P(n^{-1})$,  as well. So, under  undersmoothing of the HAL-MLE so that $(\tilde{P}_n-P_n)D^{(1)}_{P_0}=o_P(n^{-1/2})$ this will be $O_P(n^{-1})$  (Lemma \ref{lemmaexactepsn1}. So we conclude that the only real remainder in this exact expansion is given by $R^{(1)}(\tilde{P}_n^{(1)}\tilde{P}_n^{(2)}(P_n^0),\tilde{P}_n)-R^{(1)}(\tilde{P}_n^{(1)}\tilde{P}_n^{(2)}(P_0),\tilde{P}_n)$, whose third order structure is shown next. 

\subsection{Second order inference based on second order TMLE}
Let $D_n=D^{(1)}_{\tilde{P}_n^{(1)}\tilde{P}_n^{(2)}(P_n^0)}+D^{(2)}_{\tilde{P}_n^{(2)}(P_n^0)}$ and $\sigma^2_n=P_n D_n^2$ be the sample variance of this estimated influence curve $D_n$.
Then, $\Psi(P_n^{2,(1),*})\pm 1.96 \sigma_n/n^{1/2}$ is an $0.95$-confidence interval that takes into account the second order expansion and thereby only ignores a third order contribution.

\subsection{The exact third order remainder in the expansion of the second order TMLE}

We will now  derive a representation of the third order  remainder $R_n^{(1)}(\tilde{P}_n^{(1)}P_n^{(2)}(P),\tilde{P}_n)$.
Let  $P_n^*=\tilde{P}_n^{(2)}(P)$.
Recall 
${R}^{(1)}(\tilde{P}_n^{(1)}(P_n^*),\tilde{P}_n)=-\int (\tilde{p}_n^{(1)}(p_n^*)-\tilde{p}_n)^2 dx$. 
Since $R_n^{(1)}(P,\tilde{P}_n)$ is a second order polynomial difference in $P$ and $\tilde{P}_n$, Lemma \ref{lemma5} shows that 
\[
R^{(1)}(\tilde{P}_n^{(1)}(P_n^*),\tilde{P}_n)=R^{(1),+}(\tilde{P}_n^{(1)}(P_n^*),\tilde{P}_n)+1/2 \tilde{P}_n D^{(2)}_{n,P_n^*},\]
so that, by our definition of  $\tilde{P}_n^{(2)}(P)$ with $\tilde{P}_n D^{(2)}_{n,P_n^*}=0$, we have
\[
R^{(1)}(\tilde{P}_n^{(1)}(P_n^*),\tilde{P}_n)=R^{(1),+}(\tilde{P}_n^{(1)}(P_n^*),\tilde{P}_n).\]

The representation of $R^{(1),+}(\tilde{P}_n^{(1)}(P_n^*),\tilde{P}_n)$ is derived by assuming
 $\frac{d}{dp_n^*}\int (\tilde{p}_n^{(1)}(p_n^*)-\tilde{p}_n^{(1)}(\tilde{p}_n))^2 dx (p_n^*-\tilde{p}_n)=0$.
That is, one assumes for the sake of derivation of $R^{(1),+}$
\[
2\int (\tilde{p}_n^{(1)}(p_n^*)-\tilde{p}_n^{(1)}(\tilde{p}_n))\frac{d}{dp_n^*}\tilde{p}_n^{(1)}(p_n^*)(p_n^*-\tilde{p}_n) dx =0.\]
As in Lemma \ref{lemma5}, we have
\[
\begin{array}{l}
\int (\tilde{p}_n^{(1)}(\tilde{p}_n)-\tilde{p}_n^{(1)}(p_n^*)^2d\mu \\
=
\int (\tilde{p}_n^{(1)}(\tilde{p}_n)-\tilde{p}_n^{(1)}(p_n^*)\left(\tilde{p}_n^{(1)}(\tilde{p}_n)-\tilde{p}_n^{(1)}(p_n^*)-\frac{d}{dp_n^*}\tilde{p}_n^{(1)}(p_n^*)(\tilde{p}_n-p_n^*)\right) dx\\
\hfill +\int (\tilde{p}_n^{(1)}(\tilde{p}_n)-\tilde{p}_n^{(1)}(p_n^*))\frac{d}{dp_n^*}\tilde{p}_n^{(1)}(p_n^*)(\tilde{p}_n-p_n^*) dx\\
\equiv \int (\tilde{p}_n^{(1)}(\tilde{p}_n)-\tilde{p}_n^{(1)}(p_n^*))R_{2,p_n^*,\tilde{p}_n^{(1)}()}(\tilde{p}_n,p_n^*) dx.
\end{array}
\]
Thus, we have
\[
R^{(1),+}(\tilde{P}_n^{(1)}(P_n^*),\tilde{P}_n)=-\int (\tilde{p}_n^{(1)}(\tilde{p}_n)-\tilde{p}_n^{(1)}(p_n^*))R_{2,p_n^*,\tilde{p}_n^{(1)}()}(\tilde{p}_n,p_n^*) dx,\]
which already shows the third order difference of this term. 
More concretely, we can determine the type of third order terms included in this expression, as presented in the next lemma.

\begin{lemma}\label{lemmaexample1} For a $P_n^*=\tilde{P}_n^{(2)}(P)$, we have
\[
R^{(1)}(\tilde{P}_n^{(1)}(P_n^*),\tilde{P}_n)=-\int (\tilde{p}_n^{(1)}(\tilde{p}_n)-\tilde{p}_n^{(1)}(p_n^*))R_{2,p_n^*,\tilde{p}_n^{(1)}()}(\tilde{p}_n,p_n^*) dx.\]
Let 
\[
R_{2,p,\epsilon()}(\tilde{p}_n,p)=\epsilon_n^{(1)}(\tilde{p}_n)-\epsilon_n^{(1)}(p)-\frac{d}{dp}\epsilon_n^{(1)}(p)(\tilde{p}_n-p).\]
More concretely,
\[
\begin{array}{l}
-{R}^{(1)}(\tilde{p}_n^{(1)}(p_n^*),\tilde{p}_n)=
\int (\tilde{p}_n^{(1)}(\tilde{p}_n)-\tilde{p}_n^{(1)}(p_n^*))^2dx \\
=\int (\tilde{p}_n^{(1)}(\tilde{p}_n)-\tilde{p}_n^{(1)}(p_n^*))R_{2,p_n^*,\tilde{p}_n^{(1)}()}(\tilde{p}_n,p_n^*) dx\\
=\int (\tilde{p}_n^{(1)}(\tilde{p}_n)-\tilde{p}_n^{(1)}(p_n^*))(\epsilon_n^{(1)}(\tilde{p}_n)D^{(1)}_{\tilde{p}_n}-\epsilon_n^{(1)}(p)D^{(1)}_p)(\tilde{p}_n-p)dx\\
+\int (\tilde{p}_n^{(1)}(\tilde{p}_n)-\tilde{p}_n^{(1)}(p_n^*))D^{(1)}_p p R_{2,\epsilon_n^{(1)},p}(\tilde{p}_n,p)dx\\
-2\epsilon_n^{(1)}(p) \int (\tilde{p}_n-p)^2 dx\int (\tilde{p}_n^{(1)}(\tilde{p}_n)-\tilde{p}_n^{(1)}(p_n^*)) p dx\\
\end{array}
\]
The first two terms are third order differences and the last term is a fourth order difference. 
\end{lemma}
We can also work out the second order difference $R_{2,p,\epsilon()}(\tilde{p}_n,p)$, but the above terms are completely representative of the extra terms generated by this term. 

\subsection{Concluding remark}
We conclude that the square remainder for the first order TMLE of $\int p^2d\mu$ has been replaced by a sum of third order differences that also involve a lot of cancelation due to  integration of products of differences. That is, we did not only improve the order of difference from second to third order, but we also made the terms smaller due to involving cancelations. In fact, neither term is a square as $R^{(1)}(\tilde{P}_n^{(1)(P)},P_0)$, and thereby has a predictable sign.
Therefore, the second order TMLE can be expected to have significantly better finite sample performance, and allows for much better inference than the first order TMLE whose inference suffers from a  pure negative bias term $R^{(1)}(\tilde{P}_n^{(1)}(P_n^0),P_0)$ that would be ignored by a Wald-type confidence interval.

\section{Example II: Nonparametric estimation of the treatment specific mean}\label{section11}

\subsection{Formulation of estimation problem.}
Let
$O=(W,A,Y)\sim P_0$ and let the statistical model ${\cal M}$ be nonparametric.
Suppose that $A$ and $Y$ are binary. Let $\Psi(P)=E_PE_P(Y\mid W,A=1)$ be the target parameter.
Let $q_W$, $q$ and $g$ be the marginal density of $W$, conditional density of $Y$, given $W,A$, and conditional density of $A$, given $W$, respectively. Let $\bar{q}(W,A)=E_P(Y\mid W,A)$ and $\bar{g}(W)=g(1\mid W)$. Let $\bar{q}_1(W)=E_P(Y\mid W,A=1)$.
Let $q_{W,n}$ be  the empirical distribution of $W_1,\ldots,W_n$, and the marginal distribution $\tilde{q}_{W,n}$ of $\tilde{P}_n$ is also set equal to this empirical distribution $q_{W,n}$.  Since the marginal distribution of $W$ under $P_n^0$ and $\tilde{P}_n$ are equal, any higher order TMLE updates would not result in an update of $q_{W,n}$. 
We can decompose \[
\Psi(q_{W_n},\bar{q})-\Psi(q_{W,0},\bar{q}_0)=\Psi(q_{W,0},\bar{q})-\Psi(q_{W,0},\bar{q}_0)+
(P_n-P_0)\bar{q}_1.\]
 Therefore, for the sake of computing the necessary canonical gradients and exact remainders for the second order TMLE, we can focus on the target parameter $\Psi(\bar{q})=E_{P_0}\bar{q}(W,1)$ treating the expectation over $W$ as given. The resulting second order TMLE $\bar{q}_n^{2,(1),*}$ of $\bar{q}_0$ combined with the empirical $q_{W,n}$ provides then the desired second order TMLE $\Psi(q_{W,n},\bar{q}_n^{(2,(1),*})$ of $\Psi(P_0)=\Psi(q_{W,0},\bar{q}_0)$.

The canonical gradient of $\Psi(\bar{q})=E_{P_0}\bar{q}(1,W)$ at $P$ is given by $D^{(1)}_P(O)=A/\bar{g}(W)(Y-\bar{q}(W,A))$. The exact remainder for $\Psi$ is given by $R^{(1)}(P,P_0)\equiv \Psi_0(P)-\Psi_0(P_0)+P_0 D^{(1)}_P=P_0 (\bar{q}_1-\bar{q}_{10})\frac{\bar{g}-\bar{g}_0}{\bar{g}}$.

{\bf Class of paths:}
As paths through the conditional density $q(Y\mid W,A)$, we consider $\mbox{Logit}q_{\delta,C_1}(1\mid W,A)=\mbox{Logit}q(1\mid W,A)+\delta C_1(W,A)$ for a function $C_1$. This path has  a score $h_1(W,A,Y)=\frac{d}{d\delta_0}\log q_{\delta_0,C_1}=
C_1(W,A)(Y-q(1\mid W,A))$. Similarly, we define paths through $g$ by $\mbox{Logit}g_{\delta,C_2}(1\mid W)=\mbox{Logit}g(1\mid W)+\delta C_2$ for a function $C_2(W)$. This path has score $h_2(W,A)=\frac{d}{d\delta_0}\log g_{\delta_0,C_2}=C_2(W)(A-g(1\mid W))$.
Let $h=(h_1,h_2)$. These two paths imply a path through the data density  $p$ given by  $p_{\delta,h}=q_{W,0}g_{\delta,h_1}q_{\delta,h_2}$. The score of this path through $P$ is given by $S_{h_1,h_2}=h_1+h_2$:
\[
S_{h_1,h_2}=
C_1(W,A)(Y-q(1\mid W,A))+C_2(W)(A-g(1\mid W)).\]

\subsection{HAL-MLE}
Let $\tilde{p}_n=(\tilde{q}_{W,n},\tilde{q}_{n},\tilde{g}_n)$ be an HAL-MLE of $p_0=(q_{W,0},\bar{q}_0,\bar{g}_0)$. The HAL-MLE $\tilde{q}_{n}^C$ is computed with linear logistic regression of $Y$ on the zero order splines $I((W,A)>u_j)$ indexed by knot-points $u_j$ and constraining the $L^1$-norm of the coefficient vector by $C$. Similarly, $\tilde{g}_n^C$ is computed with logistic linear regression of $A$ on zero order splines $I(W>u_j)$, and constraining the $L^1$-norm of the coefficient vector by $C$. The $L_1$-norm bounds $(C_y,C_a)$ for these two HAL-MLEs can be data adaptively selected, but are chosen to be larger than their respective cross-validation selectors $C_{y,n,cv}$ and $C_{a,n,cv}$. Both HAL-MLEs can be computed with the R-function HAL9001 based on glmnet.
  
  \subsection{First order TMLE-update}
  We will interchange $p=(q,g)$.
Let $\tilde{p}^{(1)}(p,\epsilon)=q_W g \tilde{q}^{(1)}(p,\epsilon)$, where
\[
\mbox{Logit}\tilde{q}^{(1)}(q,g,\epsilon)(1\mid W,A)=\mbox{Logit}q(1\mid W,A)+\epsilon C_g(W,A)\]
and $C_g\equiv A/\bar{g}(W)$.
Let $\epsilon_n^{(1)}(p)=\arg\max_{\epsilon}\tilde{P}_n \log \tilde{q}^{(1)}(q,g,\epsilon)$ be the MLE of $\epsilon$ under the HAL-MLE $\tilde{P}_n$. This defines the first order TMLE-update $\tilde{q}_{n}^{(1)}(p)=\tilde{q}^{(1)}(p,\epsilon_n^{(1)}(p))$ of $q$.
The corresponding first order TMLE-update of $p$ is given by $\tilde{p}_n^{(1)}(p)=\tilde{p}^{(1)}(p,\epsilon_n^{(1)}(p))=q_Wg\tilde{q}_n^{(1)}(p)$. Notice that only $q$ is updated, while $g$ is not. 
 Let $\bar{q}_n^{(1)}(p)=\tilde{q}_n^{(1)}(p)(1\mid W,A)$ be the corresponding first order TMLE-update of $\bar{q}$. We use notation $\bar{q}_n^{(1)}(p)(1)$ for $\bar{q}_n^{(1)}(p)(W,1)$.
 Notice that $\tilde{P}_n D^{(1)}_{\tilde{P}_n^{(1)}(P)}=0$ exactly.

This now defines the first order target parameter $\Psi_n^{(1)}(P)=\Psi(\tilde{P}_n^{(1)}(P))$ defined by
\[
\Psi_n^{(1)}(P)=\Psi(\tilde{q}_n^{(1)}(P))=E_{P_0}\bar{q}_n^{(1)}(p).\]

\subsection{Second order canonical gradient.}

We need to compute the score operator $A_{n,P}^{(1)}(h)=\frac{d}{\delta_0}\log \tilde{p}_n^{(1)}(p_{\delta_0,h})=\frac{d}{d\delta_0}\tilde{p}_n^{(1)}(p_{\delta_0,h})/\tilde{p}_n^{(1)}(p)$. The numerator equals $\frac{d}{dp}\tilde{p}_n^{(1)}(p)(\frac{d}{d\delta_0}p_{\delta_0,h})=\frac{d}{dp}\tilde{p}_n^{(1)}(p)(p \bar{h})$, where $\bar{h}=h_1+h_2$.
This shows that the score operator is linear mapping in $\bar{h}=h_1+h_2\in L^2_0(P)$.
Due to factorization of $\tilde{p}_n^{(1)}(p)=q_Wg\tilde{q}_n^{(1)}(q,g)$ it follows that $A_{n,P}^{(1)}(\bar{h})=A_{n,p,1}^{(1)}(h_1)+A_{n,p,2}^{(1)}(h_2)$, where $A_{n,p,1}^{(1)}:L^2_0(q)\rightarrow L^2_0(\tilde{P}_n^{(1)}(p))$ is defined by 
\[
A_{n,p,1}^{(1)}(h_1)=\frac{\frac{d}{\delta_0}\tilde{q}_n^{(1)}(q_{\delta_0,h_1},g)}{\tilde{q}_n^{(1)}(q,g)},\]
 and $A_{n,p,2}:L^2_0(g)\rightarrow L^2_0(\tilde{P}_n^{(1)}(P))$ is given by 
 \[
 A_{n,p,2}^{(1)}(h_2)=h_2+\frac{\frac{d}{d\delta_0}\tilde{q}_n^{(1)}(q,g_{\delta_0,h_2})}{\tilde{q}_n^{(1)}(q,g)}.\]
 The following lemma establishes the precise forms of these two score operators. 
 
 \begin{lemma}\label{lemmaexample2a}
 Define \[
\begin{array}{l}
c_P^{(1)}=
\tilde{P}_n C_g^2\bar{q}_n^{(1)}(p) (1-\bar{q}_n^{(1)}(p))=P \frac{\tilde{g}_n}{\bar{g}^2} \bar{q}_n^{(1)}(p)(1-\bar{q}_n^{(1)}(p))(1).
\end{array}
\]
 We have
\[
\begin{array}{l}
A_{n,p,1}(h_1)=\frac{d}{d\delta_0}\log \tilde{q}_n^{(1)}(q_{\delta_0,h_1},g)\\
=C_1(W,A)(Y-\bar{q}_n^{(1)}(q,g)(W,A))\\
-c_P^{(1),-1}\left\{\tilde{P}_n C_g C_1 \bar{q}_n^{(1)}(q,g)(1-\bar{q}_n^{(1)}(q,g))\right\}
C_g(Y-\bar{q}_n^{(1)}(q,g) ),
\end{array}
\]
and
\[
\begin{array}{l}
A_{n,p,2}(h_2)=h_2+\frac{d}{d\delta_0}\log \tilde{q}_n^{(1)}(q,g_{\delta_0})\\
=h_2-\epsilon_n^{(1)}(q,g)\frac{A}{\bar{g}}h_{2}(W,A)  (Y-\bar{q}_n^{(1)}(q,g))\\
-c_P^{(1),-1}\left\{ \tilde{P}_n\frac{A}{\bar{g}}C_2(1-\bar{g}) \left(Y-\bar{q}_n^{(1)}(q,g)\right)\right\}
C_g(Y-\bar{q}_n^{(1)}(q,g))\\
+c_P^{(1),-1}\left\{ \tilde{P}_n C_g\epsilon_n^{(1)}(q,g) \frac{A}{\bar{g}}C_2(1-\bar{g}) \bar{q}_n^{(1)}(q,q)(1-\bar{q}_n^{(1)}(q,g))\right\} C_g(Y-\bar{q}_n^{(1)}(q,g)).
\end{array}
\]
  \end{lemma}
 
 Subsequently, as shown in Appendix, we determine the adjoints of these two score operators and apply them to $D^{(1)}_{\tilde{P}_n^{(1)}(P)}$ to obtain the desired second order canonical gradient, and, specifically, its two components corresponding with the dependence on $P$ through $q(P)$ and $g(P)$, respectively.
 \begin{lemma}\label{lemmaexample2b}
 Define
 \[
\begin{array}{l}
\tilde{c}_{\tilde{P}_n,P}^{(1)}=c_P^{(1),-1}P\frac{\bar{q}_n^{(1)}(1-\bar{q}_n^{(1)})(p)(1)}{\bar{g}}
=\left\{P \frac{\tilde{g}_n}{\bar{g}^2} \bar{q}_n^{(1)}(1-\bar{q}_n^{(1)})(p)(1) \right\}^{-1}
P\frac{\bar{q}_n^{(1)}(1-\bar{q}_n^{(1)})(p)(1)}{\bar{g}}.
\end{array}
\]
We have
\[
\begin{array}{l}
D^{(2)}_{n,P,1}=A_{n,p,1}^{\top}(D^{(1)}_{\tilde{P}_n^{(1)}(p)})\\
=\frac{A}{\bar{g}}\frac{\bar{q}_n^{(1)}(p)(1-\bar{q}_n^{(1)}(p))}{\bar{q}(1-\bar{q})}\left\{1-\tilde{c}_{\tilde{P}_n,P}^{(1)}\frac{\tilde{g}_n}{\bar{g}}\right\}(Y-\bar{q}).
\end{array}
\]
This represents the canonical gradient of $q\rightarrow \Psi_n^{(1)}(q,g)$. 
The second order canonical gradient $D^{(2)}_{n,P}$ is given by $D^{(2)}_{n,P}=D^{(2)}_{n,P,1}+D^{(2)}_{n,P,2}$, where $D^{(2)}_{n,P,2}$ is the canonical gradient of $g\rightarrow \Psi_n^{(1)}(q,g)$:
\[
\begin{array}{l}
D^{(2)}_{n,P,2}=A_{n,p,2}^{\top}(D^{(1)}_{\tilde{P}_n^{(1)}(p)})\\
=-\epsilon_n^{(1)}(p)\frac{\bar{q}_n^{(1)}(1-\bar{q}_n^{(1)})(p)}{\bar{g}^2}\left\{1-\tilde{c}_{\tilde{P}_n,P}^{(1)}\frac{\tilde{g}_n}{\bar{g}}\right\}(A-\bar{g})\\
-  \tilde{c}_{\tilde{P}_n,P}^{(1)}
 \frac{\tilde{g}_n}{\bar{g}^2} (\tilde{q}_n-\bar{q}_n^{(1)}(p)(1) )(A-\bar{g} ).
\end{array}
\]
Thus, we have
\[
\begin{array}{l}
D^{(2)}_{n,P}=\frac{A}{\bar{g}}\frac{\bar{q}_n^{(1)}(p)(1-\bar{q}_n^{(1)}(p))}{\bar{q}(1-\bar{q})}\left\{1-\tilde{c}_{\tilde{P}_n,P}^{(1)}\frac{\tilde{g}_n}{\bar{g}}\right\}(Y-\bar{q})\\
-\epsilon_n^{(1)}(p)\frac{\bar{q}_n^{(1)}(1-\bar{q}_n^{(1)})(p)}{\bar{g}^2}\left\{1-\tilde{c}_{\tilde{P}_n,P}^{(1)}\frac{\tilde{g}_n}{\bar{g}}\right\}(A-\bar{g})\\
-  \tilde{c}_{\tilde{P}_n,P}^{(1)}
 \frac{\tilde{g}_n}{\bar{g}^2} (\tilde{q}_n-\bar{q}_n^{(1)}(p)(1) )(A-\bar{g} ).
 \end{array}
 \]
 \end{lemma}
 We note that at $P=\tilde{P}_n$ we have $\bar{q}_n^{(1)}(p)=\bar{q}$; $\tilde{g}_n=\bar{g}$; $\epsilon_n^{(1)}(p)=0$; $C_{\tilde{P}_n}^{(1)}=P\bar{q}(1-\bar{q})/\bar{g}$, where the latter cancels with   $P\{\frac{1}{\bar{g}}
\bar{q}_n^{(1)}(1-\bar{q}_n^{(1)})(p)(1) \}$, and thereby 
\[
D^{(2)}_{n,\tilde{P}_n,1}=\frac{A}{\bar{g}}(Y-\bar{q})-C_g (Y-\bar{q})=0.\] 
 Note also that $D^{(2)}_{n,\tilde{P}_n,2}=0$. Thus, we indeed have that $D^{(2)}_{n,\tilde{P}_n}=0$.
 
 \subsection{The second order TMLE update of initial estimator}
 
Define
\[
\begin{array}{l}
C^y_{n,P}\equiv \frac{A}{\bar{g}}\frac{\bar{q}_n^{(1)}(p)(1-\bar{q}_n^{(1)}(p))}{\bar{q}(1-\bar{q})}\left\{1-\tilde{c}_{\tilde{P}_n,P}^{(1)}\frac{\tilde{g}_n}{\bar{g}}\right\}\\
C^a_{n,P}\equiv -\epsilon_n^{(1)}(p)\frac{\bar{q}_n^{(1)}(1-\bar{q}_n^{(1)})(p)(1)}{\bar{g}^2}\left\{1-\tilde{c}_{\tilde{P}_n,P}^{(1)}\frac{\tilde{g}_n}{\bar{g}}\right\}-  \tilde{c}_{\tilde{P}_n,P}^{(1)}
 \frac{\tilde{g}_n}{\bar{g}^2} (\tilde{q}_n-\bar{q}_n^{(1)}(p)(1) ).
\end{array}
\]
 Notice that
 \[
 D^{(2)}_{n,P}=C^y_{n,P}(W,A)(Y-\bar{q}(W,A))+C^a_{n,P}(A-\bar{g}(W)).\]
 We define a least favorable path $\{\tilde{P}_n^{(2)}(P,\epsilon):\epsilon=(\epsilon_1,\epsilon_2)\}$ through $q$ and $g$ targeting $\Psi_n^{(1)}(P_0)$:
 \[
 \begin{array}{l}
 \mbox{Logit}\bar{q}_n^{(2)}(p,\epsilon)=\mbox{Logit}\bar{q}+\epsilon_1 C^y_{n,P}\\
 \mbox{Logit}\bar{g}_n^{(2)}(p,\epsilon)=\mbox{Logit}\bar{g}+\epsilon_2 C^a_{n,P}.
 \end{array}
 \]
 Note that the resulting density \[
 \tilde{p}_n^{(2)}(p,\epsilon)(W,A,Y)=q_{W,n}\tilde{g}_n^{(2)}(p,\epsilon_2)(A\mid W)\tilde{q}_n^{(2)}(p,\epsilon_1)(Y\mid W,A)\]
 is a submodel through $p$ at $\epsilon =0$ with scores $C^y_{n,P}(W,A)(Y-\bar{q}(1,W))$ and $C^a_{n,P}(A-\bar{g}(W))$, which thus spans $D^{(2)}_{n,P}$.

\subsection{The second order TMLE for the treatment specific mean}

\begin{description}
\item[Initial estimator:] Let $p_n^0=(q_{W,n},\bar{g}_n^0,\bar{q}_n^0)$ be an initial density estimator with $q_{W,n}$ being the empirical probability density with mass $1/n$ on each $W_i$, $i=1,\ldots,n$. $\bar{g}_n^0$ and $\bar{q}_n^0$ are  initial estimators of $P(A=1\mid W)$ and $P(Y=1\mid W,A)$, respectively, such as a logistic regression based HAL-MLE using cross-validation to select the $L^1$-norm of its coefficient vector or a super learner including such HAL-MLEs.
\item[HAL-MLE:] Recall the HAL-MLE $\tilde{P}_n^C=(q_{W,n},\tilde{g}_n^{C_a},\tilde{q}_n^{C_y})$ indexed by $L_1$-norm bounds on the coefficient vectors. Let $\tilde{P}_n$ be the HAL-MLE corresponding with a data adaptive selector $C_n$.
\item[Computing clever covariates for second order TMLE update:]
The estimators $\tilde{P}_n$ and $P_n^0$ map into $C^y_{n,P_n^0}$ and $C^a_{n,P_n^0}$.  Specifically, the evaluation of these clever covariates relies on $\bar{g}_n^0$; $\tilde{g}_n$; $\bar{q}_n^0$; and the first order TMLE update $\bar{q}_n^{(1)}(\bar{g}_n^0,\bar{q}_n^0)$, which involves computing the HAL-regularized MLE $\epsilon_n^{(1)}(P_n^0)=\arg\min_{\epsilon}\tilde{P}_n L(\tilde{q}_n^{(1)}(p_n^0,\epsilon))$.
\item[Computing the first order MLE:]
Notice that $\epsilon_n^{(1)}(P_n^0)$ solves \[
\begin{array}{l}
0=\tilde{P}_n C_{g_n^0}(Y-\bar{q}_n^{(1)}(q_n^0,g_n^0,\epsilon))\\
=
P_n \frac{\tilde{g}_n}{\bar{g}_n^0} \left(\tilde{q}_n(1)-
\frac{1}{1+\exp\left (-\log \bar{q}_n^0/(1-\bar{q}_n^0)-\epsilon C_{g_n^0}\right)}\right).\end{array}
\]
This could be computed with a grid search. One might also compute $\epsilon_n^{(1)}(P_n^0)$ with   univariate logistic regression using as  outcome (proportion-valued) $\tilde{q}_n^{(1)}(W_i,1)=\tilde{P}_n(Y=1\mid W=W_i,A=1)$, off-set $\log \bar{q}_n^0/(1-\bar{q}_n^0)(W_i,1)$, and single covariate $C_{g_n^0}(W_i,1)=1/\bar{g}_n^0(W_i)$, using weights $\tilde{g}_n/\bar{g}_n^0(W_i)=\tilde{P}_n(A=1\mid W=W_i)/P_n^0(A=1\mid W=W_i)$ (our experience is that logistic regression can be unstable when using it in this manner, so that a grid search is preferable).

\item[Computing the second order TMLE update:]
 We can now compute $\tilde{P}_n^{(2)}(P_n^0,\epsilon_n^{(2)}(P_n^0))$ as
 \[
 \begin{array}{l}
 \mbox{Logit}\bar{q}_n^{(2)}(p_n^0,\epsilon_{1n}^{(2)})=\mbox{Logit}\bar{q}+\epsilon_{1n}^{(2)} C^y_{n,P_n^0}\\
 \mbox{Logit}\bar{g}_n^{(2)}(p_n^0,\epsilon_{2n}^{(2)})=\mbox{Logit}\bar{g}+\epsilon_{2n}^{(2)} C^a_{n,P_n^0},
 \end{array}
 \]
where $\epsilon_{1n}^{(2)}=\arg\max_{\epsilon} \tilde{P}_n  \log \tilde{q}_n^{(2)}(p_n^0,\epsilon)$ and
$\epsilon_{2n}^{(2)}=\arg\max_{\epsilon}\tilde{P}_n \log \tilde{g}_n^{(2)}(p_n^0),\epsilon)$.
If we use the empirical TMLE update, then the coefficient $\epsilon_{1n}^{(2)}$  is fitted with standard univariate logistic regression of $Y_i$ on $(W_i,A_i)$ 
using as off-set $\log \bar{q}_n^0/(1-\bar{q}_n^0)(W_i,A_i)$, and univariate covariate $C^y_{n,P_n^0}(W_i,A_i)$. Similarly, in that case the coefficient $\epsilon_{2n}^{(2)}$  is fitted with standard univariate logistic regression of $A_i$ on $W_i$ 
using as off-set $\log \bar{g}_n^0/(1-\bar{g}_n^0)(W_i)$, and univariate covariate $C^a_{n,P_n^0}(W_i)$. Similarly, the regularized TMLE updates can be obtained by solving the score equations, as discussed above.
This provides us with $\tilde{q}_n^{(2)}(p_n^0)$ and $\tilde{g}_n^{(2)}(p_n^0)$ while the empirical measure $q_{W,n}$ for $W$ remains unchanged. 

We could also use a iterative TMLE or universal least favorable path TMLE,  which then exactly solves $\tilde{P}_n D^{(2)}_{P_n^{(2)}(P_n^0)}=0$. The latter was carried out in our simulation study.

\item[Computing the first order TMLE update:]  
We now compute the first order TMLE $\tilde{P}_n^{(1)}\tilde{P}_n^{(2)}(P_n^0)$ targeting $\Psi(P_0)$ using $\tilde{P}_n^{(2)}(P_n^0)$ as initial estimator and using $\tilde{P}_n$ to compute the MLE update.  Let $\bar{g}_n^{(2),*}=\bar{g}_n^{(2)}(P_n^0)$ and $\bar{q}_n^{(2),*}=\bar{q}_n^{(2)}(P_n^0)$. 
This first order TMLE does not update $\bar{g}_n^{(2)}(P_n^0)$ and just computes 
$\bar{q}_n^{(1)}(\bar{g}_n^{(2)}(P_n^0),\bar{q}_n^{(2)}(P_n^0) )$ obtained by maximizing the $\tilde{P}_n$-logistic regression log-likelihood with outcome $Y$, using as off-set $log (\bar{q}_n^{(2),*}/(1-\bar{q}_n^{(2),*})(W,A)$, and single covariate $C_{\bar{g}_n^{(2),*}}=A/\bar{g}_n^{(2),*}(W)$. Let's denote this update with $\bar{q}_n^{(1)}(p_n^{(2),*})$.
One can compute $\epsilon_n^{(1)}(\bar{q}_n^{(2),*},\bar{g}_n^{(2),*})$ by solving  \[\begin{array}{l}0=\tilde{P}_n C_{\bar{g}_n^{(2),*}}(Y-\bar{q}_n^{(1)}(\bar{q}_n^{(2),*},\bar{g}_n^{(2),*},\epsilon))\\=P_n \frac{\tilde{g}_n}{\bar{g}_n^{(2),*}} \left(\tilde{q}_n(1)-\frac{1}{1+\exp\left (-\log \bar{q}_n^{(2),*}/(1-\bar{q}_n^{(2),*})-\epsilon C_{\bar{g}_n^{(2),*}}\right)}\right).\end{array}\]
\item[Final second order TMLE of target estimand:]
Let $P_n^{2,(1),*}=\tilde{P}_n^{(1)}(P_n^{(2),*})$ be this second order TMLE of $P_0$.
Our second order TMLE of $\Psi(P_0)$ is given by $\Psi(\tilde{P}_n^{(1)}(P_n^{(2),*})=
\frac{1}{n}\sum_{i=1}^n \bar{q}_n^{(1)}(p_n^{(2),*})(W_i,1)$. 

Note that we always have $\tilde{P}_n D^{(1)}_{P_n^{2,(1),*}}=0$ and $\tilde{P}_nD^{(2)}_{n,P_n^{(2),*}}=0$. 
We want these equations also  solved under $P_n$ at  $o_P(n^{-1/2})$, such as making its absolute value smaller than $\sigma_{1n}/(n^{1/2}\log n)$. Some undersmoothing of $\tilde{P}_n$ might be needed to achieve this threshold.

\item[Iterative second order TMLE:] If we use a single step so that $\tilde{P}_n D^{(2)}_{\tilde{P}_n^{(2)}(P_n^0)}$ will not be zero exactly, then we can iterate the second order TMLE a few times (i.e. replace $P_n^0$ by $P_n^{2,(1),*}$ and iterate) so that we will still have $\tilde{P}_n D^{(2)}_{n,P_n^{(2),*}}\approx 0$ at user supplied precision, beyond the guaranteed $\tilde{P}_n D^{(1)}_{\tilde{P}_n^{(1)}(P_n^{(2),*}}=0$.
\end{description}

\subsection{Targeted selection of amount of undersmoothing for HAL-MLEs}
Recall that the second order TMLE $P_n^{2,(1),*}$ is indexed by the $L_1$-norm bound $C$ in the  HAL-MLE $\tilde{P}_n$. Thus, we can denote it with $P_{n,C}^{2,(1),*}$, where $C=(C_1,C_2)$ is the $L_1$-norm bound for the HAL-MLEs $\tilde{q}_n$ and $\tilde{g}_n $ of $\bar{q}_0=E_0(Y\mid A=1,W)$ and $\bar{g}_0=E_0(A\mid W)$, respectively. 
Let $\psi_n(C)\equiv \Psi(P_{n,C}^{2,(1),*})$ be the resulting plug-in second order TMLE of $\Psi(P_0)$.
For each choice of $C$ we can also estimate the variance of $D_n=D^{(1)}_n+D^{(2)}_n$, the natural plug-in estimate of the  sum of the first and second order efficient influence curve $D^{(1)}_{P_0}$ and $D^{(2)}_{n,P_0}$. Let's denote this variance with $\sigma^2_n(C)$, which will generally be increasing in $C$. Let $\psi_n(C)\pm 1.96 \sigma_n(C)/n^{1/2}$ be a $C$-specific confidence interval. We now search among  $C>C_{n,cv}$ for a local optimum of either the lower or the upper bound of this confidence band: if $\psi_n(C)$ is increasing, then we maximize the lower bounds and if $\psi_n(C)$ is decreasing, then we minimize the upper bounds. This corresponds with solving $
\frac{d}{dC}\psi_n(C)\pm 1.96\frac{d}{dC}\sigma_n(C)/n^{1/2}\approx 0$. This is essentially Lepski's method trading off bias and variance for $\psi_n(C)$.
One can think of this as a method that keeps increasing $C>C_{n,cv}$ till $\psi_n(C)$ reaches a plateau, so that its changes are washed out by the noise.

\subsection{Exact expansion for the second order TMLE and second order inference}
We have $R^{(1)}(\tilde{P}_n^{(1)}(P_0),P_0)=0$ since $g_0$ is not updated by $\tilde{P}_n^{(1)}(P_0)$ and $R^{(1)}(P,P_0)=0$ if $g(P)=g(P_0)$.
Since $\tilde{P}_n^{(2)}(P_0)$ is a standard MLE in a 2-dimensional correctly specified parametric model, and we undersmooth so that $(\tilde{P}_n-P_n)D^{(2)}_{n,P_0}=o_P(n^{-1/2})$, by the analogue of Lemma \ref{lemmaexactepsn1} for $\epsilon_n^{(2)}(P_0)$, we have $R_n^{(2)}(\tilde{P}_n^{(2)}(P_0),P_0)=O_P(n^{-1})$.

By our general exact expansion for the second order TMLE of Theorem \ref{theorem2} we have
\[
\begin{array}{l}
\Psi(\tilde{P}_n^{(1)}\tilde{P}_n^{(2)}(P_n^0))-\Psi(P_0)=
(P_n-P_0)D^{(1)}_{\tilde{P}_n^{(1)}(P_0)}+R^{(1)}(\tilde{P}_n^{(1)}(P_0),P_0)\\
+(P_n-P_0)D^{(2)}_{\tilde{P}_n^{(2)}(P_0)}+R_n^{(2)}(\tilde{P}_n^{(2)}(P_0),P_0)\\
+R_n^{(2)}(P_n^{(2)}(P_n^0),\tilde{P}_n)-R_n^{(2)}(P_n^{(2)}(P_0),\tilde{P}_n)\\
+(\tilde{P}_n-P_n) D^{(1)}_{\tilde{P}_n^{(1)}(P_0)}-P_n D^{(2)}_{n,\tilde{P}_n^{(2)}(P_0)}\\
=(P_n-P_0)D^{(1)}_{\tilde{P}_n^{(1)}(P_0)}+ (P_n-P_0)D^{(2)}_{\tilde{P}_n^{(2)}(P_0)}\\
+R_n^{(2)}(P_n^{(2)}(P_n^0),\tilde{P}_n)-R_n^{(2)}(P_n^{(2)}(P_0),\tilde{P}_n)\\
+(\tilde{P}_n-P_n) D^{(1)}_{\tilde{P}_n^{(1)}(P_0)}+(\tilde{P}_n-P_n) D^{(2)}_{n,\tilde{P}_n^{(2)}(P_0)}+O_P(n^{-1}).\\
\end{array}
\]
Due to $D^{(2)}$ already being  a first order difference, we will have
$(\tilde{P}_n-P_n) D^{(2)}_{n,\tilde{P}_n^{(2)}(P_0)}\approx O_P(n^{-1})$ under undersmoothing of $\tilde{P}_n$, so that we can declare this term as negligible.
Assuming we used a TMLE update so that $\tilde{P}_n D^{(2)}_{n,\tilde{P}_n^{(2)}(P)}=0$,
by Lemma \ref{lemma2}, we also have that 
\[
\begin{array}{l}
R_n^{(2)}(\tilde{P}_n^{(2)}(P),\tilde{P}_n)=R^{(1)}(\tilde{P}_n^{(1)}\tilde{P}_n^{(2)}(P),\tilde{P}_n)-\tilde{P}_n D^{(1)}_{\tilde{P}_n^{(1)}\tilde{P}_n^{(2)}(P)}+\tilde{P}_n D^{(2)}_{n,\tilde{P}_n^{(2)}(P)}\\
=R^{(1)}(\tilde{P}_n^{(1)}\tilde{P}_n^{(2)}(P),\tilde{P}_n).\end{array}
\]

So in the case that $\tilde{P}_n D^{(2)}_{n,\tilde{P}_n^{(2)}(P)}=0$, we have
\[
\begin{array}{l}
\Psi(\tilde{P}_n^{(1)}\tilde{P}_n^{(2)}(P_n^0))-\Psi(P_0)=
(P_n-P_0)D^{(1)}_{\tilde{P}_n^{(1)}(P_0)}
+(P_n-P_0)D^{(2)}_{n,\tilde{P}_n^{(2)}(P_0)}\\
+R^{(1)}(\tilde{P}_n^{(1)}\tilde{P}_n^{(2)}(P_n^0),\tilde{P}_n)-R^{(1)}(\tilde{P}_n^{(1)}\tilde{P}_n^{(2)}(P_0),\tilde{P}_n)\\
+(\tilde{P}_n-P_n)  D^{(1)}_{\tilde{P}_n^{(1)}(P_0)}+O_P(n^{-1}).\end{array}
\]
The term  $(\tilde{P}_n-P_n)D^{(1)}_{\tilde{P}_n^{(1)}(P_0)}$ is bounded by $(\tilde{P}_n-P_n)D^{(1)}_{P_0}$ and  is well behaved and   $O_P(n^{-2/3})$ under our level of undersmoothing of $\tilde{P}_n$ (Lemma \ref{lemmaexactepsn1}). 

So we conclude that the remainder is  driven by the $R^{(1)}$-difference, which will be a third order term by our general theory, and as shown below for this example.

\subsection{Second order inference based on the second order TMLE}
Let $D_n=D^{(1)}_{\tilde{P}_n^{(1)}\tilde{P}_n^{(2)}(P_n^0)}+D^{(2)}_{n,\tilde{P}_n^{(2)}(P_n^0)}$ and $\sigma^2_n=P_n D_n^2$ be the sample variance of this estimated influence curve $D_n$.
Then, $\Psi(P_n^{2,(1),*})\pm 1.96 \sigma_n/n^{1/2}$ is an $0.95$-confidence interval that takes into account the second order expansion and thereby only ignores a third order contribution coming from the difference of the $R^{(1)}(\cdot,\tilde{P}_n)$-remainder.
Notice
\[
D_n=C^y_{n,P_n^{(2),*} } (Y-\bar{q}_n^{(2),*})+C^a_{n,P_n^{(2),*}}(A-\bar{g}_n^{(2),*})+
C_{\bar{g}_n^{(2),*}}(Y-\bar{q}_n^{(1)}(p_n^{(2),*})),\]
where $P_n^{(2),*}$ is the TMLE $P_n^{(2)}(P_n^0)$ targeting $\Psi_n^{(1)}(P_0)$ and $\tilde{q}_n^{(1)}(p_n^{(2),*})$ is the TMLE of $\bar{q}_0$ targeting $\Psi(P_0)$ using as initial estimator this $P_n^{(2),*}$.

\subsection{Determining the third order remainder for the second order TMLE.}
By Lemma \ref{lemma2}, given that $R^{(1)}(\tilde{P}_n^{(1)}(P),\tilde{P}_n)$ is a second order difference, it allows a representation $R^{(1),+}(\tilde{P}_n^{(1)}(P),\tilde{P}_n)+1/2\frac{d}{dP}R^{(1)}(\tilde{P}_n^{(1)}(P),\tilde{P}_n)(P-\tilde{P}_n)$, which will be applied to $P=\tilde{P}_n^{(2)}(P)$, so that the directional derivative equals $- \tilde{P}_n D^{(2)}_{\tilde{P}_n^{(2)}(P)}$.  The term $R^{(1),+}$ is then a generalized third order difference. Depending on the definition of $\tilde{P}_n
^{(2)}(P)$, we either have $\tilde{P}_n D^{(2)}_{\tilde{P}_n^{(2)}(P)}=0$ or it equals $(P_n-\tilde{P}_n)D^{(2)}_{\tilde{P}_n^{(2)}(P)}$. Our third order remainder in the exact expansion equals the difference of this representation at $P=P_n^0$ and $P=P_0$.
So our job is to determine the form of this $R^{(1),+}(\tilde{P}_n^{(1)}\tilde{P}_n^{(2)}(P),\tilde{P}_n)$ in this example. 
The following lemma provides this result.
\begin{lemma}\label{lemmaexample2c}
Let $P_n^*=\tilde{P}_n^{(2)}(P)$ for some $P\in {\cal M}$.
Let
\[
\begin{array}{l}
R_{2,\tilde{q}_n^{(1)},p_n^*}(p_n^*,\tilde{p}_n)\equiv 
\tilde{q}_n^{(1)}(p_n^*)-\tilde{q}_n^{(1)}(\tilde{p}_n)-\frac{d}{dq_n^*}\tilde{q}_n^{(1)}(q_n^*,g_n^*)(q_n^*-\tilde{q}_n)\\
\hfill -
\frac{d}{dg_n^*}\tilde{q}_n^{(1)}(q_n^*,g_n^*)(g_n^*-\tilde{g}_n)\\
{R}^{(1),p,+,a}(\tilde{P}_n^{(1)}(P_n^*),\tilde{P}_n))
=
1/2 \tilde{P}_n R_{2,\tilde{q}_n^{(1)},p_n^*}(p_n^*,\tilde{p}_n)
\frac{\bar{g}_n^*-\tilde{g}_n}{\tilde{g}_n} \\
R^{(1),+,b}(\tilde{P}_n^{(1)}(P_n^*),\tilde{P}_n)=
1/2\tilde{P}_n \frac{d}{dp_n^*}\bar{q}_n^{(1)}(p_n^*)(p_n^*-\tilde{p}_n) \frac{(\bar{g}_n^*-\tilde{g}_n)^2}{\bar{g}_n^*\tilde{g}_n}\\
\hfill -\tilde{P}_n  (\bar{q}_n^{(1)}(p_n^*)-\tilde{q}_n)\frac{1}{\tilde{g}_n(\bar{g}_n^*)^2}(\bar{g}_n^*-\tilde{g}_n)^3.
\end{array}
\]

Then, 
\[
\begin{array}{l}
R^{(1),+}(\tilde{P}_n^{(1)}(P_n^*),\tilde{P}_n)={R}^{(1),p,+,a}(\tilde{P}_n^{(1)}(P_n^*),\tilde{P}_n^{(1)}(\tilde{P}_n))+R^{(1),+,b} (\tilde{P}_n^{(1)}(P_n^*),\tilde{P}_n),\end{array}
\]
and
\[
R^{(1)}(\tilde{P}_n^{(1)}(P_n^*),\tilde{P}_n)= R^{(1),+}(\tilde{P}_n^{(1)}(P_n^*),\tilde{P}_n)-
1/2\tilde{P}_n D^{(2)}_{n,P_n^*}.\]
\end{lemma}
The proof of this lemma is presented in the Appendix \ref{AppendixH}.
Thus, $R^{(1),+}(\tilde{P}_n^{(1)}\tilde{P}_n^{(2)}(P),\tilde{P}_n)$  is a sum of three third and fourth order differences.
This sum could be written itself as a third order difference, by combining the last two terms in $R^{(1),+,b}$ in a single third order difference. Alternatively, it can be represented as a sum of a  third order polynomial difference and fourth order difference.

\subsection{The trade-off in selecting the variation norm for the HAL-MLE}
The HAL-regularized second order TMLE is indexed by the $L^1$-norm $C$ used in $\tilde{P}_n$, possibly, a separate $C_1$ for the HAL-MLE $\tilde{g}_n$ of $\bar{g}_0$  and $C_2$ for $\tilde{q}_n$ of $\bar{q}_0$. There is a minor trade-off when selecting $C\geq C_{n,cv}$: one wants $\tilde{P}_n$ to be a good estimator of $P_0$ under the restriction that $(\tilde{P}_n-P_n)D^{(1)}_{\tilde{P}_n^{(1)}(P_0)}\approx 0$. 

Specifically, this is seen as follows. Firstly, $D^{(2)}_{n,\tilde{P}_n^{(2)}(P_0)}$ involves ratios $\tilde{g}_n/\bar{g}_n^{(2)}(p_0)$ and $\tilde{q}_n-\bar{q}_n^{(1)}(\tilde{p}_n^{(2)}(p_0))$, so that one wants $\tilde{g}_n$ and $\tilde{q}_n$ to be good estimators of  $\tilde{g}_n^{(2)}(p_0)\approx \bar{g}_0$ and $\tilde{q}_n^{(1)}(\tilde{p}_n^{(2)}(p_0))\approx \bar{q}_0$. One also wants that if $\bar{g}_0\approx 0$, then $\tilde{g}_n\approx 0$ so that the ratio $\tilde{g}_n/\tilde{g}_n^{(2)}(p_0)$ is well behaved.  Similarly, $R^{(1)}(\tilde{P}_n^{(1)}\tilde{P}_n^{(2)}(P_0),\tilde{P}_n)$ will be rewarded by $\tilde{p}_n$ being a good estimator of $p_0$. The term $ R^{(1),+}(\tilde{P}_n^{(1)}\tilde{P}_n^{(2)}(P_n^0),\tilde{P}_n)$ is rewarded by $P_n^0$ not being different from $\tilde{P}_n$. So all these terms are only concerned with  $\tilde{P}_n$ being a good estimator  of $P_0$, suggesting $C\approx C_{n,cv}$. The only term that suggests some undersmoothing of $\tilde{P}_n$ might be needed is that $(\tilde{P}_n-P_n)  D^{(1)}_{\tilde{P}_n^{(1)}(P_0)}$ has to be small (and this will generally also take care that $(\tilde{P}_n-P_n)D^{(2)}_{\tilde{P}_n^{(2)}(P_0)}$ is small as well). This is generally a small well behaved term, but to make it $o_P(n^{-1/2})$  a weak amount of undersmoothing might be needed. 

In this subsection we provide more details of how undersmoothing the HAL-MLE reduces this term.
By Lemma \ref{lemmaexactepsn1}, we have $(\tilde{P}_n-P_n) D^{(1)}_{\tilde{P}_n^{(1)}(P_0)}$ is bounded by $(\tilde{P}_n-P_n) A/\bar{g}_0(Y-\bar{q}_0)$.

{\bf Contribution of HAL-MLE $\tilde{q}_n$:}
An HAL-MLE $\tilde{q}_n$ approximately solves $P_n \phi_j(Y-\tilde{q}_n)\approx 0$ for the basis functions $\phi_j$ that have non-zero coefficients in the HAL-MLE fit $\tilde{q}_n$: it solves all scores along paths that keep the $L_1$-norm identical, so there is only a single constraint obstructing the HAL-MLE to solve these score equations exactly. Let $\phi_j(W,1)$ be the spline basis functions evaluated at $A=1$. However, $P_n \phi_j(W,1) \tilde{q}_n(W,1)=
\tilde{P}_n \phi_j Y$, using that the distribution of $W$ under $\tilde{P}_n$ equals  the empirical measure. Therefore $P_n \phi_j(Y-\tilde{q}_n)=P_n \phi_j Y-\tilde{P}_n \phi_j Y=(P_n-\tilde{P}_n)\phi_j Y$. 
Thus an HAL-MLE $\tilde{q}_n$ approximately solves $(P_n-\tilde{P}_n )\phi_j Y=0$ across all its selected basis functions, and therefore $(P_n-\tilde{P}_n)\phi_n(W,A)Y\approx 0$ for any function $\phi_n$ in the linear span of these basis functions $\phi_j$. We need to solve $(P_n-\tilde{P}_n)A/\bar{g}_0(W)(Y-\bar{q}_0)$. Thus, the score equations solved by the HAL-MLE $\tilde{q}_n$ take care that we approximately solve $(P_n-\tilde{P}_n)A/\bar{g}_0Y$. 

{\bf Contribution of HAL-MLE $\tilde{g}_n$:}
We also need to solve $(P_n-\tilde{P}_n)(A/\bar{g}_0) \bar{q}_0=(P_n-\tilde{P}_n ) \frac{\bar{q}_0}{\bar{g}_0}A$.
Now we can use that $\tilde{g}_n$ is an HAL-MLE that approximately solves $P_n \phi_j(A-\tilde{g}_n)=0$ where $\phi_j$ are the basis functions with non-zero coefficients in $\tilde{g}_n$. 
As above it follows that $(P_n-\tilde{P}_n)\phi_j A\approx 0$ and thus $(P_n-\tilde{P}_n)\phi_n A\approx 0$ for $\phi_n$ in the linear span of these basis functions $\phi_j$. This linear span approximates $ \frac{\bar{q}_0(W,1)}{\bar{g}_0(W)}$, showing that also $(P_n-\tilde{P}_n )\frac{\bar{q}_0}{\bar{g}_0} A\approx 0$. 

Therefore, we conclude that the combined set of score equations solved by the  two HAL-MLEs $\tilde{g}_n$ and $\tilde{q}_n$ together approximately span $(P_n-\tilde{P}_n ) A/\bar{g}_0(Y-\bar{q}_0)$. This proves that indeed the bound $C=(C_1,C_2)$ can be tuned to make this term $(P_n-\tilde{P}_n) D^{(1)}_{\tilde{P}_n^{(1)}(P_0)}\approx 0$ and negligible relative to the third order contributions. 

\subsection{Second order TMLE for continuous outcome.}
In the Appendix \ref{AppendixJ} we derive the second order canonical gradient for the case that $Y$ is a continuous outcome; we use a linear  regression least favorable path $\bar{q}^{(1)}(p,\epsilon)=\bar{q}+\epsilon A/\bar{g}$; least squares regression $\epsilon_n^{(1)}(p)=\arg\min_{\epsilon}E_{\tilde{P}_n}(Y-\bar{q}_n^{(1)}(p,\epsilon))^2$ as MLE-step, and $\bar{q}_n^{(1)}(p)=\bar{q}_n^{(1)}(p,\epsilon_n^{(1)}(p))$.  Here one can decide to move  $1/\bar{g}$ in the weight of the squared error loss and use $\bar{q}_n^{(1)}(p)=\bar{q}+\epsilon A$ instead. 
Lemma \ref{lemmactsoutcome} in Appendix \ref{AppendixJ} shows that  the canonical gradient $D^{(2)}_{n,P}=D^{(2)}_{n,P,1}+D^{(2)}_{n,P,2}$
 for  $\Psi_n^{(1)}(P)=E_0\bar{q}_n^{(1)}(P)$  is now given by
\[
\begin{array}{l}
D^{(2)}_{n,P}=
C_g\frac{\bar{g}-\tilde{g}_n}{\bar{g}}   (Y-\bar{q}) \\
+ \epsilon_n^{(1)}(p)\frac{\tilde{g}_n-\bar{g}}{\bar{g}^3}(A-\bar{g})\\
- \frac{\tilde{g}_n}{\bar{g}^2} (\tilde{q}_n-\bar{q}_n^{(1)}(p) )(A-\bar{g}) ,
 \end{array}
 \]
and the first term represent $D^{(2)}_{n,P,1}$. Notice that this formula corresponds with the formula for the binary outcome by setting $\tilde{c}_{\tilde{P}_n,P}^{(1)}$ and $(\bar{q}_n^{(1)}(1-\bar{q}_n^{(1)}(p))/(\bar{q}(1-\bar{q}))$ equal to 1.
The second order TMLE is now defined completely analogue to above.

\section{Simulations for the HAL-regularized second order TMLE of integrated square of density} \label{sec:ave_dens}

In this section, we focus on the integrated square of density for univariate discrete variables. Consider observed $n$ i.i.d. copies of $O \sim P_0 \in \mathcal{M}$. The p.m.f. function is $p = dP/d\mu$ with support points $\{x_i: i\in \mathcal{I}\}$, where $\mu$ is the counting measure. In this case, the integrated square of density, or the average density value, $\Psi: \mathcal{M} \to \openr$, is a summation over the supports, $\Psi(P) = \sum_{i\in\mathcal{I}} p(x_i)^2$. The negative bias, $R^{(1)}(P, P_0) = \Psi(P) - \Psi(P_0) + P_0D^{(1)}_P = - \sum_{i\in\mathcal{I}} (p(x_i) - p_0(x_i))^2$, now takes the form of a negative sum of squares, and larger biases can be created by setting $P$ to be more distant from $P_0$.

\begin{figure}[ht!]
\includegraphics[scale = .5]{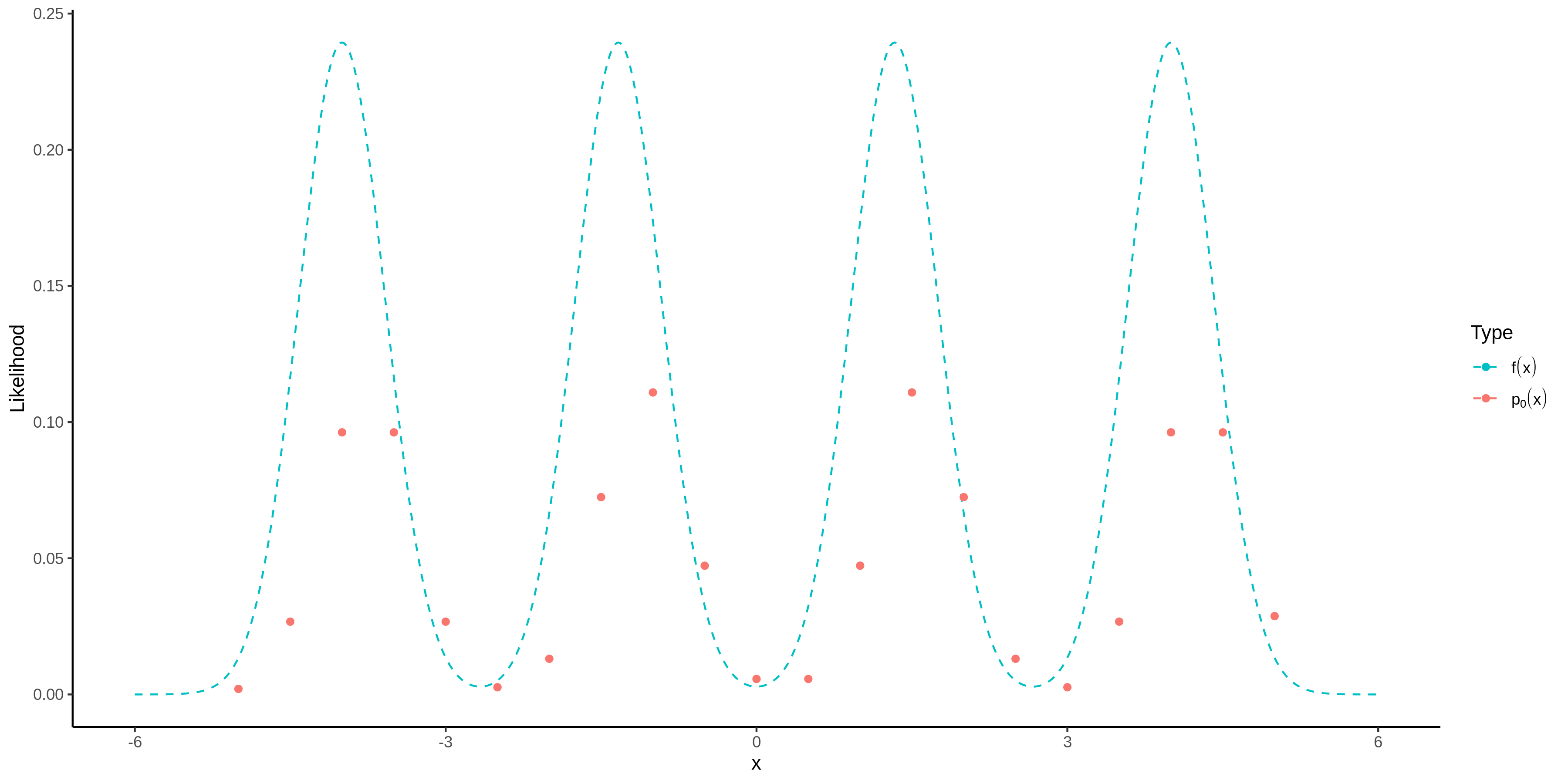}
  \caption{The discretized $p_0$ from the Gaussian mixture density $f(x) = \sum_{k=1}^4 g_k(x)$. $g_k$ is the density of $N(\mu_k, \sigma_K)$. $\mu_k$s are evenly placed across the interval $[-4, 4]$. $\sigma_K = 10/K/6$. $K = 4$. Discrete supports: -5, -4.5, \ldots, 5. } 
\label{fig:p0}
\end{figure}

In the simulation, we define $P_0$ (Figure \ref{fig:p0}) by discretizing Gaussian mixtures that have densities $f(x) = \frac{1}{K}\sum_{k=1}^K g_k(x)$, $g_k(x) = \frac{1}{\sqrt{2\pi} \sigma_K} \exp[ -\frac{1}{2} (x - \mu_k)^2 / \sigma_K^2] $. For a given $K$, $\mu_k$s are evenly placed across the interval $[-4, 4]$. $\sigma_K = 10/K/6$. As for discretization, the support points $x_1, \ldots, x_I$ are chosen to be evenly placed across the interval $[-5, 5]$ with $I=21$, and $p_0(x_i)$ is defined as $\int_{x_{i-1}}^{x_i} f(u) du$ except for that $p_0(x_1) = \int_{-\infty}^{x_1} f(u) du$ and $p_0(x_I) = \int_{x_{I-1}}^{\infty} f(u) du$. 

The potentially biased initial estimator, $p_n^0$, is created by  adding a same mass to all the support points of the empirical pmf $p_n$ and scaling to sum $1$. 

\subsection{Performance of $\tilde P_n$-based second order TMLE}

Suppose that $\tilde p_n$ is a discrete HAL-MLE whose knot points equal the pmf supports, $x_1, \ldots, x_I$. The first order HAL-regularized TMLE is implemented with a locally least favorable path $\tilde p_n^{(1)}(p) = \tilde p^{(1)}(p, \epsilon_n^{(1)}(P))$ where $\tilde p^{(1)}(p, \epsilon) = (1 + \epsilon D^{(1)}_P)p$ and $\epsilon_n^{(1)}(P) = \argmax{\epsilon} \tilde P_n \log \tilde p^{(1)}(p, \epsilon)$. The maximum first order step size is restricted to ensured that all the updates remain as valid pmf functions. The second order target is $\Psi_n^{(1)}: \mathcal{M} \to \openr$, $\Psi_n^{(1)}(P) = \Psi(\tilde P_n^{(1)}(P))$, where the canonical gradient $D_{n, P}^{(2)}$ is calculated as Lemma \ref{lemmaex1a}. Let $\tilde p^{(2)}_n(p, \epsilon) = (1+\epsilon D^{(2)}_{n, P})p$ follow a similar locally least favorable path. 

In the simulation, we restrict both $\epsilon^{(1)}_n$ and $\epsilon^{(2)}_n$ with maximum step size $0.1$, and iterate by letting $\tilde P_n^{(1)}\tilde P_n^{(2)}(P_n^0)$ be the new initial $P_n^0$. Second order updates are skipped, i.e. $\epsilon_n^{(2)}(P_n^0)=0$, if $|\tilde P_n D^{(2)}_{n, P_n^0}| < 1/n$ at any iteration, and similarly for first order steps $\epsilon_n^{(1)}(\tilde P_n^{(2)}(P_n^0))=0$ if $|\tilde P_n D^{(1)}_{\tilde P_n^{(2)}(P_n^0)}| < 1/n$. The algorithm stops when both $|\tilde P_n D^{(1)}_{\tilde P_n^{(1)}\tilde P_n^{(2)}(P_n^*)}| < 1/n$ and $|\tilde P_n D^{(2)}_{n, \tilde P_n^{(2)}(P_n^*)}| < 1/n$ are reached at the final update $P_n^*$. The iterative second order TMLE of $\Psi(P_0)$ is given by $\Psi(\tilde P_n^{(1)}\tilde P_n^{(2)}(P_n^*))$. 
Results show that the second order TMLE provides additional bias control over the first order TMLE, and remains consistent even with the increasing bias of the initial $P_n^0$ and the squared remainders (Figure \ref{fig:int_sq_dens}; Table \ref{tab:int_sq_dens}).

\subsection{Effect of undersmoothing HAL}

Undersmoothing of the HAL-MLE is conducted by selecting a larger variation norm than the cross validation selector $C_{n, cv}$. In the simulation, we use \texttt{glmnet} for the HAL fits, so that undersmoothing is equivalently achieved by selecting a smaller \texttt{glmnet lambda} value than the \texttt{cv.glmnet} selector. Specifically, we calculate iterative second order TMLE updates $P_n^*(\texttt{lambda})$ as a function of $10$ positive and equidistant \texttt{lambda} candidates that are smaller than or equal to the cross validation selector. The largest candidate $\texttt{lambda}_n$ such that $|P_n D^{(1)}_{\tilde P_n^{(1)}\tilde P_n^{(2)}(P_n^*(\texttt{lambda}_n))}| < 1/n$ and $|(P_n - \tilde P_n) D^{(1)}_{\tilde P_n^{(1)}\tilde P_n^{(2)}(P_n^*(\texttt{lambda}_n))}| < 1/n$ is used to calculated the undersmoothed second order TMLE, $\Psi(\tilde P_n^{(1)}\tilde P_n^{(2)}(P_n^*(\texttt{lambda}_n)))$. Undersmoothing ensures that both $\tilde P_n$ and $P_n$-based estimating equations are solved to the desired precision, and provides finite sample advantage in bias control (Figure \ref{fig:int_sq_dens}; Table \ref{tab:int_sq_dens}). 

 \subsection{Effect of empirical TMLE-updates}
\begin{figure}[ht!]
\includegraphics[scale = .45]{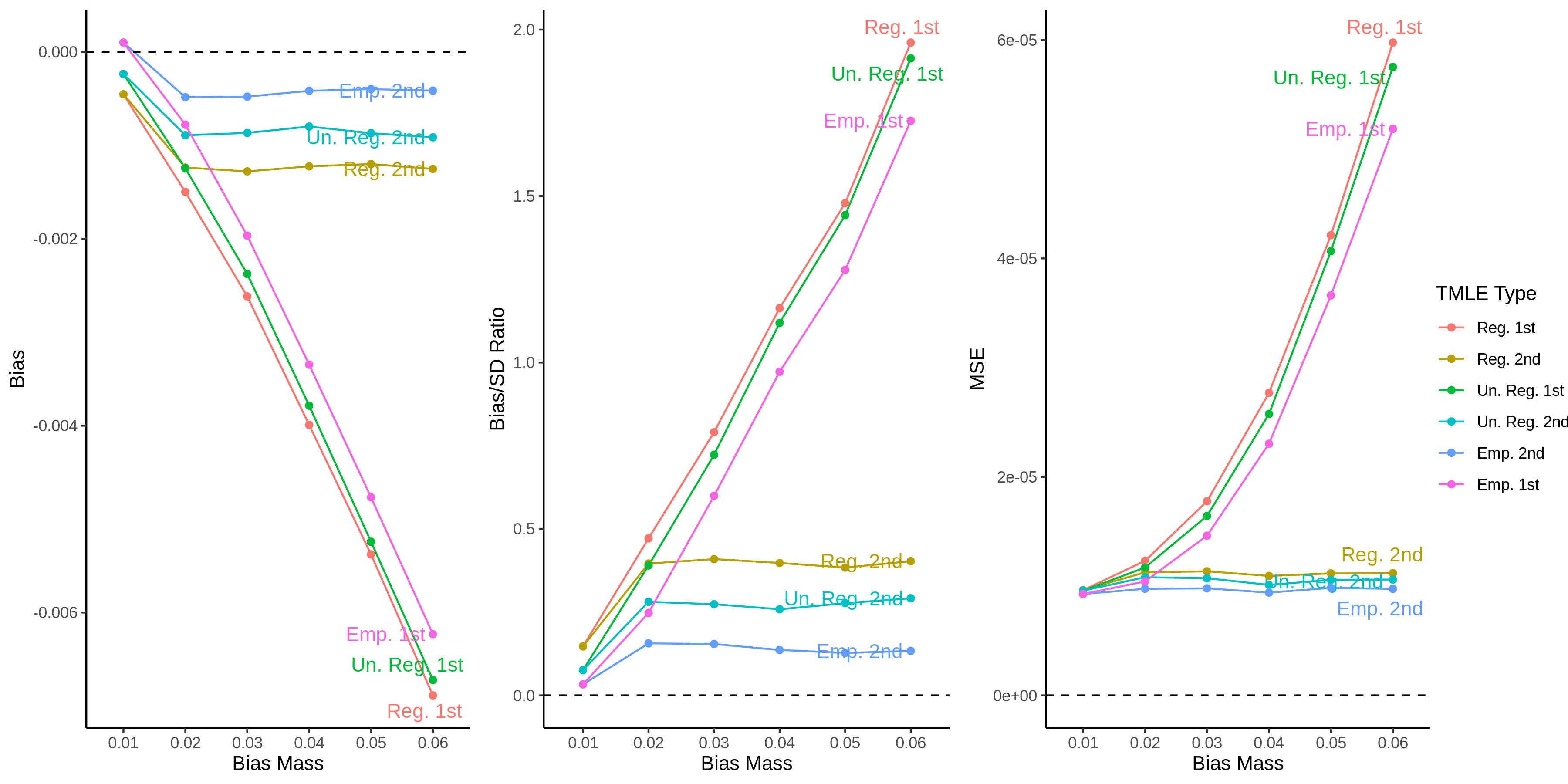}
  \caption{Bias, Bias/SD ratio, and MSE performance of HAL-regularized first and second order TMLE (Reg. 1st, Reg. 2nd), compared with empirical second order TMLE (Emp. 2nd), in the integrated square of density simulation, with or without undersmoothing (Un.). n=500. Biased initial estimators $P_n^0$ are created by adding a point mass to randomly selected five of the support points of the empirical $P_n$. Second order TMLEs provide additional total remainder control over first order TMLEs following likelihood guidance in all scenarios. Empirical updates for second order TMLE provide the similar level of bias reduction with or without undersmoothing. }
\label{fig:int_sq_dens}
\end{figure}

 Replace HAL-regularized TMLE-updates $\tilde{P}_n^{(1)}\tilde{P}_n^{(2)}(P_n^0)$ with empirical TMLE-updates $P_n^{(1)}P_n^{(2)}(P_n^0)$, where $p_n^{(1)}(p) = \tilde{p}^{(1)}(p, \epsilon_{P_n}^{(1)}(P))$, 
 $\epsilon_{P_n}^{(1)}(P) = \argmax{\epsilon} P_n \log \tilde p^{(1)}(p, \epsilon)$, and
 $p_n^{(2)}(p) = \tilde{p}^{(2)}_n(p, \epsilon_{P_n}^{(2)}(P))$, 
 $\epsilon_{P_n}^{(2)}(P) = \argmax{\epsilon} P_n \log \tilde p^{(2)}_n(p, \epsilon)$, following the same paths $\tilde p^{(1)}(p, \epsilon)$ and $\tilde p_n^{(2)}(p, \epsilon)$ as above. Figure \ref{fig:int_sq_dens} demonstrates that empirical updates achieve the similar level of bias reduction without undersmoothing. 

\begin{table}[ht!] 
\caption{HAL-regularized (Reg.) first and second order TMLE, and empirical (Emp.) second order TMLE, in the integrated square of density simulation, with or without undersmoothing. Initial $P_n^0$ is adding a bias mass to randomly selected five of the supports of the empirical $P_n$, and then scaling to sum 1. Second order TMLE controls the exact total remainder from first order TMLE. Undersmoothing controls $(\tilde P_n - P_n)D^{(1)}_{\tilde P_n^{(1)}\tilde P_n^{(2)}(P_n^*)}$ at the final update $P_n^*$. Empirical updates achieve bias control without undersmoothing. }
\centering
\begin{tabular}{r|rrr}
\hline
 & n=500 &  &  \\ 
  \hline
Bias mass: 0.02 & Bias & SD & MSE \\ 
  \hline
Reg. 1st order & -1.498E-03 & 3.176E-03 & 1.232E-05 \\ 
  Undersmoothed Reg. 1st order & -1.245E-03 & 3.188E-03 & 1.170E-05 \\ 
  Emp. 1st order & -7.776E-04 & 3.136E-03 & 1.043E-05 \\ 
  Reg. 2nd order & -1.237E-03 & 3.122E-03 & 1.127E-05 \\ 
  Undersmoothed Reg. 2nd order & -8.904E-04 & 3.168E-03 & 1.082E-05 \\ 
  Emp. 2nd order & -4.829E-04 & 3.088E-03 & 9.759E-06 \\ 
  Undersmoothed Emp. 2nd order & -4.832E-04 & 3.088E-03 & 9.758E-06 \\ 
   \hline
    Bias mass: 0.04 & Bias & SD & MSE \\ 
  \hline
Reg. 1st order & -3.991E-03 & 3.431E-03 & 2.769E-05 \\ 
  Undersmoothed Reg. 1st order & -3.785E-03 & 3.382E-03 & 2.576E-05 \\ 
  Emp. 1st order & -3.346E-03 & 3.443E-03 & 2.304E-05 \\ 
  Reg. 2nd order & -1.223E-03 & 3.074E-03 & 1.094E-05 \\ 
  Undersmoothed Reg. 2nd order & -7.972E-04 & 3.081E-03 & 1.012E-05 \\ 
  Emp. 2nd order & -4.150E-04 & 3.042E-03 & 9.414E-06 \\ 
  Undersmoothed Emp. 2nd order & -4.166E-04 & 3.041E-03 & 9.411E-06 \\ 
   \hline
   Bias mass: 0.06 & Bias & SD & MSE \\ 
  \hline
Reg. 1st order & -6.887E-03 & 3.512E-03 & 5.976E-05 \\ 
  Undersmoothed Reg. 1st order & -6.722E-03 & 3.513E-03 & 5.752E-05 \\ 
  Emp. 1st order & -6.231E-03 & 3.610E-03 & 5.185E-05 \\ 
  Reg. 2nd order & -1.251E-03 & 3.104E-03 & 1.119E-05 \\ 
  Undersmoothed Reg. 2nd order & -9.129E-04 & 3.128E-03 & 1.061E-05 \\ 
  Emp. 2nd order & -4.139E-04 & 3.097E-03 & 9.752E-06 \\ 
  Undersmoothed Emp. 2nd order & -4.157E-04 & 3.095E-03 & 9.743E-06 \\ 
   \hline
\end{tabular}
\label{tab:int_sq_dens}
\end{table}

\section{Simulations for the iterative HAL-regularized second order TMLE of the treatment specific mean outcome}\label{section11b}
Recall our observed data structure is $(W,A,Y)$ where $W$ is a vector of baseline covariates, $A$ is a binary treatment, and $Y$ is a binary outcome, the model is nonparametric and the target parameter is $EY_1=E_PE_P(Y\mid A=1,W)$.
We implemented the iterative second order TMLE algorithm that also iteratively targets the HAL-MLE as presented in previous section.
That is, given $P_n^0=(Q_{W,n},\bar{g}_n^0,\bar{Q}_n^0)$, we replace an initial HAL-MLE $\tilde{P}_n=(Q_{W,n},\tilde{g}_n,\tilde{Q}_n)$ by its TMLE $(\tilde{P}_n^*=(Q_{W,n},\tilde{g}_n^*,\tilde{Q}_n^*)$, using as initial estimator $\tilde{P}_n$, and targeted so that  $(P_n-\tilde{P}_n^*)D^{(1)}_{\tilde{P}_n^{(2)}(P_n^0)}\approx 0$ and $(P_n-\tilde{P}_n^*)D^{(2)}_{n,\tilde{P}_n^{(2)}(P_n^0)}\approx 0$, and we iterate this process by replacing the initial estimator by the previous second order TMLE full update $\tilde{P}_n^{(1)}\tilde{P}_n^{(2)}(P_n^0)$. The algorithm stops when $(P_n-\tilde{P}_n^*)D^{(1)}_{\tilde{P}_n^{(1)}\tilde{P}_n^{(2)}(P_n^0)}\approx 0$ and $(P_n-\tilde{P}_n^*)D^{(2)}_{n,\tilde{P}_n^{(2)}(P_n^0)}\approx 0$. In our simulations, it converged in one step so that the HAL-MLE was only targeted once: this might be partly due to the fact that we used a somewhat undersmoothed HAL-MLE $\tilde{P}_n$.   Recall also that  in this algorithm we  use the regular TMLE updates $\tilde{P}_n^{(1)}(P)$ and $\tilde{P}_n^{(2)}(P)$  based on $P_n$, while the targeting of the HAL-MLE arranged that this regular TMLE update behaves the same as if we would have used the targeted  HAL-MLE $\tilde{P}_n^*$ in our TMLE-updates. Given the relevant initial $P$, the TMLE $\tilde{P}_n^{(2)}(P)$ was implemented with a universal least favorable path, while $\tilde{P}_n^{(1)}(P)$ is the simple closed form logistic regression update. 
In our simulations we report the bias and standard error scaled by $n^{1/2}$ and the MSE scaled by $n$.
The goal of the simulations is to evaluate if indeed the second order TMLE is able to achieve a strong bias reduction relative to the first order TMLE in the case that the first order TMLE has a non-negligible bias, while also making sure that it behaves as the first order TMLE otherwise.

\subsection{Simulation 1: $n^{-1/4}$-consistent $g_n^0$; inconsistent $Q_n^0$, targeted undersmoothed HAL-MLE}
The data is generated as follows: $W$  uniform $[-1,1]$ distributed; $A$, given $W$, is Bernoulli with probability $\bar{g}(W)=  \mbox{Expit}(2W - W^2 )$; and $Y$ is Bernoulli with probability $Q(W,A) =   \mbox{Expit}(W + A/2 )$.
The initial estimates $\bar{g}_n^0$ and $Q_n^0$ inconsistent by adding $n^{-1/4}$ bias to a misspecified version of the true $\bar{g}_0$ and $Q_0$. Specifically, we define
\begin{eqnarray*} 
\bar{g}_n^0(W) &=&   \mbox{Expit}(2W-W^2 ) +  \frac{0.1 + 2|W|}{2 n^{0.25}}\\
Q_n^0(A, W) &=&   \mbox{Expit}(2W + 2  A + AW/2 ) +  \frac{| 0.1 + 2|W| -A  |}{3 n^{0.25}}.
\end{eqnarray*}
The initial HAL-MLEs $\tilde g_n$ and $\tilde Q_n$ are given by undersmoothed Highly Adaptive Lasso logistic regressions, where the lasso tuning parameter lambda is chosen to be 10 lambdas above the cross-validation selected lambda with respect to glmnet's default grid of 100 lambdas.

\begin{table}
\centering
\begin{tabular}{r|r|r|r|r|r|r} 
\hline
n &   bias 1-st &   bias 2-nd &   se 1-st &   se 2-nd & mse 1-st & mse 2-nd \\  
\hline 400 &   -0.720 &   0.078 &   0.815 &   1.175 & 1.087 & 1.178\\
\hline 750 &   -0.996 &   0.029 &   0.800 &   1.102 &  1.278 & 1.102\\ 
\hline 1000 &   -1.258 &   -0.062 &   0.786 &   1.066 & 1.483 & 1.068\\ 
\hline 1200 &   -1.345 &   0.022 &   0.809 &   1.028 &  1.570 & 1.028\\ 
\hline 1600 &   -1.549 &   -0.019 &   0.818 &   1.055 & 1.752  & 1.055\\
\hline 2500 &   -2.066 &   -0.094 &   0.819 &   0.999 & 2.222 & 1.003\\
\hline 
\end{tabular} 
\caption{Simulation I:  $g_n^0$ is $n^{-1/4}$-consistent, while $Q_n^0$ is inconsistent.  The HAL-MLE $\tilde{P}_n$ is targeted and undersmoothed. The first order TMLE should have $n^{1/2}$-scaled bias that increases with $n$  while the second order TMLE has a $n^{1/2}$-bias that should be constant in $n$. We observe that the second order TMLE has a negligible bias and thereby still provides valid inference.}
\end{table}  

\subsection{Simulation 2: $n^{-1/4}$-non-random consistent $g_n^0$ and $Q_n^0$, targeted undersmoothed HAL-MLE}
As a second simulation, $\overline{g}_n^0$ and $Q_n^0$ are obtained by adding a fixed $n^{-1/4}$ bias to the true functions $\overline{g}_0$ and $Q_0$. The initial estimators $\overline{g}_n^0$ and $Q_n^0$ are given by
\begin{eqnarray*}
\overline{g}_n^0(W) &=&    \mbox{Expit}(2W-W^2 ) +  \frac{0.1 + 2|W|}{2 n^{0.25}}\\
Q_n^0(A, W)& = &  \mbox{Expit}(W + A/2 ) +  \frac{|0.1 +2|W|+  A/2 |}{3 n^{0.25}}.
\end{eqnarray*}
$\tilde{P}_n^*$ is given by the same targeted undersmoothed Highly Adaptive lasso estimator as in the previous simulation.

\begin{table}
\centering
\begin{tabular}{r|r|r|r|r|r|r} 
\hline
n &   bias 1-st &   bias 2-nd &   se 1-st &   se 2-nd & mse 1-st & mse 2-nd \\  
\hline
500 &   -0.193 &   0.079 &   0.858 &   1.062 & 0.879 & 1.065\\
\hline 
1000 &   -0.226 &   0.041 &   0.942 &   1.126 &  0.968 & 1.126\\ 
\hline 
1500 &   -0.273 &   -0.022 &   0.887 &   1.000 & 0.928 & 1.000\\ 
\hline
2500 &   -0.244 &   0.027 &   0.888 &   0.955 & 0.920 & 0.955\\ 
\hline 
4000 &   -0.256 &   0.077 &   0.892 &   0.940 & 0.928  & 0.943\\
\hline 
\end{tabular} 
\caption{Simulation II: $g_n^0$  and $Q_n^0$ are both  $n^{-1/4}$-consistent.  The HAL-MLE is targeted and undersmoothed. The first order TMLE should have $n^{1/2}$-scaled bias that does not converge to zero (but is constant in $n$), while the second order TMLE should have a  $n^{-1/2}$-scaled bias that converges to zero at rate $n^{-1/4}$. We indeed observe that the second order TMLE has a negligible bias (bias/SE $<10$), and thereby still provides valid inference.}

\end{table}

\subsection{Simulation 3: $n^{-1/4}$-consistent $g_n^0$ and $Q_n^0$, targeted HAL-MLE} The third simulation is identical to the previous simulation except that the $\tilde{P}_n^*$
uses as initial estimator the HAL-MLEs  of $g_0$ and $Q_0$ at the cross-validation selected lambda. In this way, we aimed to evaluate if the extra undersmoothing of $\tilde{P}_n^*$ changes the finite sample performance of the second order TMLE.

\begin{table}
\centering
\begin{tabular}{r|r|r|r|r|r|r} 
\hline
n &   bias 1-st &   bias 2-nd &   se 1-st &   se 2-nd & mse 1-st & mse 2-nd \\  
\hline
500 &   -0.138 &   0.045 &   0.881 &   1.056 & 0.891 & 1.057\\
\hline 
1000 &   -0.197 &   0.037 &   0.875 &   1.005 &  0.897 & 1.006\\ 
\hline 
1500 &   -0.260 &   -0.017 &   0.885 &   0.969 & 0.923 & 0.969\\ 
\hline
2500 &   -0.278 &   0.013 &   0.895 &   0.991 &  0.937 & 0.991\\ 
\hline 
4000 &   -0.268 &   0.084 &   0.867 &   0.934 & 0.908  & 0.938\\
\hline 
\end{tabular} 
\caption{Simulation III: $g_n^0$  and $Q_n^0$ are both  $n^{-1/4}$-consistent. The HAL-MLE is targeted but not undersmoothed. The first order TMLE should have $n^{1/2}$-scaled bias that does not converge to zero (but is constant in $n$), while the second order TMLE should have a  $n^{-1/2}$-scaled bias that converges to zero at rate $n^{-1/4}$. We indeed observe that the second order TMLE has a negligible bias (bias/SE $<10$), and thereby still provides valid inference.}
\end{table} 

\subsection{Simulation 4: using HAL-MLE for initial estimators, targeted  undersmoothed HAL-MLE}

For simulation 4, the same simulation design is used. However, $g_n^0$ and $Q_n^0$ are estimated nonparametrically with Highly adaptive lasso.
$\tilde{P}_n^*$ is a targeted HAL-MLE that uses as initial $\tilde g_n$ and $\tilde Q_n$  obtained with an undersmoothed Highly Adaptive Lasso, where the lasso tuning parameter lambda is chosen to be 5 lambdas above the cross-validation selected lambda with respect to glmnet's default grid of 100 lambdas.

\begin{table}
\centering
\begin{tabular}{r|r|r|r|r|r|r} 
\hline
n &   bias 1-st &   bias 2-nd &   se 1-st &   se 2-nd & mse 1-st & mse 2-nd \\  
\hline 
1000 &   -0.032 &   -0.037 &   0.952 &   0.997 &  0.953 & 0.997\\ 
\hline 
1500 &   0.013 &   0.004 &   0.947 &   0.977 & 0.947 & 0.977\\ 
\hline
2500 &   -0.003 &   0.007 &   0.938 &   0.954 &  0.938 & 0.954\\ 
\hline 
4000 &   -0.009 &   0.002 &   0.991 &   1.005 & 0.991  & 1.005\\
\hline 
\end{tabular} 
\caption{Simulation IV: $g_n^0$  and $Q_n^0$ are HAL-MLE using a cross-validation selector for $\lambda$ (converging at rate $n^{-1/3}$). The HAL-MLE $\tilde{P}_n^*$ is targeted and undersmoothed. Both TMLEs should have  $n^{1/2}$-scaled bias converging to zero. We indeed observe that both TMLEs have  negligible bias (bias/SE $<10$), and thereby  provide valid inference.}

\end{table}


\section{Discussion}\label{section12}
The HAL-regularized higher order TMLE  $P_n^{k,(1),*}=\tilde{P}_n^{(1)}\ldots\tilde{P}_n^{(k)}(P_n^0)$ of $P_n^0$, with $\tilde{P}_n^{(j)}(P)$ solving $\tilde{P}_n D^{(j)}_{n,\tilde{P}_n^{(j)}(P)}=0$, satisfies an exact expansion 
\[
\Psi(P_n^{k,(1),*})-\Psi(P_0)=\sum_{j=1}^k\left\{ (\tilde{P}_n-P_0)D^{(j)}_{n,\tilde{P}_n^{(j)}(P_0)}+R^{(j)}_n(\tilde{P}_n^{(j)}(P_0),P_0)\right\}+\tilde{R}_n^k \]
where $\tilde{R}_n^{(k)}\equiv R_n^{(k)}(\tilde{P}_n^{(k)}(P_n^0),\tilde{P}_n)-R_n^{(k)}(\tilde{P}_n^{(k)}(P_0),\tilde{P}_n)$.
In addition, $R_n^{(k)}(\tilde{P}_n^{(k)}(P),\tilde{P}_n)$ is a generalized $k+1$-th order difference of $P$ and $\tilde{P}_n$, thereby only relying on  $\pl \tilde{p}_n-p_0\pl^{k+1}$ and $\pl p_n^0-p_0\pl^{k+1}$ to be small. The leading sum is controlled by the HAL-regularized parametric MLE $\tilde{p}_n^{(j)}(p_0)$ being a good estimator of $p_0$ and $(\tilde{P}_n-P_n)D^{(j)}_{n,\tilde{P}_n^{(j)}(P_0)}$ to be small (i.e. we should also solve the empirical score equation for this MLE).
This is where some undersmoothing of the HAL-MLE $\tilde{P}_n$ or targeting might be needed to make the regularized HAL-MLE $\tilde{p}_n^{(j)}(p_0,\epsilon_n^{(j)}(p_0))$ for a correctly specified one-dimensional parametric model $\tilde{P}_n^{(j)}(P_0,\epsilon)$ behave as the standard MLE
$\tilde{p}_n^{(j)}(p_0,\epsilon_{P_n}(p_0))$. As we showed, this comes down to controlling $(\tilde{P}_n-P_n)D^{(1)}_{P_0}$.
On the other hand, if we use the empirical higher order TMLE then we showed that the regularization bias is essentially negligible allowing one to use the cross-validated HAL-MLE. The empirical higher order TMLE is our recommended choice even though the slightly undersmoothed HAL-regularized higher order TMLE is completely competitive (but more sensitive to choice of undersmoothing). 

Therefore, this exact expansion essentially achieves a level of inference as one would have with a correctly specified parametric model, {\em up till the size of $\tilde{R}_n^{(k)}$}.
We expect that the typical first order pathwise differentiable target parameters will have that $\Psi_n^{(j)}(P)$ is pathwise differentiable up till arbitrarily large $k$. Therefore, we will, in principle, be able to compute $k$-th order TMLE for  large $k$. This should then provide a highly effective tool to fight the {\em finite sample} curse of dimensionality. As $k$ increases, it might be the case that $R_n^{(k)}$ becomes a sum of many terms or that the terms become less stable due to inverse weighting by higher powers of a nuisance  parameter estimator that might be close to zero. Therefore,  for a fixed sample size, at some $k$, there will be no further benefit to increase $k$ further. On the other hand, $\tilde{R}_n^{(k)}$ involves a lot of cancelation due to being a difference of the same sum of terms, and that it consist of terms that themselves involve cancelation of positive and negative contributions to the integrals that make up these terms.
Either way, within the TMLE framework, we can also make $k$ itself a tuning parameter, just as the variation norm bound $C$ in the HAL-MLE $\tilde{P}_n$, and thereby select a correct $k$. Our sense is that w.r.t. finite sample bias a second order TMLE is already a large advance relative to the first order TMLE 

Another important implication of our exact expansion is that it allows us to base inference on $(P_n-P_0)\sum_j D^{(j)}_{n,\tilde{P}_n^{(j)}(P_0)}$, which will still behave as a normal distribution, thereby incorporating the higher order behavior of the $k$-th order TMLE.
Therefore, we expect that, relative to the first order TMLE, the higher order TMLE will not only improve the bias and MSE of the estimator of the target estimand, but will also drastically improve coverage of the confidence intervals. There is also the potential that our exact expansion for the higher order TMLE allows one to analyze problems where the first order canonical gradient equals zero at the true $P_0$, using our second order expansion that includes an empirical mean of the second order efficient influence curve to obtain a  limit distribution. Section \ref{sectioncvhotmle} on cross-validated higher order TMLE shows progress on this front by being able to establish a normal limit distribution for the CV-HOTMLE under a condition allowing the first order influence function at the true $P_0$ being equal to zero.

The (automated) computation of the first and higher order canonical gradients represents an important area of research. Fortunately, it is able to  follow a constructive sequential approach. 
The computation of our higher order canonical gradients $D^{(k+1)}_{n,P}$  corresponds with computing pathwise derivatives of a sequential composition of functions. 
The representation of the next canonical gradient presented in Lemma \ref{lemmaformulacangradient} relies on computing the adjoint of the score operator $A_{n,P}^{(j)}$ and applying it to previous canonical gradient $D^{(j)}_{n,\tilde{P}_n^{(j)}(P)}$. Appendices \ref{sectionfirstgradient} and \ref{section9b} provide the basis of doing these computations with least squares regression or symmetric matrix inversion, thereby opening up the computation of higher order TMLEs with standard machinery, avoiding delicate analytics needed to determine closed forms.

\subsection*{Acknowledgement}

Research reported in this publication was supported by the National Institute Of Allergy And Infectious Diseases of the National Institutes of Health under Award Number R01AI074345. The content is solely the responsibility of the authors and does not necessarily represent the official views of the National Institutes of Health. 

\bibliographystyle{plain}
\bibliography{combined,combined1,combined2,CCbib,TMLEivan,gencens2003,gencens2003logfluct,TMLE1_TLB,TMLE_TLB,gencens2003a,gencens2003b,gencens2003c,gencens2003d,gencens2003e,Grant,higherordertmle,TLB2}


\clearpage

\appendix 
\section*{Appendix}\label{Appendix}
Appendix \ref{AppendixA} establishes that the existence of $d\tilde{P}_n/dP$ implies that $\Psi_n^{(j)}$ is pathwise differentiable at $P$. Appendix \ref{AppendixB} proves the exact expansion Theorem \ref{theorem1}. In the second subsection of Appendix \ref{AppendixB} we  show what happens to the remainder as one increases the order of the TMLE.
In  Appendix \ref{section7} we study till what degree the HAL-regularized higher order TMLE solves the empirical higher order efficient influence curve equations.
In Appendix \ref{Asectionemphotmle} we establish the exact expansion for the empirical $k$-th  order TMLE, and specifically demonstrates its reduction of the HAL-regularization bias term. 
Appendix \ref{AppendixC} establishes the proofs of the Lemmas \ref{lemma5} and \ref{lemma6} concerning the $R_n=R_n^++\dot{R}_n$ decomposition.
In Appendix \ref{AppendixD} we show that higher order canonical gradients are themselves higher order differences. 
Appendix \ref{AppendixF} proves the representation of the $j+1$-th order canonical gradient  in terms of the $j$-th  canonical gradient, as presented in Lemma \ref{lemmaformulacangradient}.
Appendix \ref{AppendixG}  provides the proof of the key lemmas for the second order TMLE  of the integrated square density. Appendix \ref{AppendixH} provide the proofs of the Lemmas for the second order TMLE of the  treatment specific mean outcome.
Appendix \ref{AppendixI} shows how we can target the HAL-MLE in the second order TMLE for the treatment specific mean example to make it behave as the empirical distribution for the desired first and second order efficient score equations (as in Section \ref{sectionthotmle}).
Appendix \ref{AppendixJ} derives the second order TMLE for the treatment specific mean outcome when the outcome is continuous instead of binary.
Appendix \ref{sectionfirstgradient} shows how one can approximate a canonical gradient as a linear combination of  (HAL)-basis functions spanning the tangent space. 
Appendix \ref{section9b} generalizes this to sequential computation of higher order canonical gradients in terms of linear combination of basis functions spanning the tangent space.
Appendix \ref{sectionthotmle} presents an iterative HAL-regularized higher order TMLE involving iteratively targeting the HAL-MLE, thereby making it equivalent with an empirical higher order TMLE. 
Finally, Appendix \ref{sectioncvhotmle} presents a cross-validated higher order TMLE.


\section{$\epsilon_n^{(j)}(P)$ is pathwise differentiable due to using HAL-MLE instead of empirical measure}\label{AppendixA}
Consider $j=1$ first.
Due to using $\tilde{P}_n$, one can show that $\epsilon_n^{(1)}(P)=\arg\min_{\epsilon}\tilde{P}_n L(\tilde{P}^{(1)}(P,\epsilon))$ is  generally a  pathwise differentiable target parameter in $P$ under a weak condition that $
d\tilde{P}_n/dP$ exists for $P\in {\cal M}$. Here $\tilde{P}_n^{(1)}(P,\epsilon)$ is a local least favorable path through $P$ with canonical gradient $D^{(1)}_{P}$ at $\epsilon =0$. For example, $\tilde{p}_n$ might be a Lebesgue density that is bounded away from 0 on the support of $O$, so that this holds.
Since $\epsilon_n^{(1)}(P)$ is the only data dependent part of $\Psi_n^{(1)}(P)=\Psi(\tilde{P}^{(1)}(P,\epsilon_n^{(1)}(P))$, and the remaining part $\Psi^{(1)}_{0}(P)=\Psi(\tilde{P}^{(1)}(P,\epsilon_{P_0}^{(1)}(P))$ (i.e., $\tilde{P}_n$ is replaced by $P_0$)  is generally pathwise differentiable, this will then make $\Psi_n^{(1)}(P)$ pathwise differentiable as well.  The pathwise differentiability of $\epsilon_n^{(1)}(P)$ is understood as follows.

The parameter $\epsilon_n^{(1)}(P)$ can be represented as a solution of the score equation $\tilde{P}_n D^{(1)}_{P,\epsilon_n^{(1)}(P)}=0$, where
$D^{(1)}_{P,\epsilon}= \frac{d}{d\epsilon}L(\tilde{P}^{(1)}(P,\epsilon))$ is the score for $\epsilon$ for this  least favorable path. Consider a path $\{P_{\delta,h}:\delta\}\subset{\cal M}$ through $P$ with score $h$ at $\delta=0$. 
Let $\delta_0=0$.
By the implicit function theorem we have
\[
\frac{d}{d\delta_0}\epsilon_n^{(1)}(P_{\delta_0,h})=-\left\{ \frac{d}{d\epsilon_n^{(1)}(P)}\tilde{P}_nD^{(1)}_{P,\epsilon_n^{(1)}(P)}\right\}^{-1} \frac{d}{d\delta_0}\tilde{P}_n D^{(1)}_{P_{\delta_0,h},\epsilon_n^{(1)}(P)}.\]
In the next formula let $\epsilon_n^{(1)}=\epsilon_n^{(1)}(p)$ and we view $D^{(1)}_{P,\epsilon}$ as a function of the density $p$ when taking the derivative.
The latter factor equals 
\begin{eqnarray*}
\tilde{P}_n\frac{d}{dp}D^{(1)}_{p,\epsilon_n^{(1)}}\left( \frac{d}{d\delta_0}p_{\delta_0,h}\right)&=&
\tilde{P}_n\frac{d}{dp}D^{(1)}_{n,p,\epsilon_n^{(1)}}(ph) \\
&=&P \frac{d\tilde{P}_n}{dP} \frac{d}{dp}D^{(1)}_{n,p,\epsilon_n^{(1)}}(ph).
\end{eqnarray*}
The operator $h\rightarrow A_{n,p}(h)\equiv   \frac{d}{dp}D^{(1)}_{p,\epsilon_n^{(1)}}(ph)$ is an operator on the tangent space $T_P({\cal M})$: $A_{n,p}:T_P({\cal M})\rightarrow L^2(P)$.
Let $A_{n,p}^{\top}:L^2(P)\rightarrow T_P({\cal M})$ be the adjoint of $A_{n,p}$, so that the last expression equals
\[
P A_{n,p}^{\top}\left( \frac{d\tilde{P}_n}{dP} \right) h,\]
which then proves that $\epsilon_n^{(1)}(P)$ is pathwise differentiable at $P$ w.r.t tangent space $T_P({\cal M})$ with canonical gradient
\[-\left\{ \frac{d}{d\epsilon_n^{(1)}(P)}\tilde{P}_nD^{(1)}_{P,\epsilon_n^{(1)}(P)}\right\}^{-1} 
A_{n,p}^{\top}\left( \frac{d\tilde{P}_n}{dP})\right).
\]
The key that makes it work is that $d\tilde{P}_n/dP$ exists as an element in $L^2(P)$: i.e., we need that $\int (\tilde{p}_n/p)^2 dP<\infty$. Note that if $\tilde{P}_n$ is replaced by $P_n$, then $\epsilon_n^{(1)}(p)$ would not be pathwise differentiable.

Consider now a general $j\in \{1,2,\ldots\}$ and assume that we already established that $\Psi_n^{(j-1)}(P)$ is pathwise differentiable at $P$. We have $\Psi_n^{(j)}(P)=\Psi_n^{(j-1)}(\tilde{P}_n^{(j)}(P,\epsilon_n^{(j)}(P)))$. Since $\Psi_n^{(j-1)}(\tilde{P}_n^{(j)}(P,\epsilon))$ at a fixed $\epsilon$ is pathwise differentiable due to pathwise differentiability of $\Psi_n^{(j-1)}(P)$ and of the smoothness of $P\rightarrow \tilde{P}_n^{(j)}(P,\epsilon)$, it follows that it suffices to show that $\epsilon_n^{(j)}(P)$ is pathwise differentiable at $P$. This is proved by copying the proof above replacing $^{(1)}$ by $^{(j)}$ throughout.

\section{Proof of Theorem \ref{theorem1} and understanding of the exact expansion.}\label{AppendixB}
Regarding the claimed expansion, 
consider first the case $k=2$. We have
\[
\begin{array}{l}
\Psi(\tilde{P}_n^{(1)}\tilde{P}_n^{(2)}(P_n^0))-\Psi(P_0)=\{\Psi(\tilde{P}_n^{(1)}(P_0))-\Psi(P_0)\}
\\
+
\{\Psi(\tilde{P}_n^{(1)}\tilde{P}_n^{(2)}(P_n^0))-\Psi(\tilde{P}_n^{(1)}(P_0))\}\\
=\Psi(\tilde{P}_n^{(1)}(P_0))-\Psi(P_0)+\Psi_n^{(1)}(\tilde{P}_n^{(2)}(P_n^0))-\Psi_n^{(1)}(P_0) \\.
\end{array}
\]
For $k=3$, we have
\[
\begin{array}{l}
\Psi(\tilde{P}_n^{(1)}\tilde{P}_n^{(2)}\tilde{P}_n^{(3)}(P_n^0))-\Psi(P_0)=
\Psi(\tilde{P}_n^{(1)}(P_0)-\Psi(P_0)\\
+
\Psi(\tilde{P}_n^{(1)}\tilde{P}_n^{(2)}(P_0))-\Psi(\tilde{P}_n^{(1)}(P_0))
+\Psi(\tilde{P}_n^{(1)}\tilde{P}_n^{(2)}\tilde{P}_n^{(3)}(P_n^0))-\Psi(\tilde{P}_n^{(1)}\tilde{P}_n^{(2)}(P_0))\\
=
\Psi(\tilde{P}_n^{(1)}(P_0))-\Psi(P_0)+
\Psi_n^{(1)}(\tilde{P}_n^{(2)}(P_0))-\Psi_n^{(1)}(P_0)
+\Psi_n^{(2)}(\tilde{P}_n^{(3)}(P_n^0))-\Psi_n^{(2)}(P_0).
\end{array}
\]
For general $k$ we have
\[
\begin{array}{l}
\Psi(\tilde{P}_n^{(1)}\ldots \tilde{P}_n^{(k-1)}\tilde{P}_n^{(k)}(P_n^0))-\Psi(P_0)
=
\sum_{j=1}^{k-1}\Psi_n^{(j-1)}(\tilde{P}_n^{(j)}(P_0))-\Psi_n^{(j-1)}(P_0)\\
+
\Psi_n^{(k-1)}(\tilde{P}_n^{(k)}(P_n^0))-\Psi_n^{(k-1)}(P_0).
\end{array}
\]
Thus, this decomposition decomposes the $k$-th order TMLE as a sum of $k-1$ plug-in TMLEs $\Psi_n^{(j-1)}(\tilde{P}_n^{(j)}(P_0))$ of $\Psi_n^{(j-1)}(P_0)$ using initial estimator $\tilde{P}_n^{(j)}(P_0)$ in the TMLE update $\tilde{P}_n^{(j)}(P)$, $j=1,\ldots,k-1$, and a standard (i.e., no HAL-MLE smoothing) $k$-th plug-in TMLE $\Psi_n^{(k-1)}(\tilde{P}_n^{(k)}(P_n^0))$ of $\Psi_n^{(k-1})(P_0)$ using as initial estimator $P_n^0$ in the TMLE update $P_n^{(k)}(P)$. For any TMLE $\Psi_n^{(j-1)}(\tilde{P}_n^{(j)}(P_0))$ of $\Psi_n^{(j-1)}(P_0)$ we have the exact expansion implied by the definition of $R_n^{(j)}$: for $j=1,\ldots,k-1$
\[
\begin{array}{l}
\Psi_n^{(j-1)}(\tilde{P}_n^{(j)}(P_0))-\Psi_n^{(j-1)}(P_0)=(P_n-P_0)D^{(j)}_{n,\tilde{P}_n^{(j)}(P_0)}+R_n^{(j)}(\tilde{P}_n^{(j)}(P_0),P_0)\\
\hfill -P_n D^{(j)}_{n,\tilde{P}_n^{(j)}(P_0)}\end{array}
\] 
Similarly for $\Psi_n^{(k-1)}(\tilde{P}_n^{(k)}(P_n^0))-\Psi_n^{(k-1)}(P_0)$ we have
\[
\begin{array}{l}
\Psi_n^{(k-1)}(\tilde{P}_n^{(k)}(P_n^0))-\Psi_n^{(k-1)}(P_0)=(P_n-P_0)D^{(k)}_{n,\tilde{P}_n^{(k)}(P_n^0)}\\
\hfill +R_n^{(k)}(\tilde{P}_n^{(k-1)}(P_n^0),P_0)-
P_n D^{(k)}_{n,\tilde{P}_n^{(k)}(P_n^0)}.\end{array}
\]
Plugging this in the above expansion yields the result as stated in the theorem. $\Box$

\subsection{Understanding that the $j+1$-th order TMLE reduces the remainder relative to the $j$-th order TMLE}
Suppose that $\tilde{P}_nD^{(1)}_{\tilde{P}_n^{(1)}(P^0)}=0$. Then,
\[
\begin{array}{l}
\Psi(\tilde{P}_n^{(1)}(P^0))-\Psi(P_0)=(\tilde{P}_n-P_0)D^{(1)}_{\tilde{P}_n^{(1)}(P^0)}+R^{(1)}(\tilde{P}_n^{(1)}(P^0),P_0)\\
\equiv \tilde{P}_n D^{(1)}_{\tilde{P}_n^{(1)}(P_0)}+\bar{R}^{(1)}(\tilde{P}_n^{(1)}(P^0),P_0),
\end{array}
\]
where we defined the total remainder 
\[
\bar{R}^{(1)}(\tilde{P}_n^{(1)}P^0,P_0)=(\tilde{P}_n-P_0)\{D^{(1)}_{\tilde{P}_n^{(1)}P^0}-D^{(1)}_{\tilde{P}_n^{(1)}(P_0)}\}+R^{(1)}(\tilde{P}_n^{(1)}P^0,P_0).\]
We have
\[
\Psi_n^{(1)}(\tilde{P}_n^{(2)}(P^0))-\Psi_n^{(1)}(P_0)=
\bar{R}^{(1)}(\tilde{P}_n^{(1)}(P_n^{2,*}),P_0)-\bar{R}^{(1)}(\tilde{P}_n^{(1)}(P_0),P_0),\]
where  $P_n^{2,*}=\tilde{P}_n^{(2)}(P^0)$.  Assume again that $\tilde{P}_nD^{(2)}_{n,P_n^{2,*}}=0$.
That is, $\Psi_n^{(1)}(\tilde{P}_n^{(2)}(P^0))$ acts as the TMLE of $\bar{R}^{(1)}(\tilde{P}_n^1(P),P_0)$ at $P=P_0$.
So, for the second order TMLE we have
\[
\begin{array}{l}
\Psi(\tilde{P}_n^{(1)}P_n^{2,*})-\Psi(P_0)\\
=
\tilde{P}_n D^{(1)}_{\tilde{P}_n^{(1)}(P_0)}+
\Psi_n^{(1)}(P_n^{(2,*})-\Psi_n^{(1)}(P_0)+\bar{R}^{(1)}(\tilde{P}_n^{(1)}(P_0),P_0)\\
=\tilde{P}_nD^{(1)}_{\tilde{P}_n^{(1)}(P_0)}+
(\tilde{P}_n-P_0)D^{(2)}_{n,P_n^{2,*}}+\bar{R}^{(1)}(\tilde{P}_n^{(1)}(P_0),P_0)+
R_n^{(2)}(P_n^{2,*},P_0).
\end{array}
\]
So by replacing $P^0$ in the first order TMLE by the TMLE $P_n^{2,*}$ we reduced the exact total remainder $\bar{R}^{(1)}(\tilde{P}_n^{(1)}(P^0),P_0)$ to
$(\tilde{P}_n-P_0)D^{(2)}_{n,P_n^{2,*}}+R_n^{(2)}(P_n^{2,*},P_0)$, plus a perfect remainder $\bar{R}^{(1)}(\tilde{P}_n^{(1)}(P_0),P_0)$.
In general, by replacing $P^0$ in $\tilde{P}_n^{(j)}(P^0)$ targeting $\Psi_n^{(j-1)}(P_0)$  by the TMLE $P_n^{j+1,*}$ targeting $\Psi_n^{(j)}(P_0)$ we reduced the exact total remainder term $\bar{R}_{n}^{(j)}(\tilde{P}_n^{(j)}(P^0),P_0)$ by
$(\tilde{P}_n-P_0)D^{(j+1)}_{n,P_n^{j+1,*}}+R_n^{(j+1)}(P_n^{j+1,*},P_0)$, plus a perfect remainder $\bar{R}_{n}^{(j)}(\tilde{P}_n^{(j)}(P_0),P_0)$.
Clearly, each step replacing $P^0$ by the corresponding TMLE represents an advance in reducing the remainder $\bar{R}^{(1)}(\tilde{P}_n^{(1)}(P^0),P_0)$, given that we showed that $R_n^{(j+1)}(P_n^{j+1,*},P_0)$ is a higher order difference than $R_n^{(j)}(P_n^{j,*},P_0)$  (Lemmas \ref{lemma3}, \ref{lemma2} and \ref{lemmaa2}).

The reduction is due to $\bar{R}_{n}^{(j)}(\tilde{P}_n^{(j)}(P_n^{j+1,*}),P_0)$  being a TMLE of $\bar{R}_{n}(\tilde{P}_n^{(j)}(P_0),P_0)$ instead of using the non-targeted $\bar{R}_{n}^{(j)}(\tilde{P}_n^{(j)}(P^0),P_0)$. By $P_n^{j+1,*}$ being a TMLE, the directional derivative  of $\bar{R}_n^{(j)}(\tilde{P}_n^{(j)}(P),P_0)$ at $P=P_n^{j+1,*}$ in the direction of $P_n^{j+1,*}-P_0$ becomes $o_P(n^{-1/2})$-noise, which implies it becomes a higher order difference as shown in Lemma \ref{lemmaa2}. Moreover, in essence, as we showcase in detail in a Section \ref{section3}, the TMLE minimizes the exact total remainder under the constraint that it can move in direction of $P_0$ (i.e., it is better than solving the derivative equation, which is just an implication of being an optimizer).

 \section{Solving the empirical higher order efficient score equations for the HAL-regularized $k$-th order TMLE.}
 \label{section7}

Assume that the $k$-th order TMLE is defined to satisfy $\tilde{P}_n D^{(j)}_{n,\tilde{P}_n^{(j)}(P)}=0$ exactly, for all $j=1,\ldots,k$. (The more general case is discussed in the last subsection  of this section)  By using local least favorable paths and defining $\epsilon_n^{(j)}(P)$ as the solution of $\tilde{P}_n D^{(j)}_{n,\tilde{P}_n^{(j)}(P,\epsilon)}=0$, we have that these $\tilde{P}_n$-equations are exactly solved. 
 Then, we have the corresponding exact expansion as presented in Theorem \ref{theorem2} involving 1) perfect remainders $R_n^{(j)}(\tilde{P}_n^{(j)}(P_0),P_0)$; 2) nice HAL-empirical means $(\tilde{P}_n-P_0)D^{(j)}_{n,\tilde{P}_n^{(j)}(P_0)}$; as well as a difference of the $k+1$-th order remainder $R_n^{(k)}(\tilde{P}_n^{(k)}(P),\tilde{P}_n)$ at $P=P_0$ and $P=P_n^0$. For the sake of statistical inference, we will  need that  $(\tilde{P}_n-P_n )D^{(1)}_{\tilde{P}_n^{(1)}(P_0)}=o_P(n^{-1/2})$. For the sake of bringing the HAL-empirical means of the higher order efficient influence curves into the statistical inference (see next section), one also wants $(\tilde{P}_n-P_n )D^{(j)}_{n,\tilde{P}_n^{(j)}(P_0)}=o_P(n^{-1/2})$ for $j=2,\ldots,k$. 
 Since $D^{(j)}_P$ for $j\geq 2$ involves a $j-1$-th order difference of $P$ and $\tilde{P}_n$, in practice it mostly comes down to having to solve $P_n D^{(1)}_{n,\tilde{P}_n^{(1)}(P_0)}\approx 0$, since the higher order TMLE score values will typically even be smaller.  
 By Lemma \ref{lemmaexactepsn1} we know that the latter comes down to controlling $(\tilde{P}_n-P_n)D^{(1)}_{P_0}$.



{\bf Tuning the TMLE score values by the $L_1$-norm in the HAL-MLE:}
To start with we point out that a sensible method for undersmoothing the HAL-MLE  would be to  increase the variation norm of the HAL-MLE (i.e, the $L_1$-norm of the coefficient vector) till 
\[
\frac{\max_{j\in \{1,\ldots,k-1\}}\pl  P_n D^{(j)}_{n,\tilde{P}_n^{(j)}\ldots \tilde{P}_n^{(k)}(P_n^0)}\pl }{\sigma_{jn}}=r(n)\]
 for a user supplied number $r(n)$ such as $n^{-1/2}(\log n)^{-1}$, where $\sigma^2_{jn}$ is a standard estimator of the variance $P_0\{D^{(j)}_{n,\tilde{P}_n^{(j)}(P_0)}\}^2$. One would only consider $L_1$-norms larger than the cross-validation selector. The cut-off is chosen so that the size is negligible for purpose of estimation and inference. As one increases the variation norm, the span of the scores solved by the HAL-MLE grows till all scores, so that at some large enough value the desired cut-off is achieved. This would be  conservative (therefore  undersmooth more than actually required) in the sense that we are only requiring that $P_nD^{(j)}_{n,\tilde{P}_n^{(j)}(P_0)}\approx 0$ at the initial (and thus fixed) $P_0$. The latter is easier since it concerns the TMLE for a correctly specified one dimensional least favorable model using as initial/off-set $P_0$ itself. Therefore, a better method  for selecting the $L_1$-norm $C$ of the HAL-MLE might be to  increase the $L_1$ norm, starting at its cross-validation selector $C_{n,cv}$, till   the $k$-th order plug-in TMLE $\Psi(P_{n,C}^{k,(1),*})$ reaches a plateau, as a function of the $L_1$-norm $C$ of the HAL-MLE $\tilde{P}_{n,C}$ used in each of the TMLEs $\tilde{P}_n^{(j)}$. Formal results for such a plateau selector can be obtained (see e.g., \citep{Cai&vanderLaan20})
 
 \subsection{Targeting the HAL-MLE to assist in solving the desired score equations}
 We wish to arrange that $(\tilde{P}_n-P_n)D^{(1)}_{\tilde{P}_n^{(1)}(P_0)}$ is small. This can be arranged by undersmoothing $\tilde{P}_n$. We could also target $\tilde{P}_n$ to solve
 $(\tilde{P}_n-P_n)D^{(1)}_{P^0}=0$ for a particular estimator $P^0$ (e.g. $\tilde{P}_n^{(1)}(P_n^0)$). This also assist in making $(\tilde{P}_n-P_n)D^{(1)}_{\tilde{P}_n^{(1)}(P_0)}$ small, even without undersmoothing $\tilde{P}_n$. We can still index this targeted version of $\tilde{P}_n$ with its $L_1$-norm and tune it as above with the plateau selector. Since $(\tilde{P}_n-P_n)D^{(1)}_{P^0}=P_n\{D^{(1)}_{P^0}-\tilde{P}_n D^{(1)}_{P^0}\}$ is just a score equation at $\tilde{P}_n$, this just requires targeting $\tilde{P}_n$ so that it solves this score equation (e.g., with an TMLE-update or targeted HAL-MLE). This targeting of the HAL-MLE is presented in Section \ref{sectionthotmle} and we  demonstrate such an iterative  higher order TMLE, involving iteratively targeting the HAL-MLE, in the treatment specific mean example (see Appendix \ref{AppendixJ}).
 
\subsection{Theoretical understanding why the empirical means of efficient scores  are negligible, when  undersmoothing the HAL-MLE.}
After having pointed out that this TMLE score values  are easily controlled by the user, we now proceed with analyzing conditions under which these score equations are solved at the desired precision (i.e., negligible for practical performance of the $k$-th order TMLE).
We focus on the dominating $j=1$-term, as explained above.
The asymptotic efficiency of undersmoothed plug-in HAL-MLE of target features of the data distribution, analogue to our presentation above, is presented in detail in \citep{vanderLaan19efficient,vanderLaan&Benkeser&Cai19}. This is directly relevant since we want to achieve that $\epsilon_{\tilde{P}_n}^{(1)}(P_0)$ is an asymptotically efficient estimator of $\epsilon_{P_0}^{(1)}(P_0)$, so that it is asymptotically equivalent with $\epsilon_{P_n}^{(1)}(P_0)$. 
This viewpoint is discussed in detail in Section \ref{sectionthotmle}, where it provides the rational and motivation for targeting the HAL-MLE $\tilde{P}_n$ in the higher order TMLE method.


{\bf $(\tilde{P}_n-P_n)D^{(1)}_{P_0}$ is a second order difference in discrepancy of HAL-MLE and empirical measure:}
An HAL-MLE $\tilde{P}_n^C=P_{\theta_n^C}$ solves a class of scores $P_n S_{h,\tilde{P}_n^C}=0$, where $S_{h,\tilde{P}_n^C}=\frac{d}{d\delta_0} L(P_{\theta_{n,\delta_0,h}^C})$, along all paths $\{P_{\theta_{\delta,h}}:\delta\}\subset {\cal M}$ that preserve the variation norm $C$ (which represents a single constraint). Given the HAL-MLE fit $\theta_n=\sum_j\beta_n(j)\phi_{u_j}$, this includes the scores $S_{h,\tilde{P}_n^C}$ generated by paths $(1+\epsilon h(j))\beta_n(j)$ for any vector $h$ with 
$\sum_j h(j)\mid \beta_n(j)\mid =0$. 
These scores approximate the unconstraint score equations  of the non-zero coefficients $\beta_n$. The linear span of all these scores $S_{h,\tilde{P}_n^C}$ approximate the whole tangent space at $\tilde{P}_n$ as $n$ increases, or, more generally, as the number of non-zero coefficients grows to infinity. One could define $\tilde{D}^{(1)}_n$ as the minimizer of $P_n \{D^{(1)}_{P_0}-S_{h,\tilde{P}_n^C}\}^2 $ over all the scores $S_{h,\tilde{P}_n^C}$. 
Since $P_n \tilde{D}^{(1)}_n=\tilde{P}_n^C \tilde{D}^{(1)}_n=0$, we have \begin{equation}\label{scoreeqna}
(\tilde{P}_n^C-P_n)D^{(1)}_{P_0} =(\tilde{P}_n^C-P_n)\{D^{(1)}_{P_0}-\tilde{D}^{(1)}_n\}.\end{equation}
As $C$ increases, both $(\tilde{P}_n^C-P_n)$ decreases as well as the oracle approximation $\tilde{D}^{(1)}_n$ of $D^{(1)}_{P_0}$ improves quickly due to the dimension of the linear span of scores growing towards $n$. 
Importantly, both approximations are not driven by how well $\tilde{p}_n$ estimates $p_0$ (i.e., the Kullback-Leibler divergence $d_0(\tilde{P}_n,P_0)$) but about resembling the empirical measure. In addition, the number of basis functions with non-zero coefficients in an HAL-MLE fit will be a proportion of $n$ (and maximal $n-1$) even as the $L_1$-norm remains bounded, so there is no need to let the $L_1$-norm go to infinity to obtain the full benefit of undersmoothing. 
Therefore,  this error $P_n D^{(1)}_{P_0}$ appears to be only {\em weakly} affected by the curse of dimensionality. 
Finally, as pointed out above, we can verify on data how far to undersmooth to achieve the goal of being similar to the simple well behaved $\epsilon_{P_n}^{(1)}(P_0)$.

{\bf Even without undersmoothing, it will be $O_P(n^{-2/3})$:}
Even without undersmoothing this term (\ref{scoreeqna}) will generally be $O_P(n^{-2/3}(\log n)^d)$, but this bound is not taking advantage of the above bound that it is not about estimating $p_0$, allowing us to control the size by undersmoothing $\tilde{P}_n$, as argued above.
This is seen by noting that (suppressing $C$)
\[
(\tilde{P}_n-P_n)(D^{(1)}_{P_0}-\tilde{D}^{(1)}_n)=
(\tilde{P}_n-P_0)(D^{(1)}_{P_0}-\tilde{D}^{(1)}_n)-(P_n-P_0)(D^{(1)}_{P_0}-\tilde{D}^{(1)}_n).\]
The first term can be bounded by Cauchy-Schwarz inequality by $\pl \tilde{p}_n-p_0\pl_{\mu}\pl f_n\pl_{\mu}$, where $f_n\equiv (D^{(1)}_{P_0}-\tilde{D}^{(1)}_n) $.
Due to the rate of convergence of HAL-MLE we have $\pl \tilde{p}_n-p_0\pl_{\mu}=O_P(n^{-1/3})(\log n)^{d/2})$, under the assumption that $p_0>\delta>0$ for some $\delta>0$.
The second term is $O_P(n^{-2/3}(\log n)^d)$ if the $L^2(P_0)$-norm of $f_n\equiv (D^{(1)}_{P_0}-\tilde{D}^{(1)}_n) $ is $O_P(n^{-1/3}(\log n)^{d/2})$ and $f_n$ is a multivariate cadlag function with a universal bound on its variation norm, and, either way, even when we only have $\pl f_n\pl_{P_0}\rightarrow_p 0$, it is $o_P(n^{-1/2})$. Due to being able to optimize the choice $\tilde{D}^{(1)}_n$ (i.e., select $S_{h,\tilde{P}_n}$ closest to $D^{(1)}_{P_0}$), one expects that $\pl f_n\pl_{P_0}$ converges  faster to zero than $\tilde{p}_n-\tilde{p}_0$ since the linear span of the scores, which are a linear transformation of the spline-basis functions with non-zero coefficient in the HAL-MLE fit, will generally approximate this function $D^{(1)}_{P_0}$  faster than the HAL-MLE $\tilde{p}_n$ approximates $p_0$ (HAL involves estimation, while $f_n$ is defined by an oracle fit). However, if somehow $o\rightarrow D^{(1)}_{P_0}(o)$ is a more complex function than $\theta_0$ (e.g., $d\theta_0(u)=0$ for a set of knot points $u$ while these are required to approximate $D^{(1)}_{P_0}$), then some undersmoothing may be needed. 

So we conclude that without undersmoothing one can already bound  the desired $(\tilde{P}_n-P_n)\{D^{(1)}_{P_0}-\tilde{D}^{(1)}_n)\}$ by $O_P(\pl \tilde{p}_n-p_0\pl_{\mu}\pl f_n\pl_{\mu})=O_P(n^{-1/3}(\log n)^{d/2}\pl f_n\pl_{\mu})$, plus a term that is guaranteed $o_P(n^{-1/2})$ (and generally even smaller). In addition, $\pl f_n\pl_{\mu}$ will generally have at minimal the same rate of convergence as the HAL-MLE, possibly relying on some undersmoothing of the HAL-MLE. Moreover, the norm $\pl f_n\pl_{\mu}$ is decreasing in choice of $L_1$-norm $C$ and can thereby be controlled to be finite sample small.


{\bf Targeting the HAL-MLE to further reduce it to a third order difference:}
 In Appendix \ref{sectionthotmle} we present an iterative algorithm for a higher order TMLE that also involves targeting of the HAL-MLE $\tilde{P}_n$ so that it becomes a TMLE itself tailored so that $P_n D^{(1)}_{P_n^{(k),(1),*}}=0$, while we already have $\tilde{P}_n^* D^{(1)}_{P_n^{(k),(1),*}}=0$. 
  We also point out that one might use a simple external targeting of $\tilde{P}_n$ so that
$P_n D^{(1)}_{\tilde{P}_n^*}=0$, while we also have $\tilde{P}_n^* D^{(1)}_{\tilde{P}_n^*}=0$.
Let's consider the latter case for simplicity, but the same argument applies to the iterative targeting of HAL-MLE.
We now have
 \[
 (\tilde{P}_n-P_n)D^{(1)}_{P_0}=(\tilde{P}_n-P_n)\{D^{(1)}_{P_0}-D^{(1)}_{\tilde{P}_n^*}\}.\]
  If the targeting $\tilde{P}_n^*$ is done to not affect the score equations solved by the HAL-MLE $\tilde{P}_n$ (by doing a constrained HAL-MLE, so that $\tilde{P}_n^*$ is still an MLE), we can then still subtract $f_n$ defined as a best approximation in the linear span of scores solved by $\tilde{P}_n^*$ of $D^{(1)}_{P_0}-D^{(1)}_{\tilde{P}_n^*}$, so that we obtain
 \[
 (\tilde{P}_n-P_n)D^{(1)}_{P_0}= (\tilde{P}_n-P_n)\{D^{(1)}_{P_0}-D^{(1)}_{\tilde{P}_n^*}-f_n\}.\]
 This is now a third order difference still of the nice type that for two of the differences it is about $\tilde{P}_n$ approximating $P_n$ and only the difference $\tilde{P}_n^*$ versus $P_0$ concerns 
 $d_0(\tilde{P}_n^*,P_0)$. 
 

\subsection{Behavior of the empirical mean of the higher order efficient influence curve when using a single step local least favorable TMLE update}
Even when we use a $k$-th order TMLE based on single step MLE-updates according to  local least favorable paths, thereby not guaranteeing that the $\tilde{P}_n$ score equations are solved, it is easy to show that $\tilde{P}_n D^{(1)}_{\tilde{P}_n^{(1)}(P_0)}\approx O_P(\epsilon_n^{(1),2}(P_0))$, so that Lemma \ref{lemmaexactepsn1} show that   this will still be bounded by the maximum of $(\tilde{P}_n-P_n)^2 D^{(1)}_{P_0}$
and an $O_P(n^{-1})$-term.
Our analysis above shows then again $P_n D^{(1)}_{\tilde{P}_n^{(1)}(P_0)}$ is bounded by $(\tilde{P}_n-P_n)D^{(1)}_{P_0}$.

\section{Establishing exact expansion for the higher order TMLE that uses empirical TMLE-updates}\label{Asectionemphotmle}
Consider the case that $\tilde{P}_n^{(j)}(P)=\tilde{P}(P,\tilde{\epsilon}_n^{(j)}(P))$, where $\tilde{\epsilon}_n^{(j)}(P)$ is the solution of $\tilde{P}_n D^{(j)}_{n,\tilde{P}_n^{(j)}(P,\epsilon)}=0$, $j=1,\ldots,k$. 
Let $P_n^{(j)}(P)=\tilde{P}(P,\epsilon_n^{(j)}(P))$, where $\epsilon_n^{(j)}(P)$ is the solution of ${P}_n D^{(j)}_{n,\tilde{P}_n^{(j)}(P,\epsilon)}=0$, $j=1,\ldots,k$. In other words, $P_n^{(j)}(P)$ represents the regular empirical TMLE-updates, still using the canonical gradients $D^{(j)}_{n,P}$ that depend on the HAL-MLE $\tilde{P}_n$. 

In this section we consider the empirical $k$-th order TMLE defined by $\Psi(P_n^{(1)}\ldots P_n^{(k)}(P_n^0))$, and contrast it to the HAL-regularized $k$-th order TMLE $\Psi(\tilde{P}_n^{(1)}\ldots\tilde{P}_n^{(k)}(P_n^0))$. Specifically, we will show that it improves the undersmoothing term by an order of difference, thereby further minimizing the need for undersmoothing the HAL-MLE.In particular, we present the proof of Theorem  \ref{theoremempiricalhotmle} that provides the exact expansion for the empirical $k$-th order TMLE.

\subsection{Exact expansion for the empirical higher order TMLE}
Note that $\tilde{P}_n^{(1)}\ldots\tilde{P}_n^{(k)}(P_n^0)$ is a function of $P_n^0$, a vector of MLE-fits of the fluctuation parameters along least favorable paths, where the least favorable paths depend on $\tilde{P}_n$ (due to their scores depending on $\tilde{P}_n$).
Therefore we can represent this HAL-regularized  $k$-th order TMLE as a function  $\Psi_n(\tilde{\epsilon}_n)$ of a vector $\tilde{\epsilon}_n=(\tilde{\epsilon}_n^{(1)},\ldots,\tilde{\epsilon}_n^{(k)})$, where $\Psi_n$ is indexed by $(P_n^0,\tilde{P}_n)$. 
Then, the empirical $k$-th order TMLE is given by $\Psi_n(\epsilon_n)$: that is, it only differs from the HAL-regularized $k$-th order TMLE by its values of the fluctuation parameters, since it uses the same least favorable paths and same initial estimator $P_n^0$.
We already provided an exact expansion for $\Psi_n(\tilde{\epsilon}_n)$ and we now want to obtain an exact expansion for $\Psi_n(\epsilon_n)$. 
 We have
 \begin{eqnarray*}
 \Psi_n(\epsilon_n)-\Psi(P_0)&=&\Psi_n(\tilde{\epsilon}_n)-\Psi(P_0)+
 \Psi_n(\epsilon_n)-\Psi_n(\tilde{\epsilon}_n)\\
&& =\sum_{j=1}^k (\tilde{P}_n-P_0)D^{(j)}_{n,\tilde{P}_n^{(j)}(P_0)}+R^{(j)}_n(\tilde{P}_n^{(j)}(P_0),P_0)\\
 &&+ R^{(k)}_n(\tilde{P}_n^{(k)}(P_n^0),\tilde{P}_n)-R^{(k)}_n(\tilde{P}_n^{(k)}(P_0),\tilde{P}_n)\\
&&+ \Psi_n(\epsilon_n)-\Psi_n(\tilde{\epsilon}_n)\\
&=&\sum_{j=1}^k ({P}_n-P_0)D^{(j)}_{n,\tilde{P}_n^{(j)}(P_0)}+R^{(j)}_n(\tilde{P}_n^{(j)}(P_0),P_0)\\
 &&+ R^{(k)}_n(\tilde{P}_n^{(k)}(P_n^0),\tilde{P}_n)-R^{(k)}_n(\tilde{P}_n^{(k)}(P_0),\tilde{P}_n)\\
&&+ \Psi_n(\epsilon_n)-\Psi_n(\tilde{\epsilon}_n)-\sum_{j=1}^k(P_n-\tilde{P}_n)D^{(j)}_{n,\tilde{P}_n^{(j)}(P_0)}.
\end{eqnarray*}
So we obtain the same exact empirical process expansion with a $k+1$-th order difference
as with the HAL-regularized $k$-th order TMLE but a HAL-regularization bias term that is the difference of
 $\Psi_n(\epsilon_n)-\Psi_n(\tilde{\epsilon}_n)$ and the HAL-regularization term 
 $\sum_{j=1}^k(\tilde{P}_n-P_n)D^{(j)}_{n,\tilde{P}_n^{(j)}(P_0)}$ from the HAL-regularized $k$-th order TMLE.

 However, this additional $\Psi_n(\epsilon_n)-\Psi_n(\tilde{\epsilon}_n)$ 
  actually counteracts the other bias term. In fact, one can view the resulting HAL-regularization bias term  as $\Psi_n(\epsilon_n) -\Psi_n(\tilde{\epsilon}_n)$ minus its first order Tailor approximation, thereby making it essentially a  second order difference  in $\epsilon_n-\tilde{\epsilon}_n$. Thus the extra term $\Psi_n(\epsilon_n)-\Psi_n(\tilde{\epsilon_n})$ can be viewed as a bias reduction of the regularized higher order TMLE. 
 One way to understand this bias reduction is the following. 
 If we start the regularized higher order TMLE with initial $P_n^0=\tilde{P}_n$, then it returns $\tilde{P}_n$ itself.
This shows that the the HAL-MLE $\tilde{P}_n$ is itself a  regularized $k$-th order TMLE
satisfying  this exact expansion with the HAL-regularization bias term $\sum_{j=1}^k (\tilde{P}_n-P_n)D^{(j)}_{n,\tilde{P}_n^{(j)}(P_0)}$, whose dominating term is given by $(\tilde{P}_n-P_n)D^{(1)}_{\tilde{P}_n^{(1)}(P_0)}$. On the other hand, if we plug $\tilde{P}_n$ as initial in the empirical higher order TMLE then it will sequentially target $\tilde{P}_n$ (and thus removes bias at each step) towards the target parameters $\Psi_n^{(k-1)}(P_0)$, $\ldots$, $\Psi(P_0)$, respectively. This must thus mean that it is removing bias from  $\sum_j (\tilde{P}_n-P_n)D^{(j)}_{n,\tilde{P}_n^{(j)}(P_0)}$. 
 
The above explanation of sequential bias reduction of the empirical higher order TMLE versus regularized higher order TMLE corresponds precisely with the following decomposition of this bias reduction, where each term represents the bias reduction due to  replacing the $j$-th regularized MLE by the $j$-th empirical MLE:
\[
\begin{array}{l}
\Psi(P_n^{(1)}\ldots P_n^{(k)}(P_n^0))-\Psi(\tilde{P}_n^{(1)}\ldots\tilde{P}_n^{(k)}(P_n^0))\\
=
\sum_{j=0}^{k-1} \Psi(\tilde{P}_n^{(1)}\ldots\tilde{P}_n^{(j)}P_n^{(j+1)}\ldots P_n^{(k)}(P_n^0))-
\Psi(\tilde{P}_n^{(1)}\ldots\tilde{P}_n^{j+1}P_n^{(j+2)}\ldots P_n^{(k)}(P_n^0) )\\
=\sum_{j=0}^{k-1}\Psi_n^{(j)}(P_n^{(j+1)}\ldots P_n^{(k)}(P_n^0))-\Psi_n^{(j)}(\tilde{P}_n^{(j+1)}P_n^{(j+2)}\ldots P_n^{(k)}(P_n^0))\\
=\sum_{j=0}^{k-1}\{ \Psi_n^{(j)}(P_n^{(j+1)}\ldots P_n^{(k)}(P_n^0))-\Psi_n^{(j)}(P_0)\}\\
-\sum_{j=0}^{k-1}
\{ \Psi_n^{(j)}(\tilde{P}_n^{(j+1)}P_n^{(j+2)}\ldots P_n^{(k)}(P_n^0))-\Psi_n^{(j)}(P_0)\}\\
=\sum_{j=1}^{k}\left\{ (P_n-P_0)D^{(j)}_{n,P_n^{(j)}\ldots P_n^{(k)}(P_n^0)}+R_n^{(j)}(P_n^{(j)}\ldots P_n^{(k)}(P_n^0),P_0)\right\}\\
-\sum_{j=1}^k\left\{ (\tilde{P}_n-P_0) D^{(j)}_{n,\tilde{P}_n^{(j)}P_n^{(j+1)}\ldots P_n^{(k)}(P_n^0)}\}+R_n^{(j)}(\tilde{P}_n^{(j)}P_n^{(j+1)}\ldots P_n^{(k)}(P_n^0),P_0)\right\} .
\end{array}
\]
Adding the last expression to the exact expansion for the regularized higher order TMLE yields the following exact expansion for the empirical $k$-th order TMLE:
\[
\begin{array}{l}
\Psi(P_n^{(1)}\ldots P_n^{(k)}(P_n^0))-\Psi(P_0)\\
=
\sum_{j=1}^k (\tilde{P}_n-P_0)D^{(j)}_{n,\tilde{P}_n^{(j)}(P_0)}+R^{(j)}_n(\tilde{P}_n^{(j)}(P_0),P_0)\\
 + R^{(k)}_n(\tilde{P}_n^{(k)}(P_n^0),\tilde{P}_n)-R^{(k)}_n(\tilde{P}_n^{(k)}(P_0),\tilde{P}_n)\\
 +\sum_{j=1}^{k}\left\{ (P_n-P_0)D^{(j)}_{n,P_n^{(j)}\ldots P_n^{(k)}(P_n^0)}+R_n^{(j)}(P_n^{(j)}\ldots P_n^{(k)}P_0)\right\}\\
-\sum_{j=1}^k\left\{ (\tilde{P}_n-P_0) D^{(j)}_{n,\tilde{P}_n^{(j)}P_n^{(j+1)}\ldots P_n^{(k)}(P_n^0)}\}+R_n^{(j)}(\tilde{P}_n^{(j)}P_n^{(j+1)}\ldots P_n^{(k)}(P_n^0),P_0)\right\} \\
= \sum_{j=1}^{k}\left\{ (P_n-P_0)D^{(j)}_{n,\tilde{P}_n^{(j)}(P_0)}+R^{(j)}_n(\tilde{P}_n^{(j)}(P_0),P_0)\right\}\\
+ R^{(k)}_n(\tilde{P}_n^{(k)}(P_n^0),\tilde{P}_n)-R^{(k)}_n(\tilde{P}_n^{(k)}(P_0),\tilde{P}_n)\\
+\sum_{j=1}^k(\tilde{P}_n-P_n)\{D^{(j)}_{n,\tilde{P}_n^{(j)}(P_0)} -D^{(j)}_{n,\tilde{P}_n^{(j)}P_n^{(j+1)}\ldots P_n^{(k)}(P_n^0)}\}\\
+\sum_{j=1}^k(P_n-P_0)\{ D^{(j)}_{n,P_n^{(j)}\ldots P_n^{(k)}(P_n^0)}-D^{(j)}_{n,\tilde{P}_n^{(j)}P_n^{(j+1)}\ldots P_n^{(k)}(P_n^0)}\}\\
+\sum_{j=1}^k\left\{  R_n^{(j)}(P_n^{(j)}\ldots P_n^{(k)}(P_n^0),P_0)-R_n^{(j)}(\tilde{P}_n^{(j)}P_n^{(j+1)}\ldots P_n^{(k)}(P_n^0),P_0)\right\}.
\end{array}
\]
Let's discuss the last three terms that need to be compared with the HAL-regularization bias term $\sum_{j=1}^k(\tilde{P}_n-P_n)D^{(j)}_{n,\tilde{P}_n^{(j)}(P_0)}$ for the $k$-th order regularized TMLE.
The leading term  of the form $(\tilde{P}_n-P_n)\cdot$ is dominated by its first $j=1$-term and that one is already a 
third order difference, while for the regularized $k$-th order TMLE this first term is a second order difference. The second term in this sum  of three terms is of the form $(P_n-P_0)$ and is a very nice generally negligible empirical process term that has nothing to do with HAL-regularization bias. The
final term does represents a HAL-regularization bias term concerning a difference of $\tilde{P}_n$ and $P_n$. The worst case contribution from the latter sum  would come from the first term $R^{(1)}(P_n^{(1)}P,P_0)-R^{(1)}(\tilde{P}_n^{(1)}P,P_0)$, for $P=P_n^{(1)}\ldots P_n^{(k)}(P_n^0)$, which equals 
\[
\begin{array}{l}
R^{(1)}(\tilde{P}(P,\epsilon_n^{(1)}(P_n^0)),P_0)-R^{(1)}(\tilde{P}(P,\tilde{\epsilon}_n^{(1)}(P)),P_0)\\
=\frac{d}{d\xi}R^{(1)}(\tilde{P}(P,\xi),P_0)(\epsilon_n^{(1)}(P)-\tilde{\epsilon}_n^{(1)}(P)\end{array}
\]
at some intermediate point $\xi$ between $\epsilon_n^{(1)}(P)$ and $\tilde{\epsilon}_n^{(1)}(P)$.  This behaves as $d_0(\tilde{P}_n^1(P),P_0)^{1/2}(\tilde{P}_n-P_n)D^{(1)}_{\tilde{P}_n^{(1)}(P)}$, which shows  that also  this undersmoothing term represents  a third order difference instead of second order difference as with the regularized $k$-th order TMLE. This does not take into account that $R_n^{(1)}(\tilde{P}_n^{(1)}(P),P_0)$ at $P=P_n^{(2)}\ldots P_n^{(k)}(P_n^0)$ should be a better behaved remainder than at (e.g.)  $P=P_n^0$, due to the sequential targeting targeting this remainder. Therefore, in practice we expect an additional reduction due to using a $k$-th order TMLE with $k\geq 2$ with increasing benefit as $k$ increases. 

\subsection{Improved representation of exact expansion  for empirical $k$-th order TMLE}
In the above exact expansion for the empirical $k$-th order TMLE  the canonical gradients and remainders $R_n^{(j)}$, including the final $k+1$-th order remainders $R_n^{(k)}$, are  evaluated at the regularized $\tilde{P}_n^{(j)}$. In this subsection we succeed in rewriting the above expression so that it evaluates these leading terms  at the empirical TMLE update $P_n^{(j)}$ instead.

Recall that for $j=1,\ldots,k$, we have
 \[
\begin{array}{l}
R_n^{(j)}(\tilde{P}_n^{(j)}(P),P_0)=R_n^{(j)}(\tilde{P}_n^{(j)}(P),\tilde{P}_n)-R_n^{(j)}(\tilde{P}_n^{(j)}(P_0),\tilde{P}_n)\\
+R_n^{(j)}(\tilde{P}_n^{(j)}(P_0),P_0)
+(\tilde{P}_n-P_0) \{D^{(j)}_{n,\tilde{P}_n^{(j)}(P_0)}-D^{(j)}_{\tilde{P}_n^{(j)}(P)} \}.
\end{array}
\]
This expression applies to  $\tilde{P}_n^{(j)}(P)$ replaced by $P_n^{(j)}(P)$.
Thus,
\[
\begin{array}{l}
\sum_{j=1}^k R_n^{(j)}(P_n^{(j)}(P),P_0)-R_n^{(j)}(\tilde{P}_n^{(j)}(P),P_0)\\
=
\sum_{j=1}^k \{R_n^{(j)}(P_n^{(j)}(P),\tilde{P}_n)-R_n^{(j)}({P}_n^{(j)}(P_0),\tilde{P}_n)\}\\
-\sum_{j=1}^k \{R_n^{(j)}(\tilde{P}_n^{(j)}(P),\tilde{P}_n)-R_n^{(j)}(\tilde{P}_n^{(j)}(P_0),\tilde{P}_n)\} \\
+\sum_{j=1}^k \{R_n^{(j)}(P_n^{(j)}(P_0),P_0)-
-R_n^{(j)}(\tilde{P}_n^{(j)}(P_0),P_0)\}\\
+\sum_{j=1}^k (\tilde{P}_n-P_0) \{D^{(j)}_{n,{P}_n^{(j)}(P_0)}-D^{(j)}_{P_n^{(j)}(P)} \}- 
 (\tilde{P}_n-P_0) \{D^{(j)}_{n,\tilde{P}_n^{(j)}(P_0)}-D^{(j)}_{\tilde{P}_n^{(j)}(P)} \}.
 \end{array}
 \]
Substituting this in the above expression for the empirical $k$-th order TMLE  results in various replacements/swaps  of $\tilde{P}_n^{(j)}$ with $P_n^{(j)}$.  The following steps result in our desired final expression:
\[
\begin{array}{l}
\Psi(P_n^{(1)}\ldots P_n^{(k)}(P_n^0))-\Psi(P_0)= \sum_{j=1}^{k}\left\{ (P_n-P_0)D^{(j)}_{n,\tilde{P}_n^{(j)}(P_0)}+R^{(j)}_n({\bf P}_n^{(j)}(P_0),P_0)\right\}\\
+\sum_{j=1}^k(\tilde{P}_n-P_n)\{D^{(j)}_{n,\tilde{P}_n^{(j)}(P_0)} -D^{(j)}_{n,\tilde{P}_n^{(j)}P_n^{(j+1)}\ldots P_n^{(k)}(P_n^0)}\}\\
+\sum_{j=1}^k(P_n-P_0)\{ D^{(j)}_{n,P_n^{(j)}\ldots P_n^{(k)}(P_n^0)}-D^{(j)}_{n,\tilde{P}_n^{(j)}P_n^{(j+1)}\ldots P_n^{(k)}(P_n^0)}\}\\
+\sum_{j=1}^k \left\{R_n^{(j)}(P_n^{(j)}P_n^{(j+1)}\ldots P_n^{(k)}(P_n^0),\tilde{P}_n)-R_n^{(j)}({P}_n^{(j)}(P_0),\tilde{P}_n)\right\}\\
-\sum_{j=1}^k\left\{R_n^{(j)}(\tilde{P}_n^{(j)}P_n^{(j+1)}\ldots P_n^{(k)}P_n^0,\tilde{P}_n)-R_n^{(j)}(\tilde{P}_n^{(j)}(P_0),\tilde{P}_n)\right\}\\
+\sum_{j=1}^k (\tilde{P}_n-P_n+(P_n-P_0) )\{D^{(j)}_{n,{P}_n^{(j)}(P_0)}-D^{(j)}_{P_n^{(j)}P_n^{(j+1)}\ldots P_n^{(k)}(P_n^0)} \}\\
- 
\sum_{j=1}^k  (\tilde{P}_n-P_n+(P_n-P_0)) \{D^{(j)}_{n,\tilde{P}_n^{(j)}(P_0)}-D^{(j)}_{\tilde{P}_n^{(j)}P_n^{(j+1)}\ldots P_n^{(k)}(P_n^0)} \}\\
+ R^{(k)}_n(\tilde{P}_n^{(k)}(P_n^0),\tilde{P}_n)-R^{(k)}_n(\tilde{P}_n^{(k)}(P_0),\tilde{P}_n)\\
= \sum_{j=1}^{k}\left\{ (P_n-P_0)D^{(j)}_{n,\tilde{P}_n^{(j)}(P_0)}+R^{(j)}_n({P}_n^{(j)}(P_0),P_0)\right\}\\
+\sum_{j=1}^k(\tilde{P}_n-P_n)\{D^{(j)}_{n,{\bf P}_n^{(j)}(P_0)} -D^{(j)}_{n,{\bf P}_n^{(j)}P_n^{(j+1)}\ldots P_n^{(k)}(P_n^0)}\}\\
+\sum_{j=1}^k(P_n-P_0)\{ D^{(j)}_{n,P_n^{(j)}\ldots P_n^{(k)}(P_n^0)}-D^{(j)}_{n,\tilde{P}_n^{(j)}P_n^{(j+1)}\ldots P_n^{(k)}(P_n^0)}\}\\
+\sum_{j=1}^k \{R_n^{(j)}(P_n^{(j)}\ldots P_n^{(k)}P_n^0,\tilde{P}_n)-R_n^{(j)}({P}_n^{(j)}(P_0),\tilde{P}_n)\}\\
-\sum_{j=1}^k \{R_n^{(j)}(\tilde{P}_n^{(j)}P_n^{(j+1)}\ldots P_n^{(k)}P_n^0,\tilde{P}_n)-R_n^{(j)}(\tilde{P}_n^{(j)}(P_0),\tilde{P}_n)\}\\
+ \sum_{j=1}^k (P_n-P_0) \{D^{(j)}_{n,{P}_n^{(j)}(P_0)}-D^{(j)}_{P_n^{(j)}\ldots P_n^{(k)}P_n^0} \}\\
- 
\sum_{j=1}^k (P_n-P_0) \{D^{(j)}_{n,\tilde{P}_n^{(j)}(P_0)}-D^{(j)}_{\tilde{P}_n^{(j)}P_n^{(j+1)}\ldots P_n^{(k)}P_n^0} \}\\
+ R^{(k)}_n(\tilde{P}_n^{(k)}(P_n^0),\tilde{P}_n)-R^{(k)}_n(\tilde{P}_n^{(k)}(P_0),\tilde{P}_n)\\
=
\sum_{j=1}^{k}\left\{ (P_n-P_0)D^{(j)}_{n,{\bf P}_n^{(j)}(P_0)}+R^{(j)}_n({P}_n^{(j)}(P_0),P_0)\right\}\\
+\sum_{j=1}^k(\tilde{P}_n-P_n)\{D^{(j)}_{n,{P}_n^{(j)}(P_0)} -D^{(j)}_{n,{P}_n^{(j)}P_n^{(j+1)}\ldots P_n^{(k)}(P_n^0)}\}\\
+\sum_{j=1}^k(P_n-P_0)\{ D^{(j)}_{n,P_n^{(j)}\ldots P_n^{(k)}(P_n^0)}-D^{(j)}_{n,\tilde{P}_n^{(j)}P_n^{(j+1)}\ldots P_n^{(k)}(P_n^0)}\}\\
+\sum_{j=1}^k \{R_n^{(j)}(P_n^{(j)}\ldots P_n^{(k)}P_n^0,\tilde{P}_n)-R_n^{(j)}({P}_n^{(j)}(P_0),\tilde{P}_n)\}\\
-\sum_{j=1}^k\{R_n^{(j)}(\tilde{P}_n^{(j)}P_n^{(j+1)}\ldots P_n^{(k)}P_n^0,\tilde{P}_n)-R_n^{(j)}(\tilde{P}_n^{(j)}(P_0),\tilde{P}_n)\}\\
+ {\bf \sum_{j=1}^k (P_n-P_0)\{D^{(j)}_{\tilde{P}_n^{(j)}P_n^{(j+1)}\ldots P_n^{(k)}(P_n^0)}-D^{(j)}_{P_n^{(j)}P_n^{(j+1}\ldots P_n^{(k)}P_n^0 } \}} \\
+ R^{(k)}_n(\tilde{P}_n^{(k)}(P_n^0),\tilde{P}_n)-R^{(k)}_n(\tilde{P}_n^{(k)}(P_0),\tilde{P}_n)\\
\end{array}
\]
\newpage
\[
\begin{array}{l}
=\sum_{j=1}^{k}\left\{ (P_n-P_0)D^{(j)}_{n,{P}_n^{(j)}(P_0)}+R^{(j)}_n({P}_n^{(j)}(P_0),P_0)\right\}
\\
+ R^{(k)}_n({\bf P}_n^{(k)}(P_n^0),\tilde{P}_n)-R^{(k)}_n({\bf P}_n^{(k)}(P_0),\tilde{P}_n)\\
+\sum_{j=1}^k(\tilde{P}_n-P_n)\{D^{(j)}_{n,{P}_n^{(j)}(P_0)} -D^{(j)}_{n,{P}_n^{(j)}P_n^{(j+1)}\ldots P_n^{(k)}(P_n^0)}\}\\
+\sum_{j=1}^{{\bf k-1}} \left\{R_n^{(j)}(P_n^{(j)}\ldots P_n^{(k)}P_n^0,\tilde{P}_n)-R_n^{(j)}({P}_n^{(j)}(P_0),\tilde{P}_n)\right\}\\
-\sum_{j=1}^{{\bf k-1}}\left\{R_n^{(j)}(\tilde{P}_n^{(j)}P_n^{(j+1)}\ldots P_n^{(k)}P_n^0,\tilde{P}_n)-R_n^{(j)}(\tilde{P}_n^{(j)}(P_0),\tilde{P}_n)\right\}.
\end{array}
\]
Comparison with the previous expansion above for this empirical $k$-th order TMLE shows that is has been improved by replacing in first, second and last formulas the regularized TMLEs by the empirical TMLEs. In addition, it replaced the $R_n^{(j)}$-single difference by a second order difference which should therefore result in additional  cancelation. The final three terms represent the HAL-regularization bias term concerning a difference of $\tilde{P}_n$ and $P_n$, and as discussed above it represents a higher order difference relative to the HAL-regularization bias term of the regularized $k$-th order TMLE.

This exact expansion for the empirical $k$-th order TMLE might be  viewed as the most natural expansion and is the expression presented in Theorem \ref{theoremempiricalhotmle}.


{\bf Comparison with exact expansion of HAL-regularized $k$-th order TMLE:}
Comparing this exact expansion with the exact expansion for the regularized $k$-th order TMLE, it follows that the first two lines represent the same leading expansion but with $\tilde{P}_n^{(j)}$ replaced by the empirical TMLE updates $P_n^{(j)}$. This represents an improvement since the exact remainders $R^{(j)}_n({P}_n^{(j)}(P_0),P_0)$ are now $O_P(n^{-1})$ without need for undersmoothing: i.e., $P_n^{(j)}(P_0)$ is a simple parametric standard/empirical MLE according to a correctly specified one-dimensional model.

The last three rows in the exact expansion for the empirical $k$-th order TMLE need to  be compared with the undersmoothing term of the regularized $k$-th order TMLE given by $\sum_{j=1}^k (\tilde{P}_n-P_n)D^{(j)}_{n,\tilde{P}_n^{(j)}(P_0)}$. The first row of the last three rows is directly an improvement of the latter term, by subtracting from $D^{(j)}_{n,\tilde{P}_n^{(j)}(P_0)}$ its estimate. 
We already argued in previous subsection that 
\[
\sum_{j=1}^{k-1} \left\{R_n^{(j)}(P_n^{(j)}\ldots P_n^{(k)}P_n^0,\tilde{P}_n)-R_n^{(j)}(\tilde{P}_n^{(j)}P_n^{(j+1)}\ldots P_n^{(k)}P_n^0,\tilde{P}_n)-\right\}\]
 is at minimal a third order difference. The last two rows represent this term minus 
\[
\sum_j\{
R_n^{(j)}(\tilde{P}_n^{(j)}(P_0),\tilde{P}_n)-R_n^{(j)}(P_n^{(j)}(P_0),\tilde{P}_n)\},\] suggesting that it represents an even smaller term. 
Therefore, we conclude that the exact expansion for the empirical $k$-th order TMLE is superior to the exact expansion for the regularized $k$-th order TMLE. Based on these theoretical considerations, and our practical experience in simulations, contrary to the HAL-regularized $k$-th order TLME, we conclude that  the empirical $k$-th order TMLE does not require any undersmoothing for $\tilde{P}_n$ anymore.

\subsection{Understanding that HAL-regularization bias term for is third order difference.}
The leading term in the undersmoothing term of the exact expansion for the $k$-th order TMLE is given by 
$(\tilde{P}_n-P_n)\{D^{(1)}_{n,{P}_n^{(1)}(P_0)} -D^{(1)}_{n,{P}_n^{(1)}P_n^{(2)}\ldots P_n^{(k)}(P_n^0)}\}$, which is of the  form $(\tilde{P}_n-P_n)f_n$, where $\pl f_n\pl_{P_0}=O_P(d_0^{1/2}(P_n^0,P_0))$. In this subsection we demonstrate that this is a third order difference. 
Let $\tilde{f}_n=f_n/\pl f_n\pl_{P_0}$ which is a function with a bounded $\pl \cdot\pl_{P_0}$-norm. Then,  
\[
(\tilde{P}_n-P_n)f_n=\pl f_n\pl_{P_0} (\tilde{P}_n-P_n)S_{\tilde{P}_n,\tilde{f}_n},\]
 where $S_{\tilde{P}_n,\tilde{f}_n}=\tilde{f}_n-\tilde{P}_n \tilde{f}_n$ is a score at $\tilde{P}_n$. Since $\tilde{P}_n$ is an MLE it solves a set of score equations $P_n S_{\tilde{P}_n,j}=0$ for $j=1,\ldots,J$, where this set increases as the $L_1$-norm of the HAL-MLE increases. Let $\sum_j \beta_n(j)S_{\tilde{P}_n,j}$ be the projection of $S_{\tilde{P}_n,\tilde{f}_n}$ onto the linear span of these scores in $L^2(P_0)$. 
 We have
 \[ (\tilde{P}_n-P_n)f_n=\pl f_n\pl_{P_0} (\tilde{P}_n-P_n)\left\{ S_{\tilde{P}_n,\tilde{f}_n} 
 -\sum_j\beta_n(j)S_{\tilde{P}_n,j}\right\}.\]
 One can further write $(\tilde{P}_n-P_n)=(\tilde{P}_n-P_0)-(P_n-P_0)$.
 The contribution from $(\tilde{P}_n-P_0)$ is a third order difference involving the $L^2(P_0)$ rate of $f_n$; the rate $d_0(\tilde{P}_n,P_0)^{1/2}$ and the rate of the oracle approximation of $S_{\tilde{P}_n,\tilde{f}_n}$ by the linear combination of scores solved by $\tilde{P}_n$. 
 If we set $P_n^0$ equal to an HAL-MLE, then this will be $O_P(n^{-2/3}(\log n)^{k_1})$ times the $L^2(P_0)$-norm of the oracle approximation of $S_{\tilde{P}_n,\tilde{f}_n}$. 
Thus this term can be expected to converge to zero almost as fast as $n^{-1}$. 
In fact, it is better than this third order difference since $(\tilde{P}_n-P_n)$ will behave as a nicer less biased term than $d_0^{1/2}(\tilde{P}_n,P_0)$.

\section{Proof of Lemmas \ref{lemma5b} and  \ref{lemma6}}\label{AppendixC}

\subsection{Proof of Lemma \ref{lemma5b}}
See the definition of $R_n^+(P,P_0)$ in Lemma \ref{lemma5} and consider the assumptions of Lemma \ref{lemma5b}.
Let  $f_{j,k+1}(x,y)\equiv \prod_{i\not =j}(H_{1,i}(P)-H_{1,i}(P_0))(x)\prod_{l=1}^2 (H_{1,l,j}(P)-H_{1,l,j}(P_0))(x,y)$.
We   have
\[
\begin{array}{l}
R_n^+(P,P_0)=\frac{1}{k}\sum_{j=1}^k \int \prod_{i\not =j}(H_{1,i}(P)-H_{1,i}(P_0))R_{2,P,H_{1,j}()}(P,P_0)(x) d\mu_{1,n,P_0}(x)\\
=\frac{1}{k}\sum_{j=1}^k \int f_{j,k+1}(x,x) C_{1,j,P_0}(x)\bar{H}_{1,j,P,P_0}(x,x) d\mu_{1,n,P_0}(x)\\
+
\frac{1}{k}\sum_{j=1}^k \int \int f_{j,k+1}(x,y)  \bar{H}_{1,j,P,P_0}(x,y) d\mu_{j,P_0}^c(y)  d\mu_{1,n,P_0}(x)\\
=\frac{1}{k}\sum_{j=1}^k \int f_{j,k+1}(x,x) \bar{H}_{1,j,P_0,P_0}(x,x)C_{1,j,P_0}(x) d\mu_{1,n,P_0}(x)\\
+\frac{1}{k}\sum_{j=1}^k \int f_{j,k+1}(x,x) (\bar{H}_{1,j,P,P_0}-\bar{H}_{1,j,P_0,P_0})(x,x)C_{1,j,P_0}(x) d\mu_{1,n,P_0}(x)\\
+\frac{1}{k}\sum_{j=1}^k \int \int  f_{j,k+1}(x,y)  \bar{H}_{1,j,P_0,P_0}(x,y) d\mu_{j,P_0}^c(y)  d\mu_{1,n,P_0}(x)\\
+\frac{1}{k}\sum_{j=1}^k \int \int f_{j,k+1}(x,y) (\bar{H}_{1j,P,P_0}- \bar{H}_{1j,P_0,P_0})(x,y) d\mu_{j,P_0}^c(y)  d\mu_{1,n,P_0}(x).
\end{array}
\]

Define $H_{1,j,i}^a(P)=H_{1,i}(P)$ if $i\not =j$; $H_{1,j,j}^a(P)=H_{1,1,j}(P)$; $H_{1,j,k+1}(P)=H_{1,2,j}(P)$ and $d\mu^a_{1,j,n,P_0}=C_{1,j,P_0}(x)H_{1,j,P_0,P_0}(x,x)d\mu_{1,n,P_0}(x)$.
Then, the first term equals
\[
R_n^{+,a}(P,P_0)=\frac{1}{k}\sum_{j=1}^k \int \prod_{i=1}^{k+1}(H_{1,j,i}^a(P)-H_{1,j,i}^a(P_0))(x)
d\mu^a_{1,j,n,P_0}(x).
\]
This is a $k+1$-th order polynomial difference.
Define $H_{1,j,i}^b(P)(x,y)=H_{1,i}(P)(x)$ if $i\not =j$; $H_{1,j,j}^b(x,y)=H_{1,1,j}(P)(x,y)$;
$H_{1,j,k+1}^b(P)(x,y)=H_{1,2,j}(P)(x,y)$ and $d\mu^b_{1,j,n,P_0}=H_{1,j,P_0,P_0}(x,y) d\mu_{j,P_0}^c(y)  d\mu_{1,n,P_0}(x)$. Then, the 
third term equals
\[
R_n^{+,b}(P,P_0)=\frac{1}{k}\sum_{j=1}^k \int_x \int_y \prod_{i=1}^{k+1}(H_{1,j,i}^b(P)-H_{1,j,i}^b(P_0))(x,y)
d\mu^n_{1,j,n,P_0}(x,y) .\]
This is a $k+1$-th order polynomial difference.
Define $H_{1,j,i}^c(P)=H_{1,i}(P)(x)$ if $i\not =j$; $H_{1,j,j}^c(P)=H_{1,1,j}(P)(x,x)$; $H_{1,j,k+1}(P)=H_{1,2,j}(P)(x,x)$; $H_{1,j,k+2}^c(P)=H_{1,j,P,P_0}(x,x)$ and $d\mu^c_{1,j,n,P_0}=C_{1,j,P_0}(x)d\mu_{1,n,P_0}(x)$.
Then, the second term equals
\[
R_n^{+,c}(P,P_0)=\frac{1}{k}\sum_{j=1}^k \int_x\prod_{i=1}^{k+1}(H^c_{1,j,i}(P)-H^c_{1,j,i}(P_0)(x)
d\mu^c_{1,j,n,P_0}(x).\]
This is a $k+2$-th order polynomial difference.
Define $H_{1,j,i}^d(P)(x,y)=H_{1,i}(P)(x)$ if $i\not=j$; $H_{1,j,j}^d(P)(x,y)=H_{1,1,j}(P)(x,y)$; 
$H_{1,j,k+1}^d(P)(x,y)=H_{1,2,j}(P)(x,y)$; $H_{1,j,k+2}^d(P)=H_{1,j,P,P_0}(x,y)$, and $d\mu^d_{1,j,n,P_0}=d\mu_{j,P_0}^c(y)d\mu_{1,n,P_0}(x)$. Then, the fourth term equals
\[
R_n^{+,d}(P,P_0)=
\frac{1}{k}\sum_{j=1}^k \int_x \int_y \prod_{i=1}^{k+2}(H_{1,j,i}^d(P)-H_{1,j,i}^d(P_0))(x,y)
d\mu^d_{1,j,n,P_0}(x,y)  .\]
This is a $k+2$-th order polynomial difference.
So this shows that
\[
\tilde{R}_n^+(P,P_0)=R_n^{+,a}(P,P_0)+R_n^{+,b}(P,P_0)+R_n^{+,c}(P,P_0)+R_n^{+,d}(P,P_0)\]
is a sum of a $k+1$-th order and $k+2$-th order polynomial difference.
This completes the proof. $\Box$

\subsection{The directional derivative of a $k$-th order polynomial difference in direction $P-P_0$ is a sum of $k$-th and $k+1$-th order polynomial difference.}
One statement in Lemma \ref{lemma6} is given by the following lemma proven here.
\begin{lemma}\label{lemma6b}
Suppose that $R_n:{\cal M}^2\rightarrow \openr$ is a $k$-th order polynomial difference.
Then, under weak regularity conditions,  the directional derivative $\frac{d}{dP}R_n(P,P_0)(P-P_0)$ of $P\rightarrow R_n(P,P_0)$ at $P$ in the direction $P-P_0$ is a sum of a $k$-th order and $k+1$-th order polynomial difference in $P$ and $P_0$.

In particular,  given $R_n(P,P_0)$ is defined in terms of $H_{1,j}$ and $H_{2,j}$ as in our definition of a $k$-th order polynomial difference, it suffices to assume that $\frac{d}{dP}H_{1,j}(P)(P-P_0)(x)=K_{1,j,P,P_0}(x)(p-p_0)(x)+\int K_{1,j,P,P_0,x}(o)(p-p_0)(o)d\mu(o)$ for certain functions $K_{1,j,P,P_0}(x)$ and kernel $(x,o)\rightarrow K_{1,j,P,P_0}(x,o)$, and similarly for $H_{2,j}$.
\end{lemma}
{\bf Proof:}
To demonstrate this lemma, suppose that 
\[
R_n(P,P_0)=\int \prod_{j=1}^{k}(H_{1,j}(P)-H_{1,j}(P_0)) d\mu_{1,n,P_0}.\]
Then, by the product rule of differentiation
\[
\begin{array}{l}
\frac{d}{dP}R_n(P,P_0)(P-P_0)=
\sum_{j=1}^k \int \prod_{l\not =j}^k(H_{1,l}(P)-H_{1,l}(P_0))\frac{d}{dP}H_{1,j}(P)(P-P_0) d\mu_{1,n,P_0}.
\end{array}
\]
If we replace $d/dPH_{1,j}(P)$ by $d/dP_0H_{1,j}(P_0)$, then it is a $k$-th order polynomial difference. 
Suppose $\frac{d}{dP}H_{1,j}(P)(P-P_0)(x)=K_{1,j,P,P_0}(x)(p-p_0)(x)+\int K_{1,j, P,P_0,x}(o)(p-p_0)(o)d\mu(o)$ for certain functions $K_{1,j,P,P_0}(x)$ and kernel $(x,o)\rightarrow K_{1,j,P,P_0}(x,o)$.
Then,
\[ 
\begin{array}{l}
\frac{d}{dP}H_{1,j}(P)(P-P_0) = K_{1,j,P,P_0}(x)(p-p_0)(x)+\int K_{1, j,P,P_0}(x,o)(p-p_0)(o)d\mu(o)\\
=K_{1,j,P_0,P_0}(x)(p-p_0)(x)+\int K_{1, j,P_0,P_0}(x,o)(p-p_0)(o)d\mu(o)\\
+(K_{1,j,P,P_0}-K_{1,j,P_0,P_0})(x)(p-p_0)(x)+\int (K_{1,j,P,P_0}-K_{1,j,P_0,P_0})(x,o)(p-p_0)(o)d\mu(o).
\end{array}
\]
Let $d\bar{\mu}(o,x)\equiv d\mu(o)d\mu_{1,n,P_0}(x)$.
Thus,
\[
\begin{array}{l}
\frac{d}{dP}R_n(P,P_0)(P-P_0)=
\sum_{j=1}^k \int \prod_{l\not =j}^k(H_{1,l}(P)-H_{1,l}(P_0)) (p-p_0)  K_{1,j,P_0,P_0} d\mu_{1,n,P_0}\\
+\sum_{j=1}^k \int \prod_{l\not =j}^k(H_{1,l}(P)-H_{1,l}(P_0))  \int K_{1,j,P_0,P_0})(x,o)(p-p_0)(o)d\mu(o) d\mu_{1,n,P_0}(x)\\
+\sum_{j=1}^k \int \prod_{l\not =j}^k(H_{1,l}(P)-H_{1,l}(P_0)) (K_{1,j,P,P_0}- K_{1,j,P_0,P_0})(p-p_0)d\mu_{1,n,P_0}\\
+\sum_{j=1}^k \int \prod_{l\not =j}^k(H_{1,l}(P)-H_{1,l}(P_0))  \int (K_{1,j,P,P_0}-K_{1,j,P_0,P_0})(x,o)(p-p_0)(o)d\bar{\mu}(o,x).
\end{array}
\]
Define $H_{1,j,l}^a(P)=H_{1,l}(P)$ if $j\not =l$; $H_{1,j,j}^a(P)=p$, and let $d\mu_{1,j,n,P_0}=K_{1,j,P_0,P_0}d\mu_{1,n,P_0}$.
Then, the first term equals
\[ R_n^a(P,P_0)=\sum_{j=1}^k \int \prod_{l=1}^k(H_{1,j,l}^a(P)-H_{1,j,l}^a(P_0))   d\mu_{1,j,n,P_0}.
\]
This is a $k$-th order polynomial difference.
The second term can be written as:\[
\sum_{j=1}^k \int_x \int_o \prod_{l\not =j}^k(H_{1,l}(P)-H_{1,l}(P_0))(p-p_0)   
K_{1,j,P_0,P_0})(x,o)(d\mu(o) d\mu_{1,n,P_0}(x).\]
Define $H_{1,j,l}^b(P)=H_{1,l}(P)$ for $l\not =j$; $H_{1,j,j}^b(P)=p$, and  \newline $d\mu^b_{1,j,n,P_0}=
K_{1,j,P_0,P_0}(x,o)(d\mu(o) d\mu_{1,n,P_0}(x)$.
Then, the second term equals
\[ R_n^b(P,P_0)=\sum_{j=1}^k \int\int \prod_{l=1}^k(H_{1,j,l}^b(P)-H_{1,j,l}^b(P_0)) d\mu^b_{1,j,n,P_0}.\]
This is a $k$-th order polynomial difference.
Define $H_{1,j,l}^c(P)=H_{1,l}(P)$ for $l\not =j$; $H_{1,j,j}^c(P)=K_{1,j,P,P_0}$, and $H_{1,j,k+1}^c(P)=p$. Then, the third term equals
\[
R_n^c(P,P_0)=\sum_{j=1}^k \int \prod_{l=1}^{k+1}(H_{1,j,l}^c(P)-H_{1,j,l}^c(P_0))d\mu_{1,n,P_0}.\]
This is a $k+1$-th order polynomial difference.
Define $H_{1,j,l}^d(P)=H_{1,l}(P)$ if $l\not =j$; $H_{1,j,j}^d(P)=K_{1,j,P,P_0}$; $H_{1,j,k+1}^d(P)=p$, and $d\mu^d_{1,n,P_0}=d\mu(o)d\mu_{1,n,P_0}$.
Then, the fourth term equals
\[
R_n^d(P,P_0)=\sum_{j=1}^k \int\int \prod_{l=1}^{k+1}(H_{1,j,l}^d(P)-H_{1,j,l}^d(P_0))(x)  d\mu^c_{1,n,P_0}(x,o).
\]
This is a $k+1$-th order polynomial difference.
So we have shown that 
\[
\frac{d}{dP}R_n(P,P_0)(P-P_0)=R_n^a(P,P_0)+R_n^b(P,P_0)+R_n^c(P,P_0)+R_n^d(P,P_0),\]
where $R_n^a,R_n^b$ are $k$-th order polynomial differences, and $R_n^c,R_n^c$ are $k+1$-th order polynomial differences. 
This completes the proof. $\Box$

\subsection{A generalized $k$-th order difference $R_n$ equals a generalized $k+1$-th order difference $R_n^+$ plus its scaled directional derivative}
The remaining statements in Lemma \ref{lemma6} follow from the next lemma proven here.

\begin{lemma}\label{lemma6c}
Let $R_n(P,P_0)=R_{n,k}(P,P_0)+R_{n,k+1}(P,P_0)$ be a sum of a $k$ and $k+1$-th order polynomial difference.
Then,
\[
\begin{array}{l}
R_n(P,P_0)=R_{n,k}^{+}(P,P_0)+R_{n,k+1}(P,P_0)\\
+
k^{-1}\frac{d}{dP}R_n(P,P_0)(P-P_0)\\
\hfill -k^{-1}\frac{d}{d P}R_{n,k+1}(P,P_0)(P-P_0).\end{array}
\]

By Lemma \ref{lemma5b}, under specified regularity conditions, $R_{n,k}^+(P,P_0)$ is (globally)  a sum of a $k+1$ and $k+2$-th order polynomial difference. Also, by Lemma \ref{lemma6b}, under specified regularity conditions, $k^{-1}d/dP R_{n,k+1}(P,P_0)(P-P_0)$ is (globally) a sum of a $k+1$ and $k+2$-th order polynomial difference. Thus, $R_n(P,P_0)$ is (globally) a generalized $k+1$-th order difference plus $k^{-1}\frac{d}{dP}R_{n}(P,P_0)(P-P_0)$.

If  $\frac{d}{d\tilde{P}_n(P)}R_n(\tilde{P}_n(P),P_0)(\tilde{P}_n(P)-P_0)=0$, this globally $k$-th order difference $R_n:{\cal M}^2\rightarrow \openr$, satisfies that $R_n(\tilde{P}_n(P),P_0)$ is a generalized $k+1$-th order difference in $P$ and $P_0$.

In general, if $R_n(P,P_0)$ is a generalized $k$-th order difference, then 
$R_n(P,P_0)$ is a sum of a generalized $k+1$-th order difference $R_n^+(P,P_0)$  and $k^{-1}\frac{d}{dP}R_n(P,P_0)(P-P_0)$. 

That is, given a $k$-th order generalized difference $R_n:{\cal M}^2\rightarrow\openr$, we can define a $k+1$-th order generalized difference  $R_n^+:{\cal M}^2\rightarrow\openr$, and derivative operator $\dot{R}_n:{\cal M}\rightarrow\openr$ defined by 
$\dot{R}_n(P,P_0)=\frac{d}{dP}R_n(P,P_0)(P-P_0)$,  so that $R_n=R_n^++\dot{R}_n$.
By applying this global representation of $R_n$ to a particular $(\tilde{P}_n(P),P_0)$ at which $\dot{R}_n(\tilde{P}_n(P),P_0)= 0$, we obtain that $R_n(\tilde{P}_n(P),P_0)$ equals $R_n^+(\tilde{P}_n(P),P_0)$.

\end{lemma}
{\bf Proof:} 
We have
\[
\frac{d}{dP}R_{n,k}(P,P_0)(P-P_0)=\frac{d}{dP}R_n(P,P_0)(P-P_0)-\frac{d}{dP}R_{n,k+1}(P,P_0)(P-P_0).\]
We showed that $\frac{d}{dP}R_{n,k+1}(P,P_0)(P-P_0)$ can be decomposed into a sum
$dR_{n,k+1}^{1}(P,P_0)+dR_{n,k+2}^{2}(P,P_0)$ of a $k+1$-th and $k+2$-th order polynomial difference:
\[\frac{d}{dP}R_{n,k+1}(P,P_0)(P-P_0)=dR_{n,k+1}^{1}(P,P_0)+dR_{n,k+2}^{2}(P,P_0).\]
Thus,
\[
\begin{array}{l}
\frac{d}{dP}R_{n,k}(P,P_0)(P-P_0)=\frac{d}{dP}R_n(P,P_0)(P-P_0)-dR_{n,k+1}^{1}(P,P_0)\\
\hfill -dR_{n,k+2}^{2}(P,P_0).\end{array}
\]
By Lemma \ref{lemma5}, we have
\[
R_{n,k}(P,P_0)=R_{n,k}^+(P,P_0)+k^{-1}\frac{d}{dP}R_{n,k}(P,P_0)(P-P_0).\]
Thus,
\[
\begin{array}{l}
R_{n,k}(P,P_0) 
=R_{n,k}^{+}(P,P_0) +k^{-1}
\frac{d}{dP}R_n(P,P_0)(P-P_0)-k^{-1}dR_{n,k+1}^{1}(P,P_0)\\
\hfill -k^{-1}dR_{n,k+2}^{2}(P,P_0).
\end{array}
\]
By Lemma \ref{lemma5b}, $R_{n,k}^+(P,P_0)$ can be decomposed into a $k+1$ and $k+2$-th order polynomial difference:
\[
{R}_{n,k}^{+}(P,P_0)=R_{n,k}^{+,a}(P,P_0)+R_{n,k}^{+,b}(P,P_0).\]
So we obtained:
\[
\begin{array}{l}
R_{n,k}(P,P_0)=R_{n,k}^{+,a}(P,P_0)+R_{n,k}^{+,b}(P,P_0)\\
+k^{-1}
\frac{d}{dP}R_n(P,P_0)(P-P_0)-k^{-1}dR_{n,k+1}^{1}(P,P_0)-k^{-1}dR_{n,k+2}^{2}(P,P_0).
\end{array}
\]
We have
$R_n(P,P_0)=R_{n,k}(P,P_0)+R_{n,k+1}(P,P_0)$ so that we obtained
\[
\begin{array}{l}
R_n(P,P_0)=
R_{n,k}^{+,a}(P,P_0)+R_{n,k}^{+,b}(P,P_0)+R_{n,k+1}(P,P_0) \\
+k^{-1}
\frac{d}{dP}R_n(P,P_0)(P-P_0)-k^{-1}dR_{n,k+1}^{1}(P,P_0)-k^{-1}dR_{n,k+2}^{2}(P,P_0)\\
\end{array}
\]
Combine all the terms beyond $k^{-1}d/dPR_n(P,P_0)(P-P_0)$ and note that they represent a generalized $k+1$-th order difference $R_n^{+}(P,P_0)$. At $P=\tilde{P}_n(P)$ for which $d/dP R_n(P,P_0)(P-P_0)=0$, this yields $R_n^+(\tilde{P}_n(P),P_0)$.
This proves the lemma. $\Box$

\section{
Understanding the order of difference of higher order canonical gradients }\label{AppendixD}

Recall
 \[
\Psi_n^{(j-1)}(\tilde{P}_n^{(j)}(P))-\Psi_n^{(j-1)}(\tilde{P}_n)=-\tilde{P}_n D^{(j)}_{n,\tilde{P}_n^{(j)}(P)}+R_n^{(j)}(\tilde{P}_n^{(j)}(P),\tilde{P}_n).\]
The left-hand side also equals $\Psi_n^{(j)}(P)-\Psi_n^{(j)}(\tilde{P}_n)$. 
So
\[
\Psi_n^{(j)}(P)-\Psi_n^{(j)}(\tilde{P}_n)= -\tilde{P}_n D^{(j)}_{n,\tilde{P}_n^{(j)}(P)}+R_n^{(j)}(\tilde{P}_n^{(j)}(P),\tilde{P}_n).\]
Let $\{P_{\delta,h}:\delta\}\subset{\cal M}$ be a path through $P$ with score at $\delta=0$ equal to a $h\in T_P({\cal M})$.
We consider paths $\tilde{P}_{n,\delta,h}$ through $\tilde{P}_n$. We have
\[
\Psi_n^{(j)}(\tilde{P}_{n,\delta,h})-\Psi_n^{(j)}(\tilde{P}_n)=
-\tilde{P}_n D^{(j)}_{n,\tilde{P}_n^{(j)}(\tilde{P}_{n,\delta,h})}+R_n^{(j)}(\tilde{P}_n^{(j)}(\tilde{P}_{n,\delta,h}),\tilde{P}_n).\]
As a consequence, we have
\[
\frac{d}{d\delta_0}\Psi_n^{(j)}(\tilde{P}_{n,\delta,h})=\frac{d}{d\delta_0} R_n^{(j)}(\tilde{P}_n^{(j)}(\tilde{P}_{n,\delta,h}),\tilde{P}_n)-\frac{d}{d\delta_0}\tilde{P}_n D^{(j)}_{n,\tilde{P}_n^{(j)}(\tilde{P}_{n,\delta,h})}.
\]
For example, if we use a local least favorable path and $\epsilon_n^{(j)}(P)$ is defined as solution of $\tilde{P}_n D^{(j)}_{n,\tilde{P}_n^{(j)}(P,\epsilon)}=0$, then we have
\[
\frac{d}{d\delta_0}\Psi_n^{(j)}(\tilde{P}_{n,\delta,h})=\frac{d}{d\delta_0} R_n^{(j)}(\tilde{P}_n^{(j)}(\tilde{P}_{n,\delta,h}),\tilde{P}_n),\]
and $R_n^{(j)}(\tilde{P}_n^{(j)}(P),\tilde{P}_n)=R_n^{(1)}(\tilde{P}_n^{(1)}\ldots\tilde{P}_n^{(j)}(P),\tilde{P}_n)$ would be exactly a $j+1$-th order difference of $P$ and $\tilde{P}_n$. The key for this result is that $\tilde{P}_n^{(j)}(P)$ is defined to solve $\tilde{P}_n D^{(j)}_{n,\tilde{P}_n^{(j)}(P)}=0$ in all $P\in {\cal M}$. 
So in that case, the pathwise derivative would be exactly equal to zero for all paths, and thereby $D^{(j)}_{n,\tilde{P}_n}=0$ for all $j=2,\ldots,k$. In our examples, we observe that even when using a local least favorable path with MLE we still have that $D^{(2)}_{n,\tilde{P}_n}=0$ exactly. 
Moreover, higher order pathwise derivatives up till order $j$ would all be equal to zero as well. 
This demonstrates that in that case $D^{(j)}_{n,P}$ involves a difference itself between $P$ and $\tilde{P}_n$ of order $j-1$. 

By solving the TMLE score equations up till  finite sample negligible precision these statements will hold for finite sample purposes. Either way, our results and exact expansions are all presented in terms of the exact TMLE score values so that the impact of not solving the scores is directly expressed.

\section{Proof of Lemma \ref{lemmaformulacangradient}.}\label{AppendixF}
In this proof we suppress the dependence on $(j)$.

{\bf Pathwise derivative of $\epsilon_n$:}
We have $\epsilon_n(P)=\arg\min_{\epsilon}\tilde{P}_n\log (1+\epsilon D_{n,P})p$.
So $\epsilon_n(P)$ solves $\tilde{P}_n\frac{D_{n,P}}{1+\epsilon_n(P)D_{n,P} }=0$.
Let $U(\epsilon,p)=\tilde{P}_n \frac{D_{n,p}}{(1+\epsilon D_{n,p})}$, and note that
$U(\epsilon_n(p),p)=0$ and $U(\epsilon_n(p_{\delta,h}),p_{\delta,h})=0$.
By the implicit function theorem we have
\[
\begin{array}{l}
A_{\epsilon_n,P}(h)=\frac{d}{d\delta_0}\epsilon_n(P_{\delta_0,h})\\
=\left\{\tilde{P}_n \frac{D_{n,P}^2}{(1+\epsilon_n(P)D_{n,P})^2}\right\}^{-1}\frac{d}{d\delta_0}U(\epsilon_n(P),p_{\delta_0,h}).\end{array}
\] 
We have
\[
\begin{array}{l}
\frac{d}{d\delta_0}U(\epsilon_n(P),p_{\delta_0,h})=
\tilde{P}_n\frac{ \frac{d}{d\delta_0}D_{n,P_{\delta_0,h} }}{(1+\epsilon_n(P)D_{n,P})}\\
-\frac{D_{n,P} }{(1+\epsilon_n(P)D_{n,P})^2}\epsilon_n(P) \frac{d}{d\delta_0}D_{n,p_{\delta_0,h}}.
\end{array}
\]
Let $A_{D_n,P}(h)=\frac{d}{d\delta_0}D_{n,p_{\delta_0,h}}$.
Then,
\[
\begin{array}{l}
A_{\epsilon_n,P}(h)=\frac{d}{d\delta_0}\epsilon_n(P_{\delta_0,h} )\\
= c_{n,P}^{-1} 
\tilde{P}_n \left\{ \frac{A_{D_n,P}(h)}{(1+\epsilon_n(P)D_{n,P})}-\frac{D_{n,P}}{(1+\epsilon_n(P)D_{n,P})^2}\epsilon_n(P) A_{D_n,P}(h) \right\} \\
=c_{n,P}^{-1}\tilde{P}_n \frac{A_{D_n,P}(h)}{(1+\epsilon_n(P)D_{n,P})^2} \\
=c_{n,P}^{-1}P \frac{d\tilde{P}_n}{dP}  \frac{A_{D_n,P}(h)}{(1+\epsilon_n(P)D_{n,P})^2} \\
=c_{n,P}^{-1} P A_{D_n,P}^{\top}\left( \frac{d\tilde{P}_n}{dP}  \frac{1}{(1+\epsilon_n(P)D_{n,P})^2} \right) h\\
=\langle D_{\epsilon_n,P},h\rangle_{P},
\end{array}
\]
where we used the definition of the adjoint $A_{D_n,P}^{\top}$ of $A_{D_n,P}$, and defined
\[
D_{\epsilon_n,P}=c_{n,P}^{-1}A_{D_n,P}^{\top}\left( \frac{d\tilde{P}_n}{dP}  \frac{1}{(1+\epsilon_n(P)D_{n,P})^2} \right)\in L^2_0(P).\]

{\bf Score operator $A_{n,P}(h)$:}
By the product rule, we have
\[
\begin{array}{l}
A_{n,P}(h)=\frac{d}{d\delta_0}\tilde{p}_n(p_{\delta_0,h})/\tilde{p}_n(p)\\
=\frac{d}{d\delta_0}(1+\epsilon_n(p_{\delta_0,h})D_{n,p_{\delta_0,h}})p_{\delta_0,h}/\tilde{p}_n(p)\\
=\{\tilde{p}_n(p)\}^{-1}\left\{ hp+A_{\epsilon_n,P}(h)D_{n,P}p+
\epsilon_n(P)A_{D_n,P}(h)p+\epsilon_n(P)D_{n,P} hp\right\} \\
=\{\tilde{p}_n(p)\}^{-1}\left\{
(1+\epsilon_n(P)D_{n,P} )hp+A_{\epsilon_n,P}(h)D_{n,P}p +\epsilon_n(P)  A_{D_n,P}(h)p\right\}\\
=h+A_{\epsilon_n,P}(h)D_{n,P}\frac{p}{\tilde{p}_n(p)}+\epsilon_n(P)  A_{D_n,P}(h)\frac{p}{\tilde{p}_n(p)}.
\end{array}
\]
Let's first confirm that this is indeed an element of $L^2_0(\tilde{P}_n(P))$?
We have
\[
\begin{array}{l}
\tilde{P}_n(P) A_{n,P}(h)=\int  \left\{ (1+\epsilon_n(P)D_{n,P} )hp+A_{\epsilon_n,P}(h)D_{n,P}p +\epsilon_n(P)  A_{D_n,P}(h)p \right\} d\mu \\
=\int   (1+\epsilon_n(P)D_{n,P} )hp d\mu+\int \epsilon_n(P)A_{D_n,P}(h)p d\mu\\
=\epsilon_n(p) \int D_{n,P} hp d\mu+\epsilon_n(P)\int A_{D_n,P}(h)p d\mu\\
=\epsilon_n(P)\{\int D_{n,P}hpd\mu+\int \frac{d}{d\delta_0}D_{n,P_{\delta_0,h}} p d\mu\}\\
=0,
\end{array}
\]
where we use that $\frac{d}{d\delta_0}E_P D_{n,p_{\delta_0,h}}=-E_P D_{n,p} h$.
So indeed the expectation equals zero, proving $A_{n,P}(h)\in L^2_0(\tilde{P}_n(P))$.

{\bf Deriving the next canonical gradient:}
We will now determine the next canonical gradient.
We have 
\[
\begin{array}{l}
\frac{d}{d\delta_0}\Psi_n(\tilde{P}_n(P_{\delta_0,h}))=
\langle D_{\tilde{P}_n(P))},A_{n,P}(h)\rangle_{\tilde{P}_n(P)}\\
=\langle D_{n,\tilde{P}_n(P)}, (1+\epsilon_n(P)D_{n,P} )hp+A_{\epsilon_n,P}(h)D_{n,P}p +\epsilon_n(P)  A_{D_n,P}(h) p\rangle_{\mu}\\
=E_P \{D_{n,\tilde{P}_n(P)}(1+\epsilon_n(P)D_{n,P} )h\}+
E_P \{D_{n,\tilde{P}_n(P)}A_{\epsilon_n,P}(h)D_{n,P} \} \\
\hfill +E_P \{D_{n,\tilde{P}_n(P)} \epsilon_n(P)  A_{D_n,P}(h) \}\\
=E_P h\{ D_{n,\tilde{P}_n(P)}(1+\epsilon_n(P)D_{n,P} )\}+
\langle D_{\epsilon_n,P},h\rangle_P  E_P \{D_{n,\tilde{P}_n(P)}D_{n,P} \}\\
\hfill 
+E_P h  \epsilon_n(P) A_{D_n,P}^{\top}\left ( D_{n,\tilde{P}_n^{(1)}(P)}\right )\\
 =\langle A_{n,P}^{\top} (D_{n,\tilde{P}_n^{(1)}(P)}),h\rangle_P\\
 =\langle D^{+}_{n,P},h\rangle_P,
 \end{array}
 \]
 where the next canonical gradient $D^{+}_{n,P}$ is thus given by
 \[
 \begin{array}{l}
 D^{+}_{n,P}=A_{n,P}^{\top}(D_{n,\tilde{P}_n(P)})=
 D_{n,\tilde{P}_n(P)}(1+\epsilon_n(P)D_{n,P} )
 +D_{\epsilon_n,P}E_P \{D_{n,\tilde{P}_n(P)}D_{n,P} \}\\
\hfill  + \epsilon_n(P) A_{D_n,P}^{\top}\left( D_{n,\tilde{P}_n(P)}\right).\end{array}
 \]
This completes the proof. $\Box$

\section{Lemmas for  Example I: Second order TMLE of integrated square of density} \label{AppendixG}
 \subsection{Proof of Lemma \ref{lemmaex1a}.}

  By Lemma \ref{lemmaformulacangradient}  
  we have
  \[
\begin{array}{l}
D^{(2)}_{n,P}=A_{n,P}^{(1),\top}(D^{(1)}_{\tilde{P}_n^{(1)}(P)} )\\
 =
 D^{(1)}_{\tilde{P}_n^{(1)}(P)}(1+\epsilon_n^{(1)}(P)D_{P}^{(1)} )
 +D_{\epsilon_n^{(1)},P}E_P \{D^{(1)}_{\tilde{P}_n^{(1)}(P)}D^{(1)}_{P} \}\\
 \hfill + \epsilon_n^{(1)}(P) A_{D^{(1)},P}^{\top}\left( D^{(1)}_{\tilde{P}_n^{(1)}(P)}\right).
 \end{array}
 \]
 We have $A_{D^{(1)},P}:L^2_0(P)\rightarrow L^2(P)$ is given by
  \[
A_{D^{(1)},P}(h)=  \frac{d}{d\delta_0}D^{(1)}_{p_{\delta_0,h}}=2 hp-4\int h p^2 d\mu .\]
Let's compute its adjoint $A_{D^{(1)},P}^{\top}:L^2(P)\rightarrow L^2_0(P)$.
For a $h\in L^2_0(P)$ and $V\in L^2(P)$ we have
\begin{eqnarray*}
E_P A_{D^{(1)},P}(h) V&=& E_P 2ph V-4\int h p dP E_P V\\
&=&\langle h, 2p (V-2 E_PV)\rangle_P.
\end{eqnarray*}
We still need to center it to have mean zero under $P$. Thus, $A_{D^{(1)},P}^{\top}:L^2(P)\rightarrow L^2_0(P)$ is given by
\[
A_{D^{(1)},P}^{\top}(V)=2p(V-2E_PV)-\int 2p(V-2E_PV) dP.\]
Let $c_{n,P}^{(1)}=\tilde{P}_n \frac{D^{(1),2}_P}{(1+\epsilon_n^{(1)}(P)D^{(1)}_P)^2}$.
We now determine $D_{\epsilon_n^{(1)},P}=c_{n,P}^{-1}A_{D^{(1)},P}^{\top}(V)$ with 
$V=(\tilde{p}_n/p) 1/(1+\epsilon_n^{(1)}(P)D^{(1)}_P)^2$.
We have
\[
\begin{array}{l}
c_{n,P}D_{\epsilon_n^{(1)},P}=
2p(V-2E_PV)-\int 2p(V-2E_PV) dP\\
=2p\left\{ (\tilde{p}_n/p) \frac{1}{(1+\epsilon_n^{(1)}(P)D^{(1)}_P)^2} -
2\tilde{P}_n \frac{1}{(1+\epsilon_n^{(1)}(P)D^{(1)}_P)^2} \right\}\\
-2 P p\left\{ (\tilde{p}_n/p) \frac{1}{(1+\epsilon_n^{(1)}(P)D^{(1)}_P)^2}-2
\tilde{P}_n  \frac{1}{(1+\epsilon_n^{(1)}(P)D^{(1)}_P)^2} \right\}\\
=2\tilde{P}_n  \frac{1}{(1+\epsilon_n^{(1)}(P)D^{(1)}_P)^2} 
-4p \tilde{P}_n \frac{1}{(1+\epsilon_n^{(1)}(P)D^{(1)}_P)^2} \\
-2P \tilde{p}_n  \frac{1}{(1+\epsilon_n^{(1)}(P)D^{(1)}_P)^2}+4\Psi(P) \tilde{P}_n  \frac{1}{(1+\epsilon_n^{(1)}(P)D^{(1)}_P)^2}.\end{array}
\]
In particular, at $P=\tilde{P}_n$ this expression reduces to  $D_{\epsilon_n^{(1)},P}=-\left\{\tilde{P}_n D^{(1),2}_{\tilde{P}_n}\right\}^{-1}
D^{(1)}_{\tilde{P}_n}$.

 We also need to determine  $A_{D^{(1)},P}^{\top}(V)$ with $V=D^{(1)}_{\tilde{P}_n^{(1)}(P)}$.
We have
\[
\begin{array}{l}
A_{D^{(1)},P}^{\top}\left( D^{(1)}_{\tilde{P}_n^{(1)}(P)}\right)=
2p(V-2E_PV)-\int 2p(V-2E_PV) dP\\
=2p\left(D^{(1)}_{\tilde{P}_n^{(1)}(P)}-2E_P D^{(1)}_{\tilde{P}_n^{(1)}(P)}\right)
-2P p\left( D^{(1)}_{\tilde{P}_n^{(1)}(P)}-2E_P D^{(1)}_{\tilde{P}_n^{(1)}(P)}\right).
\end{array}
\]
Define the scalar
\[
d_P\equiv  E_P \{D^{(1)}_{\tilde{P}_n^{(1)}(P)}D^{(1)}_{P} \} c_{n,P}^{-1}.\]
 Thus,
\[
\begin{array}{l}
D^{(2)}_{n,P}=A_{n,P}^{(1),\top}(D^{(1)}_{\tilde{P}_n^{(1)}(P)} )\\
 =
 D^{(1)}_{\tilde{P}_n^{(1)}(P)}(1+\epsilon_n^{(1)}(P)D_{P}^{(1)} )\\
 +2d_P
 \left\{ \tilde{p}_n \frac{1}{(1+\epsilon_n^{(1)}(P)D^{(1)}_P)^2} -
2p \tilde{P}_n \frac{1}{(1+\epsilon_n^{(1)}(P)D^{(1)}_P)^2} \right\}\\
-2d_P  P \left\{ \tilde{p}_n\frac{1}{(1+\epsilon_n^{(1)}(P)D^{(1)}_P)^2}\right\}\\
+4d_P \Psi(P)
\tilde{P}_n  \frac{1}{(1+\epsilon_n^{(1)}(P)D^{(1)}_P)^2} \\
  +2 \epsilon_n^{(1)}(P) p\left(D^{(1)}_{\tilde{P}_n^{(1)}(P)}-2E_P D^{(1)}_{\tilde{P}_n^{(1)}(P)}\right)\\
-2\epsilon_n^{(1)}(P) P\left\{ p\left( D^{(1)}_{\tilde{P}_n^{(1)}(P)}-2E_P D^{(1)}_{\tilde{P}_n^{(1)}(P)}\right)\right\}.  
  \end{array}
 \]
 This completes the proof.  $\Box$

\subsection{Proof of Lemma \ref{lemmaexample1}.}
We will now work out the different contributions so that we can showcase the type of third order differences in the expression for $R^{(1),+}$.

We note that  $\frac{d}{dp_n^*}\tilde{p}_n^{(1)}(p_n^*)(\tilde{p}_n-p_n^*)=A^{(1)}_{p_n^*}((\tilde{p}_n-p_n^*)/p_n^*)$ can be expressed in terms of the score operator $A^{(1)}_p(h)=\frac{d}{d\delta_0}\tilde{p}_n^{(1)}(p_{\delta_0,h})$ we derived in the Appendix as part of proof of representation $D^{(2)}_{n,P}$.  
Analogue to the calculations of the score operator $A^{(1)}_p$, we have
\[
\begin{array}{l}
\tilde{p}_n^{(1)}(\tilde{p}_n)-\tilde{p}_n^{(1)}(p)=
(1+\epsilon_n^{(1)}(\tilde{p}_n)D^{(1)}_{\tilde{p}_n})\tilde{p}_n-(1+\epsilon_n^{(1)}(p)D^{(1)}_p)p\\
=(\tilde{p}_n-p)+(\epsilon_n^{(1)}(\tilde{p}_n)-\epsilon_n^{(1)}(p))D^{(1)}_pp+
\epsilon_n^{(1)}(\tilde{p}_n)\{D^{(1)}_{\tilde{p}_n}-D^{(1)}_p\} p\\
\hfill +\epsilon_n^{1}(\tilde{p}_n)D^{(1)}_{\tilde{p}_n}(\tilde{p}_n-p)\\
=(1+\epsilon_n^{(1)}(\tilde{p}_n)D^{(1)}_{\tilde{p}_n})(\tilde{p}_n-p)+(\epsilon_n^{(1)}(\tilde{p}_n)-\epsilon_n^{(1)}(p)) D^{(1)}_p p\\
+
\epsilon_n^{(1)}(\tilde{p}_n)\{2(\tilde{p}_n-p) -2(\Psi(\tilde{p}_n)-\Psi(p))\}p\\
=(1+\epsilon_n^{(1)}(\tilde{p}_n)D^{(1)}_{\tilde{p}_n})(\tilde{p}_n-p)+
D^{(1)}_p p \left\{ \frac{d}{dp}\epsilon_n^{(1)}(p)(\tilde{p}_n-p)+R_{2,\epsilon_n^{(1)},p}(\tilde{p}_n,p)\right\}\\
+\epsilon_n^{(1)}(\tilde{p}_n)\{2(\tilde{p}_n-p) -2(\Psi(\tilde{p}_n)-\Psi(p))\}p \\
=(\tilde{p}_n-p)+(\epsilon_n^{(1)}(p)D^{(1)}_p)(\tilde{p}_n-p)+
D^{(1)}_p p \frac{d}{dp}\epsilon_n^{(1)}(p)(\tilde{p}_n-p)\\
+\epsilon_n^{(1)}(p)p \frac{d}{dp}D^{(1)}_p(\tilde{p}_n-p)\\
+(\epsilon_n^{(1)}(\tilde{p}_n)D^{(1)}_{\tilde{p}_n}-\epsilon_n^{(1)}(p)D^{(1)}_p)(\tilde{p}_n-p)\\
+
D^{(1)}_p p R_{2,\epsilon_n^{(1)},p}(\tilde{p}_n,p)\\
-2\epsilon_n^{(1)}(p)p \int (\tilde{p}_n-p)^2 dx\\
\equiv \frac{d}{dp}\tilde{p}_n^{(1)}(p)(\tilde{p}_n-p)+R_{2,p,\tilde{p}_n^{(1)}}(\tilde{p}_n,p),\end{array}
\]
where we defined $R_{2,p,\epsilon()}(\tilde{p}_n,p)=\epsilon_n^{(1)}(\tilde{p}_n)-\epsilon_n^{(1)}(p)-\frac{d}{dp}\epsilon_n^{(1)}(p)(\tilde{p}_n-p)$ and we used that
\[
\begin{array}{l}
D^{(1)}_{\tilde{p}_n}-D^{(1)}_p=
2(\tilde{p}_n-p)-2\int \tilde{p}_n^2+2\int p^2\\
=
2(\tilde{p}_n-p)-2\int (\tilde{p}_n-p)\tilde{p}_n-2\int p(\tilde{p}_n-p)=2(\tilde{p}_n-p)-2\int (\tilde{p}_n+p) (\tilde{p}_n-p)\\
=2(\tilde{p}_n-p)-2\int (2p+\tilde{p}_n-p)(\tilde{p}_n-p)\\
=2(\tilde{p}_n-p)-2\int 2p (\tilde{p}_n-p)-2\int(\tilde{p}_n-p)^2\\
\equiv \frac{d}{dp}D^{(1)}_p(\tilde{p}_n-p)-2\int (\tilde{p}_n-p)^2.
\end{array}
\]
Thus, $\frac{d}{dp}D^{(1)}_p(\tilde{p}_n-p)=2(\tilde{p}_n-p)-4\int p(\tilde{p}_n-p)dx$.
We also concluded that the exact second order remainder in first order Tailor expansion of $\tilde{p}_n^{(1)}(p)-\tilde{p}_n^{(1)}(p_n^*)$ is given by
\[
\begin{array}{l}
R_{2,p,\tilde{p}_n^{(1)}}(\tilde{p}_n,p)=
(\epsilon_n^{(1)}(\tilde{p}_n)D^{(1)}_{\tilde{p}_n}-\epsilon_n^{(1)}(p)D^{(1)}_p)(\tilde{p}_n-p)\\
+
D^{(1)}_p p R_{2,\epsilon_n^{(1)},p}(\tilde{p}_n,p)
-2\epsilon_n^{(1)}(p)p \int (\tilde{p}_n-p)^2 dx\\
\equiv \frac{d}{dp}\tilde{p}_n^{(1)}(p)(\tilde{p}_n-p)+R_{2,p,\tilde{p}_n^{(1)}}(\tilde{p}_n,p).
\end{array}
\]
Thus, we conclude that
\[
\begin{array}{l}
-{R}^{(1),+}(\tilde{p}_n^{(1)}(p_n^*),\tilde{p}_n)=
\int (\tilde{p}_n^{(1)}(\tilde{p}_n)-\tilde{p}_n^{(1)}(p_n^*))^2dx \\
=\int (\tilde{p}_n^{(1)}(\tilde{p}_n)-\tilde{p}_n^{(1)}(p_n^*))R_{2,p_n^*,\tilde{p}_n^{(1)}()}(\tilde{p}_n,p_n^*) dx\\
=\int (\tilde{p}_n^{(1)}(\tilde{p}_n)-\tilde{p}_n^{(1)}(p_n^*))(\epsilon_n^{(1)}(\tilde{p}_n)D^{(1)}_{\tilde{p}_n}-\epsilon_n^{(1)}(p)D^{(1)}_p)(\tilde{p}_n-p)dx\\
+\int (\tilde{p}_n^{(1)}(\tilde{p}_n)-\tilde{p}_n^{(1)}(p_n^*))D^{(1)}_p p R_{2,\epsilon_n^{(1)},p}(\tilde{p}_n,p)dx\\
-2\epsilon_n^{(1)}(p) \int (\tilde{p}_n-p)^2 dx\int (\tilde{p}_n^{(1)}(\tilde{p}_n)-\tilde{p}_n^{(1)}(p_n^*)) p dx\\
\end{array}
\]

\section{Proof of Lemmas for treatment specific mean example}\label{AppendixH}

\subsection{Proof of  Lemma \ref{lemmaexample2a}.}
 Note that
$\tilde{q}_n^{(1)}(q_{\delta,h_1},g)$ is a logistic regression model with $\delta C_1+\epsilon_n(q_{\delta,h},g)C_g$ in its linear form.
Therefore, we know 
\[
\begin{array}{l}
A_{n,p,1}(h_1)=\frac{d}{d\delta_0}\log \tilde{q}_n^{(1)}(q_{\delta_0},g)\\
=C_1(W,A)(Y-\tilde{q}_n^{(1)}(q,g)(1\mid W,A))\\
+\frac{d}{d\delta_0}\epsilon_n(q_{\delta_0},g)C_g(Y-\tilde{q}_n^{(1)}(q,g)(1\mid W,A)).
\end{array}
\]
Since $C_1(W,A)=h_1(W,A,Y)/(Y-\bar{q})$ this is indeed a linear mapping in $h_1$.
 

{\bf Score operator mapping score of path through initial $g$ into score of TMLE update:}
Note that $\tilde{q}_n^{(1)}(q,g_{\delta,h_2})$ is a logistic regression model with $
\epsilon_n(q,g_{\delta})C_{g_{\delta}}$ in its linear form. 
Also note that $\frac{d}{d\delta}\bar{g}_{\delta}=C_2\bar{g}(1-\bar{g})$.
We have 
\[
\begin{array}{l}
\frac{d}{d\delta_0}C_{g_{\delta_0}}=-\frac{A}{\bar{g}^2} \frac{d}{d\delta_0}\bar{g}_{\delta_0}
=-\frac{A}{\bar{g}^2}C_2\bar{g}(1-\bar{g})=-\frac{A}{\bar{g}}C_2(1-\bar{g}).
\end{array}
\]
Thus,
\[
\begin{array}{l}
A_{n,p,2}(h_2)=h_2+\frac{d}{d\delta_0}\log \tilde{q}_n^{(1)}(q,g_{\delta_0})\\
=h_2+\frac{d}{d\delta_0}\epsilon_n(q,g_{\delta_0}) C_g(Y-\tilde{q}_n^{(1)}(1\mid W,A))\\
+\epsilon_n(q,g)\frac{d}{d\delta}C_{g_{\delta_0}}(Y-\tilde{q}_n^{(1)}(1\mid W,A))\\
=h_2+\frac{d}{d\delta_0}\epsilon_n(q,g_{\delta_0}) C_g(Y-\tilde{q}_n^{(1)}(1\mid W,A))\\
-\epsilon_n(q,g)\frac{A}{\bar{g}}C_2(1-\bar{g}) (Y-\tilde{q}_n^{(1)}(1\mid W,A))\\
=h_2+\frac{d}{d\delta_0}\epsilon_n(q,g_{\delta_0}) C_g(Y-\tilde{q}_n^{(1)}(1\mid W,A))\\
-\epsilon_n(q,g)\frac{A}{\bar{g}}h_{2g}(W,A)  (Y-\tilde{q}_n^{(1)}(1\mid W,A)).
\end{array}
\]
Here we used that  $C_2(1-\bar{g})=h_2(1,W)$ so that we have succeeded to express the score operator as a linear mapping in $h_{2g}$.

 Let $A_{n,p,1}:H_1(P)\rightarrow L^2_0(\tilde{P}_n^{(1)}(P))$, where $H_1(P)=\{h_1\in L^2_0(P): E(h_1\mid W,A)=0\}$. We have $H_1(P)=\{C_1(Y-\bar{q}(W,A)):C_1\}$.
 Let $A_{n,p,2}:H_2(P)\rightarrow L^2_0(\tilde{P}_n^{(1)}(P))$, where $H_2(P)=\{h_2(W,A): E(h_2\mid W)=0\}$. We have $H_2(P)=\{C_2(A-\bar{g}(W)):C_2\}$.
Note that $H_1(P)\perp H_2(P)$ are orthogonal spaces in $L^2_0(P)$.
Above we derived the  form of these score operators $A_{n,p,1}$ and $A_{n,p,2}$ up till the pathwise derivative of $\epsilon_n^{(1)}(q,g)$:
\begin{eqnarray*}
A_{n,p,1}(h_1)&=&C_1(W,A)(Y-\tilde{q}_n^{(1)}(q,g)(1\mid W,A))\\
&&+\frac{d}{d\delta_0}\epsilon_n(q_{\delta_0},g)C_g(Y-\tilde{q}_n^{(1)}(q,g)(1\mid W,A))\\
A_{n,p,2}(h_2)&=&h_2+ \frac{d}{d\delta_0}\epsilon_n(q,g_{\delta_0}) C_g(Y-\tilde{q}_n^{(1)}(1\mid W,A))
\\
&&-\epsilon_n(q,g)\frac{A}{\bar{g}}h_{2g}(W,A)  (Y-\tilde{q}_n^{(1)}(1\mid W,A)).
\end{eqnarray*}

{\bf Representing second order canonical gradient and understanding its two components:}
Let $A_{n,p,1}^{\top}:L^2_0(\tilde{P}_n^{(1)}(P))\rightarrow H_1(P)$ and $A_{n,p,2}^{\top}:L^2_0(\tilde{P}_n^{(1)}(P))\rightarrow H_2(P))$ be the adjoints of $A_{n,p,1}$ and $A_{n,p,2}$, respectively. 
It follows that
\[
\begin{array}{l}
\frac{d}{d\delta_0}\Psi_n^{(1)}(p_{\delta_0,h})=
\tilde{P}_n^{(1)}(p) D^{(1)}_{\tilde{p}_n^{(1)}(p)}\{A_{n,p,1}(h_1)+A_{n,p,2}(h_2)\} \\
=P A_{n,p,1}^{\top}(D^{(1)}_{\tilde{p}_n^{(1)}(p)}) h_1+
P A_{n,p,2}^{\top} (D^{(1)}_{\tilde{p}_n^{(1)}(p)}) h_2\\
=P A_{n,p,1}^{\top}(D^{(1)}_{\tilde{p}_n^{(1)}(p)}) (h_1+h_2)+
P A_{n,p,2}^{\top} (D^{(1)}_{\tilde{p}_n^{(1)}(p)}) (h_1+h_2)\\
=P (h_1+h_2)\left\{A_{n,p,1}^{\top}(D^{(1)}_{\tilde{p}_n^{(1)}(p)})+A_{n,p,2}^{\top}(D^{(1)}_{\tilde{p}_n^{(1)}(p)} )\right\}.
\end{array}
\]
This shows that the second order canonical gradient can be represented as
\[
D^{(2)}_{n,P}=A_{n,p,1}^{\top}(D^{(1)}_{\tilde{p}_n^{(1)}(p)})+A_{n,p,2}^{\top}(D^{(1)}_{\tilde{p}_n^{(1)}(p)} ).\] This corresponds with our lemma in previous section: $D^{(j+1)}_{n,P}=A_{n,p}^{\top}(D^{(j)}_{n,\tilde{P}_n^{(j)}(P)})$.

{\bf Separate second order canonical gradient components for targeting $q_Y$ and $g$:}
We should also note that $A_{n,p,1}^{\top}(D^{(1)}_{\tilde{p}_n^{(1)}(p)})$ is the canonical gradient of $q\rightarrow \Psi(\tilde{P}_n^{(1)}(q,g))$ and $A_{n,p,2}^{\top}(D^{(1)}_{\tilde{p}_n^{(1)}(p)})$ is the canonical gradient of $g\rightarrow \Psi(\tilde{P}_n^{(1)}(q,g))$.
We will have computed both components $D^{(2)}_{n,P,1} =A_{n,p,1}^{\top}(D^{(1)}_{\tilde{p}_n^{(1)}(p)})$ and $D^{(2)}_{n,P,2} =A_{n,p,2}^{\top}(D^{(1)}_{\tilde{p}_n^{(1)}(p)})$
of $D^{(2)}_{n,P}=D^{(2)}_{n,P,1}+D^{(2)}_{n,P,2}$. One can be used to target $q$ and the other can be used to target $g$. Targeting these two components of the initial estimator in the first order TMLE  separately will provide additional representations of the third order remainder of the second order TMLE, which can thereby be beneficial relative to only aiming for $D^{(2)}_{n,P}$.

{\bf Determining the pathwise derivative of $\epsilon_n^{(1)}(p)$:}
We have that $\epsilon_n^{(1)}(p)$ solves
\[
\begin{array}{l}
0=\tilde{P}_n C_g(W,A)(Y-\bar{q}^{(1)}(q,g,\epsilon)(W,A)).\\
\end{array}
\]
Let's denote this equation with $U(\epsilon(p),p=(q,g))$. 
We have $U(\epsilon(p_{\delta}),p_{\delta})=0$. So, at $\epsilon =\epsilon_n^{(1)}(p)$, we have
\[
\frac{d}{d\delta_0}\epsilon_n^{(1)}(p_{\delta_0})=-\left\{\frac{d}{d\epsilon}U(\epsilon,p)\right\}^{-1}
\frac{d}{d\delta_0}U(\epsilon,p_{\delta_0}).\]
We have
\[
c_P^{(1)}\equiv -\frac{d}{d\epsilon}U(\epsilon,p)=
\tilde{P}_n C_g\frac{d}{d\epsilon}\bar{q}^{(1)}(q,g,\epsilon)(W,A).\]
Note that 
\[
\begin{array}{l}
\frac{d}{d\delta_0}\epsilon_n^{(1)}(q_{\delta_0},g)=c_P^{(1),-1}\frac{d}{d\delta_0}U(\epsilon_n^{(1)}(p),q_{\delta_0},g)\\
\frac{d}{d\delta_0}\epsilon_n^{(1)}(q,g_{\delta_0})=
c_P^{(1),-1}\frac{d}{d\delta_0}U(\epsilon_n^{(1)}(p),q,g_{\delta_0}).
\end{array}
\]

We know that
\[
\frac{d}{d\epsilon}\bar{q}(p,\epsilon)=
\bar{q}^{(1)}(p,\epsilon)(1-\bar{q}^{(1)}(p,\epsilon) )C_g.\]
So we obtain
\[
\begin{array}{l}
c_P^{(1)}=-\frac{d}{d\epsilon}U(\epsilon,p)=
\tilde{P}_n C_g^2\bar{q}^{(1)}(p,\epsilon)(1-\bar{q}^{(1)}(p,\epsilon))\\
=P \frac{\tilde{g}_n}{\bar{g}^2} \bar{q}^{(1)}(p,\epsilon)(1-\bar{q}^{(1)}(p,\epsilon)).
\end{array}
\]

We have \[
\begin{array}{l}
\frac{d}{d\delta}U(\epsilon,q_{\delta,h_1},g)=
\frac{d}{d\delta}
\tilde{P}_n C_g\left(Y-\frac{1}{1+\exp(-\log (\bar{q}_{\delta}/(1-\bar{q}_{\delta}))-\epsilon C_g)}\right) \\
=-\frac{d}{d\delta}
\tilde{P}_nC_g\frac{1}{1+\exp(-\{ \log (\bar{q}/(1-\bar{q}))+delta C_1 +\epsilon C_g\} )}\\
=-\tilde{P}_n C_g C_1 \bar{q}^{(1)}(q,g,\epsilon)(1-\bar{q}^{(1)}(q,g,\epsilon))
.
\end{array}
\]
We have  
\[
\begin{array}{l}
\frac{d}{d\delta}U(\epsilon,q,g_{\delta,h_2})=
\frac{d}{d\delta}
\tilde{P}_n C_{g_{\delta}}\left(Y-\frac{1}{1+\exp(-\{ \log (\bar{q}/(1-\bar{q}))+\epsilon C_{g_{\delta}}\}) }\right)\\
=
\tilde{P}_n\frac{d}{d\delta} C_{g_{\delta}}\left (Y-\frac{1}{1+\exp(-\{ \log (\bar{q}/(1-\bar{q}))+\epsilon C_g \}) }\right)\\
-\tilde{P}_n C_g\epsilon \frac{d}{d\delta}C_{g_{\delta}} \bar{q}^{(1)}(q,g,\epsilon)(1-\bar{q}^{(1)}(q,g,\epsilon))\\
=-\tilde{P}_n\frac{A}{\bar{g}^2}C_2\bar{g}(1-\bar{g}) \left(Y-\frac{1}{1+\exp(-\{ \log (\bar{q}/(1-\bar{q}))+\epsilon C_g \}) }\right)\\
+\tilde{P}_n C_g\epsilon \frac{A}{\bar{g}^2}C_2\bar{g}(1-\bar{g}) \bar{q}^{(1)}(q,g,\epsilon)(1-\bar{q}^{(1)}(q,g,\epsilon))\\
=-\tilde{P}_n\frac{A}{\bar{g}}C_2(1-\bar{g}) \left(Y-\bar{q}^{(1)}(q,g,\epsilon)\right) \\
+\tilde{P}_n C_g\epsilon \frac{A}{\bar{g}}C_2(1-\bar{g}) \bar{q}^{(1)}(q,g,\epsilon)(1-\bar{q}^{(1)}(q,g,\epsilon)).
\end{array}
\]

So we have shown
\[
\begin{array}{l}
\frac{d}{d\delta_0}\epsilon_n^{(1)}(q_{\delta_0,h_1},g)=
-c_P^{(1),-1}\tilde{P}_n C_g C_1 \bar{q}^{(1)}(q,g,\epsilon_n^{(1)}(p))(1-\bar{q}^{(1)}(q,g,\epsilon_n^{(1)}(p)))\\
\frac{d}{d\delta_0}\epsilon_n^{(1)}(q,g_{\delta_0,h_2})= -c_P^{(1),-1}\tilde{P}_n\frac{A}{\bar{g}}C_2(1-\bar{g}) \left(Y-\bar{q}^{(1)}(q,g,\epsilon_n^{(1)}(p))(1)\right) \\
\hfill +c_P^{(1),-1}\tilde{P}_n C_g\epsilon_n^{(1)}(p) \frac{A}{\bar{g}}C_2(1-\bar{g}) \bar{q}^{(1)}(q,g,\epsilon)(1-\bar{q}^{(1)}(q,g,\epsilon))
\end{array}
\]
{\bf Final forms of score operators $A_{n,p,1}$ and $A_{n,p,2}$:}
Thus,
\[
\begin{array}{l}
A_{n,p,1}(h_1)=\frac{d}{d\delta_0}\log \tilde{q}_n^{(1)}(q_{\delta_0},g)\\
=C_1(W,A)(Y-\bar{q}_n^{(1)}(q,g)(W,A))\\
+\frac{d}{d\delta_0}\epsilon_n^{(1)}(q_{\delta_0},g)C_g(Y-\bar{q}_n^{(1)}(q,g)(W,A))\\
=C_1(W,A)(Y-\bar{q}_n^{(1)}(q,g)(W,A))\\
-c_P^{(1),-1}\left\{\tilde{P}_n C_g C_1 \bar{q}^{(1)}(q,g,\epsilon)(1-\bar{q}_n^{(1)}(q,g))\right\}
C_g(Y-\bar{q}_n^{(1)}(q,g) ).
\end{array}
\]
\[
\begin{array}{l}
A_{n,p,2}(h_2)=h_2+\frac{d}{d\delta_0}\log \tilde{q}_n^{(1)}(q,g_{\delta_0})\\
=h_2-\epsilon_n^{(1)}(q,g)\frac{A}{\bar{g}}h_{2g}(W,A)  (Y-\bar{q}_n^{(1)}(q,g) )\\
+\frac{d}{d\delta_0}\epsilon_n^{(1)}(q,g_{\delta_0}) C_g(Y-\bar{q}_n^{(1)}(q,g))\\
=h_2-\epsilon_n^{(1)}(q,g)\frac{A}{\bar{g}}h_{2g}(W,A)  (Y-\bar{q}_n^{(1)}(q,g))\\
-c_P^{(1),-1}\left\{ \tilde{P}_n\frac{A}{\bar{g}}C_2(1-\bar{g}) \left(Y-\bar{q}_n^{(1)}(q,g)\right)\right\}
C_g(Y-\bar{q}_n^{(1)}(q,g))\\
+c_P^{(1),-1}\left\{ \tilde{P}_n C_g\epsilon_n^{(1)}(p) \frac{A}{\bar{g}}C_2(1-\bar{g}) \bar{q}_n^{(1)}(q,q)(1-\bar{q}_n^{(1)}(q,g))\right\} C_g(Y-\bar{q}_n^{(1)}(q,g)).
\end{array}
\]
This proves Lemma \ref{lemmaexample2a}.

\subsection{Proof of Lemma \ref{lemmaexample2b}.}

{\bf Determining the adjoint of $A_{n,p,1}$:}
We have
\[
\begin{array}{l}
\tilde{P}_n^{(1)}(p) D^{(1)}_{\tilde{P}_n^{(1)}(p)} A_{n,p,1}(h_1)=
\tilde{P}_n^{(1)}(p) D^{(1)}_{\tilde{P}_n^{(1)}(p)}C_1(Y-\bar{q}_n^{(1)}(p))\\
-\tilde{P}_n^{(1)}(p) D^{(1)}_{\tilde{P}_n^{(1)}(p)}C_g(Y-\bar{q}_n^{(1)}(p)) c_P^{(1),-1} \left\{\tilde{P}_n C_g C_1\bar{q}_n^{(1)}(1-\bar{q}_n^{(1)})(p)\right\}\\
=\tilde{P}_n^{(1)}(p) (C_g C_1(Y-\bar{q}_n^{(1)}(p))^2)\\
-\tilde{P}_n^{(1)}(p) C_g^2(Y-\bar{q}_n^{(1)})^2c_P^{(1),-1}\left\{\tilde{P}_n C_g C_1\bar{q}_n^{(1)}(1-\bar{q}_n^{(1)})(p)\right\}\\
=\tilde{P}_n^{(1)}(p) \{C_g C_1\bar{q}_n^{(1)}(1-\bar{q}_n^{(1)})(p) \}\\
-\tilde{P}_n^{(1)}(p) \{C_g^2(
\bar{q}_n^{(1)}(1-\bar{q}_n^{(1)})(p) \}c_P^{(1),-1}\left\{\tilde{P}_n C_g C_1\bar{q}_n^{(1)}(1-\bar{q}_n^{(1)})(p)\right\}\\
=P \{\frac{A}{\bar{g}} C_1\bar{q}_n^{(1)}(1-\bar{q}_n^{(1)})(p) \}\\
-P \{\frac{A}{\bar{g}^2} (
\bar{q}_n^{(1)}(1-\bar{q}_n^{(1)})(p) \}c_P^{(1),-1}\left\{P C_g\frac{\tilde{g}_n}{\bar{g}} C_1\bar{q}_n^{(1)}(1-\bar{q}_n^{(1)})(p)\right\}\\
=P \{\frac{A}{\bar{g}} \frac{1}{Y-\bar{q}} C_1(Y-\bar{q})\bar{q}_n^{(1)}(1-\bar{q}_n^{(1)})(p) \}\\
-P \{\frac{1}{\bar{g}}
\bar{q}_n^{(1)}(1-\bar{q}_n^{(1)})(p) \}c_P^{(1),-1}\left\{P C_g\frac{\tilde{g}_n}{\bar{g}} \frac{1}{Y-\bar{q}}C_1(Y-\bar{q})\bar{q}_n^{(1)}(1-\bar{q}_n^{(1)})(p)\right\}\\
=P C_1(Y-\bar{q})\frac{A}{\bar{g}} \bar{q}_n^{(1)}(1-\bar{q}_n^{(1)})(p)\left\{\frac{1}{Y-\bar{q}}
-\bar{q}/(1-\bar{q})+(1-\bar{q})/\bar{q}\right\}\\
-P \{\frac{1}{\bar{g}}
\bar{q}_n^{(1)}(1-\bar{q}_n^{(1)})(p)(1) \}c_P^{(1),-1}\\
\hfill \left\{P C_g C_1(Y-\bar{q}) \frac{\tilde{g}_n}{\bar{g}} \bar{q}_n^{(1)}(1-\bar{q}_n^{(1)})(p)\left\{ \frac{1}{Y-\bar{q}}-\frac{\bar{q}}{1-\bar{q}}+\frac{(1-\bar{q})}{\bar{q}} \right\}\right\}\\
=
P h_{1,q}\frac{A}{\bar{g}} \bar{q}_n^{(1)}(1-\bar{q}_n^{(1)})(p)\left\{\frac{1}{1-\bar{q}}+\frac{1}{\bar{q}}\right\}(Y-\bar{q})\\
-P \{\frac{1}{\bar{g}}
\bar{q}_n^{(1)}(1-\bar{q}_n^{(1)})(p)(1) \}c_P^{(1),-1}\left\{P C_g h_{1,q}\frac{ \tilde{g}_n}{\bar{g}} \bar{q}_n^{(1)}(1-\bar{q}_n^{(1)})(p)  \left\{\frac{1}{1-\bar{q}}+\frac{1}{\bar{q}}\right\}(Y-\bar{q})\right\}\\
=P h_{1,q}\frac{A}{\bar{g}} \bar{q}_n^{(1)}(1-\bar{q}_n^{(1)})(p)\frac{1}{\bar{q}(1-\bar{q})}(Y-\bar{q})\\
-P \{\frac{1}{\bar{g}}
\bar{q}_n^{(1)}(1-\bar{q}_n^{(1)})(p)(1) \}c_P^{(1),-1}\left\{P C_g h_{1,q} \frac{\tilde{g}_n}{\bar{g}} \bar{q}_n^{(1)}(1-\bar{q}_n^{(1)})(p) \frac{1}{\bar{q}(1-\bar{q})} (Y-\bar{q})\right\} .
\end{array}
\]
Thus, this shows that 
\[
\begin{array}{l}
D^{(2)}_{n,P,1}=A_{n,p,1}^{\top}(D^{(1)}_{\tilde{P}_n^{(1)}(p)})\\
=
\frac{A}{\bar{g}} \bar{q}_n^{(1)}(1-\bar{q}_n^{(1)})(p)\frac{1}{\bar{q}(1-\bar{q})}(Y-\bar{q})\\
-c_P^{(1),-1}P\{\frac{1}{\bar{g}}
\bar{q}_n^{(1)}(1-\bar{q}_n^{(1)})(p)(1)\} C_g \frac{ \tilde{g}_n}{\bar{g}} \frac{\bar{q}_n^{(1)}(1-\bar{q}_n^{(1)})(p)}{\bar{q}(1-\bar{q})} (Y-\bar{q}).
\end{array}
\]
We note that at $P=\tilde{P}_n$ we have $\bar{q}_n^{(1)}(p)=\bar{q}$, $\tilde{g}_n=\bar{g}$, $C_{\tilde{P}_n}^{(1)}=P\bar{q}(1-\bar{q})/\bar{g}$, so that it cancels with $P\{\frac{1}{\bar{g}}
\bar{q}_n^{(1)}(1-\bar{q}_n^{(1)})(p)(1) \}$, and thereby 
\[
D^{(2)}_{n,\tilde{P}_n,1}=\frac{A}{\bar{g}}(Y-\bar{q})-C_g (Y-\bar{q})=0.\]


{\bf Determining the adjoint of $A_{n,2,p}$:}
Recall that $A_{n,2,p}(h_2)$ is a sum of four terms. We will determine the adjoint  for each term separately. \nl
{\bf Term 1:}
\[
\begin{array}{l}
\tilde{P}_n^{(1)}(P)D^{(1)}_{\tilde{P}_n^{(1)}(P)} h_2=
\tilde{P}_n^{(1)}(p) \frac{A}{\bar{g}}(Y-\bar{q}_n^{(1)}(p)) h_{2}=0.
\end{array}
\]
Therefore, the adjoint obtains no contribution from this term. 
\nl
{\bf Term 2:}
We need minus contribution from following:
\[
\begin{array}{l}
\tilde{P}_n^{(1)}(P)D^{(1)}_{\tilde{P}_n^{(1)}(P)}  \epsilon_n^{(1)}(p)\frac{A}{\bar{g}}h_{2}  (Y-\bar{q}_n^{(1)}(p))\\
=\tilde{P}_n^{(1)}(P) \frac{A}{\bar{g}}(Y-\bar{q}_n^{(1)}(p))  \epsilon_n^{(1)}(p)\frac{A}{\bar{g}}h_{2}(W,A)  (Y-\bar{q}_n^{(1)}(p))\\
=
\epsilon_n^{(1)}(p) \tilde{P}_n^{(1)}(p) \frac{A}{\bar{g}^2}(Y-\bar{q}_n^{(1)})^2 h_2\\
=\epsilon_n^{(1)}(p)P \frac{A}{\bar{g}^2}\bar{q}_n^{(1)}(1-\bar{q}_n^{(1)}) h_2\\
=\epsilon_n^{(1)}(p)P \frac{A-\bar{g}}{\bar{g}^2}\bar{q}_n^{(1)}(1-\bar{q}_n^{(1)})(1) h_2.
\end{array}
\]
Thus, the contribution is 
\[ - \epsilon_n^{(1)}(p)\frac{A-\bar{g}}{\bar{g}^2}\bar{q}_n^{(1)}(1-\bar{q}_n^{(1)})(1,W).\]

{\bf Term 3:} We need minus sign from following: 

\[
\begin{array}{l}
\tilde{P}_n^{(1)}(P) D^{(1)}_{\tilde{P}_n^{(1)}(P)} c_P^{(1),-1}\left\{ \tilde{P}_n\frac{A}{\bar{g}}C_2(1-\bar{g}) \left(Y-\bar{q}_n^{(1)}(q,g)\right)\right\}  C_g(Y-\bar{q}_n^{(1)}(q,g))           \\
=\tilde{P}_n^{(1)}(p)\frac{A}{\bar{g}}(Y-\bar{q}_n^{(1)}(p)) 
C_g(Y-\bar{q}_n^{(1)}(p))c_P^{(1),-1}
\tilde{P}_n \frac{A}{\bar{g}}(1-\bar{g}) C_2(Y-\bar{q}_n^{(1)}(p))\\
=P A/\bar{g}^2 \bar{q}_n^{(1)}(1-\bar{q}_n^{(1)}) c_P^{(1),-1}\tilde{P}_n  A/\bar{g}(1-\bar{g})C_2(\tilde{q}_n-\bar{q}_n^{(1)}(p))\\
=c_P^{(1),-1} P\frac{\bar{q}_n^{(1)}(1-\bar{q}_n^{(1)}(1)}{\bar{g}}P A \frac{\bar{g}_n}{\bar{g}^2}
(1-\bar{g}) C_2(\tilde{q}_n-\tilde{q}_n^{(1)}(p))\\
=c_P^{(1),-1}P\frac{\bar{q}_n^{(1)}(1-\bar{q}_n^{(1)}(1)}{\bar{g}}
P\frac{\tilde{g}_n}{\bar{g}^2} (1-\bar{g}) (\tilde{q}_n-\bar{q}_n^{(1)} 
\frac{ A}{A-\bar{g}} C_2(A-\bar{g}))\\
 =
c_P^{(1),-1}P\frac{\bar{q}_n^{(1)}(1-\bar{q}_n^{(1)}(1)}{\bar{g}}
Ph_{2}\frac{ \tilde{g}_n}{\bar{g}^2} (1-\bar{g}) (\tilde{q}_n-\bar{q}_n^{(1)}\left\{\frac{ A}{A-\bar{g}} -\bar{g}/(1-\bar{g})\right\}\\
=c_P^{(1),-1}P\frac{\bar{q}_n^{(1)}(1-\bar{q}_n^{(1)}(1)}{\bar{g}}
Ph_{2} \frac{\tilde{g}_n}{\bar{g}^2} (1-\bar{g}) (\tilde{q}_n-\bar{q}_n^{(1)}
\frac{1}{(1-\bar{g})}(A-\bar{g}) \\
=P h_{2}\left\{ \left\{ c_P^{(1),-1}P\frac{\bar{q}_n^{(1)}(1-\bar{q}_n^{(1)}(1)}{\bar{g}} \right\}
 \frac{\tilde{g}_n}{\bar{g}^2} (\tilde{q}_n-\bar{q}_n^{(1)} (A-\bar{g})\right\}\\
 \end{array}
 \]
 So contribution is given by 
 \[
-  \left\{ c_P^{(1),-1}P\frac{\bar{q}_n^{(1)}(1-\bar{q}_n^{(1)}(1)}{\bar{g}} \right\}
 \frac{\tilde{g}_n}{\bar{g}^2} (\tilde{q}_n-\bar{q}_n^{(1)}(p) )(A-\bar{g} .
 \]
 {\bf Term 4:}
 \[
 \begin{array}{l}
\tilde{P}_n^{(1)}(P)D^{(1)}_{\tilde{P}_n^{(1)}(P)} 
c_P^{(1),-1}\left\{ \tilde{P}_n C_g\epsilon \frac{A}{\bar{g}}C_2(1-\bar{g}) \bar{q}_n^{(1)}(q,q)(1-\bar{q}_n^{(1)}(q,g) )\right\} C_g(Y-\bar{q}_n^{(1)}(q,g))\\ 
= 
 \tilde{P}_n^{(1)}\frac{ A}{\bar{g}}(Y-\bar{q}_n^{(1)}(p))C_g(Y-\bar{q}_n^{(1)}(p))c_P^{(1),-1} \tilde{P}_n
 C_g \frac{A}{\bar{g}} (1-\bar{g})\epsilon_n^{(1)}(p) C_2\bar{q}_n^{(1)}(1-\bar{q}_n^{(1)})(p)\\
 =c_P^{(1),-1}P\frac{A}{\bar{g}^2} \bar{q}_n^{(1)}(1-\bar{q}_n^{(1)})(p)
 \tilde{P}_n A\frac{1-\bar{g}}{\bar{g}^2}\epsilon_n^{(1)}(p)\bar{q}_n^{(1)}(1-\bar{q}_n^{(1)}) C_2\\
 \equiv d_n P\frac{\tilde{g}_n}{\bar{g}}  A\frac{1-\bar{g}}{\bar{g}^2}\epsilon_n^{(1)}(p)\bar{q}_n^{(1)}(1-\bar{q}_n^{(1)}) C_2 \\
 =d_n P \frac{\tilde{g}_n}{\bar{g}^3}(1-\bar{g})\epsilon_n^{(1)}(p)\bar{q}_n^{(1)}(1-\bar{q}_n^{(1)})\frac{A}{A-\bar{g}} h_{2}\\
 =d_n P \frac{\tilde{g}_n}{\bar{g}^3}(1-\bar{g})\epsilon_n^{(1)}(p)\bar{q}_n^{(1)}(1-\bar{q}_n^{(1)})(1)\frac{1}{1-\bar{g}}(A-\bar{g})  h_{2} \\
 =P h_2\left\{ d_n  \frac{\tilde{g}_n}{\bar{g}^3}\epsilon_n^{(1)}(p)\bar{q}_n^{(1)}(1-\bar{q}_n^{(1)})(1)(A-\bar{g}) \right\}.
\end{array}
\]
Thus, the contribution is given by:
\[ d_n  \frac{\tilde{g}_n}{\bar{g}^3}\epsilon_n^{(1)}(p)\bar{q}_n^{(1)}(1-\bar{q}_n^{(1)})(1)(A-\bar{g}) ,\]
where
$d_n= c_P^{(1),-1}P\frac{1}{\bar{g}} \bar{q}_n^{(1)}(1-\bar{q}_n^{(1)})(p)$.

Summing up the contributions yields:
\[
\begin{array}{l}
D^{(2)}_{n,P,2}=A_{n,P,2}^{\top}(D^{(1)}_{\tilde{P}_n^{(1)}(P)}\\
=-\epsilon_n^{(1)}(p)\frac{\bar{q}_n^{(1)}(1-\bar{q}_n^{(1)})(1)}{\bar{g}^2}(A-\bar{g})
\\ -  \left\{ c_P^{(1),-1}P\frac{\bar{q}_n^{(1)}(1-\bar{q}_n^{(1)}(1)}{\bar{g}} \right\}
 \frac{\tilde{g}_n}{\bar{g}^2} (\tilde{q}_n-\bar{q}_n^{(1)}(p) )(A-\bar{g} )\\
+c_P^{(1),-1}\left\{ P\frac{1}{\bar{g}} \bar{q}_n^{(1)}(1-\bar{q}_n^{(1)})(p)(1)\right\}
\frac{\tilde{g}_n}{\bar{g}^3}\epsilon_n^{(1)}(p)\bar{q}_n^{(1)}(1-\bar{q}_n^{(1)})(1)(A-\bar{g}) .
\end{array}
\]
Note that $D^{(2)}_{n,\tilde{P}_n,2}=0$.
Define
\[
\begin{array}{l}
\tilde{c}_{\tilde{P}_n,P}^{(1)}=c_P^{(1),-1}P\frac{\bar{q}_n^{(1)}(1-\bar{q}_n^{(1)})(p)}{\bar{g}}\\
=\left\{P \frac{\tilde{g}_n}{\bar{g}^2} \bar{q}_n^{(1)}(1-\bar{q}_n^{(1)})(p)(1) \right\}^{-1}
P\frac{\bar{q}_n^{(1)}(1-\bar{q}_n^{(1)})(p)(1)}{\bar{g}}.
\end{array}
\]

Thus, we have derived  $D^{(2)}_{n,P}=D^{(2)}_{n,P,1}+D^{(2)}_{n,P,2}$.
\[
\begin{array}{l}
D^{(2)}_{n,P}=
\frac{A}{\bar{g}} \bar{q}_n^{(1)}(1-\bar{q}_n^{(1)})(p)(1)\frac{1}{\bar{q}(1-\bar{q})}(Y-\bar{q})\\
-\tilde{c}_{\tilde{P}_n,P}^{(1)} C_g \frac{ \tilde{g}_n}{\bar{g}} \frac{\bar{q}_n^{(1)}(1-\bar{q}_n^{(1)})(p)}{\bar{q}(1-\bar{q})} (Y-\bar{q})\\
-\epsilon_n^{(1)}(p)\frac{\bar{q}_n^{(1)}(1-\bar{q}_n^{(1)})(p)(1)}{\bar{g}^2}(A-\bar{g})
\\ -  \tilde{c}_{\tilde{P}_n,P}^{(1)}
 \frac{\tilde{g}_n}{\bar{g}^2} (\tilde{q}_n-\bar{q}_n^{(1)}(p) )(1)(A-\bar{g} )\\
+\tilde{c}_{\tilde{P}_n,P}^{(1)}
\frac{\tilde{g}_n}{\bar{g}^3}\epsilon_n^{(1)}(p)\bar{q}_n^{(1)}(1-\bar{q}_n^{(1)})(p)(1)(A-\bar{g}) .
\end{array}
\]
The first two terms represent $D^{(2)}_{n,P,1}$.

This proves Lemma \ref{lemmaexample2b}.

\subsection{Proof of Lemma \ref{lemmaexample2c}.}

For notational convenience, let $p_n^*=\tilde{P}_n^{(2)}(P)$.
As in our general section, we first obtain a decomposition in a pure second order polynomial difference and a third order difference:
\[
\begin{array}{l}
R^{(1)}(\tilde{P}_n^{(1)}(P_n^*),\tilde{P}_n)=
\tilde{P}_n (\bar{q}_n^{(1)}(p_n^*)-\tilde{q}_n)(1)\frac{\bar{g}_n^*-\tilde{g}_n}{\bar{g}_n^*}\\
=
\tilde{P}_n (\bar{q}_n^{(1)}(p_n^*)-\tilde{q}_n)(1)\frac{\bar{g}_n^*-\tilde{g}_n}{\tilde{g}_n}
-\tilde{P}_n (\bar{q}_n^{(1)}(p_n^*)-\tilde{q}_n)(1)\frac{(\bar{g}_n^*-\tilde{g}_n)^2}{\bar{g}_n^*\tilde{g}_n}\\
=R^{(1),p}(\tilde{P}_n^{(1)}(p_n^*),\tilde{P}_n)-\tilde{P}_n (\bar{q}_n^{(1)}(p_n^*)-\tilde{q}_n)(1)\frac{(\bar{g}_n^*-\tilde{g}_n)^2}{\bar{g}_n^*\tilde{g}_n}.
\end{array}
\]
Note that we can represent $\tilde{p}_n=\tilde{p}_n^{(1)}(\tilde{p}_n)$.
Let's denote the second third order term with $R^{(1)}_3(\tilde{P}_n^{(1)}(P_n^*),\tilde{P}_n)$.
We have $\frac{d}{dP_n^*}R^{(1)}(\tilde{P}_n^{(1)}(P_n^*),\tilde{P}_n)(P_n^*-\tilde{P}_n)=-\tilde{P}_n D^{(2)}_{P_n^*}$.
Thus, \[
\frac{d}{dP_n^*}R^{(1),p}(\tilde{P}_n^{(1)}(P_n^*),\tilde{P}_n)(P_n^*-\tilde{P}_n)=-\tilde{P}_n D^{(2)}_{P_n^*}+
\frac{d}{dP_n^*}R^{(1)}_3(\tilde{P}_n^{(1)}(P_n^*),\tilde{P}_n)(P_n^*-\tilde{P}_n).\]
The latter term remains a third order term, and we will denote it with $R^{(1),+,b}(\tilde{P}_n^{(1)}(P_n^*),\tilde{P}_n)$.
Using that $p=(\bar{q},\bar{g})$ and $\frac{d}{dp}f(p)(p-\tilde{p}_n)=\frac{d}{d\bar{q}}f(\bar{q},\bar{g})(\bar{q}-\tilde{q}_n)+\frac{d}{d\bar{g}}f(\bar{q},\bar{g})(\bar{g}-\tilde{g}_n)$, it follows
that
\[
\begin{array}{l}
R^{(1),+,b}(\tilde{P}_n^{(1)}(P_n^*),\tilde{P}_n)=
\tilde{P}_n \frac{d}{dp_n^*}\bar{q}_n^{(1)}(p_n^*)(1)(p_n^*-\tilde{p}_n) \frac{(\bar{g}_n^*-\tilde{g}_n)^2}{\bar{g}_n^*\tilde{g}_n}\\
+\tilde{P}_n (\bar{q}_n^{(1)}(p_n^*)-\tilde{q}_n)(1)\frac{2}{\bar{g}_n^*\tilde{g}_n}(\bar{g}_n^*-\tilde{g}_n)^2-\tilde{P}_n  (\bar{q}_n^{(1)}(p_n^*)-\tilde{q}_n)(1)\frac{1}{\tilde{g}_n(\bar{g}_n^*)^2}(\bar{g}_n^*-\tilde{g}_n)^3.
\end{array}
\]

By our derivative reduction representation for a second order polynomial difference we have
\[
\begin{array}{l}
R^{(1),p}(\tilde{P}_n^{(1)}(P_n^*),\tilde{P}_n)=R^{(1),p,+,a}(\tilde{P}_n^{(1)}(P_n^*),\tilde{P}_n)
+1/2\frac{d}{dP_n^*}R^{(1),p}(\tilde{P}_n^{(1)}(P_n^*),\tilde{P}_n)(P_n^*-\tilde{P}_n),
\end{array}
\]
where $R^{(1),p,+,a}$ is a third order difference we can compute as follows.
We have
\[
\begin{array}{l}
R^{(1),p}(\tilde{P}_n^{(1)}(p_n^*),\tilde{P}_n^{(1)}(\tilde{p}_n))\\
=
1/2 \tilde{P}_n \left\{ \bar{q}_n^{(1)}(p_n^*)-\bar{q}_n^{(1)}(\tilde{p}_n)-\frac{d}{dp_n^*}\bar{q}_n^{(1)}(p_n^*)(p_n^*-\tilde{p}_n)\right\}
\frac{\bar{g}_n^*-\tilde{g}_n}{\tilde{g}_n}\\
+1/2\tilde{P}_n \frac{d}{dp_n^*}\bar{q}_n^{(1)}(p_n^*)(1)(p_n^*-\tilde{p}_n)
\frac{\bar{g}_n^*-\tilde{g}_n}{\tilde{g}_n}\\
+1/2 \tilde{P}_n ( \bar{q}_n^{(1)}(p_n^*)-\bar{q}_n^{(1)}(\tilde{p}_n))
\frac{\bar{g}_n^*-\tilde{g}_n-\frac{d}{dp_n^*}\bar{g}(p_n^*)(p_n^*-\tilde{p}_n)}{\tilde{g}_n}\\
+1/2\tilde{P}_n ( \bar{q}_n^{(1)}(p_n^*)-\bar{q}_n^{(1)}(\tilde{p}_n))(1)\frac{\frac{d}{dp_n^*}\bar{g}(p_n^*)(p_n^*-\tilde{p}_n)}{\tilde{g}_n}\\
=1/2 \tilde{P}_n \left\{ \bar{q}_n^{(1)}(p_n^*)-\bar{q}_n^{(1)}(\tilde{p}_n)-\frac{d}{dp_n^*}\bar{q}_n^{(1)}(p_n^*)(p_n^*-\tilde{p}_n)\right\}(1)
\frac{\bar{g}_n^*-\tilde{g}_n}{\tilde{g}_n}\\
+
1/2 \tilde{P}_n ( \bar{q}_n^{(1)}(p_n^*)-\bar{q}_n^{(1)}(\tilde{p}_n))
\frac{\bar{g}_n^*-\tilde{g}_n-\frac{d}{dp_n^*}\bar{g}(p_n^*)(p_n^*-\tilde{p}_n)}{\tilde{g}_n}\\
+1/2 \frac{d}{dP_n^*}R^{(1),p}(\tilde{P}_n^{(1)}(P_n^*),\tilde{P}_n^{(1)}(\tilde{P}_n))(P_n^*-\tilde{P}_n)\\
={R}^{(1),p,+,a}(\tilde{P}_n^{(1)}(P_n^*),\tilde{P}_n^{(1)}(\tilde{P}_n))+
1/2  \frac{d}{dP_n^*}R^{(1),p}(\tilde{P}_n^{(1)}(P_n^*),\tilde{P}_n^{(1)}(\tilde{P}_n))(P_n^*-\tilde{P}_n).
\end{array}
\]
$R^{(1),p,+,a}$ is a third order difference. 
Define
\[\begin{array}{l}
R_{2,\bar{q}_n^{(1)},p_n^*}(p_n^*,\tilde{p}_n)\equiv 
\bar{q}_n^{(1)}(p_n^*)-\bar{q}_n^{(1)}(\tilde{p}_n)-\frac{d}{dq_n^*}\bar{q}_n^{(1)}(q_n^*,g_n^*)(q_n^*-\tilde{q}_n)\\
\hfill -
\frac{d}{dg_n^*}\bar{q}_n^{(1)}(q_n^*,g_n^*)(g_n^*-\tilde{g}_n).
\end{array}
\]
Since $p=(q,g)$, and $\frac{d}{dp}\bar{g}(p)(p-\tilde{p}_n)=(\bar{g}-\bar{g}(\tilde{p}_n))$,
\[
\begin{array}{l}
{R}^{(1),p,+,a}(\tilde{P}_n^{(1)}(P_n^*),\tilde{P}_n^{(1)}(\tilde{P}_n))
=
1/2 \tilde{P}_n R_{2,\bar{q}_n^{(1)},p_n^*}(p_n^*,\tilde{p}_n)
\frac{\bar{g}_n^*-\tilde{g}_n}{\tilde{g}_n} .
\end{array}
\]

We conclude that
\[
\begin{array}{l}
R^{(1)}(\tilde{P}_n^{(1)}(P_n^*),\tilde{P}_n)=
{R}^{(1),p,+,a}(\tilde{P}_n^{(1)}(P_n^*),\tilde{P}_n^{(1)}(\tilde{P}_n))\\
-R^{(1)}_3(\tilde{P}_n^{(1)}(P_n^*),\tilde{P}_n)+
1/2\frac{d}{dP_n^*}R^{(1)}_3(\tilde{P}_n^{(1)}(P_n^*),\tilde{P}_n)(P_n^*-\tilde{P}_n)\\
-1/2\tilde{P}_n D^{(2)}_{P_n^*}\\
=R^{(1),p,+,a}-1/2 \tilde{P}_n D^{(2)}+1/2
\tilde{P}_n \frac{d}{dp_n^*}\bar{q}_n^{(1)}(p_n^*)(1)(p_n^*-\tilde{p}_n) \frac{(\bar{g}_n^*-\tilde{g}_n)^2}{\bar{g}_n^*\tilde{g}_n}\\
+1/2 \tilde{P}_n (\bar{q}_n^{(1)}(p_n^*)-\tilde{q}_n)(1)\frac{2}{\bar{g}_n^*\tilde{g}_n}(\bar{g}_n^*-\tilde{g}_n)^2-\tilde{P}_n  (\bar{q}_n^{(1)}(p_n^*)-\tilde{q}_n)(1)\frac{1}{\tilde{g}_n(\bar{g}_n^*)^2}(\bar{g}_n^*-\tilde{g}_n)^3 \\
-\tilde{P}_n (\bar{q}_n^{(1)}(p_n^*)-\tilde{q}_n)(1)\frac{(\bar{g}_n^*-\tilde{g}_n)^2}{\bar{g}_n^*\tilde{g}_n}
\\
=R^{(1),p,+,a}-1/2 \tilde{P}_n D^{(2)}+1/2
\tilde{P}_n \frac{d}{dp_n^*}\bar{q}_n^{(1)}(p_n^*)(1)(p_n^*-\tilde{p}_n) \frac{(\bar{g}_n^*-\tilde{g}_n)^2}{\bar{g}_n^*\tilde{g}_n}\\
-\tilde{P}_n  (\bar{q}_n^{(1)}(p_n^*)-\tilde{q}_n)(1)\frac{1}{\tilde{g}_n(\bar{g}_n^*)^2}(\bar{g}_n^*-\tilde{g}_n)^3 .

\end{array}
\]


This completes the proof of the lemma. 

\subsection{Verification that the directional derivative of $\Psi_n^{(1)}$ agrees with the second order canonical gradient.}
Suppose that $\tilde{P}_n D^{(2)}_{\tilde{P}_n^{(2)}(P)}=0$ and $\tilde{P}_n D^{(1)}_{\tilde{P}_n^{(1)}(P)}=0$ for all $P$. We have $\tilde{P}_n D^{(1)}_{\tilde{P}_n^{(1)}(P)}=0$ for our TMLE-update. By using $\tilde{P}_n^{(2)}(P)$ as a universal least favorable path update, iterative TMLE update or a one-step local least favorable path with $\epsilon_n^{(2)}(P)$ defined by exactly solving the equation $\tilde{P}_n D^{(2)}_{n,\tilde{P}_n^{(2)}(P,\epsilon)}=0$. 
For our formula for $D^{(2)}_{n,P}$ we  have
\[
\begin{array}{l}
\tilde{P}_n D^{(2)}_{n,P}=
P_n \frac{\tilde{g}_n}{\bar{g}} \frac{\bar{q}_n^{(1)}(p)(1-\bar{q}_n^{(1)}(p))}{\bar{q}(1-\bar{q})}\left\{1-\tilde{c}_{\tilde{P}_n,P}^{(1)}\frac{\tilde{g}_n}{\bar{g}}\right\}(\tilde{q}_n-\bar{q})\\
-\epsilon_n^{(1)}(p)P_n \frac{\bar{q}_n^{(1)}(1-\bar{q}_n^{(1)})(p)}{\bar{g}^2}\left\{1-\tilde{c}_{\tilde{P}_n,P}^{(1)}\frac{\tilde{g}_n}{\bar{g}}\right\}(\tilde{g}_n-\bar{g})\\
-  \tilde{c}_{\tilde{P}_n,P}^{(1)}P_n
 \frac{\tilde{g}_n}{\bar{g}^2} (\tilde{q}_n-\bar{q}_n^{(1)}(p)(1) )(\tilde{g}_n-\bar{g} ).
 \end{array}
 \]
We also have \[
\begin{array}{l}
\Psi_n^{(1)}(P)-\Psi_n^{(1)}(\tilde{P}_n)=\Psi(\tilde{P}_n^{(1)}(P))-\Psi(\tilde{P}_n)=
-\tilde{P}_n D^{(1)}_{\tilde{P}_n^{(1)}(P)}+R^{(1)}(\tilde{P}_n^{(1)}(P),\tilde{P}_n)\\
=R^{(1)}(\tilde{P}_n^{(1)}(P),\tilde{P}_n)\\
=P_n (\bar{q}_n^{(1)}(p)-\tilde{q}_n)\frac{\bar{g}-\tilde{g}_n}{\bar{g}}. \\
\end{array}
\]

We want to prove that indeed with our choice $D^{(2)}_{n,P}$ it follows that$(P-\tilde{P}_n )D^{(2)}_{n,P}=\frac{d}{dP}\Psi_n^{(1)}(P)(P-\tilde{P}_n)+R_{n,3}(P,\tilde{P}_n)$, where $R_{n,3}(P,\tilde{P}_n)$ is a third order difference. 
 For that  we need to show \[(P-\tilde{P}_n)D^{(2)}_{n,P}=P_n  \frac{d}{dp}\bar{q}_n^{(1)}(p)(p-\tilde{p}_n) \frac{\bar{g}-\tilde{g}_n}{\bar{g}}+P_n (\bar{q}_n^{(1)}(p)-\tilde{q}_n)\frac{\bar{g}-\tilde{g}_n}{\bar{g}}+R_{n,3}(P,\tilde{P}_n).\]



\begin{lemma}
We have
\[
\begin{array}{l}
\frac{d}{dp}\bar{q}_n^{(1)}(p)(p-\tilde{p}_n)\\
=
\frac{\bar{q}_n^{(1)}(p)(1-\bar{q}_n^{(1)}(p))}{\bar{q}(1-\bar{q})}(\bar{q}-\tilde{q}_n)\\
+\bar{q}_n^{(1)}(p)(1-\bar{q}_n^{(1)}(p))
(\bar{g})^{-1}\frac{d}{dp}\epsilon_n(p)(p-\tilde{p}_n)\\
-\bar{q}_n^{(1)}(p)(1-\bar{q}_n^{(1)}(p))\epsilon_n(p)/\bar{g}^2(\bar{g}-\tilde{g}_n)\\
\end{array}
\]
Plugging this term  into $P_n d/dp\bar{q}_n^{(1)}(p-\tilde{p}_n) (\bar{g}-\tilde{g}_n)/\bar{g}$ yields
\[\begin{array}{l}\frac{d}{dP}\Psi_n^{(1)}(P)(P-\tilde{P}_n)= \\
P_n\frac{\bar{q}_n^{(1)}(p)(1-\bar{q}_n^{(1)}(p))}{\bar{q}(1-\bar{q})}(\bar{q}-\tilde{q}_n)(\bar{g}-\tilde{g}_n)/\bar{g}\\-P_n \bar{q}_n^{(1)}(p)(1-\bar{q}_n^{(1)}(p)) \epsilon_n(p)/\bar{g}^2(\bar{g}-\tilde{g}_n)^2/\bar{g}\\+P_n\bar{q}_n^{(1)}(p)(1-\bar{q}_n^{(1)}(p))(\bar{g})^{-1}\frac{d}{dp}\epsilon_n(p)(p-\tilde{p}_n)\frac{\bar{g}-\tilde{g}_n}{\bar{g}}\\
+R^{(1)}(\tilde{P}_n^{(1)}(P),\tilde{P}_n)
.\end{array}\]
We have\[\begin{array}{l}P_n \bar{q}_n^{(1)}(p)(1-\bar{q}_n^{(1)}(p))(\bar{g})^{-1}\frac{d}{dp}\epsilon_n(p)(p-\tilde{p}_n)\frac{\bar{g}-\tilde{g}_n}{\bar{g}}\\=-(\tilde{c}_{\tilde{P}_n,P}^{(1)}-1) {P}_n \frac{\bar{q}_n^{(1)}(p)(1-\bar{q}_n^{(1)}(p))}{\bar{q}(1-\bar{q})}(\bar{q}-\tilde{q}_n),\end{array}\] up till third order.

So we obtain\[\begin{array}{l}\frac{d}{dP}\Psi_n^{(1)}(P)(P-\tilde{P}_n)=\\
P_n\frac{\bar{q}_n^{(1)}(p)(1-\bar{q}_n^{(1)}(p))}{\bar{q}(1-\bar{q})}(\bar{q}-\tilde{q}_n)(\bar{g}-\tilde{g}_n)/\bar{g}\\-P_n \bar{q}_n^{(1)}(p)(1-\bar{q}_n^{(1)}(p)) \epsilon_n(p)/\bar{g}^2(\bar{g}-\tilde{g}_n)^2/\bar{g}\\-(\tilde{c}_{\tilde{P}_n,P}^{(1)}-1) {P}_n \frac{\bar{q}_n^{(1)}(p)(1-\bar{q}_n^{(1)}(p))}{\bar{q}(1-\bar{q})}(\bar{q}-\tilde{q}_n)\\+R^{(1)}(\tilde{P}_n^{(1)}(P),\tilde{P}_n)
,\end{array}\] up till third order term.
\end{lemma}
{\bf Proof:}
So let's compute
\[
\frac{d}{dp}\bar{q}_n^{(1)}(p)(p-\tilde{p}_n)=
\frac{d}{d\bar{q}}\bar{q}_n^{(1)}(p)(\bar{q}-\tilde{q}_n)+\frac{d}{d\bar{g}}\bar{q}_n^{(1)}(p)(\bar{g}-\tilde{g}_n).\]
We have (at $A=1$)
\[
\bar{q}_n^{(1)}(p)=\frac{1}{1+\exp\left(-\log(\bar{q}/(1-\bar{q})-\epsilon_n^{(1)}(p)1/\bar{g}\right)}.\]
So
\[
\begin{array}{l}
\frac{d}{dp}\bar{q}_n^{(1)}(p)(p-\tilde{p}_n)=
-\{\bar{q}_n^{(1)}(p)\}^2\frac{d}{dp}\left\{ \frac{1-\bar{q}}{\bar{q}}\exp(-\epsilon_n(p)/\bar{g})\right\}(p-\tilde{p}_n)\\
=-\{\bar{q}_n^{(1)}(p)\}^2\exp(-\epsilon_n(p)/\bar{g})\frac{d}{d\bar{q}}(1-\bar{q})/\bar{q})(\bar{q}-\tilde{q}_n)\\
+\{\bar{q}_n^{(1)}(p)\}^2(1-\bar{q})/\bar{q} \exp(-\epsilon_n(p)/\bar{g})\frac{d}{dp}\{\epsilon_n(p)/\bar{g}\}(p-\tilde{p}_n)\\
=-\{\bar{q}_n^{(1)}(p)\}^2\exp(-\epsilon_n(p)/\bar{g})\left\{-1/\bar{q}(\bar{q}-\tilde{q}_n)-(1-\bar{q})/\bar{q}^2(\bar{q}-\tilde{q}_n)\right\}\\
+\{\bar{q}_n^{(1)}(p)\}^2(1-\bar{q})/\bar{q} \exp(-\epsilon_n(p)/\bar{g})
(\bar{g})^{-1}\frac{d}{dp}\{\epsilon_n(p)(p-\tilde{p}_n)\\
-\{\bar{q}_n^{(1)}(p)\}^2(1-\bar{q})/\bar{q} \exp(-\epsilon_n(p)/\bar{g})\epsilon_n(p)/\bar{g}^2(\bar{g}-\tilde{g}_n)\\
=+\{\bar{q}_n^{(1)}(p)\}^2\exp(-\epsilon_n(p)/\bar{g})(1/\bar{q})(\bar{q}-\tilde{q}_n)\\
+\{\bar{q}_n^{(1)}(p)\}^2\exp(-\epsilon_n(p)/\bar{g})(1-\bar{q})/\bar{q}^2(\bar{q}-\tilde{q}_n)\\
+\{\bar{q}_n^{(1)}(p)\}^2(1-\bar{q})/\bar{q} \exp(-\epsilon_n(p)/\bar{g})
(\bar{g})^{-1}\frac{d}{dp}\epsilon_n(p)(p-\tilde{p}_n)\\
-\{\bar{q}_n^{(1)}(p)\}^2(1-\bar{q})/\bar{q} \exp(-\epsilon_n(p)/\bar{g})\epsilon_n(p)/\bar{g}^2(\bar{g}-\tilde{g}_n).
\end{array}
\]
Notice that
\[(1-\bar{q})/\bar{q} \exp(-\epsilon_n(p)/\bar{g})=
\frac{1-\bar{q}_n^{(1)}(p)}{\bar{q}_n^{(1)}(p)}.\]
Substitution gives:
\[
\begin{array}{l} \frac{d}{dp}\bar{q}_n^{(1)}(p)(p-\tilde{p}_n)=
\{\bar{q}_n^{(1)}(p)\}^2(1-\bar{q})^{-1} \left\{ \exp(-\epsilon_n(p)/\bar{g})((1-\bar{q})/\bar{q})\right\}(\bar{q}-\tilde{q}_n)\\
+\{\bar{q}_n^{(1)}(p)\}^2(\bar{q})^{-1} \left\{ \exp(-\epsilon_n(p)/\bar{g})(1-\bar{q})/\bar{q} \right\} (\bar{q}-\tilde{q}_n)\\
+\{\bar{q}_n^{(1)}(p)\}^2\left\{ (1-\bar{q})/\bar{q} \exp(-\epsilon_n(p)/\bar{g})\right\}
(\bar{g})^{-1}\frac{d}{dp}\epsilon_n(p)(p-\tilde{p}_n)\\
-\{\bar{q}_n^{(1)}(p)\}^2\left\{ (1-\bar{q})/\bar{q} \exp(-\epsilon_n(p)/\bar{g})\right\} \epsilon_n(p)/\bar{g}^2(\bar{g}-\tilde{g}_n)\\
=
\bar{q}_n^{(1)}(p)(1-\bar{q}_n^{(1)}(p))(1-\bar{q})^{-1}(\bar{q}-\tilde{q}_n)\\
+\bar{q}_n^{(1)}(p)(1-\bar{q}_n^{(1)}(p)) (\bar{q})^{-1}(\bar{q}-\tilde{q}_n)\\
+\bar{q}_n^{(1)}(p)(1-\bar{q}_n^{(1)}(p))
(\bar{g})^{-1}\frac{d}{dp}\epsilon_n(p)(p-\tilde{p}_n)\\
-\bar{q}_n^{(1)}(p)(1-\bar{q}_n^{(1)}(p))\epsilon_n(p)/\bar{g}^2(\bar{g}-\tilde{g}_n)\\
=
\frac{\bar{q}_n^{(1)}(p)(1-\bar{q}_n^{(1)}(p))}{\bar{q}(1-\bar{q})}(\bar{q}-\tilde{q}_n)\\
+\bar{q}_n^{(1)}(p)(1-\bar{q}_n^{(1)}(p))
(\bar{g})^{-1}\frac{d}{dp}\epsilon_n(p)(p-\tilde{p}_n)\\
-\bar{q}_n^{(1)}(p)(1-\bar{q}_n^{(1)}(p))\epsilon_n(p)/\bar{g}^2(\bar{g}-\tilde{g}_n)\\
\end{array}
\]
This completes the proof of the first statement. 

So it remains to derive
\[
\begin{array}{l}
 +P_n \bar{q}_n^{(1)}(p)(1-\bar{q}_n^{(1)}(p))
\frac{1}{\bar{g}}\frac{d}{dp}\epsilon_n(p)(p-\tilde{p}_n)\frac{\bar{g}-\tilde{g}_n}{\bar{g}}
\end{array}
\]
Note
\[
\begin{array}{l}
P_n \bar{q}_n^{(1)}(p)(1-\bar{q}_n^{(1)}(p))
\frac{1}{\bar{g}}\frac{d}{dp}\epsilon_n(p)(p-\tilde{p}_n)\frac{\bar{g}-\tilde{g}_n}{\bar{g}}\\
=
\left\{ P_n \bar{q}_n^{(1)}(p)(1-\bar{q}_n^{(1)}(p))\frac{\bar{g}-\tilde{g}_n}{\bar{g}^2}\right\}
\frac{d}{dp}\epsilon_n(p)(p-\tilde{p}_n)\\
=\left\{ P_n \bar{q}_n^{(1)}(p)(1-\bar{q}_n^{(1)}(p))\frac{\bar{g}-\tilde{g}_n}{\bar{g}^2}\right\}
c_P^{(1),-1} \frac{d}{dp}U(\epsilon_n,p)(p-\tilde{p}_n).
\end{array}
\]
Notice that numerator of factor in front of $d/dpU$ equals
$P_n \bar{q}_n^{(1)}(1-\bar{q}_n^{(1)})/\bar{g}-c_P^{(1)}$, so that factor becomes
\[
(\tilde{c}_{\tilde{P}_n,P}^{(1)}-1).\]
Thus, we have
\[
\begin{array}{l}
P_n \bar{q}_n^{(1)}(p)(1-\bar{q}_n^{(1)}(p))
\frac{1}{\bar{g}}\frac{d}{dp}\epsilon_n(p)(p-\tilde{p}_n)\frac{\bar{g}-\tilde{g}_n}{\bar{g}}\\
=(\tilde{c}_{\tilde{P}_n,P}^{(1)}-1)\frac{d}{dp}U(\epsilon_n,p)(p-\tilde{p}_n).
\end{array}
\]
Below we show that
\[
\begin{array}{l}
\frac{d}{dp}U(\epsilon_n,p)(p-\tilde{p}_n)\\
=-P_n \frac{\tilde{g}_n}{\bar{g}^2}(\bar{g}-\tilde{g}_n)  (\tilde{q}_n-\bar{q}_n^{(1)}(p))
-\tilde{P}_n C_g\frac{d}{dp}\bar{q}^{(1)}(p,\epsilon_n)(p-\tilde{p}_n).
\end{array}
\]
Thus,
\[
\begin{array}{l}
P_n \bar{q}_n^{(1)}(p)(1-\bar{q}_n^{(1)}(p))
\frac{1}{\bar{g}}\frac{d}{dp}\epsilon_n(p)(p-\tilde{p}_n)\frac{\bar{g}-\tilde{g}_n}{\bar{g}}\\
=-(\tilde{c}_{\tilde{P}_n,P}^{(1)}-1)P_n \frac{\tilde{g}_n}{\bar{g}^2}(\bar{g}-\tilde{g}_n)  (\tilde{q}_n-\bar{q}_n^{(1)}(p))\\
-(\tilde{c}_{\tilde{P}_n,P}^{(1)}-1) \tilde{P}_n C_g\frac{d}{dp}\bar{q}_n^{(1)}(p)(p-\tilde{p}_n)\\
=-(\tilde{c}_{\tilde{P}_n,P}^{(1)}-1)P_n \frac{\tilde{g}_n}{\bar{g}^2}(\bar{g}-\tilde{g}_n)  (\tilde{q}_n-\bar{q}_n^{(1)}(p))\\
-(\tilde{c}_{\tilde{P}_n,P}^{(1)}-1) {P}_n \tilde{g}_n/\bar{g} \frac{d}{dp}\bar{q}^{(1)}(p,\epsilon_n)(p-\tilde{p}_n)
.
\end{array}
\]
The first term is a third order term already. 
So we can focus on second term. The second term we can write $\bar{g}_n/\bar{g}-1+1$ giving another third order term and then
\[
-(\tilde{c}_{\tilde{P}_n,P}^{(1)}-1) {P}_n \frac{d}{dp}\bar{q}^{(1)}(p,\epsilon_n)(p-\tilde{p}_n).\]
We now need an expression for $\frac{d}{dp}\bar{q}^{(1)}(p,\epsilon)(p-\tilde{p}_n)$.
This was shown above:
\[
\begin{array}{l}
\frac{d}{dp}\bar{q}^{(1)}(p,\epsilon_n)(p-\tilde{p}_n)\\
=\frac{\bar{q}_n^{(1)}(p)(1-\bar{q}_n^{(1)}(p))}{\bar{q}(1-\bar{q})}(\bar{q}-\tilde{q}_n)\\
-\bar{q}_n^{(1)}(p)(1-\bar{q}_n^{(1)}(p))\epsilon_n(p)/\bar{g}^2(\bar{g}-\tilde{g}_n).
\end{array}
\]
The second term of this expression is already second order so gives a third order.
So we are left with
\[
\begin{array}{l}
P_n \bar{q}_n^{(1)}(p)(1-\bar{q}_n^{(1)}(p))
\frac{1}{\bar{g}}\frac{d}{dp}\epsilon_n(p)(p-\tilde{p}_n)\frac{\bar{g}-\tilde{g}_n}{\bar{g}}\\
=-(\tilde{c}_{\tilde{P}_n,P}^{(1)}-1) {P}_n \frac{\bar{q}_n^{(1)}(p)(1-\bar{q}_n^{(1)}(p))}{\bar{q}(1-\bar{q})}(\bar{q}-\tilde{q}_n),
\end{array}
\]
up till third order term. 
This proves the second statement. The final statement follows from substitution. This completes the proof of the lemma.
$\Box$

This allows us to prove the desired result.
\begin{lemma}
We have that
\[
\begin{array}{l}
-\tilde{P}_n D^{(2)}_{n,P}=
P_n \frac{d}{dp}\bar{q}_n^{(1)}(p)(p-\tilde{p}_n) \frac{\bar{g}-\tilde{g}_n}{\bar{g}}+R^{(1)}(\tilde{P}_n^{(1)}(P),\tilde{P}_n)\\
+R_{n,3}(P,\tilde{P}_n),\end{array}
\]
where $R_{n,3}(P,\tilde{P}_n)$ is a third order difference.
This also proves  that 
\[
-\tilde{P}_n D^{(2)}_{n,P}=\frac{d}{dP}\Psi_n^{(1)}(P)(P-\tilde{P}_n)\]
up till a third order remainder.
\end{lemma}
{\bf Proof:}
Plugging in in $P_n d/dp\bar{q}_n^{(1)}(p-\tilde{p}_n) (\bar{g}-\tilde{g}_n)/\bar{g}$ gives
\[
\begin{array}{l}
P_n\frac{\bar{q}_n^{(1)}(p)(1-\bar{q}_n^{(1)}(p))}{\bar{q}(1-\bar{q})}(\bar{q}-\tilde{q}_n)(\bar{g}-\tilde{g}_n)/\bar{g}\\
-P_n \bar{q}_n^{(1)}(p)(1-\bar{q}_n^{(1)}(p)) \epsilon_n(p)/\bar{g}^2(\bar{g}-\tilde{g}_n)^2/\bar{g}\\
+P_n
\bar{q}_n^{(1)}(p)(1-\bar{q}_n^{(1)}(p))
(\bar{g})^{-1}\frac{d}{dp}\epsilon_n(p)(p-\tilde{p}_n)\frac{\bar{g}-\tilde{g}_n}{\bar{g}}.
\end{array}
\]
Compare with first and second term of $\tilde{P}_n D^{(2)}_{n,P}$:
\[
\begin{array}{l}
P_n \frac{\tilde{g}_n}{\bar{g}} \frac{\bar{q}_n^{(1)}(p)(1-\bar{q}_n^{(1)}(p))}{\bar{q}(1-\bar{q})}\left\{1-\tilde{c}_{\tilde{P}_n,P}^{(1)}\frac{\tilde{g}_n}{\bar{g}}\right\}(\tilde{q}_n-\bar{q})\\
-\epsilon_n^{(1)}(p)P_n \frac{\bar{q}_n^{(1)}(1-\bar{q}_n^{(1)})(p)}{\bar{g}^2}\left\{1-\tilde{c}_{\tilde{P}_n,P}^{(1)}\frac{\tilde{g}_n}{\bar{g}}\right\}(\tilde{g}_n-\bar{g}).
 \end{array}
 \]
 Let's ignore both second terms since they are third order. 
Adding up the two respective terms yields
\[
 \begin{array}{l}
 P_n \frac{\tilde{g}_n}{\bar{g}} \frac{\bar{q}_n^{(1)}(p)(1-\bar{q}_n^{(1)}(p))}{\bar{q}(1-\bar{q})}\left\{1-\tilde{c}_{\tilde{P}_n,P}^{(1)}\frac{\tilde{g}_n}{\bar{g}}\right\}(\tilde{q}_n-\bar{q})\\
 +P_n\frac{\bar{q}_n^{(1)}(p)(1-\bar{q}_n^{(1)}(p))}{\bar{q}(1-\bar{q})}(\bar{q}-\tilde{q}_n)(\bar{g}-\tilde{g}_n)/\bar{g}\\
 =P_n \left( \frac{\tilde{g}_n}{\bar{g}}-1\right) \frac{\bar{q}_n^{(1)}(p)(1-\bar{q}_n^{(1)}(p))}{\bar{q}(1-\bar{q})}\left\{1-\tilde{c}_{\tilde{P}_n,P}^{(1)}\frac{\tilde{g}_n}{\bar{g}}\right\}(\tilde{q}_n-\bar{q})\\
 +P_n  \frac{\bar{q}_n^{(1)}(p)(1-\bar{q}_n^{(1)}(p))}{\bar{q}(1-\bar{q})}\left\{1-\tilde{c}_{\tilde{P}_n,P}^{(1)}\frac{\tilde{g}_n}{\bar{g}}\right\}(\tilde{q}_n-\bar{q}) \\
 +P_n\frac{\bar{q}_n^{(1)}(p)(1-\bar{q}_n^{(1)}(p))}{\bar{q}(1-\bar{q})}(\bar{q}-\tilde{q}_n)(\bar{g}-\tilde{g}_n)/\bar{g}.
  \end{array}
 \] 
 Ignore first term since third order.
 \[
 \begin{array}{l}
 +P_n  \frac{\bar{q}_n^{(1)}(p)(1-\bar{q}_n^{(1)}(p))}{\bar{q}(1-\bar{q})}\frac{\tilde{g}_n}{\bar{g}}(1-c_{\tilde{P}_n,P}^{(1)})(\tilde{q}_n-\bar{q})\\
 +P_n  \frac{\bar{q}_n^{(1)}(p)(1-\bar{q}_n^{(1)}(p))}{\bar{q}(1-\bar{q})}\frac{\bar{g}-\tilde{g}_n}{\bar{g}}(\tilde{q}_n-\bar{q}) \\
 +P_n\frac{\bar{q}_n^{(1)}(p)(1-\bar{q}_n^{(1)}(p))}{\bar{q}(1-\bar{q})}(\bar{q}-\tilde{q}_n)(\bar{g}-\tilde{g}_n)/\bar{g}\\
 =P_n  \frac{\bar{q}_n^{(1)}(p)(1-\bar{q}_n^{(1)}(p))}{\bar{q}(1-\bar{q})}\frac{\tilde{g}_n}{\bar{g}}(1-c_{\tilde{P}_n,P}^{(1)})(\tilde{q}_n-\bar{q}).
 \end{array}
 \] 
 So  \[
 \begin{array}{l}
 (\tilde{P}_n-P)D^{(2)}_{n,P}+P_n d/dp\bar{q}_n^{(1)}(p-\tilde{p}_n) (\bar{g}-\tilde{g}_n)/\bar{g}\\
  = P_n  \frac{\bar{q}_n^{(1)}(p)(1-\bar{q}_n^{(1)}(p))}{\bar{q}(1-\bar{q})}\frac{\tilde{g}_n}{\bar{g}}(1-c_{\tilde{P}_n,P}^{(1)})(\tilde{q}_n-\bar{q})\\
 -  \tilde{c}_{\tilde{P}_n,P}^{(1)}P_n
 \frac{\tilde{g}_n}{\bar{g}^2} (\tilde{q}_n-\bar{q}_n^{(1)}(p)(1) )(\tilde{g}_n-\bar{g} )\\
  +P_n \bar{q}_n^{(1)}(p)(1-\bar{q}_n^{(1)}(p))
(\bar{g})^{-1}\frac{d}{dp}\epsilon_n(p)(p-\tilde{p}_n)\frac{\bar{g}-\tilde{g}_n}{\bar{g}},
\end{array}
\]
We already obtained
\[
\begin{array}{l}
P_n \bar{q}_n^{(1)}(p)(1-\bar{q}_n^{(1)}(p))
(\bar{g})^{-1}\frac{d}{dp}\epsilon_n(p)(p-\tilde{p}_n)\frac{\bar{g}-\tilde{g}_n}{\bar{g}}\\
=-(\tilde{c}_{\tilde{P}_n,P}^{(1)}-1) {P}_n \frac{\bar{q}_n^{(1)}(p)(1-\bar{q}_n^{(1)}(p))}{\bar{q}(1-\bar{q})}(\bar{q}-\tilde{q}_n),\end{array}
\]
up till third order.
So we have shown
 \[
 \begin{array}{l}
 (\tilde{P}_n-P)D^{(2)}_{n,P}+P_n d/dp\bar{q}_n^{(1)}(p-\tilde{p}_n) (\bar{g}-\tilde{g}_n)/\bar{g}\\
  = P_n  \frac{\bar{q}_n^{(1)}(p)(1-\bar{q}_n^{(1)}(p))}{\bar{q}(1-\bar{q})}\frac{\tilde{g}_n}{\bar{g}}(1-c_{\tilde{P}_n,P}^{(1)})(\tilde{q}_n-\bar{q})\\
 -  \tilde{c}_{\tilde{P}_n,P}^{(1)}P_n
 \frac{\tilde{g}_n}{\bar{g}^2} (\tilde{q}_n-\bar{q}_n^{(1)}(p)(1) )(\tilde{g}_n-\bar{g} )\\
  -(\tilde{c}_{\tilde{P}_n,P}^{(1)}-1) {P}_n \frac{\bar{q}_n^{(1)}(p)(1-\bar{q}_n^{(1)}(p))}{\bar{q}(1-\bar{q})}(\bar{q}-\tilde{q}_n)  \\
  = -  \tilde{c}_{\tilde{P}_n,P}^{(1)}P_n
 \frac{\tilde{g}_n}{\bar{g}^2} (\tilde{q}_n-\bar{q}_n^{(1)}(p)(1) )(\tilde{g}_n-\bar{g} )  \\
\end{array}
\]
up till third order term.
Finally, for this last term we have:
\[
\begin{array}{l}
-  \tilde{c}_{\tilde{P}_n,P}^{(1)}P_n
 \frac{\tilde{g}_n}{\bar{g}^2} (\tilde{q}_n-\bar{q}_n^{(1)}(p)(1) )(\tilde{g}_n-\bar{g} ) \\
 =- (\tilde{c}_{\tilde{P}_n,P}^{(1)}-1)P_n
 \frac{\tilde{g}_n}{\bar{g}^2} (\tilde{q}_n-\bar{q}_n^{(1)}(p)(1) )(\tilde{g}_n-\bar{g} ) \\
 -P_n
 \frac{\tilde{g}_n}{\bar{g}^2} (\tilde{q}_n-\bar{q}_n^{(1)}(p)(1) )(\tilde{g}_n-\bar{g} ) \\
 \approx   -P_n
 \frac{\tilde{g}_n}{\bar{g}} (\tilde{q}_n-\bar{q}_n^{(1)}(p)(1) )\frac{\tilde{g}_n-\bar{g} }{\bar{g}} \\
 = -P_n
\left( \frac{\tilde{g}_n}{\bar{g}} -1\right)(\tilde{q}_n-\bar{q}_n^{(1)}(p)(1) )\frac{\tilde{g}_n-\bar{g} }{\bar{g}}  \\
-P_n (\tilde{q}_n-\bar{q}_n^{(1)}(p)(1) )\frac{\tilde{g}_n-\bar{g} }{\bar{g}}\\
\approx  
 -P_n (\tilde{q}_n-\bar{q}_n^{(1)}(p)(1) )\frac{\tilde{g}_n-\bar{g} }{\bar{g}}\\
 =-R^{(1)}(\tilde{P}_n^{(1)}(P),\tilde{P}_n).
 \end{array}
 \]
 This proves our desired result. 
  $\Box$


\section{Targeting the HAL-MLE in the HAL-regularized second order TMLE of the treatment specific mean, in order to guarantee solving the desired empirical efficient score equations.}\label{AppendixI}
Suppose we want to target $\tilde{P}_n$ to solve $(\tilde{P}_n-P_n) D^{(1)}_{n,P}=0$ and
$(\tilde{P}_n-P_n) D^{(2)}_{n,P}=0$ at a given initial estimator $P=P^0$. By doing this, the second order TMLE that arranges that $P_n D^{(2)}_{n,\tilde{P}_n^{(2)}(P^0)}=0$ and $\tilde{P}_n D^{(1)}_{n,\tilde{P}_n^{(1)}\tilde{P}_n^{(2)}(P^0)}=0$, will then require less undersmoothing of the HAL-MLE to arrange that also $\tilde{P}_n D^{(2)}_{n,\tilde{P}_n^{(2)}(P^0)}\approx 0$ and $P_n D^{(1)}_{n,\tilde{P}_n^{(1)}\tilde{P}_n^{(2)}(P^0)}\approx 0 $ at a desired precision. One could also iterate this targeting of the HAL-MLE as in the iterative second order TMLE, and thereby obtain a second order TMLE $P_n^*$ that solves all the desired equations $P_n D^{(1)}_{P_n^*}=\tilde{P}_n D^{(1)}_{P_n^*}=P_n D^{(2)}_{n,P_n^*}=\tilde{P}_nD^{(2)}_{n,P_n^*}=0$.

{\bf Targeting $\tilde{P}_n$ to solve $(\tilde{P}_n-P_n) D^{(1)}_{n,P}=0$ at a given $P=P^0$:}
Let $(\bar{g},\bar{q})$ represent the initial $P^0$.
We want that $\tilde{P}_n \frac{A}{\bar{g}}(Y-\bar{q})=P_n \frac{A}{\bar{g}}(Y-\bar{q})$.
This can be written as:
\[
\begin{array}{l}
P_n \frac{A}{\bar{g}}(Y-\bar{q})-\tilde{P}_n \frac{A}{\bar{g}}(Y-\bar{q})\\
=P_n \frac{A}{\bar{g}}(Y-\bar{q})-P_n \frac{\tilde{g}_n}{\bar{g}}(\tilde{q}_n-\bar{q})\\
=P_n(-A/\bar{g} \bar{q}+\frac{\tilde{g}_n}{\bar{g}} \bar{q})+
P_n (A/\bar{g} Y-\tilde{g}_n/\bar{g} \tilde{q}_n)\\
=-P_n \frac{\bar{q}}{\bar{g}}(A-\tilde{g}_n)+P_n A/\bar{g}(Y-\tilde{q}_n)+P_n \tilde{q}_n/\bar{g}(A-\tilde{g}_n).
\end{array}
\]
So we want
\[
\begin{array}{l}
P_n \frac{\tilde{q}_n-\bar{q}}{\bar{g}}(A-\tilde{g}_n)=0\\
P_n \frac{A}{\bar{g}}(Y-\tilde{q}_n)=0.
\end{array}
\]
So we could update $\tilde{q}_n$ with standard univariate logistic regression with off-set $\log \tilde{q}_n/(1-\tilde{q}_n)$ and $\epsilon A/\bar{g}$ clever covariate. In addition, given this updated $\tilde{q}_n^*$, we could then update $\tilde{g}_n$ with  a standard univariate logistic regression using off-set $\log \tilde{g}_n/(1-\tilde{g}_n)$ and $\epsilon (\tilde{q}_n^*-\bar{q})/\bar{g}$. 
This results in an updated $\tilde{P}_n^*$ defined by $q_{W,n},\tilde{g}_n^*,\tilde{q}_n^*$ so that $P_n A/\bar{g}(Y-\tilde{q}_n^*)=0$ and $P_n (\tilde{q}_n^*-\bar{q})/\bar{g} (A-\tilde{g}_n^*)=0$. By above derivation, this implies then
$P_n \frac{A}{\bar{g}}(Y-\bar{q})-\tilde{P}_n^* \frac{A}{\bar{g}}(Y-\bar{q})=0$.

{\bf Targeting $\tilde{P}_n$ to solve $(\tilde{P}_n-P_n) D^{(2)}_{n,P}=0$ at a given $P=P^0$ and $\tilde{P}_n$:}
Recall $D^{(2)}_{n,P}=C^a_{n,P}(A-\bar{g})+C^u_{n,P}(Y-\bar{q})$.
We have
\[
\begin{array}{l}
\tilde{P}_n C^a_{n,P} (A-\bar{g})-P_n C^a_{n,P}(A-\bar{g})=P_n C^a_{n,P}(\tilde{g}_n-\bar{g})-P_n C^a_{n,P}(A-\bar{g})\\
=-P_n C^a_{n,P}(A-\tilde{g}_n).
\end{array}
\]
So for this we want to update $\tilde{g}_n$ with a standard univariate logistic regression using off-set $\log \tilde{g}_n/(1-\tilde{g}_n)$ and clever covariate $\epsilon C^a_{n,P}$, where we use the current $\tilde{P}_n$ in the expression for $C^a_{n,P}$.

Define
\[
\tilde{C}^y_{n,P}\equiv \frac{1}{\bar{g}}\frac{\bar{q}_n^{(1)}(p)(1-\bar{q}_n^{(1)}(p))}{\bar{q}(1-\bar{q})}\left\{1-\tilde{c}^{(1)}_{\tilde{P}_n,P} \frac{\tilde{g}_n}{\bar{g}}\right).\]
Notice that $C^y_{n,P}=A \tilde{C}^y_{n,P}$, and
\[
\tilde{P}_n C^y_{n,P}(Y-\bar{q})=P_n \tilde{C}^y_{n,P}\tilde{g}_n (\tilde{q}_n-\bar{q}).\]
Thus,
\[
\begin{array}{l}
(\tilde{P}_n-P_n) C^y_{n,P}(Y-\bar{q})=
P_n (\tilde{q}_n-\bar{q})\tilde{C}^y_{n,P} \tilde{g}_n- P_n (Y-\bar{q}) \tilde{C}^y_{n,P}A\\
=P_n (\tilde{q}_n-Y) \tilde{C}^y_{n,P}A+P_n (\tilde{q}_n-\bar{q})\tilde{C}^y_{n,P}(\tilde{g}_n-A)\\
=-P_n C^y_{n,P}(Y-\tilde{q}_n)-P_n (\tilde{q}_n-\bar{q})\tilde{C}^y_{n,P}(A-\tilde{g}_n).
\end{array}
\]
So for this we want to target $\tilde{q}_n$ with covariate $C^y_{n,P}$, and given the resulting update $\tilde{q}_n^*$, we want to target $\tilde{g}_n$ with covariate $(\tilde{q}_n^*-\bar{q})\tilde{C}^y_{n,P}$. Notice that $C^y_{n,P}$ and $C^a_{n,P}$ depend themselves on $\tilde{P}_n$ so that we would use the current $\tilde{P}_n$ for these covariates. 

{\bf Joint targeting of $(\tilde{P}_n-P_n)D^{(1)}_{P}=(\tilde{P}_n-P_n)D^{(2)}_{n,P}=0$ at a given $P=P^0$ and current $\tilde{P}_n$:}
To conclude, we want to target $\tilde{P}_n$ so that 
\begin{eqnarray*}
P_n C^y_{\tilde{P}_n,P}(Y-\tilde{q}_n)&=&0\\
P_n \frac{A}{\bar{g}}(Y-\tilde{q}_n)&=&0\\
P_n (\tilde{q}_n-\bar{q})\tilde{C}^y_{\tilde{P}_n,P}(A-\tilde{g}_n)&=&0\\
P_n C^a_{\tilde{P}_n,P}(A-\tilde{g}_n)&=&0\\
P_n \frac{\tilde{q}_n-\bar{q}}{\bar{g}}(A-\tilde{g}_n)&=&0.
\end{eqnarray*}
The targeting of $\tilde{q}_n$ is achieved by including two covariates $A/\bar{g}$ and $C^y_{n,P}$ 
Given the resulting $\tilde{q}_n^*$, 
we target $\tilde{g}_n$ by including three covariates 
$C^a_{\tilde{P}_n,P}$, $(\tilde{q}_n^*-\bar{q})/\bar{g}$ and $(\tilde{q}_n^*-\bar{q})\tilde{C}^y_{\tilde{P}_n,P}$.
The targeting can be carried out with TMLE as mentioned above, but, one could also add these  extra covariates to  the  spline basis functions of the lasso estimators $\tilde{g}_n$ and $\tilde{q}_n$, not penalizing the corresponding coefficients. One could also use a universal least favorable path TMLE of $\tilde{P}_n$ by tracking these local least favorable paths with small update steps. In that case, for a given $P$, we end up with a $\tilde{P}_n^*$ (indexed by $P$) that solves
all of the above equations with $\tilde{P}_n$ replaced by $\tilde{P}_n^*$.

In order to emphasize that this TMLE $\tilde{P}_n^*$ is indexed by the choice $P$, we might denote it with $\tilde{P}_n^*(P)$.

\subsection{Iterative targeting of the HAL-MLE in the second order TMLE}
\begin{itemize}
\item Let $\tilde{P}^{(1)}(P,\tilde{P}_n)$ be the first order TMLE-update of initial $P$ using $\tilde{P}_n$ as HAL-MLE and solving ${P}_n D^{(1)}_{\tilde{P}^{(1)}(P,\tilde{P}_n)}=0$. Here we decided to use a regular TMLE under $P_n$ as well, so that we have a single empirical log-likelihood criterion that will increase at each step, thereby providing us with a guarantee that this iterative second order TMLE algorithm will indeed converge.
\item Let $\tilde{P}^{(2)}(P,\tilde{P}_n)$ be the second order TMLE update which uses as initial $P$ and $\tilde{P}_n$ as HAL-MLE, and is targeted to solve $P_n D^{(2)}_{\tilde{P}_n,\tilde{P}_n^{(2)}(P,\tilde{P}_n)}=0$. 

\item We start with $\bar{q}_n^0$ and $\bar{g}_n^0$, while we use $q_{W,n}$ as estimator of marginal distribution of $W$. This represents the initial estimator $P_n^0$.

\item We target $\tilde{P}_n$ w.r.t. $\bar{q}_n^0,\bar{g}_n^0$ so that it solves
$(\tilde{P}_n-P_n)D^{(1)}_{P_n^0}=0$ and $(\tilde{P}_n-P_n)D^{(2)}_{\tilde{P}_n,P_n^0}$.
Let's denote this targeted HAL-MLE with $\tilde{P}_n^{*,0}\equiv \tilde{P}_n(P_n^0)$.
\item Given this $P_n^0$ and $\tilde{P}_n^{*,0}$, we now compute the second order TMLE.
Thus, firstly, we compute the second TMLE-update $\tilde{P}^{(2)}_n(P_n^0,\tilde{P}_n^{*,0})$ at the current T-HAL-MLE $\tilde{P}_n^*(P_n^0)$ and initial estimator $P_n^0$, so that we solve $P_n D^{(2)}_{\tilde{P}_n^{*,0},\tilde{P}_n^{(2)}(P_n^0,\tilde{P}_n^{*,0})}=0$.
Then, we compute $\tilde{P}^{(1)}_n(\tilde{P}^{(2)}_n(P_n^0,\tilde{P}_n^{*,0}) , \tilde{P}_n^{*,0})$
so that we solve \[
P_n D^{(1)}_{\tilde{P}^{(1)}(\tilde{P}^{(2)}_n(P_n^0),\tilde{P}_n^{*,0} ) } =0.\]

\item So we have now computed one round of second order TMLE starting with initial $P_n^0$ and using $\tilde{P}_n^{*,0}$ as HAL-MLE, and we have increased the empirical log-likelihood at both the second and first order TMLE update. Let $P^{2,(1),*}_n(P_n^0,\tilde{P}_n^{*,0} )$ be the resulting second order TMLE of $P_0$. This represents our current estimator of $\bar{g}_0$ and $\bar{q}_0$, while using $q_{W,n}$ for the marginal distribution of $W$.

\item If $(\tilde{P}_n-P_n)D^{(1)}_{\tilde{P}^{(1)}(\tilde{P}^{(2)}(P_n^0,\tilde{P}_n^{*,0}  ) )}$ or
$(\tilde{P}_n-P_n)D^{(2)}_{\tilde{P}_n^{*,0},\tilde{P}^{(2)}_n(P_n^0,\tilde{P}_n^{*,0}) }$
is not smaller than $\sigma_{1n}/(n^{1/2}\log n)$, then we proceed with next steps, otherwise, we are done.

\item Let $P_n^1=P^{2,(1),*}_n(P_n^0,\tilde{P}_n^{*,0}) $ the new initial estimator for the second round of second order TMLE. Let $\tilde{P}_n^{*,1}=\tilde{P}_n^*(P^{2,(1),*}_n(P_n^0,\tilde{P}_n^{*,0} ))$.  Thus, we retarget the HAL-MLE but now w.r..t its current targeted estimators and its current targeted HAL-MLE. We now compute the second round second order TMLE $P^{2,(1),*}(P_n^1,\tilde{P}_n^{*,1})$.
\item We iterate this process  of computing the second order TMLE: Starting at $m=1$, 
compute $P_n^{m+1}=P^{2,(1),*}(P_n^m,\tilde{P}_n^{*,m})$ and $\tilde{P}_n^{*,m+1}=\tilde{P}_n^*(P_n^{m+1})$; $m=m+1$, till we achieve  $(\tilde{P}_n-P_n)D^{(1)}_{P^{2,(1),*}_n(P_n^m,\tilde{P}_n^{*,m}) } \approx 0$ and $(\tilde{P}_n-P_n)D^{(2)}_{\tilde{P}_n^{(2)}(P_n^m,\tilde{P}_n^{*,m})}\approx 0$. 
\item At the final step, the resulting second order TMLE $P^{2,(1),*}_n$ and its targeted $\tilde{P}_n^*=\tilde{P}_n^*$  have solved
$P_n D^{(1)}_{P^{2,(1),*}_n}=0$, $\tilde{P}_n^* D^{(1)}_{P^{2,(1),*}_n}=0$,
$P_n D^{(2)}_{\tilde{P}_n^*,P^{2,(1),*}_n}=0$, and $\tilde{P}_n D^{(2)}_{\tilde{P}_n^*,P^{2,(1),*}_n}=0$, either all exactly or up till desired precision.
\end{itemize}
Notice that for this iterative second order TMLE that also targets the HAL-MLE there is no need to undersmooth the HAL-MLE beyond the targeting.

\subsection{Analyzing the iterative HAL-regularized second order TMLE that also targets the HAL-MLE.}
 We have $\Psi_{\tilde{P}_n}^{(1)}(P)=\Psi(\tilde{P}^{(1)}(P,\tilde{P}_n))$, where
$\tilde{P}^{(1)}(P,\tilde{P}_n)$ solves $\tilde{P}_n D^{(1)}_{\tilde{P}^{(1)}(P,\tilde{P}_n)}=0$.
Let $D^{(2)}_{\tilde{P}_n,P}$ be the canonical gradient of $\Psi_{\tilde{P}_n}^{(1)}(P)$ at $P$.
We have $\tilde{P}^{(2)}(P,\tilde{P}_n)$ solving $\tilde{P}_n D^{(2)}_{\tilde{P}_n,\tilde{P}^{(2)}(P,\tilde{P}_n)}=0$.
This defines the estimator $\Psi(\tilde{P}^{(1)}(\tilde{P}^{(2)}(P,\tilde{P}_n),\tilde{P}_n))$, which is  the second order TMLE using $\tilde{P}_n$ in its TMLE-updates. 
We have an exact expansion for this second order TMLE at any $(P,\tilde{P}_n)$.
Thus,
\[
\begin{array}{l}
\Psi(\tilde{P}_{\tilde{P}_n}^{(1)}\tilde{P}^{(2)}_{\tilde{P}_n}(P))-\Psi(P_0)=
\Psi(\tilde{P}_{\tilde{P}_n}^{(1)}(P_0))-\Psi(P_0)\\
+\Psi(\tilde{P}_{\tilde{P}_n}^{(1)}\tilde{P}_{\tilde{P}_n}^{(2)}(P))-\Psi(\tilde{P}_{\tilde{P}_n}^{(1)}(P_0))\\
=(\tilde{P}_n-P_0)D^{(1)}_{\tilde{P}_{\tilde{P}_n}^{(1)}(P_0)}+R^{(1)}(\tilde{P}_{\tilde{P}_n}^{(1)}(P_0),P_0)\\
+(\tilde{P}_n-P_0)D^{(2)}_{\tilde{P}_n,\tilde{P}_{\tilde{P}_n}^{(2)}(P_0)}+R^{(2)}_{\tilde{P}_n}(\tilde{P}_n^{(2)}(P_0),P_0)\\
+R^{(1)}(\tilde{P}_{\tilde{P}_n}^{(1)}\tilde{P}_{\tilde{P}_n}^{(2)}(P),\tilde{P}_n)-
R^{(1)}(\tilde{P}_{\tilde{P}_n}^{(1)}\tilde{P}_{\tilde{P}_n}^{(2)}(P_0),\tilde{P}_n).
\end{array}
\]
Importantly, the final line represents a difference of two third order differences, due to the TMLE-updates exactly solving the $\tilde{P}_n$-efficient score equations.

In this exact expansion we still have room to select $\tilde{P}_n$ and the initial estimator $P$, as long as $d\tilde{P}_n/dP$ exists so that we preserve the pathwise differentiability of $\Psi_n^{(1)}(P)$. We can target $\tilde{P}_n^*$ and accordingly update $P$ so that it solves the equations $(\tilde{P}_n-P_n)D^{(1)}_{\tilde{P}_{\tilde{P}_n}^{(1)}(P)}=0$ and $(\tilde{P}_n-P_n)D^{(2)}_{\tilde{P}_n,\tilde{P}_{\tilde{P}_n}^{(2)}(P)}=0$. 
This will also help with solving $P_n D^{(1)}_{\tilde{P}_{\tilde{P}_n}^{(1)}(P_0)}\approx 0$.  See also our discussion of this algorithm in Section \ref{sectionthotmle} formally showing that due to this extra targeting of the HAL-MLE $P_n D^{(1)}_{\tilde{P}_{\tilde{P}_n}^{(1)}(P_0)}$ is reduced to a third order term. 

In addition,  we could also replace the initial $P$ by this second order TMLE and again target $\tilde{P}_n$ w.r.t. this new initial estimator. More generally, we can iterate this process, which defines our iterative second order TMLE targeting $\tilde{P}_n$. At any number of iterations, we can still apply the above exact expansion for the corresponding second order TMLE with this new targeted $\tilde{P}_n^*$ and updated initial $P$. Consider now a final iteration resulting in a targeted HAL-MLE $\tilde{P}_n^*$ and final second order TMLE $P_n^*$, which now also solves
 $(P_n-\tilde{P}_n^*)D^{(1)}_{P_n^*}=(P_n-\tilde{P}_n^*)D^{(2)}_{\tilde{P}_n^*,P_n^*}=0$.
 We can still apply the above exact expansion to this particular pair of initial and targeted HAL-MLE. 
 The gain we have achieved is that the size of $P_n D^{(1)}_{\tilde{P}_n^{(1)}(P_0)}$ and $P_n D^{(2)}_{\tilde{P}_n^*,\tilde{P}_n^{(2)}(P_0)}$ will be negligible, due to these TMLE-updates using a targeted $\tilde{P}_n^*$ that behaves as an empirical $P_n$-TMLE. In particular, the need for undersmoothing $\tilde{P}_n$ is minimized, thereby controlling the finite sample robustness of the second order TMLE.  An additional practical benefit of this procedure is that the discrepancy between $\epsilon_{P_n}^{(j)}(P)$ and $\epsilon_{\tilde{P}_n}^{(j)}(P)$ in the procedure at the relevant $P$ is arbitrarily small so that the user can as well implement the empirical TMLE-updates $\epsilon_{P_n}^{(j)}(P)$, $j=1,2$, even though in our analysis we act as if we used $\epsilon_{\tilde{P}_n}^{(1)}(P)$ at the final choice of initial and targeted HAL-MLE. Since the empirical likelihood increases under such empirical TMLE-updates at each step, this provides a guarantee that the iterative second order TMLE will converge.


\section{Second order TMLE for treatment specific mean example with continuous outcome}
\label{AppendixJ}
Let $\tilde{q}_n^{(1)}(q_{\delta,h_1},g)$ be a linear regression model with $\delta C_1+\epsilon_n(q_{\delta,h},g)C_g$ in its linear form and normal independent error with variance $\sigma^2(q)(W,A)$, so that the log-likelihood loss for $\bar{q}_0$ equals the inverse weighted squared error loss. This results in a different definition of $\tilde{P}_n^{(1)}(P)$, and thereby $\Psi_n^{(1)}(P)$ and its canonical gradient $D^{(2)}_{n,P}$. In this section we will compute the latter canonical gradient which then results in a new second order TMLE for the treatment specific mean for a continuous outcome.
We have
\[
\begin{array}{l}
A_{n,p,1}(h_1)=\frac{d}{d\delta_0}\log \tilde{q}_n^{(1)}(q_{\delta_0},g)\\
=\sigma^{-2}(q)C_1(W,A)(Y-\tilde{q}_n^{(1)}(q,g)(1\mid W,A))\\
+\sigma^{-2}(q) \frac{d}{d\delta_0}\epsilon_n(q_{\delta_0},g)C_g(Y-\tilde{q}_n^{(1)}(q,g)(1\mid W,A)).
\end{array}
\]
 

{\bf Score operator mapping score of path through initial $g$ into score of TMLE update:}
Note that $\tilde{q}_n^{(1)}(q,g_{\delta,h_2})$ is a linear regression model with $
\epsilon_n(q,g_{\delta})C_{g_{\delta}}$ in its linear form. 
Also note that $\frac{d}{d\delta}\bar{g}_{\delta}=C_2\bar{g}(1-\bar{g})$.
We have 
\[
\begin{array}{l}
\frac{d}{d\delta_0}C_{g_{\delta_0}}=-\frac{A}{\bar{g}^2} \frac{d}{d\delta_0}\bar{g}_{\delta_0}
=-\frac{A}{\bar{g}^2}C_2\bar{g}(1-\bar{g})=-\frac{A}{\bar{g}}C_2(1-\bar{g}).
\end{array}
\]
Thus,
\[
\begin{array}{l}
A_{n,p,2}(h_2)=h_2+\frac{d}{d\delta_0}\log \tilde{q}_n^{(1)}(q,g_{\delta_0})\\
=h_2+\frac{d}{d\delta_0}\epsilon_n(q,g_{\delta_0}\sigma^{-2}(q)C_g(Y-\tilde{q}_n^{(1)}(1\mid W,A))\\
+\epsilon_n(q,g)\frac{d}{d\delta}C_{g_{\delta_0}}\sigma^{-2}(q)(Y-\tilde{q}_n^{(1)}(1\mid W,A))\\
=h_2+\frac{d}{d\delta_0}\epsilon_n(q,g_{\delta_0}) \sigma^{-2}(q)C_g(Y-\tilde{q}_n^{(1)}(1\mid W,A))\\
-\epsilon_n(q,g)\frac{A}{\bar{g}}C_2(1-\bar{g}) \sigma^{-2}(q)(Y-\tilde{q}_n^{(1)}(1\mid W,A))\\
=h_2+\frac{d}{d\delta_0}\epsilon_n(q,g_{\delta_0}) \sigma^{-2}(q)C_g(Y-\tilde{q}_n^{(1)}(1\mid W,A))\\
-\sigma^{-2}(q)\epsilon_n(q,g)\frac{A}{\bar{g}}h_{2}(W,A)  (Y-\tilde{q}_n^{(1)}(1\mid W,A)).
\end{array}
\]
Here we used that  $C_2(1-\bar{g})=h_2(1,W)$ so that we have succeeded to express the score operator as a linear mapping in $h_{2}$.

 Let $A_{n,p,1}:H_1(P)\rightarrow L^2_0(\tilde{P}_n^{(1)}(P))$, where $H_1(P)=\{h_1\in L^2_0(P): E(h_1\mid W,A)=0\}$. We have $H_1(P)=\{C_1(Y-\bar{q}(W,A)):C_1\}$.
 Let $A_{n,p,2}:H_2(P)\rightarrow L^2_0(\tilde{P}_n^{(1)}(P))$, where $H_2(P)=\{h_2(W,A): E(h_2\mid W)=0\}$. We have $H_2(P)=\{C_2(A-\bar{g}(W)):C_2\}$.
Note that $H_1(P)\perp H_2(P)$ are orthogonal spaces in $L^2_0(P)$. In Appendix \ref{AppendixH}  we derived the  form of these score operators $A_{n,p,1}$ and $A_{n,p,2}$ up till the pathwise derivative of $\epsilon_n^{(1)}(q,g)$:
\begin{eqnarray*}
A_{n,p,1}(h_1)&=&\sigma^{-2}(q)C_1(W,A)(Y-\tilde{q}_n^{(1)}(q,g)(1\mid W,A))\\
&&+\frac{d}{d\delta_0}\epsilon_n(q_{\delta_0},g)\sigma^{-2}(q)C_g(Y-\tilde{q}_n^{(1)}(q,g)(1\mid W,A))\\
A_{n,p,2}(h_2)&=&h_2+ \sigma^{-2}(q)\frac{d}{d\delta_0}\epsilon_n(q,g_{\delta_0}) C_g(Y-\tilde{q}_n^{(1)}(1\mid W,A))
\\
&&-\epsilon_n(q,g)\frac{A}{\bar{g}}h_{2}(W,A) \sigma^{-2}(q) (Y-\tilde{q}_n^{(1)}(1\mid W,A)).
\end{eqnarray*}

{\bf Representing second order canonical gradient and understanding its two components:}
Let $A_{n,p,1}^{\top}:L^2_0(\tilde{P}_n^{(1)}(P))\rightarrow H_1(P)$ and $A_{n,p,2}^{\top}:L^2_0(\tilde{P}_n^{(1)}(P))\rightarrow H_2(P))$ be the adjoints of $A_{n,p,1}$ and $A_{n,p,2}$, respectively. 
The second order canonical gradient can be represented as
\[
D^{(2)}_{n,P}=A_{n,p,1}^{\top}(D^{(1)}_{\tilde{p}_n^{(1)}(p)})+A_{n,p,2}^{\top}(D^{(1)}_{\tilde{p}_n^{(1)}(p)} ).\] 
We should also note that $A_{n,p,1}^{\top}(D^{(1)}_{\tilde{p}_n^{(1)}(p)})$ is the canonical gradient of $q\rightarrow \Psi(\tilde{P}_n^{(1)}(q,g))$ and $A_{n,p,2}^{\top}(D^{(1)}_{\tilde{p}_n^{(1)}(p)})$ is the canonical gradient of $g\rightarrow \Psi(\tilde{P}_n^{(1)}(q,g))$.
We will compute both components $D^{(2)}_{n,P,1} =A_{n,p,1}^{\top}(D^{(1)}_{\tilde{p}_n^{(1)}(p)})$ and $D^{(2)}_{n,P,2} =A_{n,p,2}^{\top}(D^{(1)}_{\tilde{p}_n^{(1)}(p)})$
of $D^{(2)}_{n,P}=D^{(2)}_{n,P,1}+D^{(2)}_{n,P,2}$. One can be used to target $q$ and the other can be used to target $g$. 


{\bf Determining the pathwise derivative of $\epsilon_n^{(1)}(p)$:}
We have that $\epsilon_n^{(1)}(p)$ solves
\[
\begin{array}{l}
0=\tilde{P}_n C_g(W,A)(Y-\bar{q}^{(1)}(q,g,\epsilon)(W,A))\\
=
\tilde{P}_n C_g(W,A)\left (Y-\bar{q}-\epsilon C_g(W,A) \right) .\end{array}
\]
Let's denote this equation with $U(\epsilon(p),p=(q,g))$. 
We have $U(\epsilon(p_{\delta}),p_{\delta})=0$. So, at $\epsilon =\epsilon_n^{(1)}(p)$, we have
\[
\frac{d}{d\delta_0}\epsilon_n^{(1)}(p_{\delta_0})=-\left\{\frac{d}{d\epsilon}U(\epsilon,p)\right\}^{-1}
\frac{d}{d\delta_0}U(\epsilon,p_{\delta_0}).\]
We have
\[
c_P^{(1)}\equiv -\frac{d}{d\epsilon}U(\epsilon,p)=
\tilde{P}_n C_g^2 .\]
Note also $\frac{d}{d\delta}\bar{q}_{\delta}=C_1$ and recall $\frac{d}{d\delta_0}C_{g_{\delta_0}}=
-A/\bar{g} C_2(1-\bar{g})$.
Note that 
\[
\begin{array}{l}
\frac{d}{d\delta_0}\epsilon_n^{(1)}(q_{\delta_0},g)=c_P^{(1),-1}\frac{d}{d\delta_0}U(\epsilon_n^{(1)}(p),q_{\delta_0},g)\\
=-c_P^{(1),-1}\tilde{P}_n C_g C_1\\
\frac{d}{d\delta_0}\epsilon_n^{(1)}(q,g_{\delta_0})=
c_P^{(1),-1}\frac{d}{d\delta_0}U(\epsilon_n^{(1)}(p),q,g_{\delta_0})\\
=c_P^{(1),-1}\tilde{P}_n \frac{d}{d\delta_0}C_{g_{\delta_0}}(Y-\bar{q}-\epsilon_n^{(1)}(p) C_g)\\
-c_P^{(1),-1}\tilde{P}_n \epsilon_n^{(1)}(p) C_g \frac{d}{d\delta_0}C_g\\
=-c_P^{(1),-1}\tilde{P}_n \frac{A}{\bar{g}}C_2(1-\bar{g})(Y-\bar{q}-\epsilon_n^{(1)}(p) C_g)\\
+c_P^{(1),-1}\tilde{P}_n \epsilon_n^{(1)}(p) C_g \frac{A}{\bar{g}}C_2(1-\bar{g}).
\end{array}
\]

{\bf Final forms of score operators $A_{n,p,1}$ and $A_{n,p,2}$:}
Thus,
\[
\begin{array}{l}
A_{n,p,1}(h_1)=\frac{d}{d\delta_0}\log \tilde{q}_n^{(1)}(q_{\delta_0},g)\\
=\sigma^{-2}(q)C_1(W,A)(Y-\bar{q}_n^{(1)}(q,g)(W,A))\\
+\sigma^{-2}(q)\frac{d}{d\delta_0}\epsilon_n^{(1)}(q_{\delta_0},g)C_g(Y-\bar{q}_n^{(1)}(q,g)(W,A))\\
=\sigma^{-2}(q)C_1(Y-\bar{q}_n^{(1)}(q,g))\\
-c_P^{(1),-1}\{\tilde{P}_n C_g C_1(Y-\bar{q})(Y-\bar{q})\sigma^{-2}(q) \} \sigma^{-2}(q)C_g(Y-\bar{q}_n^{(1)}(q,g)(W,A))\\
=\sigma^{-2}(q)C_1(Y-\bar{q}_n^{(1)}(q,g))\\
-c_P^{(1),-1}\{\tilde{P}_n C_g h_1(Y-\bar{q})\sigma^{-2}(q) \} \sigma^{-2}(q)C_g(Y-\bar{q}_n^{(1)}(q,g)(W,A))\\\end{array}
\]
and
\[
\begin{array}{l}
A_{n,p,2}(h_2)=h_2+\frac{d}{d\delta_0}\log \tilde{q}_n^{(1)}(q,g_{\delta_0})\\
=h_2-\epsilon_n^{(1)}(q,g)\frac{A}{\bar{g}}h_{2g}(W,A)\sigma^{-2}(q)  (Y-\bar{q}_n^{(1)}(q,g) )\\
+\frac{d}{d\delta_0}\epsilon_n^{(1)}(q,g_{\delta_0}) \sigma^{-2}(q)C_g(Y-\bar{q}_n^{(1)}(q,g))\\
=h_2-\sigma^{-2}(q)\epsilon_n^{(1)}(q,g)\frac{A}{\bar{g}}h_{2g}(W,A)  (Y-\bar{q}_n^{(1)}(q,g) )\\
-\sigma^{-2}(q) c_P^{(1),-1}\{\tilde{P}_n \frac{A}{\bar{g}}C_2(1-\bar{g})(Y-\bar{q}-\epsilon C_g)\} C_g(Y-\bar{q}_n^{(1)}(p))\\
+\sigma^{-2}(q) c_P^{(1),-1}\left\{ \tilde{P}_n \epsilon C_g \frac{A}{\bar{g}}C_2(1-\bar{g})\right\} C_g(Y-\bar{q}_n^{(1)}(p))\\
=h_2-\sigma^{-2}(q)\epsilon_n^{(1)}(q,g)\frac{A}{\bar{g}}h_{2g}(W,A)  (Y-\bar{q}_n^{(1)}(q,g) )\\
-\sigma^{-2}(q) c_P^{(1),-1}\{\tilde{P}_n \frac{A}{\bar{g}}h_2(Y-\bar{q}-\epsilon C_g)\} C_g(Y-\bar{q}_n^{(1)}(p))\\
+\sigma^{-2}(q) c_P^{(1),-1}\left\{ \tilde{P}_n \epsilon C_g \frac{A}{\bar{g}}h_2\right\} C_g(Y-\bar{q}_n^{(1)}(p))\\
\end{array}
\]
This establishes the representations of the score operators. 

\subsection{Determining the adjoints of score operators}

{\bf Determining the adjoint of $A_{n,p,1}$:}
We have
\[
\begin{array}{l}
\tilde{P}_n^{(1)}(p) D^{(1)}_{\tilde{P}_n^{(1)}(p)} A_{n,p,1}(h_1)=
\tilde{P}_n^{(1)}(p) D^{(1)}_{\tilde{P}_n^{(1)}(p)}
\sigma^{-2}(q)C_1(Y-\bar{q}_n^{(1)}(q,g))\\
-\tilde{P}_n^{(1)}(p) D^{(1)}_{\tilde{P}_n^{(1)}(P)} c_P^{(1),-1}\{\tilde{P}_n C_g C_1 \} \sigma^{-2}(q)C_g(Y-\bar{q}_n^{(1)}(q,g)(W,A))\\
=\tilde{P}_n^{(1)}(p) C_g(Y-\bar{q}_n^{(1)}(p)) 
\sigma^{-2}(q)C_1(Y-\bar{q}_n^{(1)}(q,g))\\
-\tilde{P}_n^{(1)}(p) C_g(Y-\bar{q}_n^{(1)}(p)) c_P^{(1),-1}\{\tilde{P}_n C_g C_1 \} \sigma^{-2}(q) C_g(Y-\bar{q}_n^{(1)}(q,g)(W,A))\\
=\tilde{P}_n^{(1)}(p) C_gC_1\sigma^{-2}(q) (Y-\bar{q}_n^{(1)}(p))^2\\
-\tilde{P}_n^{(1)}(p) C_g^2\sigma^{-2}(q) (Y-\bar{q}_n^{(1)}(p))^2 c_P^{(1),-1}\{\tilde{P}_n C_g C_1 \} \\
=\tilde{P}_n^{(1)}(p) C_gC_1 \\
-\tilde{P}_n^{(1)}(p) C_g^2   c_P^{(1),-1}\{\tilde{P}_n C_g C_1 \} \\
=P C_gC_1 \\
-P C_g^2  c_P^{(1),-1}\{P \frac{\tilde{g}_n}{\bar{g}}  C_g C_1 \} \\
=
P C_g C_1(Y-\bar{q}) (Y-\bar{q})\sigma^{-2}(q)\\
-P C_g^2  c_P^{(1),-1}\{P\frac{\tilde{g}_n}{\bar{g}} C_g C_1(Y-\bar{q})(Y-\bar{q})\sigma^{-2}(q) \} \\
=P C_g h_1 (Y-\bar{q})\\
-P C_g^2 c_P^{(1),-1}\{P\frac{\tilde{g}_n}{\bar{g}} C_g h_1 (Y-\bar{q}) \}
\\
=P h_1\left\{  C_g  (Y-\bar{q}))\right\}\\
-P C_g^2  c_P^{(1),-1}
P h_1\left\{ \frac{\tilde{g}_n}{\bar{g}} C_g (Y-\bar{q})\right\}.
\end{array}
\]
We  have $c_P^{(1)}=P_n 1/\bar{g}$.
Note also that $P C_g^2=P_n 1/\bar{g}$.
So $P C_g^2c_P^{(1),-1}=1$.
Thus,
\[
\begin{array}{l}
D^{(2)}_{n,P,1}\equiv A_{n,p,1}^{\top}(D^{(1)}_{\tilde{P}_n^{(1)}(p)})=
\left\{  C_g  (Y-\bar{q})\right\}\\
-\left\{ \frac{\tilde{g}_n}{\bar{g}} C_g (Y-\bar{q})\right\}.
\end{array}
\]
So we obtain
\[
\begin{array}{l}
D^{(2)}_{n,P,1}=C^y_P(Y-\bar{q}),
\end{array}
\]
where
\[
C^y_P=C_g\left\{1-\frac{\tilde{g}_n}{\bar{g}}\right\}.\]
We note that at $P=\tilde{P}_n$ we have $D^{(2)}_{n,\tilde{P}_n,1}=0$.


{\bf Determining the adjoint of $A_{n,2,p}$:}
Recall that $A_{n,2,p}(h_2)$ is a sum of four terms. We will determine the adjoint  for each term separately. \nl
{\bf Term 1:}
\[
\begin{array}{l}
\tilde{P}_n^{(1)}(P)D^{(1)}_{\tilde{P}_n^{(1)}(P)} h_2=
\tilde{P}_n^{(1)}(p) \frac{A}{\bar{g}}(Y-\bar{q}_n^{(1)}(p)) h_{2}=0.
\end{array}
\]
Therefore, the adjoint obtains no contribution from this term. 
\nl
{\bf Term 2:}
We need minus contribution from following:
\[
\begin{array}{l}
\tilde{P}_n^{(1)}(P)D^{(1)}_{\tilde{P}_n^{(1)}(P)}\sigma^{-2}(q)  \epsilon_n^{(1)}(p)\frac{A}{\bar{g}}h_{2}(W,A)  (Y-\bar{q}_n^{(1)}(p)(1))\\
=\tilde{P}_n^{(1)}(P) \frac{A}{\bar{g}}(Y-\bar{q}_n^{(1)}(p)) \sigma^{-2}(q) \epsilon_n^{(1)}(p)\frac{A}{\bar{g}}h_{2}(W,A)  (Y-\bar{q}_n^{(1)}(p))\\
=
\epsilon_n^{(1)}(p) \tilde{P}_n^{(1)}(p) \frac{A}{\bar{g}^2}\sigma^{-2}(q)(Y-\bar{q}_n^{(1)})^2 h_2\\
=\epsilon_n^{(1)}(p)P \frac{A}{\bar{g}^2} h_2\\
=\epsilon_n^{(1)}(p)P \frac{A-\bar{g}}{\bar{g}^2}h_2.
\end{array}
\]
Thus, the contribution is 
\[ - \epsilon_n^{(1)}(p)\frac{A-\bar{g}}{\bar{g}^2}.\]

{\bf Term 3:} We need minus sign from following: 

\[
\begin{array}{l}
\tilde{P}_n^{(1)}(P) D^{(1)}_{\tilde{P}_n^{(1)}(P)} \sigma^{-2}(q)c_P^{(1),-1}\left\{ \tilde{P}_n\frac{A}{\bar{g}}C_2(1-\bar{g}) \left(Y-\bar{q}_n^{(1)}(q,g)\right)\right\}  C_g(Y-\bar{q}_n^{(1)}(q,g))           \\
=\tilde{P}_n^{(1)}(p)\frac{A}{\bar{g}}\sigma^{-2}(q)(Y-\bar{q}_n^{(1)}(p)) 
C_g(Y-\bar{q}_n^{(1)}(p))c_P^{(1),-1}
\tilde{P}_n \frac{A}{\bar{g}}(1-\bar{g}) C_2(Y-\bar{q}_n^{(1)}(p))\\
=P A/\bar{g}^2 c_P^{(1),-1}\tilde{P}_n  A/\bar{g}(1-\bar{g})C_2(\tilde{q}_n-\bar{q}_n^{(1)}(p))\\
=c_P^{(1),-1} P\frac{1}{\bar{g}} P A \frac{\bar{g}_n}{\bar{g}^2}
(1-\bar{g}) C_2(\tilde{q}_n-\tilde{q}_n^{(1)}(p))\\
=c_P^{(1),-1}P\frac{q}{\bar{g}}
P\frac{\tilde{g}_n}{\bar{g}^2} (1-\bar{g}) (\tilde{q}_n-\bar{q}_n^{(1)} 
\frac{ A}{A-\bar{g}} C_2(A-\bar{g}))\\
 =
c_P^{(1),-1}P\frac{1}{\bar{g}}
Ph_{2}\frac{ \tilde{g}_n}{\bar{g}^2} (1-\bar{g}) (\tilde{q}_n-\bar{q}_n^{(1)}\left\{\frac{ A}{A-\bar{g}} -\bar{g}/(1-\bar{g})\right\}\\
=c_P^{(1),-1}P\frac{1}{\bar{g}}
Ph_{2} \frac{\tilde{g}_n}{\bar{g}^2} (1-\bar{g}) (\tilde{q}_n-\bar{q}_n^{(1)}
\frac{1}{(1-\bar{g})}(A-\bar{g}) \\
=P h_{2}\left\{ \left\{ c_P^{(1),-1}P\frac{1}{\bar{g}} \right\}
 \frac{\tilde{g}_n}{\bar{g}^2} (\tilde{q}_n-\bar{q}_n^{(1)} (A-\bar{g})\right\}\\
 \end{array}
 \]
 Note $c_P^{(1),-1}P1/\bar{g}=1$.
 So contribution is given by 
 \[
-   \frac{\tilde{g}_n}{\bar{g}^2} (\tilde{q}_n-\bar{q}_n^{(1)}(p) )(A-\bar{g} .
 \]
 {\bf Term 4:}
 \[
 \begin{array}{l}
\tilde{P}_n^{(1)}(P)D^{(1)}_{\tilde{P}_n^{(1)}(P)} 
c_P^{(1),-1}\sigma^{-2}(q) \left\{ \tilde{P}_n C_g\epsilon \frac{A}{\bar{g}}C_2(1-\bar{g})  \right\} C_g(Y-\bar{q}_n^{(1)}(q,g))\\ 
= 
 \tilde{P}_n^{(1)}\frac{ A}{\bar{g}}\sigma^{-2}(q) (Y-\bar{q}_n^{(1)}(p))C_g(Y-\bar{q}_n^{(1)}(p))c_P^{(1),-1} \tilde{P}_n
 C_g \frac{A}{\bar{g}} (1-\bar{g})\epsilon_n^{(1)}(p) C_2\\
 =c_P^{(1),-1}P\frac{A}{\bar{g}^2} 
 \tilde{P}_n A\frac{1-\bar{g}}{\bar{g}^2}\epsilon_n^{(1)}(p) C_2\\
 \equiv  P\frac{\tilde{g}_n}{\bar{g}}  A\frac{1-\bar{g}}{\bar{g}^2}\epsilon_n^{(1)}(p)C_2 \\
 = P \frac{\tilde{g}_n}{\bar{g}^3}(1-\bar{g})\epsilon_n^{(1)}(p)\frac{A}{A-\bar{g}} h_{2}\\
 =P \frac{\tilde{g}_n}{\bar{g}^3}(1-\bar{g})\epsilon_n^{(1)}(p)\frac{1}{1-\bar{g}}(A-\bar{g})  h_{2} \\
 =P h_2\left\{  \frac{\tilde{g}_n}{\bar{g}^3}\epsilon_n^{(1)}(p)(A-\bar{g}) \right\}.
\end{array}
\]
Thus, the contribution is given by:
\[ \frac{\tilde{g}_n}{\bar{g}^3}\epsilon_n^{(1)}(p)(A-\bar{g}) .\]

Summing up the contributions yields:
\[
\begin{array}{l}
D^{(2)}_{n,P,2}=A_{n,P,2}^{\top}(D^{(1)}_{\tilde{P}_n^{(1)}(P)})\\
=- \epsilon_n^{(1)}(p)\frac{A-\bar{g}}{\bar{g}^2}
- \frac{\tilde{g}_n}{\bar{g}^2} (\tilde{q}_n-\bar{q}_n^{(1)}(p) )(A-\bar{g} 
 + \frac{\tilde{g}_n}{\bar{g}^3}\epsilon_n^{(1)}(p)(A-\bar{g}) .
\end{array}
\]
Notice that the first and third term can be combined into a single term $\epsilon_n^{(1)}(p)\frac{\tilde{g}_n-\bar{g}}{\bar{g}^3}(A-\bar{g})$.
Note that $D^{(2)}_{n,\tilde{P}_n,2}=0$.

This proves the desired result.
\begin{lemma}\label{lemmactsoutcome}
Thus, we have derived  $D^{(2)}_{n,P}=D^{(2)}_{n,P,1}+D^{(2)}_{n,P,2}$.
\[
\begin{array}{l}
D^{(2)}_{n,P}=
C_g\frac{\bar{g}-\tilde{g}_n}{\bar{g}}   (Y-\bar{q}) 
+ \epsilon_n^{(1)}(p)\frac{\tilde{g}_n-\bar{g}}{\bar{g}^3}(A-\bar{g})
- \frac{\tilde{g}_n}{\bar{g}^2} (\tilde{q}_n-\bar{q}_n^{(1)}(p) )(A-\bar{g}) .
 \end{array}
 \]
The first term represent $D^{(2)}_{n,P,1}$.
\end{lemma}
The second order TMLE for the treatment specific mean  of a continuous outcomecan now be defined analogue to the second order TMLE for binary outcome.

\section{Computing a canonical gradient as linear combination of  (HAL)-basis functions spanning tangent space}\label{sectionfirstgradient}

 \subsection{Characterization of tangent space for nonparametric conditional independence model as linear combination of mean zero centered HAL-basis functions} 
 The tangent space at $P$ of a nonparametric conditional density of a univariate $Y$, given vector $X$, is given by 
 $T_{Y,P}\equiv \{h(Y\mid X): E_P(h(Y\mid X)\mid X)=0\}$. Let $\{\phi_j(X,Y)=I(X>c_{j,x},Y>c_{j,y}):j\}$ be an HAL-basis of all functions of $X,Y$, where the knot-points $c_j=(c_{j,y},c_{j,x})$,  could be selected based on a large sample $(X_i,Y_i)$ from $P$, as in the typical HAL-MLE.
 We have that $\phi_j^y(X,Y)=\phi_{j,1}^y(X)\phi_{j,2}^y(Y)$, where $\phi_{j,2}^y(Y)=I(Y>c_{j,y})$ and $\phi_{j,1}^x(X)=I(X>c_{j,x})$.  Therefore, we can approximate the tangent space at $P$ with the finite dimensional linear span of $\{\phi_{j,P}^y\equiv \phi_{j,1}^y(X)(\phi_{j,2}^y(Y)-\bar{F}_Y(c_{j,y}\mid X)):j\}$,
 where $\bar{F}_Y(c\mid X)=P(Y>c\mid X)$ is the survivor function of $Y$, given $X$.
 Thus, we have $T_{Y,P}\approx\{\sum_j \beta(j)\phi_{j,P}^y:\beta\}$, where the user controls the approximation error by selecting enough knot-points.
 
 This provides the building block for the tangent space of a nonparametric model defined in terms of conditional independence assumptions.
  Suppose that $O=(Y_1,\ldots,Y_{\tau})$ is a $\tau$-dimensional vector and that we factorize the density
as $\prod_{k=1}^{\tau} p_{Y_k}(y_k\mid Pa(y_k))$, where each conditional density of $Y_k$, given its parent nodes (subset of $Y_1,\ldots,Y_{k-1}$) is unspecified beyond the specification of the parent nodes. Then the tangent space $T_{P}({\cal M})$ at $P$ in this type of nonparametric statistical model ${\cal M}$  is approximated  by the orthogonal sum of the $k$-specific tangent spaces $T_{Y_k,P}=\{\sum_j\beta(j)\phi_{j,P}^{y_k}: \beta\}$ as defined above, $k=1,\ldots,\tau$.
Then, 
$T_P({\cal M})\approx \{\sum_{j,k}\beta(j,k)\phi_{j,P}^{y_k}:\beta\}$ is well approximated by the linear span of $\{\phi_{j,P}^{y_k}:j,k\}$. Let $J_k$ be the number of basis functions in $T_{Y_k,P}$, and let $J=\sum_{k=1}^{\tau} N_k$ be the total number of basis functions for $T_P({\cal M })$.  
An element of the tangent space $T_P({\cal M})$ is now identified by a vector $(\beta(j,k):j ,k)\in \openr^J$ of coefficients.

Consider now as target parameter $\Psi(P)=E_{P^*(P)}Y_{\tau}$, the mean of the final outcome $Y_{\tau}$ under density $p^*(p)=\prod_{k=1}^{\tau} p^*_{Y_k}$, where $p^*_{Y_k}$ equals a user supplied conditional density for $k\in {\cal I}$ while $p^*_{Y_k}=p_{Y_k}$ for a subset  ${\cal I}^c$.
 We refer to the nodes in ${\cal I}$ as the intervention nodes, and $p^*_{Y_k}$ the stochastic interventions, $k\in {\cal I}$. This parameter $\Psi(P)$ represents the stochastic intervention specific mean outcome under a causal inference model and therefore covers a very large class of interesting causal quantities. Therefore, in this section we want to demonstrate how one can compute the first order canonical gradient of this target parameter with  matrix inversion.

 For notational convenience, let $\phi_{k,j,P}=\phi_{j,P}^{y_k}$.
Let $\{ P_{\delta,\phi_{k_0,j_0,P}}:\delta\}\subset{\cal M}$ be a path through $P$ defined by
$p_{\delta,\phi_{k_0,j_0,P}}(Y_k\mid Pa(Y_k))=p_{Y_k}$ if $k\not =k_0$, and
it equals $(1+\delta \phi_{k_0,j_0,P})p_{Y_k}$ when $k=k_0$. In other words, it is the path that only fluctuates the conditional density of $p_{Y_{k_0}}$ with a score equal to basis function $\phi_{k_0,j_0,P}$.
We note that, for each $\delta>0$, 
\[
\frac{\Psi(P_{\delta,\phi_{k_0,j_0,P}})-\Psi(P)}{\delta}=
E_P Y_{\tau}\phi_{k_0,j_0,P}.\]
The left-hand side represents the pathwise derivative of $\Psi$ along this path.
Let $\dot{\psi}_P\equiv (E_PY_{\tau}\phi_{k_0,j_0,P}:k_0,j_0)$ be the $M$-dimensional vector, where $M=\sum_{k=1}^{\tau}J_k$ and $J_k$ is number of basis functions $\{\phi_{k,j,P}:j\}$ for the tangent space $T_{Y_k,P}$.
Let $O_1^*,\ldots,O_N^*$ be a large sample of $N$ i.i.d. observations from $P^*$.
Then, we can compute $\dot{\psi}_P$ in one round as
\[
\dot{\psi}_P\approx \left( \frac{1}{N}\sum_{i=1}^N Y_{\tau,i}^*\phi_{k_0,j_0,P}(O_i^*): k_0,j_0\right).\]
By definition of the canonical gradient $D^*_P$, this pathwise derivative $\dot{\psi}_P(k_0,j_0)$ equals $E_P D^*_P \phi_{k_0,j_0,P}$. We represent 
$D^*_P=\sum_{k,j}\beta_P(k,j)\phi_{k,j,P}$. Thus, we have
\[
\dot{\psi}_P=\left( \sum_{k,j}\beta_P(k,j)P \phi_{k,j,P}\phi_{k_0,j_0,P}: k_0,j_0\right) .\]
 Let $\Sigma_P$ be the $M\times M$ symmetric matrix defined by
 \[
 \Sigma_P(k,j,k_0,j_0)= P \phi_{k,j,P}\phi_{k_0,j_0,P}.\]
Then,
\[
\beta_P=\Sigma_P^{-1}(\dot{\psi}_P),\]
and
\[
D^*_P=\sum_{k,j}\beta_P(k,j)\phi_{k,j,P}.\]

Instead of computing $D^*_P$ directly, one can also compute each component $D^*_{P,k}=\Pi(D^*_P\mid T_{P,Y_k})$ separately, the component of $D^*_P$ in the tangent space of the conditional distribution of $Y_k$, given $Pa(Y_k)$.
Let $\dot{\psi}_{P,k_0}\equiv (E_PY_{\tau}\phi_{k_0,j_0,P}:j_0)$ be $J_{k_0}$-dimensional vector.
Let $O_1^*,\ldots,O_N^*$ be a large sample of $N$ i.i.d. observations from $P^*$.
Then, we can compute $\dot{\psi}_{P,k_0}$ as
\[
\dot{\psi}_{P,k_0}\approx \left( \frac{1}{N}\sum_{i=1}^N Y_{\tau,i}^*\phi_{k_0,j_0,P}(O_i^*): j_0\right).\]
This pathwise derivative $\dot{\psi}_P(k_0,j_0)$ equals $E_P D^*_P \phi_{k_0,j_0,P}$. 
We represent 
$D^*_P=\sum_{k,j}\beta_P(k,j)\phi_{k,j,P}$, so that $D^*_{P,k_0}=\sum_j \beta_P(k_0,j)\phi_{k_0,j,P}$. 
The tangent spaces $T_{Y_k,P}$ are orthogonal across $k$.  Therefore, 
\[
\dot{\psi}_{P,k_0,j_0}=E_P D^*_{P,k_0}\phi_{k_0,j_0,P}=\sum_{j} \beta_P(k_0,j)P \phi_{k_0,j,P}\phi_{k_0,j_0,P}.\]
 Let $\Sigma_{P,k_0}$ be the $J_{k_0}\times J_{k_0}$ symmetric matrix defined by
 \[
 \Sigma_{P,k_0}(j,j_0)= P \phi_{k_0,j,P}\phi_{k_0,j_0,P}.\]
Then,
\[
\beta_{P,k_0}=\Sigma_{P,k_0}^{-1}(\dot{\psi}_{P,k_0}),\]
and
\[
D^*_{P,k_0}=\sum_{j}\beta_{P,k_0}(j)\phi_{k_0,j,P}.\]

  \subsection{Form of matrix $\Sigma_{P,k}$.}
 Consider the above definition of $\Sigma_{P,k}$ for the basis functions spanning the tangent space of a nonparametric conditional density of $Y_k$, given $X_k\equiv Pa(Y_k)$. 
 This matrix has the following form: given the knot points $c_1=(x_1,y_1)$ and $c_2=(x_2,y_2)$ of the two centered-spline basis functions $\phi_{P,(x_1,y_1)}=I(X>x_1)(I(Y>y_1)-\bar{F}_k(y_1\mid x))$ and $\phi_{P,(x_2,y_2)}=I(X>x_2)(I(Y>y_2)-\bar{F}_k(y_2\mid x))$, we have
 \[
 \Sigma_{P,k}(c_1,c_2)=\int_{\max(x_1,x_2)}\left( \bar{F}_k(\max(y_1,y_2)\mid x)-\bar{F}_k(y_1\mid x)\bar{F}_k(y_2\mid x)\right) dP_X(x),\]
 where $\bar{F}_k(y\mid x)=P(Y_k>x\mid X_k=x)$ and $\max(x_1,x_2)(j)=\max(x_1(j),x_2(j))$ is component wise defined.

\subsection{Explicit representation of canonical gradient in terms of pathwise derivative along HAL-basis functions of tangent space}
Consider conditional density of $Y$ given $X$ and that the statistical model does not restrict this conditional density. Let $(X,Y)$ be a $d$-dimensional vector and, for simplicity, suppose $(X,Y)\in [0,1]^d$.The tangent space $T_{Y}(P)$ of this conditional density is spanned by the basis functions $\phi_{P,(y,x)}(Y,X)=(I(Y>y)-\bar{F}(y\mid X))I(X>x)$ across all $(x,y)\in [0,1]^d$.
Let $\Psi:{\cal M}\rightarrow \openr$ be a given pathwise differentiable target parameter  with canonical gradient $D^{(1)}_P$.
We want to determine the component $D^{(1)}_{P,Y}\equiv \Pi(D^*_P\mid T_Y(P))$ of the canonical gradient in the sub-tangent space $T_Y(P)$.
Let $\dot{\psi}_{P,(x,y)}=\frac{d}{d\delta_0}\Psi(P_{\delta_0,\phi_{P,(x,y)}} )$ be the pathwise derivative along the fluctuation $(1+\delta \phi_{P,(x,y)})p_{Y\mid X}$.
We can show that
\[
D^{(1)}_{P,Y}(x,y)=\frac{1}{p(x,y)}\frac{d}{dx}\frac{d}{dy}\dot{\psi}_{P,(x,y)}
-\frac{1}{p_X(x)} \int_v \frac{d}{dx}\frac{d}{dv}\dot{\psi}_{P,(x,v)} dv .\]
The second term is simply subtracting the conditional mean of the first term, given $X$.
This is shown as follows.
\[
\begin{array}{l}
\dot{\psi}_{P,(x_0,y_0)}=P D^{(1)}_{P,Y} \phi_{P,(x_0,y_0)}\\
=
\int  D^{(1)}_{P,Y}(x,y) \{I(y>y_0)-(\bar{F}(y_0\mid x)\} I(x>x_0)dP(x,y)\\
=\int D^{(1)}_{P,Y}(x,y) I(y>y_0)I(x>x_0) dP(x,y)\\
=\int_{x_0}\int_{y_0} D^{(1)}_{P,Y}(x,y) dP(x,y),
\end{array}
\]
where we used that $D^{(1)}_{P,Y}$ has conditional mean zero, given $X$.
Let $\dot{\psi}_P(x_0,y_0)=\dot{\psi}_{P,(x_0,y_0)}$ so that we view it as a function in $(x_0,y_0)$. 
Clearly, this is a $d$-variate cadlag function, so that it generates a measure on $\openr^d$.
The above relation teaches us that
\[
\dot{\psi}_P(dx,dy)=D^{(1)}_{P,Y}(x,y)dP(x,y).\]
Thus, this proves that 
\[
D^{(1)}_{P,Y}(x,y)=\frac{d \dot{\psi}_P}{dP}(x,y)\]
is the Radon-Nikodym derivative of the measure $d\dot{\psi}_P$  w.r.t. $dP$.
Since this solution has conditional mean zero, given $X$, we can represent it as 
\[
D^{(1)}_{P,Y}(x,y)=\frac{d\dot{\psi}_P}{dP}(x,y)-\int_v \frac{d\dot{\psi}_P}{dP}(x,v) dP(v\mid x) =
\frac{d\dot{\psi}_P}{dP}(x,y)- \int_v\frac{ \dot{\psi}_P(dx,v) }{dP_X(x)} d\mu(v).\]

\subsection{First order canonical gradient in terms of initial gradient}
The above formulas apply to general basis functions $\{\phi_{j,P}: j\}$ spanning $T_{P}$, and general target parameters $\Psi$, where $\dot{\psi}_{P}(j)$ represents the pathwise derivative of $\Psi$ along  a path $\{P_{\delta,j}:\delta\}\subset {\cal M}$ with score at $\delta=0$ equal to the basis function $\phi_{j,P}$.
That is, we have $\beta_P=\Sigma_P^{-1}(\dot{\psi}_P)$, and $D^*_P\approx \sum_j \beta_P(j)\phi_{j,P}$. If we have an initial gradient $D_P$, then $\dot{\psi}_P=(E_P D_P\phi_{j,P}:j)$, so that we can avoid having to compute pathwise derivatives, but instead just evaluate empirical means.

Another way to understand this formula is as follows.
We have $D^*_P=\Pi(D_P\mid T_P)$, and $T_P$ is approximated by the linear span of $\phi_{j,P}$, $j=1,\ldots,J$.
Ignoring the approximation error due to choosing a finite set of basis functions, we have $D^*_P=\sum_j\beta_P(j)\phi_{j,P}$, where \[
\beta(P)=\arg\min_{\beta}P\left\{D_P-\sum_j \beta(j)\phi_{j,P}\right\}^2.\]
Let $\Sigma_P$ be the $J\times J$-covariance matrix defined by $\Sigma_P(j_1,j_2)=P\phi_{j_1,P}\phi_{j_2,P}$, as above.
Then, by the general least squares formula $\beta=(X^{\top}X)^{-1}X^{\top}Y$ with $X$ the design matrix with $j$-th column $X(b,j)=\phi_{j,P}(O_b)$ and $Y(b)=D_P(O_b)$ across large sample $O_b\sim P$, $b=1,\ldots,B$, it follows that
\[
\beta(P)=\Sigma_P^{-1} (E_P(\phi_{j,P}D_P):j=1,\ldots,J).\]
If the tangent space is an orthogonal sum of tangent spaces spanned by disjoint  complementary subsets of these basis functions, as in the model above defined by conditional independencies, then one can separately determine this $\beta(P)$ for each sub-tangent space.

As mentioned earlier, if the number of basis functions is too large for $\Sigma_P$ to be put in memory, then one could decide to approximate $\beta(P)$ with a highly adaptive lasso least squares regression estimator instead.

 \subsection{Using the linear combination of basis functions representation of the canonical gradient in a TMLE-update}
We note that  this algorithm does not require any analytic computations, while we still generate algebraic expressions $D^{(1)}_{P}\approx \sum_j \beta^{(1)}_P(j) \phi_{j,P}$.
These calculations defining $D^{(1)}_{P}$ need to be repeated for each $P$ at which this object is needed. If the TMLE involves MLE-updates along a local least favorable path based on $D^{(1)}_{P}$, then one only needs these objects at a single initial  $P$. Even, when it involves solving $\tilde{P}_nD^{(1)}_{\tilde{P}_n^{(k)}(P,\epsilon)}=0$, a simple steepest descent algorithm would only require knowing $D^{(1)}_{P}$ at each step of the algorithm till the equation is solved at the desired level.
 
Nonetheless, the following remark is of interest, demonstrating that the linear combination representation of a canonical gradient is quite convenient in achieving any type of corresponding TMLE-update. 
As an example, suppose that we are interested in computing a universal least favorable path TMLE-update based on a local least favorable path $\tilde{P}_n^{(1)}(P,\epsilon)$ where $D^{(1)}_{n,P}=\sum_j \beta_P(j)\phi_{j,P}$.
Let $P_n^0$ be the initial estimator. We can now apply the universal least favorable path TMLE-update based on the local least favorable path that uses $P\rightarrow \sum_j \beta^{(k)}_{P_n^0}(j)\phi_{j,P}$ as gradient. That is, we fix the coefficients at what they are at the initial estimator $P_n^0$. In this manner, we can compute a fast  universal LFM-TMLE-update $P_n^1$ that solves
$P_n \sum_j \beta^{(1)}_{P_n^0}(j)\phi_{j,P_n^1}=0$. Now, we compute $\beta^{(1)}_{P_n^1}$ and thereby obtain
$\sum_j \beta^{(1)}_{P_n^1}(j)\phi_{j,P_n^1}$. Similarly, fixing $\beta^{(1)}_{P_n^1}$ in $P\rightarrow \sum_j \beta^{(1)}_{P_n^1}(j)\phi_{j,P}$, we now use the above universal LFM-TMLE update to compute a $P_n^2$ that solves $P_n \sum_j\beta^{(1)}_{P_n^1}(j)\phi_{j,P_n^2}=0$. We iterate this $m$ times times till $P_n\sum_j \beta^{(1)}_{P_n^m}(j)\phi_{j,P_n^m}\approx 0$ at the desired level. We suggest that convergence will occur in few steps, and we also know that the log-likelihood increases at each step so that convergence is guaranteed. Note that in this algorithm, we only need to evaluate the coefficients $\beta^{(1)}_P$ at $m+1$ choices $P\in \{P_n^0,\ldots,P_n^m\}$. The same idea can be applied to determine $\epsilon_n^{(1)}(P_n^0)$ that solves $P_n D^{(1)}_{\tilde{P}_n^{(1)}(P_n^0,\epsilon)}=0$.

\subsection{Note regarding impact of using finite dimensional approximation of tangent space for the TMLE-update}
Let $D^{(1)}_P$ be the canonical gradient at $P$, while  $\tilde{D}^{(1)}_P\equiv \sum_j \beta^{(1)}_P(j)\phi_{j,P}$ with 
\[
\beta^{(1)}_P=\arg\min_{\beta}P \left(D^{(1)}_P-\sum_j \beta(j)\phi_{j,P}\right)^2\]
is the main  term linear least squares approximation.
Suppose that we compute this least squares regression on a large sample of $N$ observations from $P$ with the lasso. This would then result in a fit with at most $N-1$ non-zero coefficients. 
We will use the linear span of these basis functions as approximation of the tangent space in the calculation of TMLE-updates.

 Let $r(N)\equiv \pl D^{(1)}_P-\tilde{D}^{(1)}_P\pl_P$ be the approximation error due to using a  finite dimensional lasso fit based on $N$ observations from $P$.
 By the known rate of convergence of HAL-MLE we have $r(N)=O_P(N^{-1/3}(\log N)^{d/2})$.
Consider now a TMLE-update $P_n^*$ based on using this $\tilde{D}^{(1)}_P$ approximation, so that  it solves $P_n \tilde{D}^{(1)}_{P_n^*})=0$.
We have
\[
\begin{array}{l}
P_0 \tilde{D}^{(1)}_{P_n^*}=(P_0-P_n^*)(\tilde{D}^{(1)}_{P_n^*}-D^{(1)}_{P_n^*})+P_0 D^{(1)}_{P_n^*}\\
=(P_0-P_n^*)(\tilde{D}^{(1)}_{P_n^*}-D^{(1)}_{P_n^*})+\Psi(P_0)-\Psi(P_n^*)+R^{(1)}(P_n^*,P_0).
\end{array}
\]
Thus, 
\[
\Psi(P_n^*)-\Psi(P_0)=-P_0 \tilde{D}^{(1)}_{P_n^*}+R^{(1)}(P_n^*,P_0)-(P_n^*-P_0)(\tilde{D}^{(1)}_{P_n^*}-D^{(1)}_{P_n^*}).\]
Combining with $P_n \tilde{D}^{(1)}_{P_n^*}=0$ yields
\[
\Psi(P_n^*)-\Psi(P_0)=(P_n-P_0)\tilde{D}^{(1)}_{P_n^*} +R^{(1)}(P_n^*,P_0)-(P_n^*-P_0)(\tilde{D}^{(1)}_{P_n^*}-D^{(1)}_{P_n^*}).\]
Therefore, we conclude that the analysis of a first order TMLE based on using $\tilde{D}^{(1)}_P$ instead of $D^{(1)}_P$ generates an extra term given by  the last term. This term can be bounded with Cauchy-Schwarz inequality by $\pl p_n^*-p_0\pl_{\mu}\pl \tilde{D}^{(1)}_{P_n^*}-D^{(1)}_{P_n^*}\pl_{\mu}$. Assuming $p_n^*$ is a TMLE using as initial estimator an HAL-MLE or a super-learner including an HAL-MLE, and the above bound, implies that this term is generally
$O_P(n^{-1/3}(\log n)^{d/2})N^{-1/3}(log N)^{d/2})$. Therefore, theoretically it is fine to select $N=n$, but in practice one might want to select $N$ larger than $n$ to control possible large constants.


 \section{Sequential computation of higher order canonical gradients in terms of linear combination of basis functions spanning tangent space}\label{section9b}
 
 \subsection{Sequential computation of the higher order canonical gradients}
 In the previous section we  showed the following lemma.
 \begin{lemma}
For a given target parameter $\Psi:{\cal M}\rightarrow\openr$ with canonical gradient $D_P^*$;  tangent space $T_P$ that is approximated by the linear span  of basis functions $\{\phi_{j,P}:j\}$; 
$\Sigma_P=(P \phi_{k,P}\phi_{l,P}: k,l)$; paths $\{P_{\delta,\phi_{j,P}}:\delta\}\subset{\cal M}$ with score $\phi_{j,P}$ at $\delta_0=0$;  we have
\begin{eqnarray*}
\beta_P&=&\Sigma_P^{-1}(\dot{\psi}_P)\\
D^*_P&\approx& \sum_j \beta_P(j)\phi_{j,P}\\
\dot{\psi}_P(j)&=& \frac{d}{d\delta_0}\Psi(P_{\delta_0,\phi_{j,P}})\\
\dot{\psi}_P(j)&=& E_P D_P \phi_{j,P}\mbox{ if $D_P$ is gradient at $P$}.
\end{eqnarray*}
\end{lemma}
In the previous section we used this result to compute $D^{(1)}_P$. However, we can also apply it to sequentially compute
$D^{(k+1)}_{n,P}$ of target parameter $\Psi_n^{(k)}:{\cal M}\rightarrow\openr$, $k=0,\ldots,K-1$.
Suppose we already computed $D^{(k)}_{n,P}$; $\tilde{P}_n^{(k)}(P)$;  and $D^{(k)}_{n,\tilde{P}_n^{(k)}(P)}$ for $k=1$. We want to compute $D^{(k+1)}_{n,P}$. Application of the above lemma to $\Psi_n^{(k)}$ shows that $D_{n,P}^{(k+1)}\approx \sum_j \beta_P^{(k+1)}(j)\phi_{j,P}$ with $\beta_P^{(k+1)}=\Sigma_P^{-1}(\dot{\psi}^{(k+1)}_P)$ and, for a small $\delta\approx 0$,  
\begin{eqnarray*}
\dot{\psi}^{(k+1)}_P(j_0)&=&\frac{d}{d\delta_0}\Psi_n^{(k)}(P_{\delta_0,j})\\
&=&\frac{d}{d\delta_0}\Psi_n^{(k-1)}(\tilde{P}_n^{(k)}(P_{\delta_0,j_0}))\\
&=&\tilde{P}_n^{(k)} D^{(k)}_{n,\tilde{P}_n^{(k)}(P)}\frac{\tilde{p}_n^{(k)}(p_{\delta,j_0})-\tilde{p}_n^{(k)}(p)}{\delta\tilde{p}_n^{(k)}(p)} .
\end{eqnarray*}
Consider the case that $\tilde{p}_n^{(k)}(p)= (1+\epsilon_n^{(k)}(p)D^{(k)}_{n,P})p$, so that
$\tilde{p}_n^{(k)}(p_{\delta,j_0})$ requires knowing  $D^{(k)}_{n,P_{\delta,j_0})}$. This would then require computing $D^{(k)}_{n,P_{\delta,j_0}}$ for each $j_0$, and thus require inverting $\Sigma_{P_{\delta,j_0})}$ across all $j$. We want to avoid such massive computations. We have $D^{(k)}_{n,P}=\sum_j \beta^{(k)}_P(j)\phi_{j,P}$. Instead of using 
$D^{(k)}_{n,P_{\delta,j_0}}$ we can use $\tilde{D}^{(k)}_{n,P_{\delta,j_0}}=\sum_j \beta^{(k)}_{P}(j)\phi_{j,P_{\delta_0,j_0} } $ in the definition of the TMLE-update $\tilde{p}_n^{(k)}(p_{\delta,j_0})$.
That is, we use the same $\beta^{(k)}_{P}$  in $D^{(k)}_{n,P}$, but we do update the basis functions.
It is straightforward to  show that these TMLE-updates  using $\tilde{D}^{(k)}_{n,P_{\delta,j_0}}$ minus the correct TMLE updates using $D^{(k)}_{n,P_{\delta,j_0}}$ results in a second order approximation error for the resulting TMLE-path, assuming that $\epsilon_n^{(k)}(P)\approx 0$, as it would be if $P$ represents a consistent estimator of $P_0$. In addition, for an iterative TMLE in which $P$ is iteratively replaced by the previous TMLE, there is no-approximation error.
In this manner, the computation of all the $j_0$-specific TMLE-updates $\tilde{p}_n^{(k)}(p_{\delta,j_0})$ is very doable, and thereby computation of $\dot{\psi}^{(k)}$ is computationally feasible as well.
 
 To summarize: having computed $D^{(k)}_{n,P}$ and $D^{(k)}_{\tilde{p}_n^{(k)}(P)}$, the computation of $D^{(k+1)}_{n,P}$ requires
 \begin{itemize} 
 \item calculation of $\tilde{D}^{(k)}_{n,P_{\delta,j_0}}=\sum_{j=1} \beta^{(k)}_P(j)\phi_{j,P_{\delta,j_0}}$ across all $j_0$, beyond $D^{(k)}_{n,P}=\sum_j \beta^{(k)}_P(j)\phi_{j,P}$.
 \item calculation of  corresponding TMLE updates $\tilde{p}_n^{(k)}(p_{\delta,j_0})$ using $\tilde{D}^{(k)}_{n,P_{\delta,j_0}}$ as score instead of $D^{(k)}_{n,P_{\delta,j_0}}$, beyond $\tilde{p}_n^{(k)}(p)$ based on $D^{(k)}_{n,P}$.
 \item Calculation of $D^{(k)}_{n,\tilde{P}_n^{(k)}(P)}$.
 \item Taking large sample from $\tilde{P}_n^{(k)}(P)$ and compute \newline $\dot{\psi}^{(k)}(j_0)=
 \tilde{P}_n^{(k)}(P) D^{(k)}_{n,\tilde{P}_n^{(k)}(P)}  (\tilde{p}_n^{(k)}(p_{\delta,j_0})-\tilde{p}_n^{(k)}(p))/\tilde{p}_n^{(k)}(p)$. 
 \item Computing $\beta^{(k+1)}_P=\Sigma_P^{-1}(\dot{\psi}^{(k)})$ and then $D^{(k+1)}_{n,P}\approx \sum_j \beta^{(k+1)}_P(j)\phi_{j,P}$.
 \end{itemize}
 After having computed $D^{(k+1)}_{n,P}$, one can determine $\tilde{P}_n^{(k+1)}(P)$. We now redo the above computation at $\tilde{P}_n^{(k+1)}(P)$ instead of $P$ and obtain $D^{(k+1)}_{n,\tilde{P}_n^{(k+1)}(P)}$. At this point we went through one full cycle.
 We now set $k=k+1$ and redo the above to compute $D^{(k+2)}_{n,P}$; $\tilde{P}_n^{(k+2)}(P)$; and $D^{(k+2)}_{n,\tilde{P}_n^{(k+2)}(P)}$. 
 Iterating this results in the sequential computation of $(D^{(1)}_P,\tilde{P}_n^{(1)}(P),D^{(1)}_{\tilde{P}_n^{(1)}(P)})$, $\ldots$, $(D^{(K)}_{n,P},\tilde{P}_n^{(K)}(P),D^{(K)}_{n,\tilde{P}_n^{(K)}(P)})$.
 
 \subsection{Recursive programming to compute the $K$-th order TMLE}

Suppose we have a function $P\rightarrow HOTMLE_k(P)$ that maps $P$ into 
$\tilde{P}_n^{(k)}(P)$; $\tilde{P}_n^{(k-1)}\tilde{P}_n^{(k)}(P)$; $\ldots$; $\tilde{P}_n^{(1)}\ldots\tilde{P}_n^{(k)}(P)$, and, accordingly, the scores defining these TMLE-updates given by $D^{(k)}_{n,P}$; $D^{(k-1)}_{\tilde{P}_n^{(k)}(P)}$; $\ldots$;
$D^{(1)}_{\tilde{P}_n^{(2)}\ldots\tilde{P}_n^{(k)}(P)}$, respectively.


We now want to compute the function $P\rightarrow HOTMLE_{k+1}(HOTMLE_k(),P)$ that utilizes the previous function $HOTMLE_k()$.
We know that $HOTMLE_k(\tilde{P}_n^{(k+1)}(P))$ gives the $k+1$-th order TMLE but we do not know yet $\tilde{P}_n^{(k+1)}(P)$ which relies on $D^{(k+1)}_{n,P}$. So we need to first utilize $HOTMLE_k()$ to compute $D^{(k+1)}_{n,P}$; then compute $\tilde{P}_n^{(k+1)}(P)$, and finally run $HOTMLE_k(\tilde{P}_n^{(k+1)}(P))$.

Indeed, we can program $HOTMLE_{k+1}(HOTMLE_k(),P)$ as follows:
\begin{itemize}
\item We first need to calculate $\dot{\psi}^{(k+1)}_P$.
\begin{itemize}
\item Run $HOTMLE_k(P)$. 
\item This  yields $D^{(k)}_{n,P}$.
From this we can obtain $\tilde{D}^{(k)}_{n,P_{\delta,j_0}}$, and corresponding TMLE-updates $\tilde{p}_n^{(k)}(p_{\delta,j_0})$  and thereby its score
$A^{(k)}_{n,P}(\phi_{j_0,P})$ at $\delta=0$, for all $j_0$.
\item $HOTMLE_k(P)$ also yields $\tilde{P}_n^{(k)}(P)$.
\item We now run $HOTMLE_k(\tilde{P}_n^{(k)}(P))$ giving $D^{(k)}_{n,\tilde{P}_n^{(k)}(P)}$.
\item Take a large sample from $\tilde{P}_n^{(k)}(P)$ and compute, for each $j_0$,
$\dot{\psi}^{(k+1)}_P(j_0)=\tilde{P}_n^{(k)}(P) D^{(k)}_{n,\tilde{P}_n^{(k)}(P)}A^{(k)}_{n,P}(\phi_{j_0,P})$.
\end{itemize}
\item  Compute $\beta^{(k+1)}_P=\Sigma_P^{-1}(\dot{\psi}^{(k+1)}_P)$.
\item This yields $D^{(k+1)}_{n,P}=\sum_j \beta_P^{(k+1)}(j)\phi_{j,P}$ and allows us to determine the TMLE $\tilde{P}_n^{(k+1)}(P)$ based on $D^{(k+1)}_{n,P}$.
\item Run $HOTMLE_k(\tilde{P}_n^{(k+1)}(P))$ giving the $k+1$-th order TMLE, which is nothing else than $k$-th order TMLE at $\tilde{P}_n^{(k+1)}(P)$.
Specifically, the output of $HOTMLE_k(\tilde{P}_n^{(k+1)}(P)$ yields
$\tilde{P}_n^{(k)}\tilde{P}_n^{(k+1)}(P)$; $\tilde{P}_n^{(k-1)}\tilde{P}_n^{(k)}\tilde{P}_n^{(k+1)}(P)$; $\ldots$; $\tilde{P}_n^{(1)}\ldots\tilde{P}_n^{(k)}\tilde{P}_n^{(k+1)}(P)$, and, accordingly, $D^{(k)}_{n,\tilde{P}_n^{(k+1)}(P)}$; $D^{(k-1)}_{\tilde{P}_n^{(k)}\tilde{P}_n^{(k+1)}(P)}$; $\ldots$;
$D^{(1)}_{\tilde{P}_n^{(2)}\ldots\tilde{P}_n^{(k)}\tilde{P}_n^{(k+1)}(P)}$. 
Adding to this output $\tilde{P}_n^{(k+1)}(P)$ and $D^{(k+1)}_{n,P}$ yields now the complete output as desired for $HOTMLE_{k+1}(P)$, analogue to output of $HOTMLE_k(P)$.
\end{itemize}

The above writes the function $HOTMLE_{k+1}(HOTMLE_k(),P)$, and we already have $HOTMLE_1(P)$.
Thus, provides the program for computing the $K$-th order TMLE $HOTMLE_K(P)$:  $HOTMLE_1(P)$; $HOTMLE_2(P)=HOTMLE_2(HOTMLE_1(),P)$; $HOTMLE_3(P)=HOTMLE_3(HOTMLE_2(),P)$; $\ldots$; and finally $HOTMLE_K(P)=HOTMLE_{K}(HOTMLE_{K-1}(),P)$.

\section{Iterative HAL-regularized higher order TMLE, iteratively targeting the HAL-MLE}\label{sectionthotmle}
Let's first present the proposed iterative higher order TMLE algorithm that also iteratively targets $\tilde{P}_n$. 
Subsequently, we show that our exact expansions for the higher order TMLE also applies to this estimator, and we motivate the targeting step of the algorithm.

For simplicity, we provide the algorithm for the second order TMLE.
\begin{itemize}
\item Let $\tilde{P}^{(1)}_n(P)$ be the first order TMLE-update of initial $P$ solving ${P}_n D^{(1)}_{\tilde{P}^{(1)}_n(P)}=0$. Here we decided to use a regular TMLE under $P_n$. Due to the iterative targeting of the HAL-MLE, our final estimator will be indistinguishable from a second order TMLE using $\tilde{P}^{(1)}_n(P)$ and $\tilde{P}_n^{(2)}(P)$ based on the targeted $\tilde{P}_n$. The advantage of this is that we now have a single empirical log-likelihood criterion that will increase at each step, thereby providing us with a guarantee that this iterative second order TMLE algorithm will indeed converge.
\item Let $\tilde{P}^{(2)}_n(P,\tilde{P}_n)$ be the second order TMLE update which uses as initial $P$; it uses $\tilde{P}_n$ as HAL-MLE; and it is targeted to solve $P_n D^{(2)}_{\tilde{P}_n,\tilde{P}_n^{(2)}(P,\tilde{P}_n)}=0$. 

\item We start with the initial estimator $P_n^0$.
\item We target $\tilde{P}_n$ w.r.t. $P_n
^0$ so that it solves
$(\tilde{P}_n-P_n)D^{(1)}_{P_n^0}=0$ and $(\tilde{P}_n-P_n)D^{(2)}_{\tilde{P}_n,P_n^0}=0$.
Let's denote this targeted HAL-MLE with $\tilde{P}_n^{*,0}\equiv \tilde{P}_n(P_n^0)$.
\item Given this $P_n^0$ and $\tilde{P}_n^{*,0}$, we now compute the second order TMLE.
Thus, firstly, we compute the second TMLE-update $\tilde{P}^{(2)}_n(P_n^0,\tilde{P}_n^{*,0})$ at the current T-HAL-MLE $\tilde{P}_n^*(P_n^0)$ and initial estimator $P_n^0$, so that we solve $P_n D^{(2)}_{\tilde{P}_n^{*,0},\tilde{P}_n^{(2)}(P_n^0,\tilde{P}_n^{*,0})}=0$.
Then, we compute $\tilde{P}^{(1)}_n(\tilde{P}^{(2)}_n(P_n^0,\tilde{P}_n^{*,0}))$
so that we solve \[
P_n D^{(1)}_{\tilde{P}^{(1)}_n(\tilde{P}^{(2)}_n(P_n^0,\tilde{P}_n^{*,0} ))  } =0.\]

\item So we have now computed one round of second order TMLE starting with initial $P_n^0$ and using $\tilde{P}_n^{*,0}$ as HAL-MLE, and we have increased the empirical log-likelihood at both the second and first order TMLE update. Let $P^{2,(1),*}_n(P_n^0,\tilde{P}_n^{*,0} )$ be the resulting second order TMLE of $P_0$.

\item If $(\tilde{P}_n-P_n)D^{(1)}_{\tilde{P}^{(1)}_n(\tilde{P}^{(2)}_n(P_n^0,\tilde{P}_n^{*,0}  ) )}$ or
$(\tilde{P}_n-P_n)D^{(2)}_{\tilde{P}_n^{*,0},\tilde{P}^{(2)}_n(P_n^0,\tilde{P}_n^{*,0}) }$
is not smaller than $\sigma_{1n}/(n^{1/2}\log n)$, then we proceed with next steps, otherwise, we are done.

\item Let $P_n^1=P^{2,(1),*}_n(P_n^0,\tilde{P}_n^{*,0}) $ be the new initial estimator for the second round of second order TMLE. Let $\tilde{P}_n^{*,1}=\tilde{P}_n^*(P_n^1)$ be the new targeted HAL-MLE using as off-set the $\tilde{P}_n^{*,0}$ and now targeted to solve the score equations 
$(\tilde{P}_n^{*,1}-P_n)D^{(1)}_{P_n^1}=0$ and $(\tilde{P}_n^{*,1}-P_n)D^{(2)}_{\tilde{P}_n^{*,0},P_n^0}=0$,
at the  current ''initial'' estimator $P_n^1$.   We now compute the second round second order TMLE $P^{2,(1),*}(P_n^1,\tilde{P}_n^{*,1})$ with this initial estimator $P_n^1 $ and targeted HAL-MLE $\tilde{P}_n^{*,1}$.

\item We iterate this process  of computing the second order TMLE given a current initial estimator and current targeted HAL-MLE: Starting at $m=1$, compute
$P_n^{m+1}=P^{2,(1),*}_n(P_n^m,\tilde{P}_n^{*,m})$ and $\tilde{P}_n^{*,m+1}=\tilde{P}_n^*(P_n^{m+1})$; $m=m+1$, and repeat this  till  integer value $m^*$ at which $(\tilde{P}_n^{*,m^*}-P_n)D^{(1)}_{P^{2,(1),*}_n(P_n^{m^*},\tilde{P}_n^{*,m^*} ) } \approx 0$ and $(\tilde{P}_n^{*,m*}-P_n)D^{(2)}_{\tilde{P}_n^{m^*},\tilde{P}_n^{(2)}(P_n^{m^*},\tilde{P}_n^{*,m^*})}\approx 0$.
 
\item At the final step, the resulting second order TMLE $P^{2,(1),*}_n$ and its targeted $\tilde{P}_n^*=\tilde{P}_n^{*,m*}$  have solved
$P_n D^{(1)}_{P^{2,(1),*}_n}=0$, $\tilde{P}_n^* D^{(1)}_{P^{2,(1),*}_n}=0$,
$P_n D^{(2)}_{\tilde{P}_n^*,P^{2,(1),*}_n}=0$, and $\tilde{P}_n D^{(2)}_{\tilde{P}_n^*,P^{2,(1),*}_n}=0$, either all exactly or up till desired precision.
\end{itemize}
Even though we use $P_n$ in the definition of $\tilde{P}_n^{(1)}(P)$ and $\tilde{P}_n^{(2)}(P,\tilde{P}_n)$, in the end, this iterative second order TMLE can be represented as a second order TMLE that uses $\tilde{P}_n^*$ as HAL-MLE in these two TMLE-updates, initial estimator $P^{2,(1),*}_n$, and solves the HAL-MLE score equations
$\tilde{P}_n^* D^{(1)}_{P^{2,(1),*}_n}=\tilde{P}_n^* D^{(2)}_{\tilde{P}_n^*,P^{2,(1),*}_n}=0$, as if $\tilde{P}_n^{(1)}(P)$ and $\tilde{P}_n^{(2)}(P)$ were using $\tilde{P}_n^*$ instead of $P_n$ in its definition of the TMLE-update $\epsilon_n^{(j)}(P)$, $j=1,2$.  Therefore this estimator can be analyzed accordingly: see Appendix \ref{AppendixH}.

 Therefore our exact expansion of Theorem \ref{theorem2} applies stating that $\Psi(P^{2,(1),*)}_n)-\Psi(P_0)$ equals $\sum_{j=1}^2 \left\{ (\tilde{P}_n^*-P_0)D^{(j)}_{\tilde{P}_n^*,\tilde{P}_n^{(j)}(P_0)}+R^{(j)}_{\tilde{P}_n^*}(\tilde{P}_n^{(j)}(P_0),P_0)\right\}$ plus the difference of third order remainders $R^{(2)}_{\tilde{P}_n^*}(P^{2,(1),*}_n,\tilde{P}_n^*)$ and $R^{(2)}_{\tilde{P}_n^*}(\tilde{P}_n^{(2)}(P_0),\tilde{P}_n^*)$.
In addition, the algorithm  guaranteed that $P_n D^{(j)}_{\tilde{P}_n^*,P^{2,(1),*}_n}=0$ for $j=1,2$, thereby making the term $(\tilde{P}_n^*-P_n)D^{(j)}_{\tilde{P}_n^*,\tilde{P}_n^{(j)}(P_0)}$, as needed for inference based on the leading empirical means, will be reduced, as was discussed in Section \ref{section7}, and, in more detail below. 

{\bf Universal least favorable path analogue:}
We can define an analogue algorithm based on tracking the local least favorable paths locally instead with small steps, thereby making this algorithm use minimal fitting to achieve solving these score equations while increasing the log-likelihood at each step: i..e, replace each full MLE update step $\epsilon_n^{(j)}()$ in the above algorithm by a small move $\delta$ along the local least favorable path through the current estimate in the direction of increasing the log-likelihood or solving the equation (i.e., in direction of $\epsilon_n^{(j)}$ at that initial).

\subsection{Understanding role of $\tilde{P}_n$ in the HAL-regularized higher order TMLE}
Consider a $k$-th order TMLE which solves $\tilde{P}_n D^{(j)}_{\tilde{P}_n^{(j)}(P)}=0$ for all $j=1,\ldots,k$. This $k$-th order TMLE satisfies that $\Psi(\tilde{P}_n^{(1)}\ldots\tilde{P}_n^{(k)}(P_n^0))-\Psi(P_0)$ equals $\sum_{j=1}^k (\tilde{P}_n-P_0)D^{(j)}_{\tilde{P}_n^{(j)}(P_0)}+R^{(j)}_n(\tilde{P}_n^{(j)}(P_0),P_0)$ plus a final $k+1$-th order remainder $R^{(k)}(\tilde{P}_n^{(k)}(P_n^0),\tilde{P}_n)-R^{(k)}(\tilde{P}_n^{(k)}(P_0),\tilde{P}_n)$. Let's consider the case that we define $\epsilon_{\tilde{P}_n}^{(j)}(P)$ as the solution of $\tilde{P}_n D^{(j)}_{\tilde{P}^{(j)}(P,\epsilon)}=0$.
The remaining challenge for the choice of HAL-MLE $\tilde{P}_n$ is to enforce that $(\tilde{P}_n-P_n)D^{(j)}_{\tilde{P}_n^{(j)}(P_0)}$ is negligible, so that statistical inference is based on mean zero empirical  means. Since the $j=1$-term generally dominates the other terms, let's focus only on the leading term $(\tilde{P}_n-P_n)D^{(1)}_{\tilde{P}_n^{(1)}(P_0)}$. 
In this article we have mentioned undersmoothing $\tilde{P}_n$ as a possible strategy to control these terms for the HAL-regularized higher order TMLE. Here we argue for additional (or only) targeting of $\tilde{P}_n$ and thereby also motivate the above iterative HAL-regularized  higher order TMLE algorithm.


We consider the case that we define $\epsilon_{\tilde{P}_n}^{(1)}(P)$ as the solution of $\tilde{P}_n D^{(1)}_{\tilde{P}^{(1)}(P,\epsilon)}=0$. The goal of the HAL-MLE $\tilde{P}_n$ is to make $\epsilon_{\tilde{P}_n}^{(1)}(P_0)$  a good plug-in estimator of $\epsilon_{P_0}^{(1)}(P_0)$.
The parameter $\epsilon_{P_0}^{(1)}(P)$ of $P_0$ for a given $P$ is an easy to estimate parameter in the sense that $\epsilon_{P_n}^{(1)}(P)$ is a simple MLE or score-equation solution for a one-dimensional parametric model (correctly specified if $P=P_0$ and otherwise a slightly misspecified parametric model), so that it is a nice asymptotically linear estimator with good robust finite sample performance (not affected by curse of dimensionality). However, the challenge is that we can only estimate it with a substitution estimator $\epsilon_{\tilde{P}_n}^{(1)}(P)$ with $\tilde{P}_n\in {\cal M}$. We favor using an HAL-MLE since it is an MLE itself  and can be tuned to solve any score equation and thereby makes it easier to make it asymptotically equivalent with $\epsilon_{P_n}^{(1)}(P)$.
On the other hand, we still have the option to make $\tilde{P}_n$ itself a TMLE that uses as initial estimator $\tilde{P}_n^{C}$ for some $C$ (e.g, $C=C_{n,cv}$).  

{\bf Optimizing performance of $\epsilon_{\tilde{P}_n}^{(1)}(P_0)$ so that $P_n D^{(1)}_{\tilde{P}_n^{(1)}(P_0)}\approx 0$:}
We know that 
$\epsilon_{P_n}^{(1)}(P_0)$ is an estimator for which $P_n D^{(1)}_{\tilde{P}^{(1)}_n(P_0)}=0$.
Therefore, the task at hand is that we want $\epsilon_{\tilde{P}_n}^{(1)}(P_0)$ to be an excellent substitution estimator of $\epsilon_0(P)\equiv \epsilon_{P}^{(1)}(P_0)$ at $P=P_0$, preferably asymptotically equivalent with $\epsilon_0(P_n)=\epsilon_{P_n}^{(1)}(P_0)$.
This strongly suggest that we should make $\tilde{P}_n$ a TMLE targeting this target parameter $\epsilon_0(P_0)$. So let's proceed with computing its canonical gradient and corresponding TMLE.

The pathwise derivative $\frac{d}{d\delta_0}\epsilon_0(P_{\delta_0,h})$ at $P$  is given by:
\[
\left\{ -\frac{d}{d\epsilon_0(P)}P D^{(1)}_{\tilde{P}^{(1)}(P_0,\epsilon_0(P))}\right\}^{-1} P h D^{(1)}_{\tilde{P}^{(1)}(P_0,\epsilon_0(P))}.\]
Let $I_0(P)=-\frac{d}{d\epsilon_0(P)}P D^{(1)}_{\tilde{P}^{(1)}(P_0,\epsilon_0(P))}$.

In particular, the canonical gradient at $P=P_0$ of $\epsilon_0(P)$ is given by:
\[
D_{\epsilon_0,P_0}^*=\{ P_0 \{D^{(1)}_{P_0}\}^2\}^{-1} D^{(1)}_{P_0}.\]
Even at a  $P\approx P_0$, the canonical gradient of $\epsilon_0(P)$ at $P$ is  well approximated by $D^{(1)}_P$.
In a nonparametric model, we would have that
\[
D_{\epsilon_0,P}^*=I_0(P)^{-1}\left\{ D^{(1)}_{\tilde{P}^{(1)}(P_0,\epsilon_0(P))}-P D^{(1)}_{\tilde{P}^{(1)}(P_0,\epsilon_0(P))}\right\}.\]

Thus we want $\epsilon_0(\tilde{P}_n)$ to be an asymptotically linear estimator with influence curve a constant times $D^{(1)}_{P_0}$. Thus, we could select $\tilde{P}_n$ as a first order TMLE targeting $\Psi(P_0)$, based on initial estimator $\tilde{P}_n^{C}$, where $C$ could still be a tuning parameter. Thus, we could define $\tilde{P}_n^*=\tilde{P}^{(1)}(\tilde{P}_n^C,\epsilon_{P_n}^{(1)}(\tilde{P}_n))$ at, for example, $C=C_{n,cv}$. Then,
$\epsilon_0(\tilde{P}_n^*)-\epsilon_0(P_0)=(P_n-P_0)D^{(1)}_{\tilde{P}_n^*}+R_n$ for a second order remainder $R_n=O_P(n^{-2/3}(\log n)^d)$, so that it would be asymptotically efficient and equivalent with $\epsilon_0(P_n)$. However, instead of targeting it so that $P_n D^{(1)}_{\tilde{P}_n^*}=0$, one might as well target $\tilde{P}_n$ so that $P_n D^{(1)}_{P^{k,(1),*}_n}=0$ for the second order TMLE $P^{k,(1),*}$ using $\tilde{P}_n$ in its TMLE-update steps, which is higher order efficient for $\Psi(P_0)$.  This provides the direct motivation for the   iterative higher order TMLE that targets $\tilde{P}_n$ to indeed solve $P_n D^{(1)}_{P^{k,(1),*}_n}=0$ (and for the higher order equations as well).

\subsection{Targeting makes $(\tilde{P}_n^*-P_n)D^{(1)}_{P_0}$ a third order difference (just as with the empirical higher order TMLE).}
We now want to more formally explain how the targeting of $\tilde{P}_n$ assists in reducing
$(\tilde{P}_n^*-P_n)D^{(1)}_{P_0}$. Firstly, due to $\tilde{P}_n D^{(1)}_{P_n^{(k),(1),*}}=P_n D^{(1)}_{P_n^{k,(1),*}}=0$ we have $(\tilde{P}_n^*-P_n)D^{(1)}_{P_0}=(\tilde{P}_n^*-P_n)\{D^{(1)}_{P_0}-D^{(1)}_{P_n^{k,(1),*}} \}$.
Suppose that the targeting $\tilde{P}_n^*$ is carried out by a constraint HAL-MLE with the additional constraint (beyond $L_1$-norm) that it solves the desired efficient scores. Such a targeted HAL-MLE is an HAL-MLE itself so that it still solves a large class of score equations: we proposed and analyzed this estimator in a technical report (van der Laan, 2020). 
Let $f_n$ be a best approximation of $D^{(1)}_{P_0}-D^{(1)}_{P_n^{k,(1),*}}$
among the linear span of scores $\{S_{j,\tilde{P}_n^*}:j\in {\cal J}\}$ solved by the HAL-MLE $\tilde{P}_n^*$ in the sense that $P_n S_{j,\tilde{P}_n^*}=0$, $j\in {\cal J}$.
Then, we have $P_n f_n=0$ and also $\tilde{P}_n^* f_n=0$. 
Thus, we have 
\[
(\tilde{P}_n^*-P_n)D^{(1)}_{P_0}=(\tilde{P}_n^*-P_n)\{ D^{(1)}_{P_0}-D^{(1)}_{P_n^{k,(1),*}}-f_n\}.\]
The function $D^{(1)}_{P_0}-D^{(1)}_{P_n^{k,(1),*}}$ represents a first order difference, while 
subtracting its best linear approximation makes it a second order difference. Therefore, this term is now a third order difference. This representation also  makes  clear the two different properties of $\tilde{P}_n^*$ that   controls this term: 1) being an HAL-MLE allowing to use this best linear approximation $f_n$; 2) being targeted allowing us to subtract $D^{(1)}_{P_n^{k,(1),*}}$. In addition, it  is a nice third order difference that can be further controlled by undersmoothing $\tilde{P}_n^*$. Therefore, we can conclude that the bias due to using a smooth $\tilde{P}_n^*$ instead of $P_n$ in the computation of TMLE-updates is controlled at a level that makes it small relative to the gained third order term in the exact expansion for the higher order TMLE.

\section{Cross-validated higher order TMLE}\label{sectioncvhotmle}
For simplicity, let's focus on the second order TMLE: the generalization to higher order TMLE is immediate.  Consider $V$-fold sample splitting and let
$P_{n,v}^1$ and $P_{n,v}$ be the empirical measure of the validation sample and training sample, respectively, for the $v$-th sample split, $v=1,\ldots,V$. Typically, $V$ is quite large such as $V=10$, so that the training sample is much larger than the validation sample.
As some background, the regular CV-TMLE involves computing the initial estimator on the training sample $P_{n,v}$; doing the TMLE step on the validation sample $P_{n,v}^1$, or in a pooled manner, resulting in a targeted $P_{n,v}^*$; and defining the CV-TMLE as the average across $v$ of the $v$-specific plug-in TMLE $\Psi(P_{n,v}^*)$. The important motivation for CV-TMLE is that it robustifies the targeting step by protecting agains overfitting the initial estimator, and, as a theoretical analogue, it completely avoids reliance on the Donsker class condition in its asymptotic analysis.  The straightforward generalization of CV-TMLE to the higher order TMLE is problematic since it would involve a regularized TMLE-update step of an initial estimator on the training sample w.r.t. the HAL-MLE only fitted on the validation sample (e.g. $1/10$-th of the total sample size), pooled across the $V$ splits. Using an HAL-MLE on a much smaller sample size will hurt the finite sample performance. Our solution below circumvents this problem by cross-validated targeting the HAL-MLE instead, while keeping the regularized TMLE-update steps w.r.t. (targeted)  HAL-MLE on the full training sample.

Let $\tilde{P}_{n,v}$  be an HAL-MLE and $P_{n,v}^0$ be an initial estimator, both based on the training sample $P_{n,v}^0$. Let $\tilde{P}_{n,v}^{(1)}(P)$ and $\tilde{P}_{n,v}^{(2)}(P)$ be the  two TMLE-updates solving $\tilde{P}_{n,v} D^{(1)}_{\tilde{P}_{n,v}^{(1)}(P)}=\tilde{P}_{n,v}D^{(2)}_{\tilde{P}_{n,v},\tilde{P}_{n,v}^{(2)}(P_{n,v}^0)}=0$, $v=1,\ldots,V$.  Then, $\tilde{P}_{n,v}^{(1)}\tilde{P}_{n,v}^{(2)}(P_{n,v}^0)$ is the regular second order TMLE, but applied to the training sample $P_{n,v}$. Therefore, we can apply our exact expansion for the second order TMLE:
\begin{eqnarray}
\Psi(\tilde{P}_{n,v}^{(1)}\tilde{P}_{n,v}^{(2)}(P_{n,v}^0))-\Psi(P_0)&=&
(\tilde{P}_{n,v}-P_0)D^{(1)}_{\tilde{P}_{n,v}^{(1)}(P_0)}+R^{(1)}(\tilde{P}_{n,v}^{(1)}(P_0),P_0)\nonumber \\
&&\hspace*{-1cm} +(\tilde{P}_{n,v}-P_0)D^{(2)}_{\tilde{P}_{n,v},\tilde{P}_{n,v}^{(2)}(P_0) }+R^{(2)}(\tilde{P}_{n,v}^{(2)}(P_{n,v}^0),P_0)\nonumber \\
&&\hspace*{-1cm} +R^{(2)}(\tilde{P}_{n,v}^{(2)}(P_{n,v}^0),\tilde{P}_{n,v})-R^{(2)}(\tilde{P}_{n,v}^{(2)}(P_0),\tilde{P}_{n,v}),\label{exactexphotmle}
\end{eqnarray}
where we also know that the last line represents a third order difference. 

Earlier we presented an iterative second order TMLE that involves targeting of the HAL-MLE $\tilde{P}_{n,v}$ w.r.t. a given $P$, so that its targeted version $\tilde{P}_{n,v}^*$ satisfies  $(\tilde{P}_{n,v}^*-P_{n,v})D^{(1)}_{\tilde{P}_{n,v}^{(1)}P} =(\tilde{P}_{n,v}^*-P_{n,v})D^{(2)}_{\tilde{P}_{n,v}^*,P}=0$. 
We first target $\tilde{P}_{n,v}$ w.r.t. $P=\tilde{P}_{n,v}^{(2)}(P_{n,v}^0)$ resulting in a first targeted version $\tilde{P}_{n,v}^{1}$;  then at the next iteration, after having computed the second order TMLE $P_{n,v}^{0,1}=\tilde{P}_{n,v}^{(1)}\tilde{P}_{n,v}^{(2)}(P_{n,v}^0)$ under this $\tilde{P}_{n,v}^{*,1}$, we target $\tilde{P}_{n,v}^{1}$ w.r.t. $P=\tilde{P}_{n,v}^{(2)}(P_{n,v}^{0,1})$ (updated initial estimator, and using $\tilde{P}_{n,v}^1$ in TMLE update step), and we iterate this process of targeting $\tilde{P}_{n,v}$ w.r.t. current initial estimator and computing the corresponding second order TMLE (next initial estimator) till convergence.
At the final step,  we  have a final initial estimator $P_{n,v}^{0,*}$ and final targeted HAL-MLE $\tilde{P}_{n,v}^*$ so that
\[
\begin{array}{l}
0=(\tilde{P}_{n,v}^*-P_{n,v})D^{(1)}_{\tilde{P}_{n,v}^{(1)}\tilde{P}_{n,v}^{(2)}P_{n,v}^{0,*}}\\
=-P_{n,v}D^{(1)}_{\tilde{P}_{n,v}^{(1)}\tilde{P}_{n,v}^{(2)}P_{n,v}^{0,*}},
\end{array}
\]
where, at second equality, we use that the second order TMLE-update solves the $\tilde{P}_{n,v}^*$ score equations.
Similarly, it follows that  $P_{n,v}D^{(2)}_{\tilde{P}_{n,v}^{(2)}P_{n,v}^{0,*}}=0$.
So, this iterative second order TMLE algorithm ends  up with a choice of initial estimator $P_{n,v}^{0,*}$ and targeted HAL-MLE $\tilde{P}_{n,v}^*$ solving the empirical score equations
$P_{n,v} D^{(1)}_{\tilde{P}_{n,v}^{(1)}\tilde{P}_{n,v}^{(2)}(P_{n,v}^0)}=0$ and $P_{n,v}D^{(2)}_{\tilde{P}_{n,v}^*,\tilde{P}_{n,v}^{(2)}(P_{n,v}^0)}=0$, beyond solving these equations w.r.t. $\tilde{P}_{n,v}^*$.   This enhanced our finite sample expansion for the second order TMLE due to now making the terms $P_n D^{(1)}_{\tilde{P}_n^{(1)}(P_0)}$ and $P_n D^{(2)}_{\tilde{P}_n,\tilde{P}_n^{(2)}(P_0)}$ small, without much need of undersmoothing the HAL-MLE. 

Here we propose minor variation of this iterative second order TMLE.
Instead of applying the targeting of the HAL-MLE to the same sample $P_{n,v}$ as the initial HAL-MLE $\tilde{P}_{n,v}$ was based upon, we  apply the same targeting of $\tilde{P}_{n,v}$ to the corresponding validation sample $P_{n,v}^1$ instead. ( As with the CV-TMLE, instead of applying the targeting of $\tilde{P}_{n,v}$  to the validation sample $P_{n,v}^1$, for each $v$ separately, one can carry out the same targeting with a single pooled criterion such as the cross-validated log-likelihood.)

As a consequence, with this variation,  the same iterative second order TMLE ends up with a targeted HAL-MLE $\tilde{P}_{n,v}^*$ and ''initial'' estimator $P_{n,v}^{0,*}$ satisfying
\begin{eqnarray}
\frac{1}{V}\sum_{v=1}^V P_{n,v}^1 D^{(1)}_{\tilde{P}_{n,v}^{(1)}\tilde{P}_{n,v}^{(2)}(P_{n,v}^{0,*})}&=&0\nonumber \\
\frac{1}{V}\sum_{v=1}^V P_{n,v}^1 D^{(2)}_{\tilde{P}_{n,v}^*,\tilde{P}_{n,v}^{(2)}(P_{n,v}^{0,*})}&=&0,
\label{tildePnvtarget}
\end{eqnarray}
where the TMLE-update steps are based on this targeted HAL-MLE $\tilde{P}_{n,v}^*$. One does not need to solve these equations exactly but just till the desired precision.
Our proposed CV-second order TMLE of $\Psi(P_0)$  is given by $\psi_n^*=1/V\sum_v \Psi(\tilde{P}_{n,v}^{(1)}\tilde{P}_{n,v}^{(2)}(P_{n,v}^{0,*}))$ corresponding with this final choice of initial estimator $P_{n,v}^{0,*}$ and HAL-MLE $\tilde{P}_{n,v}^*$. Notice that the TMLE-update steps themselves depend on $\tilde{P}_{n,v}^*$.

We now want to analyze this CV-second order TMLE. For that purpose we assume that it suffices to iterate  a fixed  number of steps, which is reasonable based on our experience. For simplicity, let's consider the case that the algorithm is iterated twice: our analysis only relies on the number of steps being bounded by a finite $K$ with probability 1.

Let's represent the first fluctuation of $\tilde{P}_{n,v}$ in its targeting step with a functional
$\tilde{P}_{n,v,\delta_1}=\tilde{P}^1(\tilde{P}_{n,v},P_{n,v}^0,\delta_1)$ to indicate that it represents a fluctuation of initial $\tilde{P}_{n,v}$; that it relies on current initial estimator $P_{n,v}^0$ and a fluctuation parameter $\delta_1$ chosen so that
\begin{eqnarray}
\frac{1}{V}\sum_{v=1}^V (\tilde{P}_{n,v,\delta_1}-P_{n,v}^1) D^{(1)}_{\tilde{P}_{n,v}^{(1)}\tilde{P}_{n,v}^{(2)}(P_{n,v}^{0})}&=&0\nonumber \\
\frac{1}{V}\sum_{v=1}^V (\tilde{P}_{n,v,\delta_1}-P_{n,v}^1) D^{(2)}_{\tilde{P}_{n,v},\tilde{P}_{n,v}^{(2)}(P_{n,v}^{0})}&=&0.
\end{eqnarray}
Let $\delta_{1n}$ be the solution and let $\tilde{P}_{n,v}^1=\tilde{P}^1(\tilde{P}_{n,v},P_{n,v}^0,\delta_{1n})$ be the targeted update of $\tilde{P}_{n,v}$.
Let $P_{n,v}^{0,1}=\tilde{P}_{n,v}^{(1)}\tilde{P}_{n,v}^{(2)}(P_{n,v}^0)$ (i.e. the second order TMLE using $\tilde{P}_{n,v}^1$ as HAL-MLE) be the new initial estimator.
Let's represent the second fluctuation of $\tilde{P}_{n,v}^1$, based on new initial estimator $P_{n,v}^{0,1}$, with the functional $\tilde{P}_{n,v,\delta_2}^1=\tilde{P}^2(\tilde{P}_{n,v}^1,P_{n,v}^{0,1},\delta_2)$, where $\delta_2$ is chosen so that
\begin{eqnarray}
\frac{1}{V}\sum_{v=1}^V (\tilde{P}_{n,v,\delta_2}^1-P_{n,v}^1) D^{(1)}_{\tilde{P}_{n,v}^{(1)}\tilde{P}_{n,v}^{(2)}(P_{n,v}^{0,1})}&=&0\nonumber \\
\frac{1}{V}\sum_{v=1}^V (\tilde{P}_{n,v,\delta_2}^1-P_{n,v}^1) D^{(2)}_{\tilde{P}_{n,v}^1,\tilde{P}_{n,v}^{(2)}(P_{n,v}^{0,1})}&=&0,
\end{eqnarray}
and the TMLE-update steps are using $\tilde{P}_{n,v}^1$ as HAL-MLE.
Let $\delta_{2n}$ be the solution and $\tilde{P}_{n,v}^*=\tilde{P}^2(\tilde{P}_{n,v}^1,P_{n,v}^{0,1},\delta_{2n})$ be the second targeted version of $\tilde{P}_{n,v}$. For simplicity, we assume that at this step we already achieved the desired equations (\ref{tildePnvtarget}).
Note,
\[
\tilde{P}_{n,v}^*=\tilde{P}^2(\tilde{P}^1(\tilde{P}_{n,v},P_{n,v}^0,\delta_{1n}),P_{n,v,\delta_{1n}}^{0,1},\delta_{2n}),\]
where we emphasized that the update $P_{n,v}^{0,1}$ of initial estimator $P_{n,v}^0$ is itself a second order TMLE based on HAL-MLE $\tilde{P}_{n,v}^1$ and does thus depend on $\delta_{1n}$ as well.
 Therefore, this final targeted update $\tilde{P}_{n,v}^*$ of the HAL-MLE can be represented as  a functional \[
\tilde{P}_{n,v}^*=\tilde{P}(\tilde{P}_{n,v},P_{n,v},\delta_{1n},\delta_{2n}),\] involving 2 (sequentially computed) $\delta_{1n}$ and $\delta_{2n}$. 
In addition, we have that $\tilde{P}(\tilde{P}_{n,v},P_{n,v},0,0)=\tilde{P}_{n,v}$, so that it returns the initial HAL-MLE $\tilde{P}_{n,v}$ at $\delta_1=\delta_2=0$.

After one round the initial estimator $P_{n,v}^0$  gets updated to the second order TMLE $\tilde{P}_{n,v}^{0,*}=\tilde{P}_{n,v}^{(1)}\tilde{P}_{n,v}^{(2)}(P_{n,v}^0)$, where the TMLE-updates are under $\tilde{P}_{n,v}^1$.
Therefore, the final updated initial estimator can be represented as a functional $\tilde{P}_{n,v}^{0,*}=\tilde{P}^0(P_{n,v}^0,P_{n,v},\delta_{1n})$, an update of $P_{n,v}^0$ only depending on the training sample $P_{n,v}$, beyond $\delta_{1n}$.
Moreover, when we evaluate this functional at $\delta_1=0$, then $\tilde{P}^0(P_{n,v}^0,P_{n,v},\delta_1=0)=\tilde{P}_{n,v}^{(2)}\tilde{P}_{n,v}^{(1)}(P_{n,v}^0)$, where the TMLE-updates are under $\tilde{P}_{n,v}$, thereby making it a function of the training sample $P_{n,v}$ only. 

In short notation, we write $\tilde{P}_{n,v}^*=\tilde{P}_{n,v}^*(P_{n,v},\delta_n)$ and $P_{n,v}^{0,*}=P_{n,v}^{0,*}(P_{n,v},\delta_n)$ to emphasize that they only depend on the validation sample through $\delta_n$, and we know that at $\delta=0$ they reduce to initial HAL-MLE $\tilde{P}_{n,v}$ and the regular second order TMLE (without targeting of the HAL-MLE) $\tilde{P}_{n,v}^{(1)}\tilde{P}_{n,v}^{(2)}(P_{n,v}^0)$.

We are now ready to provide an analysis of the asymptotics of this CV higher order TMLE, showing that it provides us with an asymptotic distribution under weaker conditions than needed for the regular higher order TMLE.
We can apply the exact expansion (\ref{exactexphotmle}) above to this choice  of HAL-MLE $\tilde{P}_{n,v}^*$ and initial estimator $P_{n,v}^{0,*}$, and take the average  over $v=1,\ldots,V$ on both sides.
This resulting exact expansion has the term
$\frac{1}{V}\sum_v (\tilde{P}_{n,v}^*-P_0)D^{(1)}_{\tilde{P}_{n,v}^{(1)}(P_0)}$.
We can write this term  as:
\[
\begin{array}{l}
\frac{1}{V}\sum_v (\tilde{P}_{n,v}^*-P_0)D^{(1)}_{\tilde{P}_{n,v}^{(1)}(P_0)}=
\frac{1}{V}\sum_v (\tilde{P}_{n,v}^*-P_0)\{ D^{(1)}_{\tilde{P}_{n,v}^{(1)}(P_0)}-D^{(1)}_{\tilde{P}_{n,v}^{(1)}\tilde{P}_{n,v}^{(2)}(P_{n,v}^{0,*})}\}\\
+\frac{1}{V}\sum_v (\tilde{P}_{n,v}^*-P_0)D^{(1)}_{\tilde{P}_{n,v}^{(1)}\tilde{P}_{n,v}^{(2)}(P_{n,v}^{0,*})}.
\end{array}
\]
The last term equals $-\frac{1}{V}\sum_v P_0 D^{(1)}_{\tilde{P}_{n,v}^{(1)}\tilde{P}_{n,v}^{(2)}(P_{n,v}^{0,*})}$.
Using equation (\ref{tildePnvtarget}), this becomes
\[
\frac{1}{V}\sum_v  (P_{n,v}^1-P_0) D^{(1)}_{\tilde{P}_{n,v}^{(1)}\tilde{P}_{n,v}^{(2)}(P_{n,v}^{0,*})}.\]
Similarly, we can apply this approach to the term $\frac{1}{V}\sum_v  (\tilde{P}_{n,v}^*-P_0)D^{(2)}_{\tilde{P}_{n,v}^{(2)}(P_0)}$ in the exact expansion.
So this yields the following exact expansion for the CV-second order TMLE $\psi_n^*$ characterized by the HAL-MLE $\tilde{P}_{n,v}^*$ and initial estimator $P_{n,v}^{0,*}$:
\[
\begin{array}{l}
\frac{1}{V}\sum_v \Psi(\tilde{P}_{n,v}^{(1)}\tilde{P}_{n,v}^{(2)}(P_{n,v}^{0,*})-\Psi(P_0)=
\frac{1}{V}\sum_v (P_{n,v}^1-P_0)D^{(1)}_{\tilde{P}_{n,v}^{(1)}\tilde{P}_{n,v}^{(2)}(P_{n,v}^{0,*})}\\
+\frac{1}{V}\sum_v R^{(1)}(\tilde{P}_{n,v}^{(1)}(P_0),P_0)\\
+\frac{1}{V}\sum_v (P_{n,v}^1-P_0)D^{(2)}_{\tilde{P}_{n,v}^*,\tilde{P}_{n,v}^{(2)}(P_{n,v}^{0,*}) }+\frac{1}{V}\sum_v R^{(2)}(\tilde{P}_{n,v}^{(2)}(P_{n,v}^{0,*}),P_0)\\
+\frac{1}{V}\sum_v (\tilde{P}_{n,v}^*-P_0)\{ D^{(1)}_{\tilde{P}_{n,v}^{(1)}(P_0)}-D^{(1)}_{\tilde{P}_{n,v}^{(1)}\tilde{P}_{n,v}^{(2)}(P_{n,v}^{0,*})}\}\\
+\frac{1}{V}\sum_v (\tilde{P}_{n,v}^*-P_0)\{ D^{(2)}_{\tilde{P}_{n,v}^*,\tilde{P}_{n,v}^{(2)}(P_0)}-D^{(2)}_{\tilde{P}_{n,v}^*,\tilde{P}_{n,v}^{(2)}(P_{n,v}^{0,*})}\}
\\
+\frac{1}{V}\sum_v \left\{ R^{(2)}(\tilde{P}_{n,v}^{(2)}(P_{n,v}^{0,*}),\tilde{P}_{n,v}^*)-R^{(2)}(\tilde{P}_{n,v}^{(2)}(P_0),\tilde{P}_{n,v}^*)\right\},
\end{array}
\]
It is reasonable to assume that the dominating terms in this expansion for the CV-second order TMLE $\psi_n^*$ are the cross-validated empirical process terms: \[
\begin{array}{l}
\frac{1}{V}\sum_v \Psi(\tilde{P}_{n,v}^{(1)}\tilde{P}_{n,v}^{(2)}(P_{n,v}^0)-\Psi(P_0)=
\frac{1}{V}\sum_v (P_{n,v}^1-P_0)D^{(1)}_{\tilde{P}_{n,v}^{(1)}\tilde{P}_{n,v}^{(2)}(P_{n,v}^0)}\\
+\frac{1}{V}\sum_v (P_{n,v}^1-P_0)D^{(2)}_{\tilde{P}_{n,v}^*,\tilde{P}_{n,v}^{(2)}(P_{n,v}^0) }
+o_P(\mid \psi_n^*-\psi_0\mid)
\end{array}
\]
By our representation of $\tilde{P}_{n,v}^*$ and $P_{n,v}^{0,*}$ as functions of $P_{n,v}$ and $\delta_n$, we can represent $D^{(1)}_{\tilde{P}_{n,v}^{(1)}\tilde{P}_{n,v}^{(2)}(P_{n,v}^{0,*})}$ as a function of the training sample $P_{n,v}$ and of a multi-dimensional $\delta_n$ that is based on the whole sample.  Therefore, let's write it as $D^{(1)}_{P_{n,v},\delta_n}$. In addition, we know that $D^{(1)}_{P_{n,v},\delta=0}=D^{(1)}_{\tilde{P}_{n,v}^{(1)}\tilde{P}_{n,v}^{(2)}(P_{n,v}^0)}$ equals $D^{(1)}$ at the regular second order TMLE only based on $P_{n,v}$. Similarly, we can represent $D^{(2)}_{\tilde{P}_{n,v}^*,\tilde{P}_{n,v}^{(2)}(P_{n,v}^{0,*})}$ as $D^{(2)}_{P_{n,v},\delta_n}$, and 
$D^{(2)}_{P_{n,v},\delta=0}=D^{(2)}_{\tilde{P}_{n,v},\tilde{P}_{n,v}^{(2)}(P_{n,v}^0)}$.
Suppose that $\delta\rightarrow D^{(1)}_{P_{n,v},\delta}$ is differentiable in $\delta$ at $\delta_0=0$.
It is then reasonable to assume
\[
\frac{1}{V}\sum_v (P_{n,v}^1-P_0)\{D^{(1)}_{P_{n,v},\delta_n}-D^{(1)}_{P_{n,v},0}\}
=O_P(\pl \delta_n\pl) \pl \frac{1}{V}\sum_v (P_{n,v}^1-P_0)\frac{d}{d \delta_0}D^{(1)}_{P_{n,v},\delta_0}\pl.
\]
Therefore, we will assume that, for $j=1,2$, 
\[
\begin{array}{l}
\frac{1}{V}\sum_v (P_{n,v}^1-P_0)D^{(j)}_{P_{n,v},\delta_n}=
\frac{1}{V}\sum_v (P_{n,v}^1-P_0) D^{(j)}_{P_{n,v},\delta_0}+o_P(\mid \psi_n^*-\psi_0\mid)\\
=\frac{1}{V}\sum_v (P_{n,v}^1-P_0) D^{(j)}_{\tilde{P}_{n,v}^{(1)}\tilde{P}_{n,v}^{(2)}(P_{n,v}^0)}+o_P(\mid \psi_n^*-\psi_0\mid).
\end{array}
\]
Under these assumptions we have
\[
\begin{array}{l}
\frac{1}{V}\sum_v \Psi(\tilde{P}_{n,v}^{(1)}\tilde{P}_{n,v}^{(2)}(P_{n,v}^0)-\Psi(P_0)=
\frac{1}{V}\sum_v (P_{n,v}^1-P_0) D^{(1)}_{\tilde{P}_{n,v}^{(1)}\tilde{P}_{n,v}^{(2)}(P_{n,v}^0)}\\
+\frac{1}{V}\sum_v (P_{n,v}^1-P_0)D^{(2)}_{\tilde{P}_{n,v},\tilde{P}_{n,v}^{(2)}(P_{n,v}^0)}+
o_P(\psi_n-\psi_0).
\end{array}
\]
Let $\bar{D}_{n,v}= D^{(1)}_{\tilde{P}_{n,v}^{(1)}\tilde{P}_{n,v}^{(2)}(P_{n,v}^0)}+
D^{(2)}_{\tilde{P}_{n,v},\tilde{P}_{n,v}^{(2)}(P_{n,v}^0)}$. Then, 
$\psi_n^*-\psi_0=1/V\sum_v (P_{n,v}^1-P_0)\bar{D}_{n,v}+o_P(\mid \psi_n-\psi_0\mid)$.
For each $v$, conditional on training sample $P_{n,v}$, the standardized sample mean w.r.t. $P_{n,v}^1$ converges to a standard normal distribution. However, the $V$ normals are correlated, not guaranteeing that the average over $v$ of these $v$-specific terms (which marginally converge to a standard normal) converges to  a specified limit distribution. 

Therefore, we also assume that there is a fixed sequence $\tilde{D}_{n,v}$ (i.e., non random) representing an approximation of $\bar{D}_{n,v}$ so that $P_0\{\bar{D}_{n,v}-\tilde{D}_{n,v}\}^2/P_0\{\tilde{D}_{n,,v}\}^2\rightarrow_p 0$. In the case that $D^{(1)}_{P_0}\not =0$, we can simply select $\tilde{D}_{n,v}=D^{(1)}_{P_0}$. However, we also want to handle the case that at the true $P_0$ the first order efficient influence function equals zero, so that this expansion still allows us to obtain inference under such a challenging situation. In this more general case, one might define
$\tilde{D}_{n,v}(o)=E\bar{D}_{n,v}(o)$, i.e. the expectation of $\bar{D}_{n,v}$. 
We then want that the $L^2(P_0)$-norm of $\bar{D}_{n,v}$ dominates the approximation error $\pl \bar{D}_{n,v}-\tilde{D}_{n,v}\pl_{P_0}$. For example, if $P_{n,v}^0$ is an estimator of $P_0$ for which the bias dominates the standard error, then one expects this to hold: this might naturally hold, but  could also be achieved by slightly oversmoothing a carefully tuned initial estimator such as an HAL-MLE. 

Then, 
\[
\frac{1}{V}\sum_{v=1}^V (P_{n,v}^1-P_0)\bar{D}_{n,v}=
\frac{1}{V}\sum_{v=1}^V  (P_{n,v}^1-P_0)\tilde{D}_{n,v}+\frac{1}{V}\sum_v
(P_{n,v}^1-P_0)(\bar{D}_{n,v}-\tilde{D}_{n,v}),\]
where, by our assumption on $\tilde{D}_{n,v}$ above, the second term is of smaller order than the first term. 
The leading term is now an average of $V$ sample means of mean zero  independent random variables, and, moreover, due to $\tilde{D}_{n,v}$ being a fixed function, the $V$ sample means are independent. Let $\sigma_{n,v}^2=P_0\{\tilde{D}_{n,v}\}^2$. Define $\sigma^2_n=\frac{1}{V}\sum_{v=1}^V \sigma_{n,v}^2$.

Then,
\[
\begin{array}{l}
n^{1/2}/\sigma_n \frac{1}{V}\sum_{v=1}^V  (P_{n,v}^1-P_0)\tilde{D}_{n,v}\\
=
\sigma_n^{-1} \frac{1}{V}\sum_{v=1}^V \sigma_{n,v}  V^{1/2} (n/V)^{1/2}(P_{n,v}^1-P_0)\tilde{D}_{n,v}/\sigma_{n,v}
\Rightarrow_d N(0,1).
\end{array}
\]
This proves that $n^{1/2}\sigma_n^{-1}(\psi_n^*-\psi_0)\Rightarrow_d N(0,1)$.
This implies that $\psi_n^*\pm 1.96 n^{-1/2}\sigma_n$ is an asymptotic $0.95$-confidence interval for $\psi_0$. Moreover, the variance estimator $\sigma^2_n$ also takes into account the second order influence function, thereby picking up second order behavior of the sampling distribution of $\psi_n^*$. In addition, this allows that $\sigma_n\rightarrow_p 0$.

\end{document}